\documentclass[leqno,openany,draft]{book}

\usepackage{makeidx}
\makeindex

\def\diam{\mathop{\rm diam}}
\def\dim{\mathop{\rm dim}}
\def\supp{\mathop{\rm supp}}
\def\re{\mathop{\rm Re}}

\newtheorem{theorem}{Theorem}
\newtheorem{lemma}[theorem]{Lemma}
\newtheorem{proposition}[theorem]{Proposition}
\newtheorem{definition}[theorem]{Definition}
\newtheorem{corollary}[theorem]{Corollary}

\newcommand{\begintheorem}{\addtocounter{equation}{1}\begin{theorem}}
\newcommand{\beginlemma}{\addtocounter{equation}{1}\begin{lemma}}
\newcommand{\beginproposition}{\addtocounter{equation}{1}\begin{proposition}}
\newcommand{\begindefinition}{\addtocounter{equation}{1}\begin{definition}}
\newcommand{\begincorollary}{\addtocounter{equation}{1}\begin{corollary}}

\begin{document}

\frontmatter

\title{Some aspects of analysis related to \\ $p$-adic numbers}

\author{Stephen Semmes \\
        Rice University}

\date{}

\maketitle

\chapter*{Preface}

        A field with an absolute value function is a basic type of
metric space, which includes the real and complex numbers with their
standard metrics, and ultrametrics on fields like the $p$-adic
numbers.  Here we try to give some perspectives of analysis in
situations like these.

\tableofcontents

\mainmatter

\chapter{Preliminaries}
\label{perliminaries}

\section{Metrics and ultrametrics}
\label{metrics, ultrametrics}

        As usual, a \emph{metric space}\index{metric spaces} is a 
set $M$ with a nonnegative real-valued function $d(x, y)$ defined for
$x, y \in M$ that satisfies the following three conditions: first,
\begin{equation}
\label{d(x, y) = 0 if and only if x = y}
        d(x, y) = 0 \quad\hbox{if and only if}\quad x = y;
\end{equation}
second, $d(x, y)$ is symmetric in $x$ and $y$, so that
\begin{equation}
\label{d(x, y) = d(y, x)}
        d(x, y) = d(y, x)
\end{equation}
for every $x, y \in M$; and third,
\begin{equation}
\label{d(x, z) le d(x, y) + d(y, z)}
        d(x, z) \le d(x, y) + d(y, z)
\end{equation}
for every $x, y, z \in M$, which is known as the \emph{triangle
  inequality}.\index{triangle inequality}  If
\begin{equation}
\label{d(x, z) le max(d(x, y), d(y, z))}
        d(x, z) \le \max(d(x, y), d(y, z))
\end{equation}
for every $x, y, z \in M$, then the metric $d(x, y)$ is said to be an
\emph{ultrametric}\index{ultrametrics} on $M$.  Of course, (\ref{d(x,
  z) le max(d(x, y), d(y, z))}) implies (\ref{d(x, z) le d(x, y) +
  d(y, z)}).  The \emph{discrete metric}\index{discrete metric} on a
set $M$ is defined by putting $d(x, y) = 1$ when $x \ne y$, and $d(x,
y) = 0$ when $x = y$, and is an ultrametric on $M$.

        Let $a$ be a real number such that $0 < a \le 1$, and let
$r$, $t$ be nonnegative real numbers.  Observe that
\begin{equation}
\label{max(r, t) le (r^a + t^a)^{1/a}}
        \max(r, t) \le (r^a + t^a)^{1/a},
\end{equation}
and hence that
\begin{equation}
\label{r + t le max(r, t)^{1 - a} (r^a + t^a) le ... = (r^a + t^a)^{1/a}}
 r + t \le \max(r, t)^{1 - a} \, (r^a + t^a) \le (r^a + t^a)^{((1 - a)/a) + 1}
                                              = (r^a + t^a)^{1/a}.
\end{equation}
Equivalently,
\begin{equation}
\label{(r + t)^a le r^a + t^a}
        (r + t)^a \le r^a + t^a.
\end{equation}
If $d(x, y)$ is a metric on a set $M$ and $0 < a \le 1$, then it
follows that $d(x, y)^a$ also defines a metric on $M$, which
determines the same topology on $M$ as $d(x, y)$.  Similarly, if $d(x,
y)$ is an ultrametric on $M$, then $d(x, y)^a$ is an ultrametric on
$M$ for every $a > 0$, which defines the same topology on $M$ as
$d(x, y)$.

        Remember that the absolute value $|x|$ of a real number $x$
is defined to be equal to $x$ when $x \ge 0$ and to $-x$ when $x \le
0$.  Of course,
\begin{equation}
\label{|x + y| le |x| + |y|}
        |x + y| \le |x| + |y|
\end{equation}
and
\begin{equation}
\label{|x y| = |x| |y|}
        |x \, y| = |x| \, |y|
\end{equation}
for all real numbers $x$, $y$.  The standard Euclidean metric on the
set ${\bf R}$\index{R@${\bf R}$} of real numbers is defined by
\begin{equation}
\label{d(x, y) = |x - y|}
        d(x, y) = |x - y|,
\end{equation}
and is not an ultrametric.  Note that $|x|^a$ satisfies the analogue
of (\ref{|x + y| le |x| + |y|}) when $0 < a \le 1$, by (\ref{(r + t)^a
  le r^a + t^a}), and that it satisfies the analogue of (\ref{|x y| =
  |x| |y|}) for every $a > 0$.

        Let $(M, d(x, y))$ be a metric space.  The \emph{open 
ball}\index{open balls} centered at a point $x \in M$ and with radius
        $r > 0$ is defined by
\begin{equation}
\label{B(x, r) = {z in M : d(x, z) < r}}
        B(x, r) = \{z \in M : d(x, z) < r\}.
\end{equation}
If $y \in B(x, r)$, then $t = r - d(x, y) > 0$, and one can check that
\begin{equation}
\label{B(y, t) subseteq B(x, r)}
        B(y, t) \subseteq B(x, r),
\end{equation}
using the triangle inequality.  If $d(\cdot, \cdot)$ is an ultrametric
on $M$, then it is easy to see that
\begin{equation}
\label{B(y, r) subseteq B(x, r)}
        B(y, r) \subseteq B(x, r)
\end{equation}
under these conditions.  More precisely,
\begin{equation}
\label{B(x, r) = B(y, r)}
        B(x, r) = B(y, r)
\end{equation}
when $d(x, y) < r$, since the opposite inclusion may be obtained by
reversing the roles of $x$ and $y$.

        Similarly, the \emph{closed ball}\index{closed balls} centered at
$x$ with radius $r \ge 0$ is defined by
\begin{equation}
\label{overline{B}(x, r) = {z in M : d(x, z) le r}}
        \overline{B}(x, r) = \{z \in M : d(x, z) \le r\}.
\end{equation}
If $d(\cdot, \cdot)$ is an ultrametric on $M$, and $d(x, y) \le r$, then
\begin{equation}
\label{overline{B}(y, r) subseteq overline{B}(x, r)}
        \overline{B}(y, r) \subseteq \overline{B}(x, r),
\end{equation}
and hence
\begin{equation}
\label{overline{B}(x, r) = overline{B}(y, r)}
        \overline{B}(x, r) = \overline{B}(y, r),
\end{equation}
as before.  In particular, this implies that closed balls in $M$ of
positive radius are open sets when $d(\cdot, \cdot)$ is an
ultrametric.  One can also check that open balls in $M$ are closed
sets when $d(\cdot, \cdot)$ is an ultrametric, using (\ref{B(y, r)
  subseteq B(x, r)}).  Equivalently, the complement of an open ball in
$M$ is an open set when $d(\cdot, \cdot)$ is an ultrametric, which can
also be derived from the remarks in the next paragraph.

        Let us continue to suppose that $d(\cdot, \cdot)$ is an ultrametric
on $M$.  If $x, y, z \in M$ satisfy $d(y, z) \le d(x, y)$, then
\begin{equation}
\label{d(x, z) le max(d(x, y), d(y, z)) = d(x, y)}
        d(x, z) \le \max(d(x, y), d(y, z)) = d(x, y),
\end{equation}
by the ultrametric version of the triangle inequality.  If $d(y, z) <
d(x, y)$, then
\begin{equation}
\label{d(x, y) le max(d(x, z), d(y, z))}
        d(x, y) \le \max(d(x, z), d(y, z))
\end{equation}
implies that $d(x, y) \le d(x, z)$.  Combining this with (\ref{d(x, z)
  le max(d(x, y), d(y, z)) = d(x, y)}), we get that
\begin{equation}
\label{d(x, y) = d(x, z)}
        d(x, y) = d(x, z)
\end{equation}
when $d(y, z) < d(x, y)$.

        Let $(M, d(x, y))$ be an aribtrary metric space again.
As usual, a sequence $\{x_j\}_{j = 1}^\infty$ of elements of $M$
is said to be a \emph{Cauchy sequence}\index{Cauchy sequences} if
for each $\epsilon > 0$ there is a positive integer $L(\epsilon)$
such that
\begin{equation}
\label{d(x_j, x_l) < epsilon}
        d(x_j, x_l) < \epsilon
\end{equation}
for every $j, l \ge L(\epsilon)$.  Convergent sequences are Cauchy
sequences, and $M$ is said to be \emph{complete}\index{complete metric
  spaces} if every Cauchy sequence in $M$ converges to an element of
$M$.  If $\{x_j\}_{j = 1}^\infty$ is a Cauchy sequence in $M$, then
it follows that
\begin{equation}
\label{lim_{j to infty} d(x_j, x_{j + 1}) = 0}
        \lim_{j \to \infty} d(x_j, x_{j + 1}) = 0,
\end{equation}
by taking $l = j + 1$ in (\ref{d(x_j, x_l) < epsilon}).  Conversely,
if $d(\cdot, \cdot)$ is an ultrametric, and if $\{x_j\}_{j =
  1}^\infty$ is a sequence of elements of $M$ that satisfies
(\ref{lim_{j to infty} d(x_j, x_{j + 1}) = 0}), then one can check
that $\{x_j\}_{j = 1}^\infty$ is a Cauchy sequence in $M$.

\section{Completions}
\label{completions}

        Let $(M, d(x, y))$ be a metric space, and let $\{x_j\}_{j = 1}^\infty$
and $\{y_j\}_{j = 1}^\infty$ be Cauchy sequences of elements of $M$.
If
\begin{equation}
\label{lim_{j to infty} d(x_j, y_j) = 0}
        \lim_{j \to \infty} d(x_j, y_j) = 0,
\end{equation}
then $\{x_j\}_{j = 1}^\infty$ and $\{y_j\}_{j = 1}^\infty$ are said to
be \emph{equivalent} as Cauchy sequences in $M$.  It is easy to see
that this defines an equivalence relation on the collection of Cauchy
sequences in $M$.  Of course, convergent sequences in $M$ are Cauchy
sequences, and two convergent sequences in $M$ are equivalent as
Cauchy sequences if and only if they converge to the same element of
$M$.  Similarly, if a Cauchy sequence in $M$ is equivalent to a
convergent sequence in $M$, then that Cauchy sequence converges to the
same element of $M$.  In particular, constant sequences in $M$ are
convergent, and a Cauchy sequence in $M$ converges to an element of
$M$ if and only if it is equivalent to the corresponding constant
sequence.  The \emph{completion}\index{completions} of $M$ is defined
to be the set of equivalence classes of Cauchy sequences in $M$.
There is a natural embedding of $M$ into its completion, which
associates to each $x \in M$ the equivalence class of Cauchy sequences
that contains the constant sequence $\{x_j\}_{j = 1}^\infty$ with $x_j
= x$ for each $j$.

        Suppose for the moment that $M$ is the set
${\bf Q}$\index{Q@${\bf Q}$} of rational numbers, equipped with the
standard metric (\ref{d(x, y) = |x - y|}).  It is well known that the
completion of ${\bf Q}$ with respect to the standard metric can be
identified with the real line.  More precisely, if one uses this to
construct the real numbers, then one should not use the definition of
a metric space in the previous section, which assumes that the real
numbers have already been defined.  However, the standard metric on
${\bf Q}$ still makes sense, and takes values in ${\bf Q}$.  One can
also define what it means for a sequence of rational numbers to
converge to a rational number with respect to the standard metric in
the usual way, and what it means for a sequence of rational numbers to
be a Cauchy sequence.  Thus the completion of ${\bf Q}$ can be defined
as in the previous paragraph, with the same properties as before.  We
shall return to this in a moment.

        Let $(M, d(x, y))$ be an arbitrary metric space again.
If $\{x_j\}_{j = 1}^\infty$ and $\{y_j\}_{j = 1}^\infty$ are Cauchy
sequences in $M$, then it is well known that $\{d(x_j, y_j)\}_{j =
  1}^\infty$ is a Cauchy sequence in ${\bf R}$, with respect to the
standard metric on ${\bf R}$.  It follows that $\{d(x_j, y_j)\}_{j =
  1}^\infty$ converges in ${\bf R}$, by the completeness of ${\bf R}$.
If $\{x'_j\}_{j = 1}^\infty$ and $\{y'_j\}_{j = 1}^\infty$ are Cauchy
sequences in $M$ that are equivalent to $\{x_j\}_{j = 1}^\infty$ and
$\{y_j\}_{j = 1}^\infty$, respectively, then it is easy to see that
\begin{equation}
\label{lim_{j to infty} (d(x'_j, y'_j) - d(x_j, y_j)) = 0}
        \lim_{j \to \infty} (d(x'_j, y'_j) - d(x_j, y_j)) = 0,
\end{equation}
so that $\{d(x'_j, y'_j)\}_{j = 1}^\infty$ and $\{d(x_j, y_j)\}_{j =
  1}^\infty$ have the same limit in ${\bf R}$.  Thus the limit of
$\{d(x_j, y_j)\}_{j = 1}^\infty$ leads to a well-defined distance
function on the completion of $M$, and one can check that this
distance function is a metric on the completion of $M$.  It is well
known that the completion of $M$ is complete with respect to this
metric.  The natural embedding of $M$ into its completion preserves
distances, and maps $M$ onto a dense subset of its completion.

        Let us go back to the case where $M = {\bf Q}$, with the
standard metric.  If $\{x_j\}_{j = 1}^\infty$ is a Cauchy sequence in
${\bf Q}$, then it is easy to see that $\{|x_j|\}_{j = 1}^\infty$ is a
Cauchy sequence in ${\bf Q}$.  If $\{x'_j\}_{j = 1}^\infty$ is another
Cauchy sequence in ${\bf Q}$ that is equivalent to $\{x_j\}_{j =
  1}^\infty$, then one can check that $\{|x'_j|\}_{j = 1}^\infty$ and
$\{|x_j|\}_{j = 1}^\infty$ are equivalent as Cauchy sequences in ${\bf
  Q}$.  This permits one to extend the absolute value function to a
mapping from the completion of ${\bf Q}$ into itself.  Similarly, if
$\{x_j\}_{j = 1}^\infty$ and $\{y_j\}_{j = 1}^\infty$ are Cauchy
sequences in ${\bf Q}$, then $\{|x_j - y_j|\}_{j = 1}^\infty$ is a
Cauchy sequence in ${\bf Q}$, as in the previous paragraph.  If
$\{x'_j\}_{j = 1}^\infty$ and $\{y'_j\}_{j = 1}^\infty$ are Cauchy
sequences in ${\bf Q}$ that are equivalent to $\{x_j\}_{j = 1}^\infty$
and $\{y_j\}_{j = 1}^\infty$, then $\{|x'_j - y'_j|\}_{j = 1}^\infty$
and $\{|x_j - y_j|\}_{j = 1}^\infty$ are equivalent as Cauchy
sequences in ${\bf Q}$, as before.  Thus the standard metric on ${\bf
  Q}$ extends to a function defined on pairs of elements of the
completion of ${\bf Q}$, and with values in the completion of ${\bf
  Q}$.

        If $\{x_j\}_{j = 1}^\infty$ and $\{y_j\}_{j = 1}^\infty$ are
Cauchy sequences in ${\bf Q}$, then one can check that $\{x_j +
y_j\}_{j = 1}^\infty$ and $\{x_j \, y_j\}_{j = 1}^\infty$ are Cauchy
sequences too, by standard arguments.  In the case of products, this
uses the fact that Cauchy sequences are bounded.  If $\{x'_j\}_{j =
  1}^\infty$ and $\{y'_j\}_{j = 1}^\infty$ are Cauchy sequences in
${\bf Q}$ that are equivalent to $\{x_j\}_{j = 1}^\infty$ and
$\{y_j\}_{j = 1}^\infty$, respectively, then $\{x'_j + y'_j\}_{j =
  1}^\infty$ and $\{x'_j \, y'_j\}_{j = 1}^\infty$ are equivalent as
Cauchy sequences in ${\bf Q}$ to $\{x_j + y_j\}_{j = 1}^\infty$ and
$\{x_j \, y_j\}_{j = 1}^\infty$, respectively.  Using this, one can
extend addition and multiplication to the completion of ${\bf Q}$, so
that the completion of ${\bf Q}$ becomes a commutative ring.  Note
that the extension of the standard metric on ${\bf Q}$ to its
completion mentioned in the preceding paragraph is the same as the
extension of the absolute value to the completion of ${\bf Q}$ applied
to the difference of two elements of the completion of ${\bf Q}$.

        If $\{x_j\}_{j = 1}^\infty$ is a sequence of rational numbers
that does not converge to $0$, then there is an $r \in {\bf Q}$ such
that $r > 0$ and
\begin{equation}
\label{|x_j| ge 2 r for infinitely many j}
        |x_j| \ge 2 \, r \quad\hbox{for infinitely many } j.
\end{equation}
If $\{x_j\}_{j = 1}^\infty$ is also a Cauchy sequence of elements of
${\bf Q}$, then it follows that
\begin{equation}
\label{|x_j| ge r for all but finitely many j}
        |x_j| \ge r \quad\hbox{for all but finitely many } j,
\end{equation}
and in particular $x_j \ne 0$ for all but finitely many $j$.  If $x_j
\ne 0$ for every $j$, then one can use this to show that $\{1/x_j\}_{j
  = 1}^\infty$ is a Cauchy sequence in ${\bf Q}$.  If $\{x'_j\}_{j =
  1}^\infty$ is another Cauchy sequence of nonzero elements of ${\bf
  Q}$ that is equivalent to $\{x_j\}_{j = 1}^\infty$, then
$\{1/x'_j\}_{j = 1}^\infty$ is a Cauchy sequence in ${\bf Q}$ that is
equivalent to $\{1/x_j\}_{j = 1}^\infty$.  Using these remarks, one
can extend the mapping $x \mapsto 1/x$ to the nonzero elements of the
completion of ${\bf Q}$, so that the completion of ${\bf Q}$ becomes a
field.

        Let $\{x_j\}_{j = 1}^\infty$ be a sequence of elements of ${\bf Q}$
that does not converge to $0$ again, so that (\ref{|x_j| ge 2 r for
  infinitely many j}) holds for some $r \in {\bf Q}$ with $r > 0$.
This implies that either
\begin{equation}
\label{x_j ge 2 r for infinitely many j}
        x_j \ge 2 \, r \quad\hbox{for infinitely many } j,
\end{equation}
or that
\begin{equation}
\label{x_j le - 2 r for infinitely many j}
        x_j \le - 2 \, r \quad\hbox{for infinitely many } j.
\end{equation}
If $\{x_j\}_{j = 1}^\infty$ is a Cauchy sequence in ${\bf Q}$, then it
follows that either
\begin{equation}
\label{x_j ge r for all but finitely many j}
        x_j \ge r \quad\hbox{for all but finitely many } j,
\end{equation}
or that
\begin{equation}
\label{x_j le -r for all but finitely many j}
        x_j \le -r \quad\hbox{for all but finitely many } j.
\end{equation}
Let us say that $\{x_j\}_{j = 1}^\infty$ is positive in the first
case, and negative in the second case.  If $\{x'_j\}_{j = 1}^\infty$
is another Cauchy sequence in ${\bf Q}$ that is equivalent to
$\{x_j\}_{j = 1}^\infty$, then it is easy to see that $\{x'_j\}_{j =
  1}^\infty$ is positive when $\{x_j\}_{j = 1}^\infty$ is positive,
and that $\{x'_j\}_{j = 1}^\infty$ is negative when $\{x_j\}_{j =
  1}^\infty$ is negative.  This permits one to extend the standard
ordering on ${\bf Q}$ to its completion, with the usual properties
with respect to addition and multiplication.  Once the ordering on the
completion of ${\bf Q}$ is defined, it is easy to see that the
extensions of the absolute value and distance functions to the
completion of ${\bf Q}$ satisfy the corresponding triangle
inequalities.  Thus the completion of ${\bf Q}$ basically becomes a
metric space, but where the metric also takes values in the completion
of ${\bf Q}$.  One can define convergence of sequences in the
completion of ${\bf Q}$ in the usual way, as well as Cauchy sequences,
and show that the completion of ${\bf Q}$ is complete, as before.

        Suppose now that $d(x, y)$ is an ultrametric on a set $M$.
If $\{x_j\}_{j = 1}^\infty$ and $\{y_j\}_{j = 1}^\infty$ are Cauchy sequences
of elements of $M$ that are not equivalent, then there is an $r > 0$
such that
\begin{equation}
\label{d(x_j, y_j) ge r}
        d(x_j, y_j) \ge r
\end{equation}
for infinitely many $j$, and in fact for all but finitely many $j$.
This implies that $d(x_j, y_j)$ is eventually constant in this case,
as in (\ref{d(x, y) = d(x, z)}).  Similarly, if $\{x'_j\}_{j =
  1}^\infty$ and $\{y'_j\}_{j = 1}^\infty$ are Cauchy sequences in $M$
that are equivalent to $\{x_j\}_{j = 1}^\infty$ and $\{y_j\}_{j =
  1}^\infty$, respectively, then
\begin{equation}
\label{d(x'_j, y'_j) = d(x_j, y_j)}
        d(x'_j, y'_j) = d(x_j, y_j)
\end{equation}
for all but finitely many $j$ under these conditions.  It is easy to see
that the extension of $d(\cdot, \cdot)$ to the completion of $M$ is also
an ultrametric in this situation.

\section{Continuous extensions}
\label{continuous extensions}

        Let $(M_1, d_1(x, y))$ and $(M_2, d_2(u, v))$ be metric spaces,
and suppose that $f$ is a uniformly continuous mapping from $M_1$ into
$M_2$.\index{uniform continuity} Thus for each $\epsilon > 0$ there is
a $\delta = \delta(\epsilon) > 0$ such that
\begin{equation}
\label{d_2(f(x), f(y)) < epsilon}
        d_2(f(x), f(y)) < \epsilon
\end{equation}
for every $x, y \in M_1$ with $d_1(x, y) < \delta$.  If $\{x_j\}_{j =
  1}^\infty$ and $\{x'_j\}_{j = 1}^\infty$ are sequences of elements
of $M_1$ such that
\begin{equation}
\label{lim_{j to infty} d_1(x_j, x'_j) = 0}
        \lim_{j \to \infty} d_1(x_j, x'_j) = 0,
\end{equation}
then it is easy to see that
\begin{equation}
\label{lim_{j to infty} d_2(f(x_j), f(x'_j)) = 0}
        \lim_{j \to \infty} d_2(f(x_j), f(x'_j)) = 0.
\end{equation}
Conversely, if $f$ is not uniformly continuous, then there are
sequences $\{x_j\}_{j = 1}^\infty$ and $\{x'_j\}_{j = 1}^\infty$ in
$M_1$ that satisfy (\ref{lim_{j to infty} d_1(x_j, x'_j) = 0}) and not
(\ref{lim_{j to infty} d_2(f(x_j), f(x'_j)) = 0}).  If $f$ is
uniformly continuous and $\{x_j\}_{j = 1}^\infty$ is a Cauchy sequence
of elements of $M_1$, then $\{f(x_j)\}_{j = 1}^\infty$ is a Cauchy
sequence of elements of $M_2$, which converges to an element of $M_2$
when $M_2$ is complete.

        Now let $E$ be a dense subset of $M_1$, and suppose that $f$
is a uniformly continuous mapping from $E$ into $M_2$, with respect to
the restriction of $d_1(x, y)$ to $x, y \in E$.  If $M_2$ is complete,
then it is well known that there is a unique extension of $f$ to a
uniformly continuous mapping from $M_1$ into $M_2$.  More precisely,
uniqueness only requires continuity instead of uniform continuity, and
completeness of $M_2$ is not needed.  To get the existence of such an
extension, let $x \in M_1$ be given, and let $\{x_j\}_{j = 1}^\infty$
be a sequence of elements of $E$ that converges to $x$ in $M_1$.  Thus
$\{x_j\}_{j = 1}^\infty$ is a Cauchy sequence in $E$, so that
$\{f(x_j)\}_{j = 1}^\infty$ converges in $M_2$, as before.  If
$\{x'_j\}_{j = 1}^\infty$ is another sequence of elements of $E$ that
converges to $x$ in $M_1$, then $\{x_j\}_{j = 1}^\infty$ and
$\{x'_j\}_{j = 1}^\infty$ satisfy (\ref{lim_{j to infty} d_1(x_j,
  x'_j) = 0}), and hence (\ref{lim_{j to infty} d_2(f(x_j), f(x'_j)) = 0}).
This implies that $\{f(x_j)\}_{j = 1}^\infty$ and $\{f(x'_j)\}_{j = 1}^\infty$
converge to the same element of $M_2$.  If we put $f(x)$ equal to the
limit of $\{f(x_j)\}_{j = 1}^\infty$ under these conditions, then this
agress with the original definition of $f(x)$ when $x \in E$, and it does
not depend on the choice of sequence $\{x_j\}_{j = 1}^\infty$.

        To show that this extension is uniformly continuous on $M_1$,
let $\epsilon > 0$ be given, and let $\delta = \delta(\epsilon)$ be a
positive real number such that (\ref{d_2(f(x), f(y)) < epsilon}) holds
for every $x, y \in E$ with $d_1(x, y) < \delta$.  Let $x, y \in M_1$
be given, with $d_1(x, y) < \delta$, and let $\{x_j\}_{j = 1}^\infty$
and $\{y_j\}_{j = 1}^\infty$ be sequences of elements of $E$ that converge
to $x$ and $y$ in $M_1$, respectively.  Thus $d_1(x_j, y_j) < \delta$
for all sufficiently large $j$, so that
\begin{equation}
\label{d_2(f(x_j), f(y_j)) < epsilon}
        d_2(f(x_j), f(y_j)) < \epsilon
\end{equation}
for all sufficiently large $j$.  This implies that
\begin{equation}
\label{d_2(f(x), f(y)) le epsilon}
        d_2(f(x), f(y)) \le \epsilon,
\end{equation}
since $\{f(x_j)\}_{j = 1}^\infty$ and $\{f(y_j)\}_{j = 1}^\infty$
converge to $f(x)$ and $f(y)$ in $M_2$, respectively, by construction.
Similarly, if $f : E \to M_2$ is an isometric
embedding,\index{isometric embeddings} in the sense that
\begin{equation}
\label{d_2(f(x), f(y)) = d_1(x, y)}
        d_2(f(x), f(y)) = d_1(x, y)
\end{equation}
for every $x, y \in E$, then $f$ is obviously uniformly continuous,
and this extension of $f$ to $M_1$ satisfies (\ref{d_2(f(x), f(y)) =
  d_1(x, y)}) for every $x, y \in M_1$.

        Let $(M, d(x, y))$ be a metric space, and suppose that $\phi_1$
and $\phi_2$ are isometric embeddings of $M$ into $M_1$ and $M_2$,
respectively, so that
\begin{equation}
\label{d_1(phi_1(x), phi_1(y)) = d_2(phi_2(x), phi_2(y)) = d(x, y)}
        d_1(\phi_1(x), \phi_1(y)) = d_2(\phi_2(x), \phi_2(y)) = d(x, y)
\end{equation}
for every $x, y \in M$.  Suppose also that $\phi_j(M)$ is dense in
$M_j$ for $j = 1, 2$, which can always be arranged by replacing $M_j$
with the closure of $\phi_j(M_j)$.  Put
\begin{equation}
\label{f = phi_2 circ phi_1^{-1}}
        f = \phi_2 \circ \phi_1^{-1}
\end{equation}
on $\phi_1(M)$, which is an isometric embedding of $\phi_1(M)$ into
$M_2$, with respect to the restriction of $d_1(\cdot, \cdot)$ to
$\phi_1(M)$.  If $M_2$ is complete, then $f$ has a unique extension to
an isometric embedding of $M_1$ into $M_2$, as before.  If $M_1$ is
complete, then $f(M_1)$ is complete with respect to the restriction of
$d_2(\cdot, \cdot)$ to $f(M_1)$. This implies that $f(M_1)$ is a
closed subset of $M_2$, because any sequence of elements of $f(M_1)$
that converges to an element of $M_2$ is a Cauchy sequence in
$f(M_1)$, and hence converges to an element of $f(M_1)$, by
competeness.  It follows that
\begin{equation}
\label{f(M_1) = M_2}
        f(M_1) = M_2
\end{equation}
under these conditions, because $f(M_1) = \phi_2(M)$ is dense in $M_2$,
by hypothesis.

\section{Quasimetrics}
\label{quasimetrics}

        Let $M$ be a set, and let $d(x, y)$ be a nonnegative real-valued
function defined for $x, y \in M$ such that
\begin{equation}
\label{d(x, y) = 0 if and only if x = y, 2}
        d(x, y) = 0 \quad\hbox{if and only if}\quad x = y,
\end{equation}
and
\begin{equation}
\label{d(x, y) = d(y, x), 2}
        d(x, y) = d(y, x)
\end{equation}
for every $x, y \in M$.  We say that $d(x, y)$ is a
\emph{quasimetric}\index{quasimetrics} on $M$ if
\begin{equation}
\label{d(x, z) le C (d(x, y) + d(y, z))}
        d(x, z) \le C \, (d(x, y) + d(y, z))
\end{equation}
for some $C \ge 1$ and every $x, y, z \in M$.  This is equivalent to
asking that
\begin{equation}
\label{d(x, z) le C' max(d(x, y), d(y, z))}
        d(x, z) \le C' \, \max(d(x, y), d(y, z))
\end{equation}
for some $C' \ge 1$ and every $x, y, z \in M$.  More precisely,
(\ref{d(x, z) le C' max(d(x, y), d(y, z))}) implies (\ref{d(x, z) le C
  (d(x, y) + d(y, z))}) with $C = C'$, and (\ref{d(x, z) le C (d(x, y)
  + d(y, z))}) implies (\ref{d(x, z) le C' max(d(x, y), d(y, z))})
with $C' = 2 \, C$.  Of course, (\ref{d(x, z) le C (d(x, y) + d(y,
  z))}) reduces to the ordinary triangle inequality (\ref{d(x, z) le
  d(x, y) + d(y, z)}) when $C = 1$, and (\ref{d(x, z) le C' max(d(x,
  y), d(y, z))}) reduces to the ultrametric version of the triangle
inequality (\ref{d(x, z) le max(d(x, y), d(y, z))}) when $C' = 1$.

        If $d(\cdot, \cdot)$ satisfies (\ref{d(x, z) le C' max(d(x, y), 
d(y, z))}) and $a$ is a positive real number, then
\begin{equation}
\label{d(x, z)^a le (C')^a max(d(x, y)^a, d(y, z)^a)}
        d(x, z)^a \le (C')^a \, \max(d(x, y)^a, d(y, z)^a)
\end{equation}
for every $x, y, z \in M$, so that $d(\cdot, \cdot)^a$ is also a
quasimetric on $M$.  Similarly, if $d(\cdot, \cdot)$ satisfies
(\ref{d(x, z) le C (d(x, y) + d(y, z))}) and $0 < a \le 1$, then
\begin{equation}
\label{d(x, z)^a le C^a (d(x, y) + d(y, z))^a le C^a (d(x, y)^a + d(y, z)^a)}
 d(x, z)^a \le C^a \, (d(x, y) + d(y, z))^a \le C^a \, (d(x, y)^a + d(y, z)^a)
\end{equation}
for every $x, y, z \in M$, using (\ref{(r + t)^a le r^a + t^a}) in the
second step.  If $a \ge 1$, then $f(r) = r^a$ is a convex function on
$[0, \infty)$, and hence
\begin{equation}
\label{(r + t)^a = 2^a (r/2 + t/2)^a le ... = 2^{a - 1} (r^a + t^a)}
        (r + t)^a = 2^a \, (r/2 + t/2)^a \le 2^a \, (r^a/2 + t^a/2)
                                          = 2^{a - 1} \, (r^a + t^a)
\end{equation}
for every $r, t \ge 0$.  This implies that
\begin{equation}
\label{d(x, z)^a le ... le 2^{a - 1} C^a (d(x, y)^a + d(y, z)^a)}
        d(x, z)^a \le C^a \, (d(x, y) + d(y, z))^a
                   \le 2^{a - 1} \, C^a \, (d(x, y)^a + d(y, z)^a)
\end{equation}
for every $x, y, z \in M$ when $d(\cdot, \cdot)$ satisfies (\ref{d(x,
  z) le C (d(x, y) + d(y, z))}) and $a \ge 1$.  Note that $|x - y|^a$
is not a metric on ${\bf R}$ when $a > 1$.

        Let $d(\cdot, \cdot)$ be a quasimetric on a set $M$ that
satisfies (\ref{d(x, z) le C (d(x, y) + d(y, z))}) for some $C \ge 1$.
Also let $n$ be a nonnegative integer, and let us check that
\begin{equation}
\label{d(x_0, x_{2^n}) le C^n sum_{j = 1}^{2^n} d(x_{j - 1}, x_j)}
        d(x_0, x_{2^n}) \le C^n \, \sum_{j = 1}^{2^n} d(x_{j - 1}, x_j)
\end{equation}
for any finite sequence $x_0, x_1, \ldots, x_{2^n}$ of $2^n + 1$
elements of $M$.  This is trivial when $n = 0$, and this is the same
as (\ref{d(x, z) le C (d(x, y) + d(y, z))}) when $n = 1$.  Suppose now
that (\ref{d(x_0, x_{2^n}) le C^n sum_{j = 1}^{2^n} d(x_{j - 1},
  x_j)}) holds for some $n \ge 0$, and let us verify that the
analogous statement holds for $n + 1$ as well.  If $x_0, x_1, \ldots,
x_{2^{n + 1}}$ is a finite sequence of $2^{n + 1} + 1$ elements of
$M$, then we can apply (\ref{d(x_0, x_{2^n}) le C^n sum_{j = 1}^{2^n}
  d(x_{j - 1}, x_j)}) to the first $2^n + 1$ terms $x_0, x_1, \ldots,
x_{2^n}$ of this sequence.  Similarly, we can apply the induction
hypothesis to the sequence $x_{2^n}, x_{2^n + 1}, \ldots, x_{2^{n + 1}}$
of $2^n + 1$ elements of $M$, to get that
\begin{equation}
\label{d(x_{2^n}, x_{2^{n + 1}}) le ...}
d(x_{2^n}, x_{2^{n + 1}}) \le C^n \, \sum_{j = 1}^{2^n} d(x_{2^n + j - 1}, x_{2^n + j}).
\end{equation}
It follows that
\begin{eqnarray}
\label{d(x_0, x_{2^{n + 1}}) le ...}
 \qquad d(x_0, x_{2^{n + 1}})
                    & \le & C \, (d(x_0, x_{2^n}) + d(x_{2^n}, x_{2^{n + 1}})) \\
 & \le & C^{n + 1} \, \sum_{j = 1}^{2^n} d(x_{j - 1}, x_j)
      + C^{n + 1} \, \sum_{j = 1}^{2^n} d(x_{2^n + j - 1}, x_{2^n + j}) \nonumber \\
 & = & C^{n + 1} \, \sum_{j = 1}^{2^{n + 1}} d(x_{j - 1}, x_j) \nonumber
\end{eqnarray}
using (\ref{d(x, z) le C (d(x, y) + d(y, z))}) in the first step.

        If instead $d(\cdot, \cdot)$ satisfies
(\ref{d(x, z) le C' max(d(x, y), d(y, z))}) for some $C' \ge 1$, then
\begin{equation}
\label{d(x_0, x_{2^n}) le (C')^n max {d(x_{j - 1}, x_j) : j = 1, ldots, 2^n}}
 d(x_0, x_{2^n}) \le (C')^n \, \max \{d(x_{j - 1}, x_j) : j = 1, \ldots, 2^n\}
\end{equation}
for any sequence finite $x_0, x_1, \ldots, x_{2^n}$ of $2^n + 1$
elements of $M$.  As before, this is trivial when $n = 0$, and this is
the same as (\ref{d(x, z) le C' max(d(x, y), d(y, z))}) when $n = 1$.
If (\ref{d(x_0, x_{2^n}) le (C')^n max {d(x_{j - 1}, x_j) : j = 1,
    ldots, 2^n}}) holds for some $n \ge 0$ and $x_0, x_1, \ldots,
x_{2^{n + 1}}$ is a finite sequence of $2^{n + 1} + 1$ elements of
$M$, then we can apply (\ref{d(x_0, x_{2^n}) le (C')^n max {d(x_{j -
      1}, x_j) : j = 1, ldots, 2^n}}) to the first $2^n + 1$ terms
$x_0, x_1, \ldots, x_{2^n}$ of this sequence.  We can also apply this
induction hypothesis to the sequence $x_{2^n}, x_{2^n + 1}, \ldots,
x_{2^{n + 1}}$ of $2^n + 1$ elements of $M$, to get that
\begin{equation}
\label{d(x_{2^n}, x_{2^{n + 1}}) le ..., 2}
 d(x_{2^n}, x_{2^{n + 1}}) \le (C')^n \, \max\{d(x_{2^n + j - 1}, x_{2^n + j}) :
                                               j = 1, \ldots, 2^n\}.
\end{equation}
It follows that
\begin{eqnarray}
\label{d(x_0, x_{2^{n + 1}}) le ..., 2}
d(x_0, x_{2^{n + 1}}) & \le & C' \, \max(d(x_0, x_{2^n}), d(x_{2^n}, x_{2^{n + 1}})) \\
 & \le & (C')^{n + 1} \, \max \{d(x_{j - 1}, x_j) : j = 1, \ldots, 2^{n + 1}\},
                                                                    \nonumber
\end{eqnarray}
using (\ref{d(x, z) le C' max(d(x, y), d(y, z))}) in the first step.

\section{Quasimetrics, 2}
\label{quasimetrics, 2}

        Let $d(x, y)$ be a quasimetric on a set $M$.  Thus the open
ball $B(x, r)$ in $M$ centered at a point $x \in M$ and with radius $r
> 0$ can be defined with respect to $d(x, y)$ as in (\ref{B(x, r) = {z
    in M : d(x, z) < r}}).  Let us say that a set $U \subseteq M$
is an open set if for each $x \in M$ there is an $r > 0$ such that
\begin{equation}
\label{B(x, r) subseteq U}
        B(x, r) \subseteq U,
\end{equation}
as usual.  It is easy to see that this defines a topology on $M$, in
the same way as for metric spaces.  If $d(x, y)$ is a metric on $M$,
then open balls in $M$ are open sets, as in (\ref{B(y, t) subseteq
  B(x, r)}).  This uses the ordinary version of the triangle
inequality in a significant way, and does not work for quasimetrics,
without additional hypotheses.  However, there are some substitutes
for this, as follows.

        Let us begin with some variants of (\ref{B(y, t) subseteq B(x, r)}).
Suppose that $d(x, y)$ satisfies (\ref{d(x, z) le C (d(x, y) + d(y,
  z))}) for some $C \ge 1$, and let $x \in M$ and $r > 0$ be given.
If $y \in M$ and $d(x, y) < (2 \, C)^{-1} \, r$, then one can check
that
\begin{equation}
\label{B(y, (2 C)^{-1} r) subseteq B(x, r)}
        B(y, (2 \, C)^{-1} \, r) \subseteq B(x, r).
\end{equation}
Similarly, if $d(\cdot, \cdot)$ satisfies (\ref{d(x, z) le C' max(d(x,
  y), d(y, z))}) for some $C' \ge 1$, then
\begin{equation}
\label{B(y, (C')^{-1} r) subseteq B(x, r)}
        B(y, (C')^{-1} \, r) \subseteq B(x, r)
\end{equation}
for every $y \in M$ with $d(x, y) < (C')^{-1} \, r$.

        Let $E$ be any subset of $M$, and put
\begin{equation}
\label{U = {x in M : B(x, r) subseteq E for some r > 0}}
        U = \{x \in M : B(x, r) \subseteq E \hbox{ for some } r > 0\}.
\end{equation}
Let $x \in U$ be given, and let $r$ be a positive real number such
that $B(x, r) \subseteq E$.  If $d(\cdot, \cdot)$ satisfies (\ref{d(x,
  z) le C (d(x, y) + d(y, z))}), then
\begin{equation}
\label{B(y, (2 C)^{-1} r) subseteq B(x, r) subseteq E}
        B(y, (2 \, C)^{-1} \, r) \subseteq B(x, r) \subseteq E
\end{equation}
for every $y \in M$ with $d(x, y) < (2 \, C)^{-1} \, r$, by (\ref{B(y,
  (2 C)^{-1} r) subseteq B(x, r)}).  Similarly, if $d(\cdot, \cdot)$
satsfies (\ref{B(y, (C')^{-1} r) subseteq B(x, r)}), then
\begin{equation}
\label{B(y, (C')^{-1} r) subseteq B(x, r) subseteq E}
        B(y, (C')^{-1} \, r) \subseteq B(x, r) \subseteq E
\end{equation}
for every $y \in M$ with $d(x, y) < (C')^{-1} \, r$.  In both cases,
it follows that $y \in U$, which means that $U$ contains an open ball
centered at $x$ with positive radius.  This implies that $U$ is an
open set in $M$.  Of course, any open subset of $M$ that is contained
in $E$ is also contained in $U$, so that $U$ is equal to the interior
of $E$ in $M$ with respect to this topology.

        Let $w \in M$ and $t > 0$ be given, and let us apply the previous
remarks to $E = B(w, t)$.  This leads to a set $U = U(w, t)$ as in
(\ref{U = {x in M : B(x, r) subseteq E for some r > 0}}), which is the
interior of $B(w, t)$.  Of course, $w \in U(w, t)$, by construction.
More precisely, if $d(\cdot, \cdot)$ satisfies (\ref{d(x, z) le C
  (d(x, y) + d(y, z))}), then we get that
\begin{equation}
\label{B(w, (2 C)^{-1} t) subseteq U(w, t)}
        B(w, (2 \, C)^{-1} \, t) \subseteq U(w, t),
\end{equation}
by (\ref{B(y, (2 C)^{-1} r) subseteq B(x, r)}).  Similarly, if
$d(\cdot, \cdot)$ satisfies (\ref{d(x, z) le C' max(d(x, y), d(y,
  z))}), then
\begin{equation}
\label{B(w, (C')^{-1} t) subseteq U(w, t)}
        B(w, (C')^{-1} \, t) \subseteq U(w, t),
\end{equation}
by (\ref{B(y, (C')^{-1} r) subseteq B(x, r)}).

        Now let $w, z \in M$ be given, with $w \ne z$, so that $d(w, z) > 0$.
If $d(\cdot, \cdot)$ satisfies (\ref{d(x, z) le C (d(x, y) + d(y,
  z))}), then it follows that
\begin{equation}
\label{B(w, (2 C)^{-1} d(w, z)) cap B(z, (2 C)^{-1} d(w, z)) = emptyset}
 B(w, (2 \, C)^{-1} \, d(w, z)) \cap B(z, (2 \, C)^{-1} \, d(w, z)) = \emptyset.
\end{equation}
Similarly, if $d(\cdot, \cdot)$ satisfies (\ref{d(x, z) le C' max(d(x,
  y), d(y, z))}), then
\begin{equation}
\label{B(w, (C')^{-1} d(w, z)) cap B(z, (C')^{-1} d(w, z)) = emptyset}
        B(w, (C')^{-1} \, d(w, z)) \cap B(z, (C')^{-1} \, d(w, z)) = \emptyset.
\end{equation}
This implies that $M$ is Hausdorff, since $w$ and $z$ are in the
interiors of the corresponding open balls, by the remarks in the
preceding paragraph.

        Let $\{x_j\}_{j = 1}^\infty$ be a sequence of elements of $M$,
and let $x$ be an element of $M$.  It is natural to say that
$\{x_j\}_{j = 1}^\infty$ converges to $x$ in $M$ with respect to
$d(\cdot, \cdot)$ when
\begin{equation}
\label{lim_{j to infty} d(x_j, x) = 0}
        \lim_{j \to \infty} d(x_j, x) = 0.
\end{equation}
Alternatively, one might use the definition of convergence of
sequences in a topological space, so that $\{x_j\}_{j = 1}^\infty$
converges to $x \in M$ if for every open set $U \subseteq M$ with $x
\in U$ there is an $L \ge 1$ such that
\begin{equation}
\label{x_j in U}
        x_j \in U
\end{equation}
for each $j \ge L$.  The first condition obviously implies the second
condition, by the definition of open sets in $M$.  Conversely, the
second condition implies the first condition, because every open ball
centered at $x$ contains an open set in $M$ that contains $x$ as an
element, as before.

        Similarly, one might like to say that a point $x \in M$
is a limit point of a set $E \subseteq M$ with respect to the
quasimetric $d(\cdot, \cdot)$ if for each $r > 0$ there is a $y \in E$
such that $x \ne y$ and $d(x, y) < r$.  The usual topological
definition says that $x \in M$ is a limit point of $E \subseteq M$ if
for each open set $U \subseteq M$ with $x \in U$, there is a point $y
\in E$ such that $x \ne y$ and $y \in U$.  It is easy to see that the
first definition implies the second definition in this situation, by
the definition of open sets in $M$.  Conversely, the second definition
implies the first definition, because every open ball in $M$ centered
at $x$ contains an open set that contains $x$ as an element.  As usual,
the closure of $E \subseteq M$ is the set $\overline{E}$ of $x \in M$
such that $x \in E$ or $x$ is a limit point of $E$, and is a closed
subset of $E$, and the topological characterization implies that
$\overline{E}$ is always a closed set in $M$.

        Let $x \in M$ and a positive integer $j$ be given, and let
$U_j(x)$ be the interior of $B(x, 1/j)$.  Thus $U_j(x)$ is an open
set that contains $x$ and is contained in $B(x, 1/j)$, as before.
If $U$ is any open set in $M$ that contains $x$, then $U_j(x)
\subseteq U$ for all sufficiently large $j$, by the definition of an
open set in $M$.  This shows that there is a local base for the
topology of $M$ at $x$ with only finitely or countably many elements,
as in the case of metric spaces.  In particular, this implies that
sequences can be used for many standard topological arguments
involving $M$, concerning limit points and continuity, for instance.

        Of course, one can define the closed ball $\overline{B}(x, r)$
centered at a point $x \in M$ and with radius $r \ge 0$ as in
(\ref{overline{B}(x, r) = {z in M : d(x, z) le r}}).  In a metric
space, closed balls are closed sets, but this does not work in
quasimetric spaces without additional hypotheses.  However, if
$d(\cdot, \cdot)$ satisfies (\ref{d(x, z) le C (d(x, y) + d(y, z))}),
then it is easy to see that the closure of $\overline{B}(x, r)$ is
contained in $\overline{B}(x, C \, r)$.  This uses the
characterization of limit points of subsets of $M$ in terms of
$d(\cdot, \cdot)$ mentioned earlier.  Similarly, if $d(\cdot, \cdot)$
satisfies (\ref{d(x, z) le C' max(d(x, y), d(y, z))}), then the
closure of $\overline{B}(x, r)$ is contained in $\overline{B}(x, C' \,
r)$.

        One can also define a uniform structure on $M$ corresponding
to the quasimetric $d(\cdot, \cdot)$ in essentially the same way as
for metric spaces, as in \cite{jk}.  The topology on $M$ determined by
$d(\cdot, \cdot)$ described earlier is the same as the topology
associated to this uniform structure as in \cite{jk}.  Note that the
characterization of the interior of a set $E \subseteq M$ as the set
$U$ in (\ref{U = {x in M : B(x, r) subseteq E for some r > 0}}) is the
same as Theorem 4 on p178 of \cite{jk} in this context.  Cauchy
sequences and uniform continuity can be defined for quasimetrics in
the same way as for metrics, and are determined by the corresponding
uniform structure as well.

        The metrization theorem for uniform spaces discussed in \cite{jk}
implies that there is a metric on $M$ that determines the same uniform
structure as the one associated to $d(x, y)$, and hence the same
topology.  Of course, many related properties of $M$ can be shown more
directly, as before.  Remember that $d(x, y)^a$ is a quasimetric on
$M$ for every positive real number $a$, as in the previous section.
It is easy to see that $d(x, y)^a$ determines the same uniform
structure on $M$ as $d(x, y)$ for each $a > 0$, and the same topology
on $M$ in particular.  In \cite{m-s}, it is shown that there is a
metric $d_0(x, y)$ on $M$ and a positive real number $a_0$ such that
$d(x, y)$ and $d_0(x, y)^{a_0}$ are each bounded by constant multiples
of the other.

\section{Lipschitz mappings}
\label{lipschitz mappings}

        Let $(M_1, d_1(x, y))$ and $(M_2, d_2(u, v))$ be quasimetric spaces,
so that $M_1$ and $M_2$ are sets, and $d_1(x, y)$ and $d_2(u, v)$ are
quasimetrics on them, respectively.  A mapping $f : M_1 \to M_2$ is said
to be \emph{Lipschitz of order $a > 0$}\index{Lipschitz mappings} if 
there is a nonnegative real number $C$ such that
\begin{equation}
\label{d_2(f(x), f(y)) le C d_1(x, y)^a}
        d_2(f(x), f(y)) \le C \, d_1(x, y)^a
\end{equation}
for every $x, y \in M_1$.  If $a = 1$, then one may simply say that
$f$ is Lipschitz.  Note that a Lipschitz mapping of any order is
uniformly continuous.  Of course, $f$ satisfies (\ref{d_2(f(x), f(y))
  le C d_1(x, y)^a}) with $C = 0$ if and only if $f$ is constant.

        Remember that $d_1(x, y)^a$ is also a quasimetric on $M_1$ for 
every $a > 0$, as in Section \ref{quasimetrics}.  Thus $f$ is
Lipschitz of order $a$ as a mapping from $(M_1, d_1(x, y))$ into
$(M_2, d_2(u, v))$ if and only if $f$ is Lipschitz of order $1$ as a
mapping from $(M_1, d_1(x, y)^a)$ into $(M_2, d_2(u, v))$, with the
same constant $C$.  Similarly, $f$ is Lipschitz of order $a$ with
constant $C$ as a mapping from $(M_1, d_1(x, y))$ into $(M_2, d_2(u,
v))$ if and only if $f$ is Lipschitz of order $1$ with constant
$C^{1/a}$ as a mapping from $(M_1, d_1(x, y))$ into $(M_2, d_2(u,
v)^{1/a})$.  If $d_1(x, y)$ and $d_2(u, v)$ are metrics on $M_1$ and
$M_2$, respectively, then $d_1(x, y)^a$ is a metric on $M_1$ when $0 <
a \le 1$, and $d_2(u, v)^{1/a}$ is a metric on $M_2$ when $a \ge 1$,
as in Section \ref{metrics, ultrametrics}.

        Let us now restrict our attention to the case where $M_2$ is
the real line, equipped with the standard metric.  If $f$ is a
real-valued function on $M_1$ that satisfies
\begin{equation}
\label{f(x) le f(y) + C d_1(x, y)}
        f(x) \le f(y) + C \, d_1(x, y)
\end{equation}
for some $C \ge 0$ and every $x, y \in M_1$, then we also have that
\begin{equation}
\label{f(y) le f(x) + C d_1(x, y)}
        f(y) \le f(x) + C \, d_1(x, y)
\end{equation}
for every $x, y \in M_1$, by interchanging the roles of $x$ and $y$.
This implies that
\begin{equation}
\label{|f(x) - f(y)| = max(f(x) - f(y), f(y) - f(x)) le C d_1(x, y)}
        |f(x) - f(y)| = \max(f(x) - f(y), f(y) - f(x)) \le C \, d_1(x, y)
\end{equation}
for every $x, y \in M_1$, so that $f$ is Lipschitz of order $1$ with
constant $C$.  In particular,
\begin{equation}
\label{f_p(x) = d_1(x, p)}
        f_p(x) = d_1(x, p)
\end{equation}
has this property with $C = 1$ for every $p \in M_1$ when $d_1(x, y)$
is a metric on $M_1$, by the triangle inequality.  In this case,
$d_1(x, y)^a$ also defines a metric on $M_1$ when $0 < a \le 1$, as in
Section \ref{metrics, ultrametrics}.  It follows that
\begin{equation}
\label{f_{p, a}(x) = d_1(x, p)^a}
        f_{p, a}(x) = d_1(x, p)^a
\end{equation}
defines a Lipschitz mapping from $(M_1, d_1(x, y)^a)$ into ${\bf R}$
of order $1$ with constant $C = 1$ for each $p \in M_1$, by the same
argument.  Equivalently, this means that (\ref{f_{p, a}(x) = d_1(x,
  p)^a}) is a Lipschitz mapping from $(M_1, d_1(x, y))$ into ${\bf R}$
of order $a$ with constant $C = 1$ for each $p \in M_1$ when $0 < a
\le 1$ and $d_1(\cdot, \cdot)$ is a metric on $M_1$.

          Suppose for the moment that $M_1 = {\bf R}$, and that
$f : {\bf R} \to {\bf R}$ is Lipschitz of order $a > 1$ with
respect to the standard metric on ${\bf R}$ on the domain and range.
It is easy to see that $f$ is constant on ${\bf R}$ under these
conditions, because $f'(x) = 0$ for every $x \in {\bf R}$.
Equivalently, if $M_1 = {\bf R}$ equipped with the quasimetric
\begin{equation}
\label{d_1(x, y) = |x - y|^a}
        d_1(x, y) = |x - y|^a
\end{equation}
for some $a > 1$, and if $f$ is Lipschitz of order $1$ as a mapping
from $(M_1, d_1(x, y))$ into ${\bf R}$ with the standard metric, then
$f$ is constant.  However, if $d_1(x, y)$ is any quasimetric on a set
$M_1$, then one can use metrics on $M_1$ as in \cite{m-s} to get
real-valued Lipschitz functions of positive order on $M_1$ with
respect to $d_1(x, y)$.

\section{Haar measure}
\label{haar measure}

        Let $A$ be a commutative group, with the group operations
expressed additively.  Suppose that $A$ is also equipped with a
topology, such that the group operations on $A$ are continuous.  More
precisely, this means that addition on $A$ is continuous as a mapping
from $A \times A$ into $A$, with respect to the product topology on $A
\times A$ associated to the given topology on $A$.  The mapping
\begin{equation}
\label{x mapsto -x}
        x \mapsto -x
\end{equation}
should be continuous on $A$ too, where $-x$ is the additive inverse of
$x \in A$.  In order for $A$ to be a \emph{topological
  group},\index{topological groups} it is customary to ask that
$\{0\}$ be a closed set in $A$.  It is well known that this implies
that $A$ is Hausdorff, and in fact regular as a topological space.
Note that the translation mapping
\begin{equation}
\label{x mapsto a + x}
        x \mapsto a + x
\end{equation}
is continuous on $A$ for every $a \in A$, because of continuity of
addition on $A$.  This implies that (\ref{x mapsto a + x}) is a
homeomorphism from $A$ onto itself for each $a \in A$, since the
inverse mapping corresponds to translation by $-a$.  Similarly,
(\ref{x mapsto -x}) is a homeomorphism on $A$, because it is its
own inverse mapping.

        Put
\begin{equation}
\label{-E = {-x : x in E}}
        -E = \{-x : x \in E\}
\end{equation}
for each $E \subseteq A$, and
\begin{equation}
\label{a + E = {a + x : x in E}}
        a + E = \{a + x : x \in E\}
\end{equation}
for each $a \in A$ and $E \subseteq A$.  If $E$ is an open set in $A$,
then (\ref{-E = {-x : x in E}}) is an open set in $A$ too, and (\ref{a
  + E = {a + x : x in E}}) is an open set in $A$ for every $a \in A$,
because (\ref{x mapsto -x}) and (\ref{x mapsto a + x}) are
homeomorphisms on $A$.  There are analogous statements for closed
sets, compact sets, and Borel sets.  In particular, if there is a
nonempty open subset of $A$ that is contained in a compact set, then
it follows that $A$ is locally compact as a topological space.

        If $A$ is locally compact, then a famous theorem states that
there is a well-behaved nonnegative translation-invariant Borel
measure $H$ on $A$, known as \emph{Haar measure}.\index{Haar measure}
To say that $H$ is invariant under translations on $A$ means that
\begin{equation}
\label{H(a + E) = H(E)}
        H(a + E) = H(E)
\end{equation}
for every Borel set $E \subseteq A$ and $a \in A$.  Haar measure is
also supposed to satisfy $H(U) > 0$ for every nonempty open set $U
\subseteq A$, $H(K) < \infty$ for every compact set $K \subseteq A$,
and some additional regularity properties.  It is well known that $H$
is uniquely determined up to multiplication by a positive real number
under these conditions.  Using this, one can show that
\begin{equation}
\label{H(-E) = H(E)}
        H(-E) = H(E)
\end{equation}
for every Borel set $E \subseteq A$.  Of course, any commutative group
$A$ is a locally compact topological group with respect to the
discrete topology, in which case counting measure on $A$ satisfies the
requirements of Haar measure.  The real line is a locally compact
commutative topological group with respect to addition and the
standard topology, and Lebesgue measure on ${\bf R}$ satisfies the
requirements of Haar measure.

        Let $A$ be a locally compact commutative topological group again,
and let $C_{com}(A)$ be the space of continuous real or complex-valued
functions on $A$ with compact support.  If $H$ satisfies the
requirements of Haar measure on $A$, then
\begin{equation}
\label{L(f) = int_A f dH}
        L(f) = \int_A f \, dH
\end{equation}
defines a nonnegative linear functional on $C_{com}(A)$.  More
precisely, if $f$ is a continuous nonnegative real-valued function
with compact support on $A$ such that $f(x) > 0$ for some $x \in A$,
then it is easy to see that $L(f)$ is a positive real number.  If $f
\in C_{com}(A)$ and $a \in A$, then
\begin{equation}
\label{f_a(x) = f(x + a)}
        f_a(x) = f(x + a)
\end{equation}
defines an element of $C_{com}(A)$, and
\begin{equation}
\label{L(f_a) = L(f)}
        L(f_a) = L(f),
\end{equation}
because of the translation-invariance of $H$.  A linear functional on
$C_{com}(A)$ with these properties is known as a \emph{Haar
  integral}\index{Haar integral} on $A$.  Haar measure on $A$ can also
be obtained from a Haar integral, using the Riesz representation
theorem.  A Haar integral can be defined on the real line using the
Riemann integral, for instance.

\chapter{Absolute value functions}
\label{absolute value functions}

\section{Definitions and examples}
\label{definitions, examples}

        Let $k$ be a field.  A nonnegative real-valued function $|\cdot|$
on $k$ is said to be an \emph{absolute value function}\index{absolute
value functions} on $k$ if it satisfies the following three conditions:
first,
\begin{equation}
\label{|x| = 0 if and only if x = 0}
        |x| = 0 \hbox{ if and only if } x = 0;
\end{equation}
second,
\begin{equation}
\label{|x y| = |x| |y|, 2}
        |x \, y| = |x| \, |y|
\end{equation}
for every $x, y \in k$; and third,
\begin{equation}
\label{|x + y| le |x| + |y|, 2}
        |x + y| \le |x| + |y|
\end{equation}
for every $x, y \in k$.  Of course, the standard absolute value
function on ${\bf R}$ satisfies these conditions, as in Section
\ref{metrics, ultrametrics}.  Similarly, it is well known that the
standard norm or modulus on the field ${\bf C}$ of complex numbers
satisfies these conditions.  If $k$ is any field, then the
\emph{trivial absolute value function}\index{trivial absolute value
function} on $k$ is defined by putting $|0| = 0$ and $|x| = 1$
for every $x \in k$ with $x \ne 0$, and is easily seen to satisfy
these conditions as well.

        Suppose for the moment that $|\cdot|$ is a nonnegative real-valued
function on $k$ that satisfies (\ref{|x| = 0 if and only if x = 0})
and (\ref{|x y| = |x| |y|, 2}).  Thus $|1| > 0$, since $1 \ne 0$ in
$k$, by definition of a field.  Here we use $0$ and $1$ to denote both
the additive and multiplicative identity elements in $k$ and their
counterparts in ${\bf R}$, and it should always be clear from the
context which is being considered in any given instance.  We also have
that $|1|^2 = |1^2| = |1|$, by (\ref{|x y| = |x| |y|, 2}), which
implies that
\begin{equation}
\label{|1| = 1}
        |1| = 1.
\end{equation}
Similarly, if $x \in k$ satisfies $x^n = 1$ for some positive integer
$n$, then
\begin{equation}
\label{|x|^n = |x^n| = |1| = 1}
        |x|^n = |x^n| = |1| = 1,
\end{equation}
and hence $|x| = 1$.  In particular, $(-1)^2 = 1$ in $k$, so that
\begin{equation}
\label{|-1| = 1}
        |-1| = 1.
\end{equation}
If $x \in k$ and $x \ne 0$, then $x$ has a multiplicative inverse $x^{-1}$
in $k$, which satisfies
\begin{equation}
\label{|x| |x^{-1}| = |1| = 1}
        |x| \, |x^{-1}| = |1| = 1,
\end{equation}
and hence
\begin{equation}
\label{|x^{-1}| = |x|^{-1}}
        |x^{-1}| = |x|^{-1}.
\end{equation}

        If $|\cdot|$ is an absolute value function on $k$, then it
follows from (\ref{|-1| = 1}) that
\begin{equation}
\label{d(x, y) = |x - y|, 2}
        d(x, y) = |x - y|
\end{equation}
is symmetric in $x$ and $y$.  Thus (\ref{d(x, y) = |x - y|, 2})
defines a metric on $k$.  Let us say that $|\cdot|$ is an
\emph{ultrametric absolute value function}\index{ultrametric absolute
  value functions} on $k$ if
\begin{equation}
\label{|x + y| le max(|x|, |y|)}
        |x + y| \le \max(|x|, |y|)
\end{equation}
for every $x, y \in k$.  This implies that the associated metric
(\ref{d(x, y) = |x - y|, 2}) is an ultrametric on $k$.  The trivial
absolute value function on any field $k$ is an ultrametric absolute
value function, for which the associated metric (\ref{d(x, y) = |x -
  y|, 2}) is the same as the discrete metric on $k$.

        The \emph{$p$-adic absolute value}\index{p-adic absolute
value@$p$-adic absolute value} $|\cdot|_p$ is defined on ${\bf Q}$
for each prime number $p$ as follows.  Let $x \in {\bf Q}$ be given,
and put $|x|_p = 0$ when $x = 0$.  Otherwise, if $x \ne 0$, then $x$
can be expressed as $p^j \, a / b$ for some integers $a$, $b$, and $j$
such that neither $a$ nor $b$ is an integer multiple of $p$, including
$0$.  In this case, we put
\begin{equation}
\label{|x|_p = p^{-j}}
        |x|_p = p^{-j},
\end{equation}
which does not depend on the particular choices of $a$ and $b$.  One
can check that $|\cdot|_p$ is an ultrametric absolute value function
on ${\bf Q}$, and the corresponding ultrametric
\begin{equation}
\label{d_p(x, y) = |x - y|_p}
        d_p(x, y) = |x - y|_p,
\end{equation}
is known as the \emph{$p$-adic metric}\index{p-adic metric@$p$-adic
metric} on ${\bf Q}$.

        Let $k$ be an arbitrary field again, and let $|\cdot|$ be a
nonnegative real-valued function on $k$ that satisfies (\ref{|x| = 0
  if and only if x = 0}) and (\ref{|x y| = |x| |y|, 2}).  This implies
that (\ref{|1| = 1}) and (\ref{|-1| = 1}) still hold, for the same
reasons as before.  Let us say that $|\cdot|$ is a \emph{quasimetric
absolute value function}\index{quasimetric absolute value functions}
if there is a real number $C \ge 1$ such that
\begin{equation}
\label{|x + y| le C (|x| + |y|)}
        |x + y| \le C \, (|x| + |y|)
\end{equation}
for every $x, y \in k$.  This means that (\ref{d(x, y) = |x - y|, 2})
satisfies (\ref{d(x, z) le C (d(x, y) + d(y, z))}), and hence is a
quasimetric on $k$.  Equivalently, $|\cdot|$ is a quasimetric absolute
value function on $k$ if there is a real number $C' \ge 1$ such that
\begin{equation}
\label{|x + y| le C' max(|x|, |y|)}
        |x + y| \le C' \, \max(|x|, |y|)
\end{equation}
for every $x, y \in k$, in which case (\ref{d(x, y) = |x - y|, 2})
satisfies (\ref{d(x, z) le C' max(d(x, y), d(y, z))}) on $k$.  As
before, (\ref{|x + y| le C' max(|x|, |y|)}) implies (\ref{|x + y| le C
  (|x| + |y|)}) with $C = C'$, and (\ref{|x + y| le C (|x| + |y|)})
implies (\ref{|x + y| le C' max(|x|, |y|)}) with $C' = 2 \, C$.  Of
course, (\ref{|x + y| le C (|x| + |y|)}) reduces to (\ref{|x + y| le
  |x| + |y|, 2}) when $C = 1$, and (\ref{|x + y| le C' max(|x|, |y|)})
reduces to (\ref{|x + y| le max(|x|, |y|)}) when $C' = 1$.

        If $|\cdot|$ satisfies (\ref{|x + y| le C (|x| + |y|)}) and
$0 < a \le 1$, then
\begin{equation}
\label{|x + y|^a le C^a (|x| + |y|)^a le C^a (|x|^a + |y|^a)}
        |x + y|^a \le C^a \, (|x| + |y|)^a \le C^a \, (|x|^a + |y|^a)
\end{equation}
for every $x, y \in k$, by (\ref{(r + t)^a le r^a + t^a}).  In
particular, if $|x|$ is an absolute value function on $k$, then
$|x|^a$ is also an absolute value function on $k$ when $0 < a \le 1$.
If $|\cdot|$ satisfies (\ref{|x + y| le C (|x| + |y|)}) and $a \ge 1$,
then
\begin{equation}
\label{|x + y|^a le C^a (|x| + |y|)^a le 2^{a - 1} C^a (|x|^a + |y|^a)}
 |x + y|^a \le C^a \, (|x| + |y|)^a \le 2^{a - 1} \, C^a \, (|x|^a + |y|^a)
\end{equation}
for every $x, y \in k$, by (\ref{(r + t)^a = 2^a (r/2 + t/2)^a le
  ... = 2^{a - 1} (r^a + t^a)}).  Similarly, if $|\cdot|$ satisfies
(\ref{|x + y| le C' max(|x|, |y|)}) and $a > 0$, then
\begin{equation}
\label{|x + y|^a le (C')^a max(|x|^a, |y|^a)}
        |x + y|^a \le (C')^a \, \max(|x|^a, |y|^a)
\end{equation}
for every $x, y \in k$.  It follows that $|x|^a$ is an ultrametric
absolute value function on $k$ for every $a > 0$ when $|x|$ is an
ultrametric absolute value function on $k$, and that $|x|^a$ is a
quasimetric absolute value function on $k$ for every $a > 0$ when
$|x|$ is a quasimetric absolute value function on $k$.

        Let $|\cdot|$ be a nonnegative real-valued function on $k$
that satisfies (\ref{|x| = 0 if and only if x = 0}) and (\ref{|x y| =
  |x| |y|, 2}) again, and hence (\ref{|1| = 1}) and (\ref{|-1| = 1}).
If $|\cdot|$ also satisfies (\ref{|x + y| le C' max(|x|, |y|)}) for
some $C' \ge 1$, then
\begin{equation}
\label{|1 + z| le C' for every z in k with |z| le 1}
        |1 + z| \le C' \quad\hbox{for every } z \in k \hbox{ with } |z| \le 1,
\end{equation}
by (\ref{|1| = 1}).  Conversely, suppose that $|\cdot|$ satisfies
(\ref{|1 + z| le C' for every z in k with |z| le 1}) for some $C' \ge
1$, and let us check that (\ref{|x + y| le C' max(|x|, |y|)}) holds
for every $x, y \in k$.  We may as well restrict our attention to the
case where $|y| \le |x|$, since otherwise we can interchange the roles
of $x$ and $y$.  If $x = 0$, then (\ref{|x + y| le C' max(|x|, |y|)})
is trivial, and so we can suppose that $x \ne 0$ too.  Thus we can put
$z = y/x$, so that $|z| = |y|/|x| \le 1$, by (\ref{|x^{-1}| =
  |x|^{-1}}).  This permits us to use (\ref{|1 + z| le C' for every z
  in k with |z| le 1}) to get that
\begin{equation}
\label{|x + y| = |1 + z| |x| le C' max(|x|, |y|)}
        |x + y| = |1 + z| \, |x| \le C' \, \max(|x|, |y|),
\end{equation}
as desired.

        The definition of a quasimetric absolute value function in terms
of (\ref{|1 + z| le C' for every z in k with |z| le 1}) corresponds to
Definition 1.1 on p12 of \cite{cas}, but with different terminology.
The definition of an ordinary absolute value function corresponds to
Definition 2.1.1 on p21-2 of \cite{fg}.  The relationship between
ordinary absolute value functions and quasimetric absolute value
functions will be clarified in the next section, as in Lemma 1.2 on
p13-4 of \cite{cas}.  Ultrametric absolute value functions are also called
\emph{non-archimedian}\index{non-archimedian absolute value functions},
and we shall return to this in Section \ref{some more refinements}.

\section{Some refinements}
\label{some refinements}

        Let $|\cdot|$ be a quasimetric absolute value function on a
field $k$, that satisfies (\ref{|x + y| le C' max(|x|, |y|)}) for some
$C' \ge 1$.  Also let $n$ be a nonnegative integer, and let us check
that
\begin{equation}
\label{|sum_{j = 1}^{2^n} z_j| le (C')^n max {|z_j| : j = 1, ldots, 2^n}}
        \biggl|\sum_{j = 1}^{2^n} z_j\biggr|
                \le (C')^n \, \max \{|z_j| : j = 1, \ldots, 2^n\}
\end{equation}
for any finite sequence $z_1, \ldots, z_{2^n}$ of $2^n$ elements of
$k$.  This is trivial when $n = 0$, and this is the same as (\ref{|x +
  y| le C' max(|x|, |y|)}) when $n = 1$.  Suppose now that
(\ref{|sum_{j = 1}^{2^n} z_j| le (C')^n max {|z_j| : j = 1, ldots,
    2^n}}) holds for some $n \ge 1$, and let us verify that the
analogous statement holds for $n + 1$.  Let $z_1, \ldots, z_{2^{n +
    1}}$ be a finite sequence of $2^{n + 1}$ elements of $k$, so that
(\ref{|sum_{j = 1}^{2^n} z_j| le (C')^n max {|z_j| : j = 1, ldots,
    2^n}}) can be applied to the first $2^n$ terms $z_1, \ldots,
z_{2^n}$ of this sequence.  Similarly, we can apply the induction
hypothesis to the last $2^n$ terms $z_{2^n + 1}, \ldots, z_{2^{n +
    1}}$ of this sequence, to get that
\begin{equation}
\label{|sum_{j = 1}^{2^n} z_{2^n + j}| le ...}
        \biggl|\sum_{j = 1}^{2^n} z_{2^n + j}\biggr|
                \le (C')^n \, \max\{|z_{2^n + j}| : j = 1, \ldots, 2^n\}.
\end{equation}
It follows that
\begin{eqnarray}
\label{|sum_{j = 1}^{2^{n + 1}} z_j| le ...}
 \biggl|\sum_{j = 1}^{2^{n + 1}} z_j\biggr|
          & \le & C' \, \max\bigg\{\biggl|\sum_{j = 1}^{2^n} z_j\biggr|,
                           \biggl|\sum_{j = 1}^{2^n} z_{2^n + j}\biggr|\bigg\} \\
 & \le & (C')^{n + 1} \, \max\{|z_j| : j = 1, \ldots, 2^{n + 1}\}, \nonumber
\end{eqnarray}
using (\ref{|x + y| le C' max(|x|, |y|)}) in the first step.  Note
that (\ref{|sum_{j = 1}^{2^n} z_j| le (C')^n max {|z_j| : j = 1,
    ldots, 2^n}}) is equivalent to (\ref{d(x_0, x_{2^{n + 1}}) le ...,
  2}) in this setting, where $d(\cdot, \cdot)$ is the quasimetric
(\ref{d(x, y) = |x - y|, 2}) corresponding to $|\cdot|$.  More
precisely, (\ref{|sum_{j = 1}^{2^n} z_j| le (C')^n max {|z_j| : j = 1,
    ldots, 2^n}}) follows from (\ref{d(x_0, x_{2^n}) le (C')^n max
  {d(x_{j - 1}, x_j) : j = 1, ldots, 2^n}}) with $x_0 = 0$ and $x_l =
\sum_{j = 1}^l z_j$ when $l \ge 1$.  Conversely, (\ref{d(x_0, x_{2^n})
  le (C')^n max {d(x_{j - 1}, x_j) : j = 1, ldots, 2^n}}) follows from
(\ref{|sum_{j = 1}^{2^n} z_j| le (C')^n max {|z_j| : j = 1, ldots,
    2^n}}) applied to $z_j = x_j - x_{j - 1}$ in this situation.

        Suppose that $|\cdot|$ is a quasimetric absolute value function
on $k$ that satisfies (\ref{|x + y| le C' max(|x|, |y|)}) with $C' =
2$.  Let $N$ be a positive integer, and let $n$ be the smallest
nonnegative integer such that
\begin{equation}
\label{N le 2^n}
        N \le 2^n,
\end{equation}
so that $2^{n - 1} < N$.  Also let $z_1, \ldots, z_N$ be a finite
sequence of $N$ elements of $k$, and put $z_j = 0$ when $N < j \le
2^n$.  Applying (\ref{|sum_{j = 1}^{2^n} z_j| le (C')^n max {|z_j| : j
    = 1, ldots, 2^n}}) with $C' = 2$, we get that
\begin{eqnarray}
\label{|sum_{j = 1}^N z_j| le ... le 2 N max{|z_j| : j = 1, ldots, N}}
        \biggl|\sum_{j = 1}^N z_j\biggr|
                 & \le & 2^n \, \max\{|z_j| : j = 1, \ldots, N\} \\
                & \le & 2 \, N \, \max\{|z_j| : j = 1, \ldots, N\}. \nonumber
\end{eqnarray}
Let $z \in k$ be given, and let $N \cdot z$ be the sum of $N$ $z$'s in
$k$.  Observe that
\begin{equation}
\label{|N cdot z| le 2 N |z|}
        |N \cdot z| \le 2 \, N \, |z|,
\end{equation}
by applying (\ref{|sum_{j = 1}^N z_j| le ... le 2 N max{|z_j| : j = 1,
    ldots, N}}) with $z_j = z$ for each $j = 1, \ldots, N$.  In
particular,
\begin{equation}
\label{|N cdot 1| le 2 N}
        |N \cdot 1| \le 2 \, N,
\end{equation}
by (\ref{|1| = 1}).

        Let $x, y \in k$ be given, and let $r$ be a positive integer.
The binomial theorem implies that
\begin{equation}
\label{(x + y)^r = sum_{j = 0}^r {r choose j} cdot x^j y^{r - j}}
        (x + y)^r = \sum_{j = 0}^r {r \choose j} \cdot x^j \, y^{r - j},
\end{equation}
where ${r \choose j}$ is the usual binomial coefficient, which is a
positive integer, and $x^0$, $y^0$ are both interpreted as being equal
to $1$ in $k$.  If $|\cdot|$ is a quasimetric absolute value function
on $k$ that satisfies (\ref{|x + y| le C' max(|x|, |y|)}) with $C' = 2$,
then it follows that
\begin{equation}
\label{|(x + y)^r| le ...}
 |(x + y)^r| \le 2 \, (r + 1) \,
   \max \bigg\{\biggl|{r \choose j} \cdot x^j \, y^{r - j}\biggr| :
                                               j = 0, \ldots, r\bigg\},
\end{equation}
by (\ref{|sum_{j = 1}^N z_j| le ... le 2 N max{|z_j| : j = 1, ldots, N}})
with $N = r + 1$.  This implies that
\begin{equation}
\label{|(x + y)^r| le ..., 2}
 |(x + y)^r| \le 2 \, (r + 1) \,
   \max \bigg\{2 \, {r \choose j} \, |x|^j \, |y|^{r - j} :
                                                 j = 0, \ldots, r\bigg\},
\end{equation}
using (\ref{|N cdot z| le 2 N |z|}) with $N = {r \choose j}$.  Of course,
\begin{eqnarray}
\label{max {{r choose j} |x|^j |y|^{r - j} : j = 0, ldots, r} le ...}
 \max \bigg\{{r \choose j} \, |x|^j \, |y|^{r - j} : j = 0, \ldots, r\bigg\}
 & \le & \sum_{j = 0}^\infty {r \choose j} \, |x|^j \, |y|^{r - j} \\
  & = & (|x| + |y|)^r,                                     \nonumber
\end{eqnarray}
using the binomial theorem again in the second step.  Combining this
with (\ref{|(x + y)^r| le ..., 2}), we get that
\begin{equation}
\label{|x + y|^r = |(x + y)^r| le 4 (r + 1) (|x| + |y|)^r}
        |x + y|^r = |(x + y)^r| \le 4 \, (r + 1) \, (|x| + |y|)^r
\end{equation}
for each positive integer $r$ under these conditions.  Thus
\begin{equation}
\label{|x + y| le (4 (r + 1))^{1/r} (|x| + |y|)}
        |x + y| \le (4 \, (r + 1))^{1/r} \, (|x| + |y|)
\end{equation}
for every $x, y \in k$ and $r \ge 1$, which implies that
\begin{equation}
\label{|x + y| le |x| + |y|, 3}
        |x + y| \le |x| + |y|,
\end{equation}
by taking the limit as $r \to \infty$.

        This shows that a quasimetric absolute value function $|\cdot|$
on $k$ that satisfies (\ref{|x + y| le C' max(|x|, |y|)}) with $C' =
2$ is actually an absolute value function on $k$.  Suppose that
$|\cdot|$ is any quasimetric absolute value function on $k$, so that
$|\cdot|$ satisfies (\ref{|x + y| le C' max(|x|, |y|)}) for some $C'
\ge 1$.  As in the previous section, $|\cdot|^a$ is also a quasimetric
absolute value function on $k$ for each $a > 0$, which satisfies
(\ref{|x + y|^a le (C')^a max(|x|^a, |y|^a)}).  If $a$ is sufficiently
small, then
\begin{equation}
\label{(C')^a le 2}
        (C')^a \le 2,
\end{equation}
which implies that $|\cdot|^a$ is an absolute value function on $k$,
by the preceding argument.

\section{Some more refinements}
\label{some more refinements}

        Let $|\cdot|$ be an absolute value function on a field $k$,
and suppose that there is a real number $A \ge 1$ such that
\begin{equation}
\label{|N cdot 1| le A}
        |N \cdot 1| \le A
\end{equation}
for every positive integer $N$.  Observe that
\begin{equation}
\label{N cdot z = (N cdot 1) z}
        N \cdot z = (N \cdot 1) \, z
\end{equation}
for every $z \in k$ and positive integer $N$, so that
\begin{equation}
\label{|N cdot z| = |N cdot 1| |z| le A |z|}
        |N \cdot z| = |N \cdot 1| \, |z| \le A \, |z|
\end{equation}
under these conditions.  Let $x, y \in k$ and a positive integer $r$
be given, so that $(x + y)^r$ can be expressed as in (\ref{(x + y)^r =
  sum_{j = 0}^r {r choose j} cdot x^j y^{r - j}}), using the binomial
theorem.  This implies that
\begin{eqnarray}
\label{|(x + y)^r| le ... le A (r + 1) max(|x|, |y|)^r}
 |(x + y)^r| & \le &
      \sum_{j = 0}^r \biggl|{r \choose j} \cdot x^j \, y^{r - j}\biggr| \\
             & \le & A \, \sum_{j = 0}^r |x|^j \, |y|^{r - j}
             \le A \, (r + 1) \, \max(|x|, |y|)^r, \nonumber
\end{eqnarray}
using (\ref{|N cdot z| = |N cdot 1| |z| le A |z|}) in the second step.
Equivalently,
\begin{equation}
\label{|x + y|^r = |(x + y)^r| le A (r + 1) max(|x|, |y|)^r}
        |x + y|^r = |(x + y)^r| \le A \, (r + 1) \, \max(|x|, |y|)^r
\end{equation}
for each $r \ge 1$, and hence
\begin{equation}
\label{|x + y| le (A (r + 1))^{1/r} max(|x|, |y|)}
        |x + y| \le (A \, (r + 1))^{1/r} \, \max(|x|, |y|).
\end{equation}
Taking the limit as $r \to \infty$, we get that
\begin{equation}
\label{|x + y| le max(|x|, |y|), 2}
        |x + y| \le \max(|x|, |y|)
\end{equation}
for every $x, y \in k$, so that $|\cdot|$ is an ultrametric absolute
value function on $k$.  If $|\cdot|$ is a quasimetric absolute value
function on $k$ that satisfies (\ref{|N cdot 1| le A}) for some $A \ge
1$ and every positive integer $N$, then
\begin{equation}
\label{|N cdot 1|^a le A^a}
        |N \cdot 1|^a \le A^a
\end{equation}
for every positive real number $a$ and positive integer $N$, and we
know from the previous section that $|\cdot|^a$ is an absolute value
function on $k$ when $a$ is sufficiently small.  It follows from the
preceding argument that $|\cdot|^a$ is an ultrametric absolute value
function on $k$ when $a$ is sufficiently small, which implies that
$|\cdot|$ is an ultrametric absolute value function on $k$, as in
Section \ref{definitions, examples}.  Alternatively, one could extend
the preceding argument directly to quasimetric absolute value
functions, using (\ref{|sum_{j = 1}^{2^n} z_j| le (C')^n max {|z_j| :
    j = 1, ldots, 2^n}}).

        Of course, if $|\cdot|$ is an ultrametric absolute value function
on $k$, then
\begin{equation}
\label{|N cdot 1| le 1}
        |N \cdot 1| \le 1
\end{equation}
for every positive integer $N$, by (\ref{|1| = 1}).  Note that
\begin{equation}
\label{(N_1 N_2) cdot 1 = N_1 cdot (N_2 cdot 1) = (N_1 cdot 1) (N_2 cdot 1)}
 (N_1 \, N_2) \cdot 1 = N_1 \cdot (N_2 \cdot 1) = (N_1 \cdot 1) \, (N_2 \cdot 1)
\end{equation}
for any two positive integers $N_1$, $N_2$, and hence that
\begin{equation}
\label{|N^j cdot 1| = |(N cdot 1)^j| = |N cdot 1|^j}
        |N^j \cdot 1| = |(N \cdot 1)^j| = |N \cdot 1|^j
\end{equation}
for all positive integers $N$ and $j$.  If $|N \cdot 1| > 1$ for some
positive integer $N$, then it follows that $|N^j \cdot 1|$ is
unbounded, so that (\ref{|N cdot 1| le 1}) can also be derived
directly from (\ref{|N cdot 1| le A}).  A quasimetric absolute value
function $|\cdot|$ on $k$ is said to be
\emph{archimedian}\index{archimedian absolute value functions} if $|N
\cdot 1|$ has no finite upper bound for all positive integers $N$, 
and otherwise $|\cdot|$ is said to be
\emph{non-archimedian}.\index{non-archimedian absolute value functions}
Thus $|\cdot|$ is non-archimedean if and only if it is an ultrametric
absolute value function.

        Suppose that $|\cdot|$ is an ultrametric absolute value function
on $k$.  If $x, y \in k$ satisfy $|y| \le |x|$, then
\begin{equation}
\label{|x + y| le max(|x|, |y|) = |x|}
        |x + y| \le \max(|x|, |y|) = |x|.
\end{equation}
We also have that
\begin{equation}
\label{|x| = |(x + y) - y| le max(|x + y|, |y|)}
        |x| = |(x + y) - y| \le \max(|x + y|, |y|),
\end{equation}
by (\ref{|-1| = 1}), which implies that $|x| \le |x + y|$ when $|y| <
|x|$.  It follows that
\begin{equation}
\label{|x + y| = |x|}
        |x + y| = |x|
\end{equation}
when $|y| < |x|$, which also corresponds to (\ref{d(x, y) = d(x, z)})
in Section \ref{metrics, ultrametrics}.

        If $|\cdot|$ is a nontrivial quasimetric absolute value function
on ${\bf Q}$, then a famous theorem of Ostrowki\index{Ostrowski's theorems}
states that $|\cdot|$ is either a positive power of the standard
absolute value function on ${\bf Q}$, or a positive power of the
$p$-adic absolute value function on ${\bf Q}$ for some prime number
$p$.  More precisely, if $|\cdot|$ is archimedian, then $|N \cdot 1| >
1$ for some positive integer $N$, and one can show that $|\cdot|$ is a
positive power of the standard absolute value function on ${\bf Q}$.
Otherwise, if $|\cdot|$ is non-achimedian, then $|N \cdot 1| \le 1$
for every positive integer $N$, and $|N \cdot 1| < 1$ for some
positive integer $N$, because $|\cdot|$ is nontrivial.  If $p$ is the
smallest positive integer such that $|p \cdot 1| < 1$, then one can
show that $p$ is a prime number, and that $|\cdot|$ is a positive
power of the $p$-adic absolute value on ${\bf Q}$.  See Theorem 2.1 on
p16 of \cite{cas} or Theorem 3.1.3 on p44 of \cite{fg} for more
details.

\section{Some topological properties}
\label{some topological properties}

        Let $k$ be a field, and let $|\cdot|$ be a quasimetric absolute
value function on $k$, with the associated quasimetric $d(\cdot,
\cdot)$ on $k$, as in (\ref{d(x, y) = |x - y|, 2}).  By construction,
$d(\cdot, \cdot)$ is invariant under translations on $k$, in the sense
that
\begin{equation}
\label{d(x + z, y + z) = d(x, y)}
        d(x + z, y + z) = d(x, y)
\end{equation}
for every $x, y, z \in k$.  If $a$ is a positive real number, then
$|\cdot|^a$ is also a quasimetric absolute value function on $k$, as
in Section \ref{definitions, examples}, for which the corresponding
quasimetric on $k$ is equal to
\begin{equation}
\label{d(x, y)^a = |x - y|^a}
        d(x, y)^a = |x - y|^a.
\end{equation}
Each of these quasimetrics determines the same topology on $k$, and in
fact the same uniform structure.

        Remember that $|\cdot|^a$ is an absolute value function on $k$
when $a$ is sufficiently small, as in Section \ref{some refinements},
in which case (\ref{d(x, y)^a = |x - y|^a}) is a metric on $k$.  This
implies that open balls in $k$ with respect to (\ref{d(x, y)^a = |x -
  y|^a}) are open sets with respect to the corresponding topology when
$a$ is sufficiently small, and that closed balls in $k$ with respect
to (\ref{d(x, y)^a = |x - y|^a}) are closed sets.  It is easy to see
that an open or closed ball in $k$ with respect to $d(\cdot, \cdot)$
centered at a point $x \in k$ and with radius $r$ is the same as the
open or closed ball in $k$ with respect to (\ref{d(x, y)^a = |x -
  y|^a}) centered at $x$ with radius $r^a$, for each $a > 0$.  It
follows that open balls in $k$ with respect to $d(\cdot, \cdot)$ are
open sets, and that closed balls in $k$ with respect to $d(\cdot,
\cdot)$ are closed sets, even when $|\cdot|$ is a quasimetric absolute
value function on $k$.

        If $a > 0$ is sufficiently small so that $|\cdot|^a$ is an
absolute value function on $k$, then $|x|^a$ is continuous with
respect to the corresponding metric (\ref{d(x, y)^a = |x - y|^a}) on
$k$.  More precisely, this means that $|x|^a$ is continuous as a
mapping from $k$ into ${\bf R}$, with respect to the standard topology
on ${\bf R}$.  In fact, $|x|^a$ is Lipschitz of order $1$ with
constant $C = 1$ with respect to the metric (\ref{d(x, y)^a = |x -
  y|^a}) on $k$ and the standard metric on ${\bf R}$, as in Section
\ref{lipschitz mappings}.  Equivalently, this means that $|x|^a$ is
Lipschitz of order $a$ with constant $C = 1$ with respect to the
associated quasimetric (\ref{d(x, y) = |x - y|, 2}) on $k$ and the
standard metric on ${\bf R}$.  In particular, $|x|^a$ is continuous
with respect to the topology determined by the associated quasimetric
(\ref{d(x, y) = |x - y|, 2}) on $k$, so that $|x|$ is continuous with
respect to this topology as well.

        Of course, addition and multiplication on $k$ correspond to
mappings from the Cartesian product $k \times k$ of $k$ with itself
into $k$.  Using the topology on $k$ determined by the quasimetric
$d(\cdot, \cdot)$ associated to the quasimetric $|\cdot|$, one can
define the corresponding product topology on $k \times k$.  With
respect to this topology, addition and multiplication on $k$
correspond to continuous mappings from $k \times k$ into $k$.  This
can be verified in essentially the same way as for real or complex
numbers.  If $x, y \in k$ and $x, y \ne 0$, then
\begin{equation}
\label{x^{-1} - y^{-1} = (y - x) x^{-1} y^{-1}}
        x^{-1} - y^{-1} = (y - x) \, x^{-1} \, y^{-1},
\end{equation}
and hence
\begin{equation}
\label{|x^{-1} - y^{-1}| = ... = |y - x| |x|^{-1} |y|^{-1}}
        |x^{-1} - y^{-1}| = |y - x| \, |x^{-1}| \, |y^{-1}|
                         = |y - x| \, |x|^{-1} \, |y|^{-1},
\end{equation}
using (\ref{|x^{-1}| = |x|^{-1}}) in the last step.  If $x, y \in k$,
$x \ne 0$, and $y$ is sufficiently close to $x$, then one can also
check that there is a uniform positive lower bound for $|y|$ in terms
of $|x|$, using the quasimetric version of the triangle inequality.
This permits one to show that $x \mapsto x^{-1}$ is continuous as a
mapping from $k \setminus \{0\}$ into itself, in essentially the same
way as for real or complex numbers.  It follows that $k$ is a
topological field with respect to the topology determined by $d(\cdot,
\cdot)$.

\section{Completions, 2}
\label{completions, 2}

        Let $k$ be a field, and let $|\cdot|$ be an absolute value function
on $k$, with the corresponding metric $d(\cdot, \cdot)$ as in
(\ref{d(x, y) = |x - y|, 2}).  It is convenient to restrict our
attention here to absolute value functions instead of quasimetric
absolute value functions, in order to follow the usual discussion for
metric spaces, as in Section \ref{completions}.  One could also start
with a quasimetric absolute value function on $k$, and then reduce to
the case of ordinary absolute value functions, as in Section \ref{some
  refinements}.  In the context of arbitrary quasimetric spaces, one
has to be a bit careful about continuity properties of the
quasimetric, or work with a metric that determines the same uniform
structure.  In this situation, quasimetric absolute value functions
and their associated quasimetrics already have nice continuity
properties, because they are related to ordinary absolute value
functions and their associated metrics as in Section \ref{some
  refinements}, and one may as well work directly with the latter.

        Before considering the completion of $k$, let us look at some
properties of Cauchy sequences of elements of $k$.  If $\{x_j\}_{j = 1}^\infty$
is a Cauchy sequence of elements of $k$, then it is easy to see that
$\{|x_j|\}_{j = 1}^\infty$ is a Cauchy sequence in ${\bf R}$ with respect
to the standard metric, which thus converges to a nonnegative real number.
This uses the triangle inequality, and may be considered as a special
case of an analogous statement about distances between Cauchy sequences
in Section \ref{completions}, because $|x_j|$ is the same as the distance
between $x_j$ and $0$.  If $|\cdot|$ is an ultrametric absolute value function
on $k$, and if $\{x_j\}_{j = 1}^\infty$ does not converge to $0$, then
$\{|x_j|\}_{j = 1}^\infty$ is eventually constant.  This can be derived from
(\ref{|x + y| = |x|}), and may also be considered as a special case of
an analogous statement for distances in ultrametric spaces, as in
Section \ref{completions}.

        Suppose that $\{x_j\}_{j = 1}^\infty$, $\{x'_j\}_{j = 1}^\infty$
are equivalent Cauchy sequences of elements of $M$, in the sense
discussed in Section \ref{completions}.  In this case, one can check that
\begin{equation}
\label{lim_{j to infty} (|x_j| - |x'_j|) = 0}
        \lim_{j \to \infty} (|x_j| - |x'_j|) = 0
\end{equation}
in ${\bf R}$, so that $\{|x_j|\}_{j = 1}^\infty$ and $\{|x'_j|\}_{j =
  1}^\infty$ have the same limit in ${\bf R}$.  This permits one to
extend the absolute value function on $k$ to a nonnegative real-valued
function on the completion of $k$.  If $|\cdot|$ is an ultrametric
absolute value function on $k$, and $\{x_j\}_{j = 1}^\infty$ or
$\{x'_j\}_{j = 1}^\infty$ does not converge to $0$ in $k$, then $|x_j|
= |x'_j|$ for all sufficiently large $j$.  As before, these statements
may be considered as special cases of analogous statements for
distances, as in Section \ref{completions}.

        If $\{x_j\}_{j = 1}^\infty$ and $\{y_j\}_{j = 1}^\infty$ are
Cauchy sequences of elements of $k$, then one can check that $\{x_j +
y_j\}_{j = 1}^\infty$ and $\{x_j \, y_j\}_{j = 1}^\infty$ are Cauchy
sequences of elements of $k$ as well.  In the case of products, this
also uses the fact that Cauchy sequences in $k$ are bounded.  If
$\{x'_j\}_{j = 1}^\infty$ and $\{y'_j\}_{j = 1}^\infty$ are Cauchy
sequences of elements of $k$ that are equivalent to $\{x_j\}_{j =
  1}^\infty$ and $\{y_j\}_{j = 1}^\infty$, respectively, then $\{x'_j
+ y'_j\}_{j = 1}^\infty$ and $\{x'_j \, y'_j\}_{j = 1}^\infty$ are
equivalent as Cauchy sequences in $k$ to $\{x_j + y_j\}_{j =
  1}^\infty$ and $\{x_j \, y_j\}_{j = 1}^\infty$, respectively.  This
permits one to extend addition and multiplication on $k$ to the
completion of $k$, so that the completion of $k$ becomes a commutative
ring.  The extension of the absolute value function to the completion
of $k$ satisfies the same type of properties as on $k$, and the
extension of the associated metric on $k$ to the completion of $k$ as
in Section \ref{completions} corresponds to the extension of the
absolute value function to the completion of $k$ as in (\ref{d(x, y) =
  |x - y|, 2}).

        If $\{x_j\}_{j = 1}^\infty$ is a sequence of elements of $k$
that does not converge to $0$, then there is an $r > 0$ such that
\begin{equation}
\label{|x_j| ge 2 r for infinitely many j, 2}
        |x_j| \ge 2 \, r \quad\hbox{for infinitely many } j.
\end{equation}
If $\{x_j\}_{j = 1}^\infty$ is also a Cauchy sequence of elements of
$k$, then it follows that
\begin{equation}
\label{|x_j| ge r for all but finitely many j, 2}
        |x_j| \ge r \quad\hbox{for all but finitely many } j,
\end{equation}
and in particular $x_j \ne 0$ for all but finitely many $j$.  Of
course, it is not necessary to switch between $2 \, r$ in (\ref{|x_j|
  ge 2 r for infinitely many j, 2}) and $r$ in (\ref{|x_j| ge r for
  all but finitely many j, 2}) when $|\cdot|$ is an ultrametric
absolute value function.  If $\{x_j\}_{j = 1}^\infty$ is a Cauchy
sequence of nonzero elements of $k$ that does not converge to $0$,
then one can check that $\{1/x_j\}_{j = 1}^\infty$ is a Cauchy
sequence in $k$ too, using (\ref{|x^{-1} - y^{-1}| = ... = |y - x|
  |x|^{-1} |y|^{-1}}).  Similarly, if $\{x'_j\}_{j = 1}^\infty$ is
another Cauchy sequence of nonzero elements of $k$ that does not
converge to $0$ and which is equivalent to $\{x_j\}_{j = 1}^\infty$,
then one can verify that $\{1/x'_j\}_{j = 1}^\infty$ is equivalent to
$\{1/x_j\}_{j = 1}^\infty$, using (\ref{|x^{-1} - y^{-1}| = ... = |y -
  x| |x|^{-1} |y|^{-1}}) again.  This permits one to extend the
mapping $x \mapsto 1/x$ to the nonzero elements of the completion of
$k$, so that the completion of $k$ becomes a field.  Note that the
extension of the absolute value function to the completion of $k$
still satisfies (\ref{|x^{-1}| = |x|^{-1}}).

        The extension of the absolute value function to the completion
of $k$ is an absolute value function on the completion of $k$.  If
$|\cdot|$ is an ultrametric absolute value function on $k$, then its
extension to the completion of $k$ is an ultrametric absolute value
function too.  As in Section \ref{completions}, there is a natural
embedding of $k$ into its completion, which associates to each $x \in
k$ the equivalence class of Cauchy sequences that contains the
constant sequence $\{x_j\}_{j = 1}^\infty$ with $x_j = x$ for each
$j$.  By construction, this embedding is a field isomorphism from $k$
onto its image in the completion of $k$, which preserves absolute
values and hence distance.  It is customary to identify $k$ with
its image in the completion under this embedding, which is a dense
subset of the completion.

        Alternatively, suppose that we start with a completion of $k$
as a metric space, which is to say an isometric embedding of $k$ onto
a dense subset of a complete metric space.  Note that addition on $k$
may be considered as a uniformly continuous mapping from $k \times k$
into $k$, with respect to a suitable product metric on $k \times k$
corresponding to the metric on $k$ associated to the absolute value
function.  One can then extend addition on $k$ to the completion of
$k$ as in Section \ref{continuous extensions}.  Similarly,
multiplication on $k$ corresponds to a mapping from $k \times k$ into
$k$ that is uniformly continuous on bounded subsets of $k \times k$,
which is sufficient for this type of extension argument.  If $r$ is
any positive real number, then $x \mapsto 1/x$ is uniformly continuous
as a mapping from
\begin{equation}
\label{{x in k : |x| ge r}}
        \{x \in k : |x| \ge r\}
\end{equation}
into $k$, which is again sufficient for this type of extension
argument.  This gives another way to look at the extension of the
field operations to the completion of $k$.  The absolute value
function can also be considered as a uniformly continuous mapping from
$k$ into ${\bf R}$, but its extension to the completion of $k$ is
already implicitly included in the metric, since it is the same as the
distance to $0$ in $k$.

        Suppose now that $k_1$ and $k_2$ are fields equipped with
absolute value functions $|\cdot|_1$ and $|\cdot|_2$, respectively,
and let $d_1(\cdot, \cdot)$ and $d_2(\cdot, \cdot)$ be the
corresponding metrics on $k_1$ and $k_2$, as in (\ref{d(x, y) = |x -
  y|, 2}).  Also let $\phi_1$ and $\phi_2$ be field isomorphisms from
$k$ onto subsets of $k_1$ and $k_2$ that are dense with respect to the
corresponding metrics.  If $\phi_1$ and $\phi_2$ preserve the
appropriate absolute value functions, in the sense that
\begin{equation}
\label{|phi_1(x)|_1 = |phi_2(x)|_2 = |x|}
        |\phi_1(x)|_1 = |\phi_2(x)|_2 = |x|
\end{equation}
for every $x \in k$, then $\phi_1$ and $\phi_2$ are isometric
embeddings of $k$ into $k_1$ and $k_2$, respectively.  If $k_1$ and
$k_2$ are also complete with respect to their corresponding metrics,
then it follows that there is an isometry $f$ from $k_1$ onto $k_2$
such that
\begin{equation}
\label{f circ phi_1 = phi_2}
        f \circ \phi_1 = \phi_2,
\end{equation}
by the remarks at the end of Section \ref{continuous extensions}.
Under these conditions, it is easy to see that $f$ is a field
isomomorphism, and that
\begin{equation}
\label{|f(z)|_2 = |z|_1}
        |f(z)|_2 = |z|_1
\end{equation}
for every $z \in k_1$.

        Let $k$ be a field, and let $|\cdot|$ be an ultrametric
absolute value function on $k$.  As before, if $\{x_j\}_{j = 1}^\infty$
is a Cauchy sequence of elements of $k$ that does not converge to $0$,
then $\{|x_j|\}_{j = 1}^\infty$ is eventually constant.  In particular,
this implies that the natural extension of $|\cdot|$ to the completion
of $k$ takes values in the same set of nonnegative real numbers
as $|\cdot|$ does on $k$.

        Now let $k = {\bf Q}$, let $p$ be a prime number, and let
$|\cdot|_p$ be the $p$-adic absolute value function on ${\bf Q}$,
as in Section \ref{definitions, examples}.  The completion of ${\bf
  Q}$ with respect to $|\cdot|_p$ is known as the field ${\bf
  Q}_p$\index{Q_p@${\bf Q}_p$} of \emph{$p$-adic
  numbers}.\index{p-adic numbers@$p$-adic numbers} The natural
extension of $|\cdot|_p$ to ${\bf Q}_p$ is known as the $p$-adic
absolute value function on ${\bf Q}_p$,\index{p-adic absolute
value@$p$-adic absolute value} and is also denoted $|\cdot|_p$.
Similarly, the natural extension of the $p$-adic metric $d_p(\cdot,
\cdot)$ on ${\bf Q}$ to ${\bf Q}_p$ is known as the $p$-adic metric on
${\bf Q}_p$,\index{p-adic metric@$p$-adic metric} and is denoted
$d_p(\cdot, \cdot)$ as well.  The remark in the preceding paragraph
implies that $|x|_p$ is either equal to $0$ or to an integer power of
$p$ for every $x \in {\bf Q}_p$, and $d_p(\cdot, \cdot)$ has the
same property on ${\bf Q}_p$.

\section{Infinite series}
\label{infinite series}

        Let $k$ be a field, let $|\cdot|$ be a quasimetric absolute
value function on $k$, and let $d(\cdot, \cdot)$ be the associated
quasimetric on $k$, as in (\ref{d(x, y) = |x - y|, 2}).  Also let
$\sum_{j = 1}^\infty a_j$ be an infinite series with $a_j \in k$ for
each $j$.  As usual, this series $\sum_{j = 1}^\infty a_j$ is said to
\emph{converge}\index{convergence of infinite series} in $k$ if the
corresponding sequence of partial sums
\begin{equation}
\label{s_n = sum_{j = 1}^n a_j}
        s_n = \sum_{j = 1}^n a_j
\end{equation}
converges to an element of $k$ with respect to $d(\cdot, \cdot)$, in
which case the value of the sum $\sum_{j = 1}^\infty a_j$ is equal to
the limit of the sequence $\{s_n\}_{n = 1}^\infty$ of partial sums.
If $\sum_{j = 1}^\infty a_j$ converges in $k$ and $c \in k$, then it
is easy to see that $\sum_{j = 1}^\infty c \, a_j$ converges in $k$
too, with
\begin{equation}
\label{sum_{j = 1}^infty c a_j = c sum_{j = 1}^infty a_j}
        \sum_{j = 1}^\infty c \, a_j = c \, \sum_{j = 1}^\infty a_j.
\end{equation}
Similarly, if $\sum_{j = 1}^\infty a_j$ and $\sum_{j = 1}^\infty b_j$
are convergent series with terms in $k$, then $\sum_{j = 1}^\infty
(a_j + b_j)$ converges in $k$, with
\begin{equation}
\label{sum_{j = 1}^infty (a_j + b_j) = ...}
 \sum_{j = 1}^\infty (a_j + b_j) = \sum_{j = 1}^\infty a_j + \sum_{j = 1}^\infty b_j.
\end{equation}

        Note that $\{s_n\}_{n = 1}^\infty$ is a Cauchy sequence in $k$
with respect to $d(\cdot, \cdot)$ if and only if for each $\epsilon > 0$
there is an $L \ge 1$ such that
\begin{equation}
\label{|s_n - s_l| = |sum_{j = l + 1}^n a_j| < epsilon}
        |s_n - s_l| = \biggl|\sum_{j = l + 1}^n a_j\biggr| < \epsilon
\end{equation}
for every $n > l \ge L$.  In particular, this implies that $\{a_j\}_{j
  = 1}^\infty$ converges to $0$ in $k$, by taking $n = l + 1$.  Of
course, the convergence of $\sum_{j = 1}^\infty a_j$ implies that
$\{s_n\}_{n = 1}^\infty$ is a Cauchy sequence in $k$, and the converse
holds when $k$ is complete with respect to $d(\cdot, \cdot)$.

        If $|\cdot|$ is an absolute value function on $k$, then
we say that $\sum_{j = 1}^\infty a_j$ \emph{converges absolutely}\index{absolute
convergence} when
\begin{equation}
\label{sum_{j = 1}^infty |a_j|}
        \sum_{j = 1}^\infty |a_j|
\end{equation}
converges as an infinite series of nonnegative real numbers.  This
implies that $\{s_n\}_{n = 1}^\infty$ is a Cauchy sequence in $k$,
because
\begin{equation}
\label{|sum_{j = l + 1}^n a_j| le sum_{j = l + 1}^n |a_j|}
        \biggl|\sum_{j = l + 1}^n a_j\biggr| \le \sum_{j = l + 1}^n |a_j|
\end{equation}
for every $n > l \ge 1$.  It follows that $\sum_{j = 1}^\infty a_j$
converges in $k$ when $k$ is complete, in which case we also have that
\begin{equation}
\label{|sum_{j = 1}^infty a_j| le sum_{j = 1}^infty |a_j|}
        \biggl|\sum_{j = 1}^\infty a_j\biggr| \le \sum_{j = 1}^\infty |a_j|.
\end{equation}

        If $|\cdot|$ is an ultrametric absolute value function on $k$, then
\begin{equation}
\label{|sum_{j = l + 1}^n a_j| le max_{l + 1 le j le n} |a_j|}
        \biggl|\sum_{j = l + 1}^n a_j\biggr| \le \max_{l + 1 \le j \le n} |a_j|
\end{equation}
for every $n > l \ge 1$.  This implies that $\{s_n\}_{n = 1}^\infty$
is a Cauchy sequence in $k$ when $\{a_j\}_{j = 1}^\infty$ converges to
$0$.  If $k$ is complete, then $\sum_{j = 1}^\infty a_j$ converges in
$k$, and satisfies
\begin{equation}
\label{|sum_{j = 1}^infty a_j| le max_{j ge 1} |a_j|}
        \biggl|\sum_{j = 1}^\infty a_j\biggr| \le \max_{j \ge 1} |a_j|.
\end{equation}
More precisely, it is easy to see that the maximum on the right side
of (\ref{|sum_{j = 1}^infty a_j| le max_{j ge 1} |a_j|}) exists under
these conditions, because $\{|a_j|\}_{j = 1}^\infty$ is a sequence of
nonnegative real numbers that converges to $0$.

        Of course, we can also consider series that start with $j = 0$.
If $x \in k$ and $n$ is a nonnegative integer, then it is well known
that
\begin{equation}
\label{(1 - x) sum_{j = 0}^n x^j = 1 - x^{n + 1}}
        (1 - x) \, \sum_{j = 0}^n x^j = 1 - x^{n + 1},
\end{equation}
where $x^0$ is interpreted as being equal to $1$, as usual.  Thus
\begin{equation}
\label{sum_{j = 0}^n x^j = frac{1 - x^{n + 1}}{1 - x}}
        \sum_{j = 0}^n x^j = \frac{1 - x^{n + 1}}{1 - x}
\end{equation}
for every $n \ge 0$ when $x \ne 1$, so that
\begin{equation}
\label{lim_{n to infty} sum_{j = 0}^n x^j = frac{1}{1 - x}}
        \lim_{n \to \infty} \sum_{j = 0}^n x^j = \frac{1}{1 - x}
\end{equation}
when $|x| < 1$, because $|x^{n + 1}| = |x|^{n + 1} \to 0$ as $n \to
\infty$.

        Let $\sum_{j = 0}^\infty a_j$ and $\sum_{l = 0}^\infty b_l$
be infinite series with terms in $k$.  The \emph{Cauchy
product}\index{Cauchy products} of these two series is the infinite
series $\sum_{n = 0}^\infty c_n$, where
\begin{equation}
\label{c_n = sum_{j = 0}^n a_j b_{n - j}}
        c_n = \sum_{j = 0}^n a_j \, b_{n - j}
\end{equation}
for each nonnegative integer $n$.  It is easy to see that
\begin{equation}
\label{sum_{n = 0}^infty c_n = (sum_{j = 0}^infty a_j) (sum_{l = 0}^infty b_l)}
        \sum_{n = 0}^\infty c_n = \Big(\sum_{j = 0}^\infty a_j\Big) \,
                                 \Big(\sum_{l = 0}^\infty b_l\Big)
\end{equation}
formally.  In particular, if $a_j = 0$ for all but finitely many $j$,
and if $b_l = 0$ for all but finitely many $l$, then $c_n = 0$ for all
but finitely many $n$, and (\ref{sum_{n = 0}^infty c_n = (sum_{j =
    0}^infty a_j) (sum_{l = 0}^infty b_l)}) holds.

        If $|\cdot|$ is an absolute value function on $k$, then
\begin{equation}
\label{|c_n| le sum_{j = 0}^n |a_j| |b_{n - j}|}
        |c_n| \le \sum_{j = 0}^n |a_j| \, |b_{n - j}|
\end{equation}
for each $n \ge 0$, where the right side corresponds to the Cauchy
product of $\sum_{j = 0}^\infty |a_j|$ and $\sum_{l = 0}^\infty
|b_l|$.  If $\sum_{j = 0}^\infty a_j$ and $\sum_{l = 0}^\infty b_l$
converge absolutely, then one can check that $\sum_{n = 0}^\infty c_n$
converges absolutely, with
\begin{equation}
\label{sum_{n = 0}^infty |c_n| le ...}
        \sum_{n = 0}^\infty |c_n| \le \Big(\sum_{j = 0}^\infty |a_j|\Big) \,
                                     \Big(\sum_{l = 0}^\infty |b_l|\Big).
\end{equation}
If $k$ is complete, then it follows that these series converge in $k$,
and one can check that (\ref{sum_{n = 0}^infty c_n = (sum_{j =
    0}^infty a_j) (sum_{l = 0}^infty b_l)}) holds.  The main point is
to approximate the infinite series by finite sums, using (\ref{sum_{n
    = 0}^infty |c_n| le ...}) to estimate the errors.

        Similarly, if $|\cdot|$ is an ultrametric absolute value function
on $k$, then
\begin{equation}
\label{|c_n| le max_{0 le j le n} (|a_j| |b_{n - j}|)}
        |c_n| \le \max_{0 \le j \le n} (|a_j| \, |b_{n - j}|)
\end{equation}
for each $n \ge 0$.  If $\{a_j\}_{j = 0}^\infty$ and $\{b_l\}_{l =
  0}^\infty$ both converge to $0$, then it is easy to see that
$\{c_n\}_{n = 0}^\infty$ onverges to $0$ as well.  This implies that
the series $\sum_{j = 0}^\infty a_j$, $\sum_{l = 0}^\infty b_l$, and
$\sum_{n = 0}^\infty c_n$ converge in $k$ when $k$ is complete, in
which case one can check that (\ref{sum_{n = 0}^infty c_n = (sum_{j =
    0}^infty a_j) (sum_{l = 0}^infty b_l)}) holds again.  As before,
the main point is to approximate the infinite series by finite sums,
but now using the ultrametric version of the triangle inequality to
estimate the errors.

\section{Topological equivalence}
\label{topological equivalence}

        Let $k$ be a field, and let $|\cdot|_1$ be a quasimetric
absolute value function on $k$.  As before, this leads to a
quasimetric $d_1(\cdot, \cdot)$ on $k$ as in (\ref{d(x, y) = |x - y|,
  2}), and thus a topology on $k$.  If $x \in k$ and $|x| < 1$, then
\begin{equation}
\label{|x^j|_1 = |x|_1^j to 0 as j to infty}
        |x^j|_1 = |x|_1^j \to 0 \quad\hbox{as } j \to \infty
\end{equation}
as a sequence of nonnegative real numbers, which implies that
$\{x^j\}_{j = 1}^\infty$ converges to $0$ in $k$ with respect to the
topology corresponding to $|\cdot|_1$.  Conversely, if $\{x^j\}_{j =
  1}^\infty$ converges to $0$ in $k$ with respect to the topology
corresponding to $|\cdot|_1$, then
\begin{equation}
\label{|x|_1^j = |x^j|_1 < 1}
        |x|_1^j = |x^j|_1 < 1
\end{equation}
for all but finitely many $j$, and hence $|x|_1 < 1$.  This shows that
the open unit ball in $k$ with respect to $|\cdot|_1$ is uniquely
determined by the topology on $k$ that corresponds to $|\cdot|_1$.
Note that $x \in k$ satisfies $|x|_1 > 1$ if and only if $x \ne 0$ and
$|x^{-1}|_1 = |x|_1^{-1} < 1$, so that the set of $x \in k$ with $|x|_1 > 1$
is also uniquely determined by the topology on $k$ that corresponds to
$|\cdot|_1$.  It follows that the set of $x \in k$ with $|x|_1 = 1$ is
uniquely determined by the topology on $k$ that corresponds to $|\cdot|_1$
as well.

        Now let $|\cdot|_1$ and $|\cdot|_2$ be quasimetric absolute value
functions on $k$, and let $d_1(\cdot, \cdot)$ and $d_2(\cdot, \cdot)$ be
the associated quasimetrics on $k$, as in (\ref{d(x, y) = |x - y|, 2}).
If
\begin{equation}
\label{|x|_2 = |x|_1^a}
        |x|_2 = |x|_1^a
\end{equation}
for some $a > 0$ and every $x \in k$, then
\begin{equation}
\label{d_2(x, y) = d_1(x, y)^a}
        d_2(x, y) = d_1(x, y)^a
\end{equation}
for every $x, y \in k$, and the corresponding topologies on $k$ are
the same, as in Section \ref{some topological properties}.
Conversely, suppose that the topologies on $k$ corresponding to
$|\cdot|_1$ and $|\cdot|_2$ are the same, and let us show that there
is an $a > 0$ such that (\ref{|x|_2 = |x|_1^a}) holds for every $x \in
k$.  The topological equivalence of $|\cdot|_1$ and $|\cdot|_2$
implies that
\begin{equation}
\label{|x|_1 < 1 if and only if |x|_2 < 1}
        |x|_1 < 1 \quad\hbox{if and only if}\quad |x|_2 < 1
\end{equation}
for every $x \in k$, by the remarks in the preceding paragraph.
Similarly,
\begin{equation}
\label{|x|_1 = 1 if and only if |x|_2 = 1}
        |x|_1 = 1 \quad\hbox{if and only if}\quad |x|_2 = 1
\end{equation}
and
\begin{equation}
\label{|x|_1 > 1 if and only if |x|_2 > 1}
        |x|_1 > 1 \quad\hbox{if and only if}\quad |x|_2 > 1
\end{equation}
for every $x \in k$.  Of course, (\ref{|x|_2 = |x|_1^a}) is trivial
when $x = 0$, and when $|x|_1 = |x|_2 = 1$.  We would like to show
that there is an $a > 0$ such that (\ref{|x|_2 = |x|_1^a}) holds for
every $x \in k$ with $x \ne 0$ and $|x|_1, |x|_2 < 1$.  This would
imply that (\ref{|x|_2 = |x|_1^a}) also holds for every $x \in k$ with
$|x|_1, |x|_2 > 1$, by applying the previous statement to $1/x$.

        Let $y, z \in k$ be given, with $y, z \ne 0$ and $|y|_1, |z|_1 < 1$,
which implies that $|y|_2, |z|_2 < 1$.  If $m$ and $n$ are positive
integers, then we can apply (\ref{|x|_1 < 1 if and only if |x|_2 < 1})
to $x = y^m / z^n$, to get that
\begin{equation}
\label{|y|_1^m < |z|_1^n if and only if |y|_2^m < |z|_2^n}
        |y|_1^m < |z|_1^n \quad\hbox{if and only if}\quad |y|_2^m < |z|_2^n.
\end{equation}
Equivalently, this means that
\begin{equation}
\label{m log |y|_1 < n log |z|_1 if and only if m log |y|_2 < n log |z|_2}
        m \, \log |y|_1 < n \, \log |z|_1 \quad\hbox{if and only if}\quad
        m \, \log |y|_2 < n \, \log |z|_2,
\end{equation}
and hence that
\begin{equation}
\label{frac{n}{m} < frac{log |y|_1}{log |z|_1} if and only if ...}
 \frac{n}{m} < \frac{\log |y|_1}{\log |z|_1} \quad\hbox{if and only if}\quad
  \frac{n}{m} < \frac{\log |y|_2}{\log |z|_2},
\end{equation}
since the logarithms are negative in this situation.  It follows that
\begin{equation}
\label{frac{log |y|_1}{log |z|_1} = frac{log |y|_2}{log |z|_2}}
        \frac{\log |y|_1}{\log |z|_1} = \frac{\log |y|_2}{\log |z|_2},
\end{equation}
so that
\begin{equation}
\label{frac{log |y|_1}{log |y|_2} = frac{log |z|_1}{log |z|_2}}
        \frac{\log |y|_1}{\log |y|_2} = \frac{\log |z|_1}{\log |z|_2}.
\end{equation}
This implies that there is an $a > 0$ such that (\ref{|x|_2 =
  |x|_1^a}) holds for every $x \in k$ with $x \ne 0$ and $|x|_1, |x|_2
< 1$, as desired.

        This follows the discussion that begins on p42 of \cite{fg},
and we shall consider a refinement in the next section, corresponding
to Lemma 3.1 on p18 of \cite{cas}.  Of course, if $|\cdot|_1$ is the
trivial absolute value function on $k$, as in Section
\ref{definitions, examples}, then the corresponding topology on $k$ is
the discrete topology.  Conversely, if the topology on $k$
corresponding to a quasimetric absolute value function $|\cdot|_1$ on
$k$ is the discrete topology, then one can check that $|\cdot|_1$ is
the trivial absolute value function on $k$, using the remarks at the
beginning of the section.

\section{Topological equivalence, 2}
\label{topological equivalence, 2}

        Let $|\cdot|_1$ and $|\cdot|_2$ be quasimetric absolute value
functions on a field $k$ again, and suppose that the topology on $k$
corresponding to $|\cdot|_1$ is at least as strong as the topology on
$k$ corresponding to $|\cdot|_2$.  More precisely, this means that
every open set in $k$ with respect to the topology corresponding to
$|\cdot|_2$ is also an open set in $k$ with respect to the topology
corresponding to $|\cdot|_1$.  In particular, if $x \in k$ and
$\{x^j\}_{j = 1}^\infty$ converges to $0$ with respect to the topology
on $k$ corresponding to $|\cdot|_1$, then it follows that $\{x^j\}_{j
  = 1}^\infty$ also converges to $0$ with respect to the topology on
$k$ corresponding to $|\cdot|_2$.  Thus
\begin{equation}
\label{|x|_1 < 1 implies that |x|_2 < 1}
        |x|_1 < 1 \quad\hbox{implies that}\quad |x|_2 < 1
\end{equation}
for every $x \in k$ under these conditions, by the remarks at the
beginning of the previous section.  If $x \ne 0$, then we can apply
(\ref{|x|_1 < 1 implies that |x|_2 < 1}) to $1/x$, to get that
\begin{equation}
\label{|x|_1 > 1 implies that |x|_2 > 1}
        |x|_1 > 1 \quad\hbox{implies that}\quad |x|_2 > 1.
\end{equation}
Equivalently,
\begin{equation}
\label{|x|_2 le 1 implies that |x|_1 le 1}
        |x|_2 \le 1 \quad\hbox{implies that}\quad |x|_1 \le 1
\end{equation}
for every $x \in k$, which is the contrapositive of (\ref{|x|_1 > 1
  implies that |x|_2 > 1}).  If $u, v \in k$ and $v \ne 0$, then we
can apply (\ref{|x|_2 le 1 implies that |x|_1 le 1}) to $x = u/v$, to
get that
\begin{equation}
\label{|u|_2 le |v|_2 implies that |u|_1 le |v|_1}
        |u|_2 \le |v|_2 \quad\hbox{implies that}\quad |u|_1 \le |v|_1.
\end{equation}

        Of course, all of these statements are trivial when $|\cdot|_1$
is the trivial absolute value function on $k$.  Suppose now that
$|\cdot|_1$ is not the trivial absolute value function on $k$, in
addition to the condition on the associated topologies in the
preceding paragraph.  This means that there is a $v_0 \in k$ such that
$v_0 \ne 0$ and $|v_0|_1 \ne 1$, and we may as well ask that $|v_0|_1
< 1$, since otherwise we could replace $v_0$ with $1/v_0$.  Note that
$|v_0|_2 < 1$ as well, by (\ref{|x|_1 < 1 implies that |x|_2 < 1}).
If we take $v = v_0$ in (\ref{|u|_2 le |v|_2 implies that |u|_1 le |v|_1}),
then we get that
\begin{equation}
\label{|u|_2 le |v_0|_2 implies that |u|_1 le |v_0|_1 < 1}
        |u|_2 \le |v_0|_2 \quad\hbox{implies that}\quad |u|_1 \le |v_0|_1 < 1
\end{equation}
for every $u \in k$.  If $x \in k$ satisfies $|x|_2 < 1$, then
\begin{equation}
\label{|x^n|_2 = |x|_2^n le |v_0|_2}
        |x^n|_2 = |x|_2^n \le |v_0|_2
\end{equation}
for all sufficiently large positive integers $n$, which implies that
$|x|_1^n = |x^n|_1 < 1$, by taking $u = x^n$ in (\ref{|u|_2 le |v_0|_2
  implies that |u|_1 le |v_0|_1 < 1}).  It follows that $|x|_1 < 1$,
so that
\begin{equation}
\label{|x|_2 < 1 implies that |x|_1 < 1}
        |x|_2 < 1 \quad\hbox{implies that}\quad |x|_1 < 1
\end{equation}
for every $x \in k$ under these conditions.  As before, one can
apply this to $1/x$ when $x \ne 0$, to get that
\begin{equation}
\label{|x|_2 > 1 implies that |x|_1 > 1}
        |x|_2 > 1 \quad\hbox{implies that}\quad |x|_1 > 1.
\end{equation}

        At this stage, we are in essentially the same situation as
in the previous section.  More precisely, the combination of
(\ref{|x|_1 < 1 implies that |x|_2 < 1}) and (\ref{|x|_2 < 1 implies
  that |x|_1 < 1}) corresponds exactly to (\ref{|x|_1 < 1 if and only
  if |x|_2 < 1}), and similarly the combination of (\ref{|x|_1 > 1
  implies that |x|_2 > 1}) and (\ref{|x|_2 > 1 implies that |x|_1 >
  1}) corresponds exactly to (\ref{|x|_1 > 1 if and only if |x|_2 >
  1}).  The third condition (\ref{|x|_1 = 1 if and only if |x|_2 = 1})
follows automatically from (\ref{|x|_1 < 1 if and only if |x|_2 < 1})
and (\ref{|x|_1 > 1 if and only if |x|_2 > 1}), and one can use this
to show that there is an $a > 0$ such that (\ref{|x|_2 = |x|_1^a})
holds for every $x \in k$, as before.  In \cite{cas}, one uses 
(\ref{|x|_1 < 1 implies that |x|_2 < 1}) and the nontriviality of
$|\cdot|_1$ to get that
\begin{equation}
\label{|x|_1 = 1 implies that |x|_2 = 1}
        |x|_1 = 1 \quad\hbox{implies that}\quad |x|_2 = 1
\end{equation}
for every $x \in k$.  The combination of (\ref{|x|_1 < 1 implies that
  |x|_2 < 1}), (\ref{|x|_1 > 1 implies that |x|_2 > 1}) and
(\ref{|x|_1 = 1 implies that |x|_2 = 1}) imply (\ref{|x|_1 < 1 if and
  only if |x|_2 < 1}), (\ref{|x|_1 = 1 if and only if |x|_2 = 1}), and
(\ref{|x|_1 > 1 if and only if |x|_2 > 1}) again, so that one can
continue as before.

        In particular, if the topology on $k$ corresponding to $|\cdot|_1$
is at least as strong as the topology on corresponding to $|\cdot|_2$,
and if $|\cdot|_1$ is nontrivial, then it follows that the topologies
on $k$ corresponding to $|\cdot|_1$ and $|\cdot|_2$ are the same.  One
can show that the topology on $k$ corresponding to $|\cdot|_2$ is at
least as strong as the topology corresponding to $|\cdot|_1$ more
directly under these conditions, using (\ref{|u|_2 le |v|_2 implies
  that |u|_1 le |v|_1}).  This also uses the fact that there is a
$v_0 \in k$ such that $v_0 \ne 0$ and $|v_0|_1 < 1$, as before, so that
$|v_0^j|_1 = |v_0|_1^j \to 0$ as $j \to \infty$.

\section{Another refinement}
\label{another refinement}

        Let $k$ be a field, and let $|\cdot|$ be a quasimetric absolute
value function on $k$.  Thus $|\cdot|$ satisfies (\ref{|x + y| le C'
  max(|x|, |y|)}) in Section \ref{definitions, examples} for some $C'
\ge 1$, which implies that $|\cdot|$ satisfies (\ref{|sum_{j =
    1}^{2^n} z_j| le (C')^n max {|z_j| : j = 1, ldots, 2^n}}) in
Section \ref{some refinements} for each nonnegative integer $n$.  Let
$r$ be a nonnegative integer, and let $n$ be the smallest nonnegative
integer such that $r + 1 \le 2^n$.  The minimality of $n$ implies that
$2^{n - 1} < r + 1$, or equivalently
\begin{equation}
\label{2^n < 2 (r + 1)}
        2^n < 2 \, (r + 1).
\end{equation}
If $z_0, z_1, \ldots, z_r$ are $r + 1$ elements of $k$, then
\begin{equation}
\label{|sum_{j = 0}^r z_j| le (C')^n max {|z_j| : j = 0, 1, ldots, r}}
        \biggl|\sum_{j = 0}^r z_j\biggr|
                \le (C')^n \, \max \{|z_j| : j = 0, 1, \ldots, r\},
\end{equation}
by (\ref{|sum_{j = 1}^{2^n} z_j| le (C')^n max {|z_j| : j = 1, ldots,
    2^n}}) and the choice of $n$.  Let us choose $b \ge 0$ so that
$2^b = C'$, and hence
\begin{equation}
\label{(C')^n = (2^b)^n = (2^n)^b le (2 (r + 1))^b = 2^b (r + 1)^b}
        (C')^n = (2^b)^n = (2^n)^b \le (2 \, (r + 1))^b = 2^b \, (r + 1)^b,
\end{equation}
by (\ref{2^n < 2 (r + 1)}).  Combining this with (\ref{|sum_{j = 0}^r
  z_j| le (C')^n max {|z_j| : j = 0, 1, ldots, r}}), we get that
\begin{equation}
\label{|sum_{j = 0}^r z_j| le 2^b (r + 1)^b max {|z_j| : j = 0, 1, ldots, r}}
        \biggl|\sum_{j = 0}^r z_j\biggr|
                \le 2^b \, (r + 1)^b \, \max \{|z_j| : j = 0, 1, \ldots, r\}
\end{equation}
for all $z_0, z_1, \ldots, z_r \in k$.

        Suppose that there is a real number $B \ge 1$ such that
\begin{equation}
\label{|N cdot 1| le B N}
        |N \cdot 1| \le B \, N
\end{equation}
for every positive integer $N$.  Of course, this implies that
\begin{equation}
\label{|N cdot z| le B N |z|}
        |N \cdot z| \le B \, N \, |z|
\end{equation}
for every $z \in k$ and positive integer $N$, using (\ref{N cdot z =
  (N cdot 1) z}) at the beginning of Section \ref{some more
  refinements}.  Let $x, y \in k$ and a positive integer $r$ be given,
so that $(x + y)^r$ can be expressed as in (\ref{(x + y)^r = sum_{j =
    0}^r {r choose j} cdot x^j y^{r - j}}) in Section \ref{some
  refinements}, by the binomial theorem.  It follows that
\begin{equation}
\label{|(x + y)^r| le ..., 3}
 |(x + y)^r| \le 2^b \, (r + 1)^b \,
    \max \bigg\{\biggl|{r \choose j} \cdot x^j \, y^{r - j}\biggr| :
                                              j = 0, 1, \ldots, r\biggr\},
\end{equation}
by (\ref{|sum_{j = 0}^r z_j| le 2^b (r + 1)^b max {|z_j| : j = 0, 1,
    ldots, r}}).  This reduces to
\begin{equation}
\label{|(x + y)^r| le ..., 4}
 \quad  |(x + y)^r| \le 2^b \, (r + 1)^b \, B \,
    \max \bigg\{{r \choose j} \, |x|^j \, |y|^{r - j} :
                                        j = 0, 1, \ldots, r\bigg\},
\end{equation}
using (\ref{|N cdot z| le B N |z|}) for each term in the maximum on
the right side.  Combining this with (\ref{max {{r choose j} |x|^j
    |y|^{r - j} : j = 0, ldots, r} le ...}) in Section \ref{some
  refinements}, we get that
\begin{equation}
\label{|x + y|^r = |(x + y)^r| le 2^b (r + 1)^b B (|x| + |y|)^r}
        |x + y|^r = |(x + y)^r| \le 2^b \, (r + 1)^b \, B \, (|x| + |y|)^r.
\end{equation}
Equivalently,
\begin{equation}
\label{|x + y| le 2^{b/r} (r + 1)^{1/r} B^{1/r} (|x| + |y|)}
        |x + y| \le 2^{b/r} \, (r + 1)^{1/r} \, B^{1/r} \, (|x| + |y|)
\end{equation}
for every $x, y \in k$ and positive integer $r$, which implies that
\begin{equation}
\label{|x + y| le |x| + |y|, 4}
        |x + y| \le |x| + |y|
\end{equation}
for every $x, y \in k$, by taking the limit as $r \to \infty$.  This
shows that $|\cdot|$ is an absolute value function on $k$ under these
conditions.

        Let $|\cdot|$ be any quasimetric absolute value function on a
field $k$ again.  If $N_0 \cdot 1 = 0$ in $k$ for some positive
integer $N_0$, then there are only finitely many elements of $k$ of
the form $N \cdot 1$, where $N$ is a positive integer.  In particular,
this implies that there is a real number $A \ge 1$ such that $|N \cdot
1| \le A$ for every positive integer $N$, which is the same as
(\ref{|N cdot 1| le A}) in Section \ref{some more refinements}.  It
follows that $|\cdot|$ is an ultrametric absolute value function on
$k$ under these conditions, as before.  Otherwise, if $N \cdot 1 \ne
0$ in $k$ for every positive integer $N$, then $k$ has characteristic
$0$, and one can define $r \cdot 1$ in $k$ for every rational number
$r$ in the usual way.  This leads to an embedding of ${\bf Q}$ into
$k$, so that
\begin{equation}
\label{|r cdot 1|}
        |r \cdot 1|
\end{equation}
defines a quasimetric absolute value function on ${\bf Q}$.  As
mentioned at the end of Section \ref{some more refinements}, a theorem
of Ostrowski implies that (\ref{|r cdot 1|}) is either a positive
power of the standard absolute value function on ${\bf Q}$, or a
positive power of the $p$-adic absolute value function on ${\bf Q}$
for some prime number $p$, or the trivial absolute value function on
${\bf Q}$.  In the second and third cases, $|N \cdot 1| \le 1$ for
every positive integer $N$, so that the discussion in Section
\ref{some more refinements} implies that $|\cdot|$ is an ultrametric
absolute value function on $k$.  In the first case, let us suppose
that (\ref{|r cdot 1|}) is equal to the standard absolute value
function on ${\bf Q}$, which can be arranged by replacing $|\cdot|$
with the appropriate positive power of itself.  This implies that
$|\cdot|$ satisfies (\ref{|N cdot 1| le B N}) with $B = 1$, and hence
that $|\cdot|$ is an absolute value function on $k$.

\section{Complex numbers}
\label{complex numbers}

        Let $|z|_0$ denote the standard absolute value function on the
complex numbers ${\bf C}$, so that
\begin{equation}
\label{|z|_0 = (x^2 + y^2)^{1/2}}
        |z|_0 = (x^2 + y^2)^{1/2}
\end{equation}
for each $z \in {\bf C}$, where $z = x + i \, y$ and $x, y \in {\bf
  R}$.  Suppose that $|z|$ is another absolute value function on ${\bf
  C}$ such that
\begin{equation}
\label{|z| = |z|_0}
        |z| = |z|_0
\end{equation}
for every $z \in {\bf R}$.  Note that
\begin{equation}
\label{|i| = 1}
        |i| = 1,
\end{equation}
as in (\ref{|x|^n = |x^n| = |1| = 1}), because $i^4 = 1$.  If $z = x +
i \, y$ with $x, y \in {\bf R}$ again, then it follows that
\begin{equation}
\label{|z| le |x| + |y| = |x|_0 + |y|_0 le 2 |z|_0}
        |z| \le |x| + |y| = |x|_0 + |y|_0 \le 2 |z|_0,
\end{equation}
where the $2$ in the last step could be replaced by a $\sqrt{2}$,
using the triangle inequality.  This implies that the standard
topology on ${\bf C}$, which corresponds to $|z|_0$, is at least as
strong as the topology corresponding to $|z|$.

        Of course, $|z|_0$ is not the trivial absolute value function
on ${\bf C}$.  Thus the discussion in Section \ref{topological
  equivalence, 2} implies that the topology on ${\bf C}$ corresponding
to $|z|$ is the same as the standard topology.  Alternatively, one
could get the same conclusion by considering ${\bf C}$ as a
two-dimensional vector space over ${\bf R}$, as in Section
\ref{finite-dimensional vector spaces}.  At any rate, the argument in
Section \ref{topological equivalence} implies that $|z|$ is equal to a
positive power of $|z|_0$.  It is easy to see that this power has to
be $1$, because of (\ref{|z| = |z|_0}), so that (\ref{|z| = |z|_0})
holds for every $z \in {\bf C}$ under these conditions.

        Now let $k$ be a field with a quasimetric absolute value function
$|\cdot|$, and suppose that $|\cdot|$ is archimedian.  This implies that
$k$ has characteristic $0$, so that there is a natural embedding of
${\bf Q}$ into $k$.  The restriction of $|\cdot|$ to the image of ${\bf Q}$
in $k$ is also archimedian, and hence the restriction of $|\cdot|$ to
the image of ${\bf Q}$ is equal to a positive power of the standard
absolute value function on ${\bf Q}$.  By replacing $|\cdot|$ with a
suitable positive power of itself, if necessary, we may as well ask
that the restriction of $|\cdot|$ to the image of ${\bf Q}$ in $k$ be
equal to the standard absolute value function on ${\bf Q}$.  The argument
in the previous section then implies that $|\cdot|$ is an absolute value
function on $k$.

        If $k$ is complete with respect to the metric corresponding
to $|\cdot|$, then the natural embedding of ${\bf Q}$ into $k$ can be
extended to an embedding of ${\bf R}$ into $k$, such that the
restriction of $|\cdot|$ to the image of ${\bf R}$ in $k$ is equal to
the standard absolute value function on ${\bf R}$.  Suppose that there
is an element $i$ of $k$ such that $i^2 = -1$.  This implies that the
embedding of ${\bf R}$ into $k$ can be extended to an embedding of
${\bf C}$ into $k$ in the obvious way.  The restriction of $|\cdot|$
to the image of ${\bf C}$ in $k$ is equal to the standard absolute
value function on ${\bf C}$, by the argument at the beginning of the
section.

        In particular, this permits one to consider $k$ as a complex
Banach algebra.  A famous theorem implies that the image of ${\bf C}$
in $k$ is equal to $k$ under these conditions, because $k$ is a field.
More precisely, suppose that $\zeta \in k$ is not in the image of
${\bf C}$ in $k$, and put
\begin{equation}
\label{f(z) = (zeta - z)^{-1}}
        f(z) = (\zeta - z)^{-1}
\end{equation}
for each $z \in {\bf C}$.  Here we identify $z$ with its image in $k$
on the right side of (\ref{f(z) = (zeta - z)^{-1}}), so that $\zeta -
z$ is considered as an element of $k$ for each $z \in {\bf C}$.  By
hypothesis, $\zeta - z \ne 0$ in $k$ for each $z \in {\bf C}$, so that
the right side of (\ref{f(z) = (zeta - z)^{-1}}) makes sense as an
element of $k$.

        The theory of holomorphic functions can be extended to functions
that take values in a complex Banach space, such as $k$.  One can show
that (\ref{f(z) = (zeta - z)^{-1}}) is holomorphic as a mapping from
${\bf C}$ into $k$ in this sense, and in fact that (\ref{f(z) = (zeta
  - z)^{-1}}) has a local power series expansion at each point.  We
also have that
\begin{equation}
\label{|zeta - z| to infty as z to infty in {bf C}}
        |\zeta - z| \to \infty \quad\hbox{as } z \to \infty \hbox{ in } {\bf C},
\end{equation}
by the triangle inequality, which implies that
\begin{equation}
\label{(zeta - z)^{-1} to 0 as z to infty in {bf C}}
 (\zeta - z)^{-1} \to 0 \quad\hbox{as } z \to \infty \hbox{ in } {\bf C},
\end{equation}
by the multiplicative property of absolute value functions.  Although
the second step does not work in the context of Banach algebras, where
the norm is submultiplicative and not necessarily multiplicative, one
can still get (\ref{(zeta - z)^{-1} to 0 as z to infty in {bf C}})
using another argument.  Of course, (\ref{f(z) = (zeta - z)^{-1}}) is
a continuous function on ${\bf C}$, which implies that it is bounded
on compact subsets of ${\bf C}$.  Thus (\ref{f(z) = (zeta - z)^{-1}})
is bounded on ${\bf C}$, by (\ref{(zeta - z)^{-1} to 0 as z to infty
  in {bf C}}).  An appropriate version of Liouville's theorem implies
that (\ref{f(z) = (zeta - z)^{-1}}) should be constant on ${\bf C}$,
and hence identically $0$, by (\ref{(zeta - z)^{-1} to 0 as z to infty
  in {bf C}}).  This is a contradiction, because (\ref{f(z) = (zeta -
  z)^{-1}}) is not supposed to be $0$ for any $z \in {\bf C}$.  It
follows that $\zeta \in k$ is in the image of ${\bf C}$, as desired.

        A more elementary approach in the setting of absolute value
functions on fields can be found on p38-9 of \cite{cas}.  Let
$\zeta \in k$ be given, and put
\begin{equation}
\label{g(z) = |zeta - z|}
        g(z) = |\zeta - z|
\end{equation}
for each $z \in {\bf C}$, where $z$ is identified with an element of
$k$, as before.  Thus $g(z)$ defines a continuous real-valued function
${\bf C}$ that satisfies (\ref{|zeta - z| to infty as z to infty in
  {bf C}}), which implies that the minimum of $g(z)$ is attained on
${\bf C}$, since closed and bounded subsets of ${\bf C}$ are compact.
If $\zeta$ is not in the image of ${\bf C}$ in $k$, then the minimum
of $g(z)$ on $k$ is positive, and it is shown in \cite{cas} that
$g(z)$ is constant on ${\bf C}$, contradicting (\ref{|zeta - z| to
  infty as z to infty in {bf C}}).  Note that minimizing $g(z)$ on
${\bf C}$ corresponds to maximizing $|f(z)|$ on ${\bf C}$, where
$f(z)$ is as in (\ref{f(z) = (zeta - z)^{-1}}).  That the maximum of
$|f(z)|$ is attained on ${\bf C}$ can be derived from (\ref{(zeta -
  z)^{-1} to 0 as z to infty in {bf C}}), which also works in the
context of Banach algebras.  One can show that $|f(z)|$ is subharmonic
on ${\bf C}$, because $f(z)$ is holomorphic on ${\bf C}$.  It follows
from well-known results about subharmonic functions that $|f(z)|$ is
constant on ${\bf C}$, since it attains its maximum.  This implies
that $|f(z)|$ is identically $0$ on ${\bf C}$, by (\ref{(zeta -
  z)^{-1} to 0 as z to infty in {bf C}}), which is a contradiction
again.  Of course, this argument is very closely related to the one
using Liouville's theorem discussed in the preceding paragraph.

\section{Local compactness}
\label{local compactness}

        Remember that a topological space $X$ is said to be
\emph{locally compact}\index{locally compact topological spaces}
if for each $x \in X$ there is an open set $U \subseteq X$ such that
$x \in U$ and $U$ is contained in compact subset $K$ of $X$.  If $X$
is Hausdorff, then $K$ is also a closed set in $X$, and hence the
closure $\overline{U}$ of $U$ in $X$ is contained in $K$.  This
implies that $\overline{U}$ is compact, since closed subsets of
compact sets are compact.  Now let $k$ be a field, and let $|\cdot|$
be a quasimetric absolute value function on $k$.  Thus $k$ is a
Hausdorff topological space with respect to the topology determined by
the corresponding quasimetric (\ref{d(x, y) = |x - y|, 2}), as usual.
If $k$ is locally compact, then there is an open set $U$ in $k$ such
that $0 \in U$ and the closure $\overline{U}$ of $U$ in $k$ is
compact.  This implies that there is an $r > 0$ such that the closed
ball $\overline{B}(0, r)$ centered at $0$ with radius $r$ in $k$ is
compact.  It follows that every closed ball in $k$ with radius $r$ is
compact, by continuity of translations.

        If $\{x_j\}_{j = 1}^\infty$ is a Cauchy sequence of elements of
$k$, then for each $\epsilon > 0$ there is an $L(\epsilon) \ge 1$ such that
\begin{equation}
\label{|x_j - x_l| < epsilon}
        |x_j - x_l| < \epsilon
\end{equation}
for every $j, l \ge L(\epsilon)$.  If $k$ is locally compact, then
there is an $r > 0$ such that every closed ball in $k$ of radius $r$
is compact, as in the previous paragraph.  This implies that all but
finitely many terms of the sequence $\{x_j\}_{j = 1}^\infty$ are
contained in a compact set, so that there is a subsequence of
$\{x_j\}_{j = 1}^\infty$ that converges to an element of this compact
set.  It follows that the whole sequence $\{x_j\}_{j = 1}^\infty$
converges to the same element of $k$, because $\{x_j\}_{j = 1}^\infty$
is a Cauchy sequence.  This shows that $k$ is complete as a
quasimetric space with respect to the quasimetric associated to
$|\cdot|$ when $k$ is locally compact.

        Suppose that $|\cdot|$ is not the trivial absolute value function
on $k$, so that there is an $w_0 \in k$ such that $w_0 \ne 0$ and
$|w_0| \ne 1$.  More precisely, either $0 < |w_0| < 1$ or $|w_0| > 1$,
and we may as well ask that $|w_0| > 1$, since otherwise we could
replace $w_0$ with $1/w_0$.  If $k$ is locally compact, then there is
an $r > 0$ such that $\overline{B}(0, r)$ is compact, as before.
Of course,
\begin{equation}
\label{x mapsto x w_0^n}
        x \mapsto x \, w_0^n
\end{equation}
is a continuous mapping on $k$ for each positive integer $n$, which
sends $\overline{B}(0, r)$ onto $\overline{B}(0, r \, |w_0|^n)$.  Thus
$\overline{B}(0, r \, |w_0|^n)$ is a compact set in $k$ for each
positive integer $n$ under these conditions.  Every bounded subset of
$k$ is contained in $\overline{B}(0, r \, |w_0|^n)$ when $n$ is
sufficiently large, because $|w_0|^n \to \infty$ as $n \to \infty$,
since $|w_0| > 1$.  It follows that closed and bounded subsets of
$k$ are compact in this situation, because closed subsets of compact
sets are compact.

        If $|\cdot|$ is the trivial absolute value function on $k$,
then the corresponding topology on $k$ is discrete, and hence locally
compact.  However, the closed unit ball in $k$ is not compact in this
case, unless $k$ has only finitely many elements.  Note that every Cauchy
sequence in $k$ is eventually constant when $|\cdot|$ is trivial, so
that $k$ is automatically complete.  If $|\cdot|$ is nontrivial on $k$
and $k$ is locally compact, then every Cauchy sequence of elements of
$k$ is contained in a compact subset of $k$, because Cauchy sequences
are bounded, and closed and bounded subsets of $k$ are compact.
This gives a slightly different way to look at the completeness of $k$
when $k$ is locally compact, by combining these two cases.

\section{An auxiliary fact}
\label{an auxiliary fact}

        Let $k$ be a field, and suppose that
\begin{equation}
\label{x^2 ne -1}
        x^2 \ne -1
\end{equation}
in $k$ for every $x \in k$.  Note that this implies that $k$ does not
have characteristic $2$, because $1 = -1$ in $k$ when $k$ has
characteristic $2$.  Also let $|\cdot|$ be an absolute value function
on $k$, and suppose that $k$ is complete with respect to the metric
corresponding to $|\cdot|$.  Under these conditions, there is a
positive real number $c$ such that
\begin{equation}
\label{|x^2 + 1| ge c}
        |x^2 + 1| \ge c
\end{equation}
for every $x \in k$.  In fact, one can take
\begin{equation}
\label{c = |4 cdot 1| (1 + |4 cdot 1|)^{-1}}
        c = |4 \cdot 1| \, (1 + |4 \cdot 1|)^{-1},
\end{equation}
as in the proof of Lemma 2.2 starting on p35 of \cite{cas}.  More
precisely, if there is an $x_1 \in k$ such that $|x_1^2 + 1|$ is less
than (\ref{c = |4 cdot 1| (1 + |4 cdot 1|)^{-1}}), then it is shown in
\cite{cas} that there is a Cauchy sequence $\{x_j\}_{j = 1}^\infty$ of
elements of $k$ such that $\{x_j^2\}_{j = 1}^\infty$ converges to
$-1$.  If $k$ is complete, then it follows that $\{x_j\}_{j =
  1}^\infty$ converges to an element $x$ of $k$ such that $x^2 = -1$,
contradicting (\ref{x^2 ne -1}).

        Of course, if $x \in k$ satisfies $|x| \ge \sqrt{2}$, then
$|x^2| \ge 2$, and hence
\begin{equation}
\label{|x^2 + 1| ge 2 - 1 = 1}
        |x^2 + 1| \ge 2 - 1 = 1,
\end{equation}
by the triangle inequality.  Thus it suffices to consider $x \in k$
with $|x| \le \sqrt{2}$ to get (\ref{|x^2 + 1| ge c}).  If the set of
$x \in k$ with $|x| \le \sqrt{2}$ is compact, then the existence of $c
> 0$ as in (\ref{|x^2 + 1| ge c}) can be derived from the extreme
value theorem.  If $k$ is locally compact, and if $|\cdot|$ is not the
trivial absolute value function on $k$, then every closed ball in $k$
is compact, as in the previous section.  If $|\cdot|$ is the trivial
absolute value function on $k$, then (\ref{|x^2 + 1| ge c}) holds with
$c = 1$ whenever $x^2 \ne -1$.

        Here is another argument that works when $k$ is complete, and
not necessarily locally compact.  Suppose that $x, y \in k$ satisfy
\begin{equation}
\label{|x^2 + 1|, |y^2 + 1| le eta}
        |x^2 + 1|, \, |y^2 + 1| \le \eta
\end{equation}
for some $\eta > 0$, so that
\begin{equation}
\label{1 - eta le |x^2|, |y^2| le 1 + eta}
        1 - \eta \le |x^2|, \, |y^2| \le 1 + \eta,
\end{equation}
by the triangle inequality.  Observe that
\begin{equation}
\label{|x^2 - y^2| = |(x^2 + 1) - (y^2 + 1)| le |x^2 + 1| + |y^2 + 1| le 2 eta}
 |x^2 - y^2| = |(x^2 + 1) - (y^2 + 1)| \le |x^2 + 1| + |y^2 + 1| \le 2 \, \eta,
\end{equation}
and hence
\begin{equation}
\label{|x - y| |x + y| = |(x - y) (x + y)| = |x^2 - y^2| le 2 eta}
        |x - y| \, |x + y| = |(x - y) \, (x + y)| = |x^2 - y^2| \le 2 \, \eta.
\end{equation}
It follows that
\begin{equation}
\label{|x - y| le sqrt{2 eta} or |x + y| le sqrt{2 eta}}
 |x - y| \le \sqrt{2 \, \eta} \quad\hbox{or}\quad |x + y| \le \sqrt{2 \, \eta}.
\end{equation}
We also have that
\begin{equation}
\label{|x + x| = |(1 + 1) x| = |2 cdot 1| |x| ge (1 - eta)^{1/2}}
        |x + x| = |(1 + 1) \, x| = |2 \cdot 1| \, |x|
                                \ge |2 \cdot 1| \, (1 - \eta)^{1/2},
\end{equation}
by (\ref{1 - eta le |x^2|, |y^2| le 1 + eta}), and thus
\begin{equation}
\label{|2 cdot 1| (1 - eta)^{1/2} le |x + x| le |x - y| + |x + y|}
        |2 \cdot 1| \, (1 - \eta)^{1/2} \le |x + x| \le |x - y| + |x + y|.
\end{equation}

        Suppose that $\{x_j\}_{j = 1}^\infty$ is a sequence of elements
of $k$ such that
\begin{equation}
\label{|x_j^2 + 1| le eta}
        |x_j^2 + 1| \le \eta
\end{equation}
for each $j$.  Of course, $-x_j$ satisfies the same condition, and so
we may suppose that
\begin{equation}
\label{|x_1 - x_j| le sqrt{2 eta}}
        |x_1 - x_j| \le \sqrt{2 \, \eta}
\end{equation}
for each $j$, by (\ref{|x - y| le sqrt{2 eta} or |x + y| le sqrt{2
    eta}}), and replacing $x_j$ with $-x_j$ whenever necessary.  This
implies that
\begin{equation}
\label{|x_j - x_l| le |x_j - x_1| + |x_1 - x_l| le 2 sqrt{2 eta}}
        |x_j - x_l| \le |x_j - x_1| + |x_1 - x_l| \le 2 \, \sqrt{2 \, \eta}
\end{equation}
for every $j, l \ge 1$, by the triangle inequality.  Combining this with
(\ref{|2 cdot 1| (1 - eta)^{1/2} le |x + x| le |x - y| + |x + y|}), we get that
\begin{equation}
\label{|2 cdot 1| (1 - eta)^{1/2} le ... le 2 sqrt{2 eta} + |x_j + x_l|}
 |2 \cdot 1| \, (1 - \eta)^{1/2} \le |x_j - x_l| + |x_j + x_l|
                                 \le 2 \, \sqrt{2 \, \eta} + |x_j + x_l|
\end{equation}
for every $j, l \ge 1$.  If $\eta$ is sufficiently small, then this gives
a positive lower bound for $|x_j + x_l|$ for every $j, l \ge 1$.

        Let us now fix $\eta > 0$ so that the preceding statement holds.
If there is no $c > 0$ such that (\ref{|x^2 + 1| ge c}) holds for every
$x \in k$, then there is a sequence $\{x_j\}_{j = 1}^\infty$ of elements
of $k$ such that $\{x_j^2\}_{j = 1}^\infty$ converges to $-1$.  In particular,
we may ask that $x_j$ satisfy (\ref{|x_j^2 + 1| le eta}) for every $j$,
and also (\ref{|x_1 - x_j| le sqrt{2 eta}}), by replacing $x_j$ with $-x_j$,
as before.  Because $\{x_j^2\}_{j = 1}^\infty$ converges to $-1$, we have that
\begin{equation}
\label{|x_j^2 - x_l^2| to 0 as j, l to infty}
        |x_j^2 - x_l^2| \to 0 \quad\hbox{as } j, l \to \infty.
\end{equation}
This implies that either $|x_j - x_l|$ or $|x_j + x_l|$ is as small as
we like when $j$ and $l$ are sufficiently large, as in (\ref{|x - y|
  |x + y| = |(x - y) (x + y)| = |x^2 - y^2| le 2 eta}) and (\ref{|x -
  y| le sqrt{2 eta} or |x + y| le sqrt{2 eta}}).  Using the positive
lower bound for $|x_j + x_l|$ mentioned in the previous paragraph, we
get that $|x_j - x_l|$ is as small as we like when $j$ and $l$ are
sufficiently small.  Thus $\{x_j\}_{j = 1}^\infty$ is a Cauchy
sequence in $k$ under these conditions, which converges to an element
$x$ of $k$ such that $x^2 = 1$ when $k$ is complete.

        This argument would also work for a quasimetric absolute value
function $|\cdot|$ on $k$, with some minor adjustments.
Alternatively, one could reduce to the case of an absolute value
function, by replacing $|\cdot|$ with a suitable positive power of
itself.  If $|\cdot|$ is an ultrametric absolute value function on
$k$, then this argument could be simplified a bit.  Otherwise, in the
archimedian case, $k$ has characteristic $0$, and one can choose
$|\cdot|$ so that it agrees with the usual absolute value function on
the natural image of ${\bf Q}$ in $k$, as in Section \ref{another
  refinement}.  In this case, $|n \cdot 1| = n$ for each positive
integer $n$, which leads to some other simplifications.

        If (\ref{|x^2 + 1| ge c}) holds for every $x \in k$, then it
is easy to see that
\begin{equation}
\label{|x^2 + y^2| ge c max(|x|^2, |y|^2)}
        |x^2 + y^2| \ge c \, \max(|x|^2, |y|^2)
\end{equation}
for every $x, y \in k$.  More precisely, this holds with $c = 1$ when
$x = 0$ or $y = 0$, and so we may as well suppose that $x, y \ne 0$.
If $y \ne 0$, then
\begin{equation}
\label{|x^2 + y^2| = |(x/y)^2 + 1| |y|^2 ge c |y|^2}
        |x^2 + y^2| = |(x/y)^2 + 1| \, |y|^2 \ge c \, |y|^2
\end{equation}
for every $x \in k$, by (\ref{|x^2 + 1| ge c}).  Using the analogous
estimate when $x \ne 0$, we get that (\ref{|x^2 + y^2| ge c max(|x|^2,
  |y|^2)}) holds for every $x, y \in k$, as desired.

\section{Complex numbers, 2}
\label{complex numbers, 2}

        Let $k$ be a field, and suppose that $x^2 \ne -1$ for every
$x \in k$.  Consider the field $k(i)$\index{k(i)@$k(i)$} obtained by
adjoining to $k$ an additional element $i$ such that $i^2 = -1$.  More
precisely, one can think of $k(i)$ initially as a two-dimensional
vector space over $k$, in which every element can be expressed as $x +
i \, y$, where $x, y \in k$ and $i \not\in k$.  It is easy to define
multiplication on $k(i)$ in such a way that $k(i)$ becomes a
commutative ring, using multiplication on $k$, and by defining $i^2$
to be $-1$ in $k$.  One can also define complex conjugation on $k(i)$
in the usual way, so that
\begin{equation}
\label{overline{(x + i y)} = x - i y}
        \overline{(x + i \, y)} = x - i \, y
\end{equation}
for every $x, y \in k$.  This determines a ring automorphism on $k$,
which satisfies
\begin{equation}
\label{(x + i y) overline{(x + i y)} = (x + i y) (x - i y) = x^2 + y^2}
        (x + i \, y) \, \overline{(x + i \, y)}
                       = (x + i \, y) \, (x - i \, y) = x^2 + y^2
\end{equation}
for every $x, y \in k$.  If $y \ne 0$, then
\begin{equation}
\label{x^2 + y^2 = ((x/y)^2 + 1) y^2 ne 0}
        x^2 + y^2 = ((x/y)^2 + 1) \, y^2 \ne 0
\end{equation}
for every $x \in k$, because $(x/y)^2 \ne -1$ by hypothesis.
Similarly, $x^2 + y^2 \ne 0$ when $x \ne 0$.  If $x + i \, y \ne 0$ in
$k(i)$, then $x \ne 0$ or $y \ne 0$, which implies that $x^2 + y^2 \ne
0$, and hence that $x^2 + y^2$ has a multiplicative inverse in $k$,
because $x^2 + y^2 \in k$ and $k$ is a field.  It follows that $x + i
\, y$ has a multiplicative inverse in $k(i)$ under these conditions,
given by
\begin{equation}
\label{(x + i y)^{-1} = (x - i y) (x^2 + y^2)^{-1}}
        (x + i \, y)^{-1} = (x - i \, y) \, (x^2 + y^2)^{-1},
\end{equation}
so that $k(i)$ is a field too.

        Let $|\cdot|$ be a quasimetric absolute value function on $k$.
It is natural to extend $|\cdot|$ to $k(i)$ by putting
\begin{equation}
\label{|x + i y| = |x^2 + y^2|^{1/2}}
        |x + i \, y| = |x^2 + y^2|^{1/2}
\end{equation}
for every $x, y \in k$, where the right side of (\ref{|x + i y| = |x^2
  + y^2|^{1/2}}) is defined using the absolute value of $x^2 + y^2$ as
an element of $k$.  If $x + i \, y \ne 0$ in $k(i)$, then $x^2 + y^2
\ne 0$ in $k$, as in the preceding paragraph.  This implies that the
right side of (\ref{|x + i y| = |x^2 + y^2|^{1/2}}) is positive when
$x + i \, y \ne 0$ in $k(i)$.  Observe also that
\begin{equation}
\label{|overline{z}| = |z|}
        |\overline{z}| = |z|
\end{equation}
for every $z \in k(i)$, by construction.

        Equivalently,
\begin{equation}
\label{|z|^2 = |z overline{z}|}
        |z|^2 = |z \, \overline{z}|
\end{equation}
for every $z \in k(i)$, by (\ref{(x + i y) overline{(x + i y)} = (x +
  i y) (x - i y) = x^2 + y^2}).  Thus
\begin{equation}
\label{|z w|^2 = ... = |z overline{z} w overline{w}|}
        |z \, w|^2 = |(z \, w) \, \overline{(z \, w)}|
                   = |z \, w \, \overline{z} \, \overline{w}|
                   = |z \, \overline{z} \, w \, \overline{w}|
\end{equation}
for every $z, w \in k(i)$, using the fact that complex conjugation is
a field automorphism in the second step, and commutativity of
multiplication on $k(i)$ in the third step.  Because $z \,
\overline{z}, w \, \overline{w} \in k$, the multiplicative property
(\ref{|x y| = |x| |y|, 2}) of $|\cdot|$ on $k$ implies that
\begin{equation}
\label{|z w|^2 = |z overline{z}| |w overline{w}| = |z|^2 |w|^2}
 |z \, w|^2 = |z \, \overline{z}| \, |w \, \overline{w}| = |z|^2 \, |w|^2
\end{equation}
for every $z, w \in k(i)$, using (\ref{|z|^2 = |z overline{z}|}) in
the second step.  This shows that this extension of $|\cdot|$ to
$k(i)$ also satisfies the multiplicative property (\ref{|x y| = |x|
  |y|, 2}).

        Note that
\begin{equation}
\label{|x + i y| = |x^2 + y^2|^{1/2} le sqrt{C'} max(|x|, |y|)}
        |x + i \, y| = |x^2 + y^2|^{1/2} \le \sqrt{C'} \, \max(|x|, |y|)
\end{equation}
for every $x, y \in k$, where $C' \ge 1$ is as in (\ref{|x + y| le C'
  max(|x|, |y|)}), and where (\ref{|x + y| le C' max(|x|, |y|)}) is
applied to the absolute value of $x^2 + y^2$ as an element of $k$ in
the second step.  Suppose now that there is a $c > 0$ such that
(\ref{|x^2 + y^2| ge c max(|x|^2, |y|^2)}) holds for every $x, y \in
k$.  This implies that
\begin{equation}
\label{|x + i y| = |x^2 + y^2|^{1/2} ge sqrt{c} max(|x|, |y|)}
        |x + i \, y| = |x^2 + y^2|^{1/2} \ge \sqrt{c} \, \max(|x|, |y|)
\end{equation}
for every $x, y \in k$.  Under these conditions, it is easy to see
that this extension of $|\cdot|$ to $k(i)$ is a quasimetric absolute
value function, using (\ref{|x + i y| = |x^2 + y^2|^{1/2} le sqrt{C'}
  max(|x|, |y|)}) and (\ref{|x + i y| = |x^2 + y^2|^{1/2} ge sqrt{c}
  max(|x|, |y|)}).

        Alternatively, if $z, w \in k(i)$, then
\begin{equation}
\label{(z + w) overline{(z + w)} = (z + w) (overline{z} + overline{w}) = ...}
 (z + w) \, \overline{(z + w)} = (z + w) \, (\overline{z} + \overline{w})
   = z \, \overline{z} + z \, \overline{w} + w \, \overline{z}
                                                    + w \, \overline{w},
\end{equation}
using the fact that complex conjugation is a field automorphism in the
first step.  It follows that
\begin{equation}
\label{|z + w|^2 = |(z + w) (overline{z} + overline{w})| = ...}
        |z + w|^2 = |(z + w) \, \overline{(z + w)}| 
 = |z \, \overline{z} + (z \, \overline{w} + w \, \overline{z})
                                    + w \, \overline{w}|
\end{equation}
for every $z, w \in k(i)$.  Of course, $z \, \overline{z}$, $w \,
\overline{w}$, and (\ref{(z + w) overline{(z + w)} = (z + w)
  (overline{z} + overline{w}) = ...}) are elements of $k$, which
implies that
\begin{equation}
\label{z overline{w} + w overline{z}}
        z \, \overline{w} + w \, \overline{z}
\end{equation}
is an element of $k$ too.  Because $z \, \overline{z}$, $w \,
\overline{w}$, and (\ref{z overline{w} + w overline{z}}) are elements
of $k$, (\ref{|z + w|^2 = |(z + w) (overline{z} + overline{w})| =
  ...}) can be estimated in terms of $|z \, \overline{z}| = |z|^2$,
$|w \, \overline{w}| = |w|^2$, and the absolute value of (\ref{z
  overline{w} + w overline{z}}), using the quasimetric version of the
triangle inequality for $|\cdot|$ on $k$.

        More precisely, $2 \cdot z \, \overline{w}$ is the sum of
(\ref{z overline{w} + w overline{z}}) and
\begin{equation}
\label{z overline{w} - w overline{z}}
      z \, \overline{w} - w \overline{z},
\end{equation}
where (\ref{z overline{w} + w overline{z}}) is in $k$, and (\ref{z
  overline{w} - w overline{z}}) is $i$ times an element of $k$.  Thus
(\ref{|x + i y| = |x^2 + y^2|^{1/2} ge sqrt{c} max(|x|, |y|)}) implies
that
\begin{equation}
        \sqrt{c} \, |z \, \overline{w} + \overline{z} \, w|
                     \le |2 \cdot z \, \overline{w}|.
\end{equation}
We also have that
\begin{equation}
\label{|2 cdot z overline{w}| = ... = |2 cdot 1| |z| |w|}
        |2 \cdot z \, \overline{w}| = |2 \cdot 1| \, |z| \, |\overline{w}|
                                    = |2 \cdot 1| \, |z| \, |w|,
\end{equation}
using the multiplicative property of $|\cdot|$ on $k(i)$ in the first
step, and (\ref{|overline{z}| = |z|}) in the second step.  This
permits one to estimate (\ref{|z + w|^2 = |(z + w) (overline{z} +
  overline{w})| = ...}) in terms of $|z|^2$, $|w|^2$, and $|z| \,
|w|$, to get that this extension of $|\cdot|$ to $k(i)$ is a
quasimetric absolute value function too.  If $|\cdot|$ is an absolute
value function on $k$ and (\ref{|x + i y| = |x^2 + y^2|^{1/2} ge
  sqrt{c} max(|x|, |y|)}) holds with $c = 1$, then one could use this
argument to show that this extension of $|\cdot|$ to $k(i)$ is an
absolute value function as well.

        Of course, $k(i)$ can be identified with $k \times k$ in an 
obvious way.  This leads to a topology on $k(i)$ that corresponds to
the product topology on $k \times k$, using the topology on $k$
determined by the quasimetric associated to $|\cdot|$.  There is also
a topology on $k(i)$ determined by the quasimetric associated to the
extension of $|\cdot|$ to $k(i)$ just defined.  It is easy to see that
these two topologies on $k(i)$ are the same, using (\ref{|x + i y| =
  |x^2 + y^2|^{1/2} le sqrt{C'} max(|x|, |y|)}) and (\ref{|x + i y| =
  |x^2 + y^2|^{1/2} ge sqrt{c} max(|x|, |y|)}).  Similarly, if $k$ is
complete with respect to the quasimetric associated to $|\cdot|$, then
$k(i)$ is complete with respect to the quasimetric associated to the
extension of $|\cdot|$ to $k(i)$ just defined.

        Let us now consider the case where $|\cdot|$ is an archimedian
absolute value function on $k$.  As before, this implies that $k$ has
characteristic $0$, and we may as well suppose that the restriction of
$|\cdot|$ to the natural image of ${\bf Q}$ in $k$ is equal to the
standard absolute value on ${\bf Q}$.  Let us also suppose that $k$ is
complete with respect to the metric associated to $|\cdot|$.  This
implies that the natural embedding of ${\bf Q}$ into $k$ extends to an
embedding of ${\bf R}$ into $k$, and that the restriction of $|\cdot|$
to the image of ${\bf R}$ in $k$ is equal to the standard absolute
value on ${\bf R}$.  If $x^2 \ne -1$ for every $x \in k$, then there
is a $c > 0$ such that (\ref{|x^2 + y^2| ge c max(|x|^2, |y|^2)})
holds for every $x, y \in k$, as in the previous section.  Thus the
extension of $|\cdot|$ to $k(i)$ is a quasimetric absolute value
function, and in fact this extension is an absolute value function on
$k(i)$, as in Section \ref{another refinement}.  Remember that $k(i)$
is complete with respect to the metric associated to the extension of
$|\cdot|$ to $k(i)$, as in the preceding paragraph.  The natural
embedding of ${\bf R}$ into $k$ leads to an embedding of ${\bf C}$
into $k(i)$.  The image of ${\bf C}$ in $k(i)$ is actually equal to
$k(i)$ under these conditions, as in Section \ref{complex numbers}.
This implies that the image of ${\bf R}$ in $k$ is equal to $k$.  The
characterization of complete fields with archimedian absolute value
functions described in Section \ref{complex numbers} and this section
is another famous theorem of Ostrowski,\index{Ostrowski's theorems}
as on p33 of \cite{cas}.

\chapter{Additional structure}
\label{additional structure}

\section{$p$-Adic integers}
\label{p-adic integers}

        Let $p$ be a prime number, and let $|\cdot|_p$ be the $p$-adic
absolute value function on ${\bf Q}$, as in Section \ref{definitions,
  examples}.  Also let ${\bf Z}$\index{Z@${\bf Z}$} be the set of integers,
and observe that
\begin{equation}
\label{|x|_p le 1}
        |x|_p \le 1
\end{equation}
for every $x \in {\bf Z}$.  Suppose now that $y \in {\bf Q}$ satisfies
$|y|_p \le 1$, and let us check that $y$ can be approximated by
integers with respect to the $p$-adic metric.  By hypothesis, $y$ can
be expressed as $a/b$, where $a$ and $b$ are integers, $b \ne 0$, and
$b$ is not an integer multiple of $p$.  Thus $b \not\equiv 0$ modulo
$p$, which implies that there is an integer $c$ such that $b \, c
\equiv 1$ modulo $p$, because $p$ is prime.  Equivalently, there is a
$w \in {\bf Z}$ such that $b \, c = 1 - p \, w$, and hence
\begin{equation}
\label{y = ... = a c lim_{n to infty} sum_{j = 0}^n (p w)^j}
        y = \frac{a}{b} = \frac{a \, c}{b \, c} = \frac{a \, c}{1 - p \, w}
                   = a \, c \, \lim_{n \to \infty} \sum_{j = 0}^n (p \, w)^j.
\end{equation}
This uses (\ref{lim_{n to infty} sum_{j = 0}^n x^j = frac{1}{1 - x}})
applied to $p \, w$, which is permissible because
\begin{equation}
\label{|p w|_p = |p|_p |w|_p le 1/p < 1}
        |p \, w|_p = |p|_p \, |w|_p \le 1/p < 1
\end{equation}
for every $w \in {\bf Z}$.  It follows from (\ref{y = ... = a c lim_{n
    to infty} sum_{j = 0}^n (p w)^j}) that $y$ can be approximated by
integers with respect to the $p$-adic metric, as desired.

        The set ${\bf Z}_p$ of \emph{$p$-adic integers}\index{p-adic 
integers@$p$-adic integers}\index{Z_p@${\bf Z}_p$} is defined by
\begin{equation}
\label{{bf Z}_p = {x in {bf Q}_p : |x|_p le 1}}
        {\bf Z}_p = \{x \in {\bf Q}_p : |x|_p \le 1\},
\end{equation}
which is the closed unit ball $\overline{B}(0, 1)$ in ${\bf Q}_p$ with
respect to the $p$-adic metric.  Let us check that
\begin{equation}
\label{{bf Z}_p = overline{bf Z}}
        {\bf Z}_p = \overline{\bf Z},
\end{equation}
where $\overline{\bf Z}$ is the closure of ${\bf Z}$ in ${\bf Q}_p$,
with respect to the $p$-adic metric.  Of course, ${\bf Z} \subseteq
{\bf Z}_p$, by (\ref{|x|_p le 1}), which implies that $\overline{\bf
  Z} \subseteq {\bf Z}_p$, because ${\bf Z}_p$ is a closed set in
${\bf Q}_p$, by construction.  Now let $x \in {\bf Z}_p$ be given, and
let us show that $x$ can be approximated by ordinary integers with
respect to the $p$-adic metric.  Remember that ${\bf Q}$ is dense in
${\bf Q}_p$, so that $x$ can be approximated by rational numbers with
respect to the $p$-adic metric.  If $y \in {\bf Q}$ and $|x - y|_p \le
1$, then
\begin{equation}
\label{|y|_p le max (|y - x|_p, |x|_p) le 1}
        |y|_p \le \max (|y - x|_p, |x|_p) \le 1,
\end{equation}
by the ultrametric version of the triangle inequality, so that $y \in
{\bf Q} \cap {\bf Z}_p$.  Thus $x$ can actually be approximated by
elements of ${\bf Q} \cap {\bf Z}_p$ with respect to the $p$-adic
metric.  The argument in the preceding paragraph says exactly that
elements of ${\bf Q} \cap {\bf Z}_p$ can be approximated by ordinary
integers with respect to the $p$-adic metric, which implies that $x$
can be approximated by ordinary integers with respect to the $p$-adic
metric, as desired.

        Note that ${\bf Z}_p$ is a subgroup of ${\bf Q}_p$ with respect
to addition, and in fact a subring of ${\bf Q}_p$ with respect to
addition and multiplication.  Put
\begin{equation}
\label{p^j {bf Z} = {p^j x : x in {bf Z}}}
        p^j \, {\bf Z} = \{p^j \, x : x \in {\bf Z}\}
\end{equation}
and
\begin{equation}
\label{p^j {bf Z}_p = {p^j x : x in {bf Z}_p} = ...}
        p^j \, {\bf Z}_p = \{p^j \, x : x \in {\bf Z}_p\}
                         = \{y \in {\bf Q}_p : |y|_p \le p^{-j}\}
\end{equation}
for each $j \in {\bf Z}$, so that $p^j \, {\bf Z}_p$ is the same as
the closure of $p^j \, {\bf Z}$ in ${\bf Q}_p$ with respect to the
$p$-adic metric, by (\ref{{bf Z}_p = overline{bf Z}}).  Of course,
$p^j \, {\bf Z}$ is a subgroup of ${\bf Q}$ with respect to addition
for every $j \in {\bf Z}$, and $p^j \, {\bf Z}$ is an ideal in ${\bf
  Z}$ when $j \ge 0$.  Similarly, $p^j \, {\bf Z}_p$ is a subgroup of
${\bf Q}_p$ with respect to addition for every $j \in {\bf Z}$, and
$p^j \, {\bf Z}_p$ is an ideal in ${\bf Z}_p$ when $j \ge 0$.  Thus
the quotients
\begin{equation}
\label{{bf Z} / p^j {bf Z}}
        {\bf Z} / p^j \, {\bf Z}
\end{equation}
and
\begin{equation}
\label{{bf Z}_p / p^j {bf Z}_p}
        {\bf Z}_p / p^j \, {\bf Z}_p
\end{equation}
are defined as commutative rings when $j \ge 0$.  The obvious
inclusion of ${\bf Z}$ into ${\bf Z}_p$ leads to a ring homomorphism
from ${\bf Z}$ into (\ref{{bf Z}_p / p^j {bf Z}_p}), whose kernel is
\begin{equation}
\label{{bf Z} cap (p^j {bf Z}_p)}
        {\bf Z} \cap (p^j \, {\bf Z}_p).
\end{equation}
This is the same as $p^j \, {\bf Z}$, because $x \in {\bf Z}$
satisfies $|x|_p \le p^{-j}$ if and only if $x \in p^j \, {\bf Z}$.
Every element of ${\bf Z}_p$ can be expressed as a sum of elements of
${\bf Z}$ and $p^j \, {\bf Z}_p$, because ${\bf Z}_p$ is the closure
of ${\bf Z}$ in ${\bf Q}_p$.  It follows that the homomorphism from
${\bf Z}$ into (\ref{{bf Z}_p / p^j {bf Z}_p}) mentioned earlier is
surjective, and leads to a ring isomorphism from (\ref{{bf Z} / p^j
  {bf Z}}) onto (\ref{{bf Z}_p / p^j {bf Z}_p}).  In particular,
(\ref{{bf Z}_p / p^j {bf Z}_p}) has exactly $p^j$ elements for each $j
\ge 0$.

        This implies that ${\bf Z}_p$ is totally bounded with respect
to the $p$-adic metric, in the sense that ${\bf Z}_p$ can be covered
by finitely many balls of arbitrarily small radius.  It follows that
${\bf Z}_p$ is compact with respect to the topology determined by the
$p$-adic metric, because ${\bf Q}_p$ is complete, and ${\bf Z}_p$ is a
closed subset of ${\bf Q}_p$.  Similarly, $p^j \, {\bf Z}_p$ is
compact for every $j \in {\bf Z}$, which can either be derived by an
analogous argument, or using continuity of multiplication on ${\bf
  Q}_p$.  If $E$ is a bounded subset of ${\bf Q}_p$ with respect to
the $p$-adic metric, then $E \subseteq p^j \, {\bf Z}_p$ when $-j$ is
sufficiently large.  If $E$ is closed and bounded, then $E$ is
compact, because closed subsets of compact sets are compact.

\section{Formal series}
\label{formal series}

        Let $k_0$ be a field, and let $T$ be an indeterminate.  In this
section, we shall be interested in formal series\index{formal series}
of the form
\begin{equation}
\label{f(T) = sum_{j = n}^infty f_j T^j}
        f(T) = \sum_{j = n}^\infty f_j \, T^j,
\end{equation}
where $n \in {\bf Z}$ and $f_j \in k_0$ for each $j \ge n$.  More
precisely, one can ask that $f_j \in k_0$ be defined for every $j \in
{\bf Z}$, with the condition that $f_j = 0$ for all but finitely many
negative integers $j$.  This permits the space
$k_0((T))$\index{k_0((T))@$k_0((T))$} of all such formal series in $T$
with coefficients in $k_0$ to be identified with the collection of
functions from ${\bf Z}$ into $k_0$ that are equal to $0$ at all but
at most finitely many negative integers.  As on p27 of \cite{cas}, it
is convenient to use the notation
\begin{equation}
\label{f(T) = sum_{j gg -infty} f_j T^j}
        f(T) = \sum_{j \gg -\infty} f_j \, T^j
\end{equation}
instead of (\ref{f(T) = sum_{j = n}^infty f_j T^j}), to indicate that
$f_j = 0$ for all but finitely many negative integers $j$, without
specifying $n \in {\bf Z}$ as in (\ref{f(T) = sum_{j = n}^infty f_j
  T^j}).

        Of course, the space of all functions from ${\bf Z}$ into $k_0$
is a vector space over $k_0$ with respect to pointwise addition and
scalar multiplication.  The space of such functions that are equal to
$0$ at all but finitely many negative integers is a linear subspace
of this vector space, and hence a vector space over $k_0$ too.  Thus
$k_0((T))$ is a vector space over $k_0$ in a natural way, where the 
vector space operations correspond to termwise addition and scalar
multiplication of formal series as in (\ref{f(T) = sum_{j = n}^infty f_j T^j}).

        It is easy to define multiplication of formal series as in
(\ref{f(T) = sum_{j = n}^infty f_j T^j}), where
\begin{equation}
\label{T^j T^l = T^{j + l}}
        T^j \, T^l = T^{j + l}
\end{equation}
for all $j, l \in {\bf Z}$.  To be more precise, suppose that $f(T)$
is as in (\ref{f(T) = sum_{j gg -infty} f_j T^j}), and that
\begin{equation}
\label{g(T) = sum_{l gg -infty} g_l T^l}
        g(T) = \sum_{l \gg -\infty} g_l \, T^l
\end{equation}
is another element of $k_0((T))$.  Under these conditions, their
product $f(T) \, g(T)$ is given by
\begin{equation}
\label{f(T) g(T) = sum_{r gg -infty} (f g)_r T^r}
        f(T) \, g(T) = \sum_{r \gg -\infty} (f \, g)_r \, T^r,
\end{equation}
where
\begin{equation}
\label{(f g)_r = sum_{j + l = r} f_j g_l}
        (f \, g)_r = \sum_{j + l = r} f_j \, g_l
\end{equation}
for each $r \in {\bf Z}$.  The sum in (\ref{(f g)_r = sum_{j + l = r}
  f_j g_l}) is taken over all $j, l \in {\bf Z}$ with $j + l = r$, and
it is easy to see that all but at most finitely many terms in this sum
are equal to $0$, because $f_j = 0$ and $g_l = 0$ for all but finitely
many negative integers $j$, $l$, respectively.  Similarly, (\ref{(f
  g)_r = sum_{j + l = r} f_j g_l}) is equal to $0$ for all but
finitely many negative integers $r$, so that (\ref{f(T) g(T) = sum_{r
    gg -infty} (f g)_r T^r}) defines an element of $k_0((T))$.  One
can check that $k_0((T))$ is a commutative ring with respect to this
definition of multiplication, and in fact an algebra over $k_0$.  One
can also identify $k_0$ with the subalgebra of $k_0((T))$ consisting
of series for which only the constant term may be nonzero.  In
particular, one can identify the multiplicative identity element $1$
of $k_0$ with $T^0$, which is the multiplicative identity element of
$k_0((T))$.

        Let $k_0[[T]]$\index{k_0[[T]]@$k_0[[T]]$} be the subset of 
$k_0((T))$ consisting of series of the form
\begin{equation}
\label{f(T) = sum_{j = 0}^infty f_j T^j}
        f(T) = \sum_{j = 0}^\infty f_j \, T^j,
\end{equation}
so that $f_j = 0$ when $j < 0$.  This is a subalgebra of $k_0((T))$
that contains $k_0$, which is the usual algebra of formal power series
with coefficients in $k_0$.  The elements of $k_0((T))$ may be
described as formal Laurent series with coefficients in $k_0$ and a
pole of finite order.  Note that every element of $k_0((T))$ can be
expressed as $T^n \, f(T)$ for some $n \in {\bf Z}$ and $f(T) \in
k_0[[T]]$.  The algebra $k_0[T]$\index{k_0[T]@$k_0[T]$} of formal
polynomials in $T$ with coefficients in $k_0$ may be identified with
the subalgebra of $k_0[[T]]$ consisting of formal power series $f(T)$
such that $f_j = 0$ for all but finitely many $j$.

        If $a(T)$ is any element of $k_0[[T]]$, then $a(T)^n \in k_0[[T]]$
for each positive integer $n$, and we interpret $a(T)^0$ as being the
series with only the constant term $1$.  This permits one to define
\begin{equation}
\label{sum_{n = 0}^infty T^n a(T)^n}
        \sum_{n = 0}^\infty T^n \, a(T)^n
\end{equation}
as an element of $k_0[[T]]$, since for each nonnegative integer $j$,
the coefficient of $T^j$ in $T^n \, a(T)^n$ is equal to $0$ when $j < n$.
Observe that
\begin{eqnarray}
\label{(1 - T a(T)) sum_{n = 0}^infty T^n a(T)^n = ...}
 (1 - T \, a(T)) \, \sum_{n = 0}^\infty T^n \, a(T)^n
          & = & \sum_{n = 0}^\infty T^n \, a(T)^n
                  - T \, a(T) \, \sum_{n = 0}^\infty T^n \, a(T)^n \nonumber \\
 & = & \sum_{n = 0}^\infty T^n \, a(T)^n - \sum_{n = 1}^\infty T^n \, a(T)^n = 1.
\end{eqnarray}
This shows that $1 - T \, a(T)$ has a multiplicative inverse in
$k_0[[T]]$ under these conditions, which is given by (\ref{sum_{n =
    0}^infty T^n a(T)^n}).

        If $f(T)$ is a nonzero element of $k_0((T))$, then $f(T)$
can be expressed as
\begin{equation}
\label{f(T) = c T^n (1 - T a(T))}
        f(T) = c \, T^n \, (1 - T \, a(T))
\end{equation}
for some $c \in k_0$ with $c \ne 0$, $n \in {\bf Z}$, and $a(T) \in
k_0[[T]]$.  It follows that $f(T)$ has a multiplicative inverse in
$k_0((T)$, which is given by
\begin{equation}
\label{f(T)^{-1} = c^{-1} T^{-n} (1 - T a(T))^{-1}}
        f(T)^{-1} = c^{-1} \, T^{-n} \, (1 - T \, a(T))^{-1},
\end{equation}
where $(1 - T \, a(T))^{-1}$ is as in the previous paragraph.  Thus
$k_0((T))$ is a field, which contains $k_0$ as a subfield.

        Let $r$ be a positive real number strictly less than $1$,
and let us use $r$ to define an absolute value function on $k_0((T))$.
Let $f(T) \in k_0((T))$ be given, and put $|f(T)| = 0$ when $f(T) =
0$.  Otherwise, if $f(T) \ne 0$, then there is a unique integer $n =
n(f)$ such that
\begin{equation}
\label{f_n ne 0 and f_j = 0 for every j < n}
        f_n \ne 0 \hbox{ and } f_j = 0 \hbox{ for every } j < n,
\end{equation}
in which case we put
\begin{equation}
\label{|f(T)| = r^n}
        |f(T)| = r^n.
\end{equation}
It is easy to see that
\begin{equation}
\label{|f(T) g(T)| = |f(T)| |g(T)|}
        |f(T) \, g(T)| = |f(T)| \, |g(T)|
\end{equation}
for every $f(T), g(T) \in k_0((T))$, because $n(f(T) \, g(T)) =
n(f(T)) + n(g(T))$ when $f(T), g(T) \ne 0$.  Similarly,
\begin{equation}
\label{|f(T) + g(T)| le max(|f(T)|, |g(T)|)}
        |f(T) + g(T)| \le \max(|f(T)|, |g(T)|)
\end{equation}
for every $f(T), g(T) \in k_0((T))$, because $n(f(T) + g(T)) \ge
\min(n(f(T)), n(g(T)))$ when $f(T), g(T) \ne 0$.

        Thus $|f(T)|$ defines an ultrametric absolute value function on
$k_0((T))$ for each $r \in (0, 1)$.  If $a$ is a positive real number,
then
\begin{equation}
\label{|f(T)|^a}
        |f(T)|^a
\end{equation}
is the same as the ultrametric absolute value function associated in
this way to $r^a$ instead of $r$.  In particular, the topology on
$k_0((T))$ determined by the ultrametric corresponding to $|f(T)|$ as
in (\ref{d(x, y) = |x - y|, 2}) does not depend on the choice of $r
\in (0, 1)$.  If one were to take $r = 1$ in the definition of
$|f(T)|$, then one would get the trivial absolute value function on
$k_0((T))$.  Note that the restriction of $|f(T)|$ to $k_0$ is the
trivial absolute value function on $k_0$ for every $r \in (0, 1)$.

        By construction,
\begin{equation}
\label{k_0[[T]] = {f(T) in k_0((T)) : |f(T)| le 1}}
        k_0[[T]] = \{f(T) \in k_0((T)) : |f(T)| \le 1\}
\end{equation}
for every $r \in (0, 1)$, which is the closed unit ball in $k_0((T))$
with respect to $|f(T)|$.  Of course, $k_0[[T]]$ can be identified
with a Cartesian product of a sequence of copies of $k_0$, indexed by
the nonnegative integers.  One can check that the topology on
$k_0[[T]]$ determined by the ultrametric (\ref{d(x, y) = |x - y|, 2})
corresponding to $|f(T)|$ is the same as the product topology
associated to the discrete topology on each copy of $k_0$ in the
Cartesian product.  In particular, this topology does not depend on
the choice of $r \in (0, 1)$, as before.  It is easy to see that
$k_0[T]$ is dense in $k_0[[T]]$ with respect to this topology.

        Similarly,
\begin{eqnarray}
\label{T^n k_0[[T]] = ... = {f(T) in k_0((T)) : |f(T)| le r^n}}
   T^n \, k_0[[T]] & = & \{T^n \, f(T) : f(T) \in k_0[[T]]\} \\
                   & = & \{f(T) \in k_0((T)) : |f(T)| \le r^n\} \nonumber
\end{eqnarray}
for each $n \in {\bf Z}$.  Thus every bounded subset of $k_0((T))$
with respect to $|f(T)|$ is contained in (\ref{T^n k_0[[T]] = ... =
  {f(T) in k_0((T)) : |f(T)| le r^n}}) for some $n \in {\bf Z}$.  In
particular, every Cauchy sequence of elements of $k_0((T))$ with
respect to the ultrametric (\ref{d(x, y) = |x - y|, 2}) corresponding
to $|f(T)|$ is contained in (\ref{T^n k_0[[T]] = ... = {f(T) in
    k_0((T)) : |f(T)| le r^n}}) for some $n \in {\bf Z}$.  One can
check that the coefficients of $T^j$ of the terms of such a Cauchy
sequence are eventually constant for each $j$, and hence that the
Cauchy sequence converges in $k_0((T))$.  This implies that $k_0((T))$
is complete with respect to the ultrametric (\ref{d(x, y) = |x - y|,
  2}) corresponding to $|f(T)|$, for each $r \in (0, 1)$.

        If $k_0$ has only finitely many elements, then $k_0[[T]]$ is
compact with respect to the topology determined by the ultrametric
(\ref{d(x, y) = |x - y|, 2}) corresponding to $|f(T)|$ for each $r \in
(0, 1)$, because $k_0[[T]]$ is topologically equivalent to a product
of finite sets.  In this case, (\ref{T^n k_0[[T]] = ... = {f(T) in
    k_0((T)) : |f(T)| le r^n}}) is also compact for each $n \in {\bf
  Z}$, either by an analogous argument, or using continuity of
multiplication on $k_0((T))$.  If $E \subseteq k_0((T))$ is bounded
with respect to $|f(T)|$, then $E$ is contained in (\ref{T^n k_0[[T]]
  = ... = {f(T) in k_0((T)) : |f(T)| le r^n}}) for some $n \in {\bf
  Z}$, as in the preceding paragraph.  If $E$ is closed and bounded,
then it follows that $E$ is compact, since closed subsets of compact
sets are compact.

\section{Discrete absolute value functions}
\label{discrete absolute value functions}

        Let $k$ be a field, and let $|\cdot|$ be a  quasimetric absolute
value function on $k$.  Thus
\begin{equation}
\label{{|x| : x in k, x ne 0}}
        \{|x| : x \in k, \, x \ne 0\}
\end{equation}
is a subgroup of the multiplicative group ${\bf R}_+$\index{R_+@${\bf
    R}_+$} of positive real numbers.  Suppose that there is a positive
real number $\rho < 1$ such that for each $x \in k$ with $|x| < 1$, we
have that
\begin{equation}
\label{|x| le rho}
        |x| \le \rho.
\end{equation}
This implies that for each $x \in k$ with $|x| > 1$, we have that
\begin{equation}
\label{|x| ge 1/rho}
        |x| \ge 1/\rho,
\end{equation}
by applying the previous statement to $1/x$.  If $y, z \in k$ satisfy
$|y| < |z|$, then we get that
\begin{equation}
\label{|y| le rho |z|}
        |y| \le \rho \, |z|,
\end{equation}
by applying (\ref{{|x| : x in k, x ne 0}}) to $x = y / z$.

        It follows that $1$ is not a limit point of 
(\ref{{|x| : x in k, x ne 0}}) with respect to the standard metric on
${\bf R}$ under these conditions.  Conversely, if $1$ is not a limit
point of (\ref{{|x| : x in k, x ne 0}}) with respect to the standard
metric on ${\bf R}$, then there is a $\rho \in (0, 1)$ with the
property described in the preceding paragraph.  In this case, one can
check that (\ref{{|x| : x in k, x ne 0}}) is a discrete subgroup of
${\bf R}_+$, in the sense that (\ref{{|x| : x in k, x ne 0}}) has no
limit points in ${\bf R}_+$.  More precisely, if $t \in {\bf R}_+$ is
a limit point of (\ref{{|x| : x in k, x ne 0}}) in ${\bf R}_+$, then
there exist $y, z \in k$ such that $|y| < |z|$ and $|y|$, $|z|$ are
arbitrarily close to $t$ with respect to the standard metric on ${\bf
  R}$.  This implies that $|y| / |z|$ is arbitrarily close to $1$,
which would contradict (\ref{|y| le rho |z|}).  Of course, $0$ is
always a limit point of (\ref{{|x| : x in k, x ne 0}}) in ${\bf R}$
when $|\cdot|$ is nontrivial on $k$.  If $|\cdot|$ satisfies the
condition described in the previous paragraph, then $|\cdot|$ is said
to be \emph{discrete}\index{discrete absolute value functions} on $k$.

        The trivial absolute value function is discrete on any field $k$.
In this case, one can take $\rho = 0$ in (\ref{|x| le rho}), and
(\ref{|x| ge 1/rho}) is vacuous.  If $p$ is a prime number, then the
$p$-adic absolute value function on $k = {\bf Q}$ is discrete, and one
can take $\rho = 1/p$ in (\ref{|x| le rho}).  However, the standard
absolute value function on $k = {\bf Q}$ is not discrete.  If $k_0$ is
a field, $T$ is an indeterminate, $k$ is the field $k_0((T))$ of
formal series in $T$ with coefficients in $k_0$ discussed in the
previous section, and $|f(T)|$ is the absolute value function on
$k_0((T))$ associated to some $r \in (0, 1)$ as in (\ref{|f(T)| =
  r^n}), then $|f(T)|$ is discrete on $k_0((T))$, and one can take
$\rho = r$ in (\ref{|x| le rho}).  If a quasimetric absolute value
function $|\cdot|$ on a field $k$ is discrete, then $|\cdot|^a$ is
discrete on $k$ for every positive real number $a$.  More precisely,
if $|\cdot|$ satisfies (\ref{|x| le rho}) for some $\rho \in (0, 1)$,
then $|\cdot|^a$ satisfies the analogous condition with $\rho^a$
instead of $\rho$.  If $|\cdot|$ is discrete on $k$, then the natural
extension of $|\cdot|$ to the corresponding completion of $k$ is also
discrete, and one can use the same value of $\rho$ in (\ref{|x| le
  rho}) on the completion of $k$.

        Suppose for the moment that $k$ has characteristic $0$, so that
there is a natural embedding of ${\bf Q}$ into $k$.  Thus a
quasimetric absolute value function $|\cdot|$ on $k$ leads to a
quasimetric absolute value function on ${\bf Q}$, as in (\ref{|r cdot
  1|}) in Section \ref{another refinement}.  If this quasimetric
absolute value function on ${\bf Q}$ is archimedian, then it is equal
to a positive power of the standard absolute value function on ${\bf
  Q}$, as mentioned in Section \ref{some more refinements}.  This
implies that $|\cdot|$ is not discrete on $k$, because any positive
power of the standard absolute value function on ${\bf Q}$ is not
discrete.

        If $|\cdot|$ is a discrete quasimetric absolute value function
on any field $k$, then $|\cdot|$ is an ultrametric absolute value
function on $k$.  More precisely, if $k$ does not have characteristic
$0$, then every quasimetric absolute value function on $k$ is an
ultrametric absolute value function.  This follows from the discussion
in Section \ref{some more refinements}, as mentioned at the end of
Section \ref{another refinement}.  If $k$ has characteristic $0$, then
there is a natural embedding of ${\bf Q}$ into $k$, which leads to a
quasimetric absolute value function on ${\bf Q}$, as before.  If this
quasimetric absolute value function on ${\bf Q}$ is non-archimedian,
then $|\cdot|$ is non-archimedian on $k$.  This implies that $|\cdot|$
is an ultrametric absolute value function on $k$ again, as in Section
\ref{some more refinements}.  The remaining possibility is that $k$
has characteristic $0$, and that the corresponding absolute value
function on ${\bf Q}$ is archimedian.  In this case, $|\cdot|$ is not
discrete on $k$, as in the preceding paragraph.

        Let $|\cdot|$ be a quasimetric absolute value function on a
field $k$ again, and put
\begin{equation}
\label{rho_1 = sup {|x| : x in k, |x| < 1}}
        \rho_1 = \sup \{|x| : x \in k, \, |x| < 1\},
\end{equation}
so that $0 \le \rho_1 \le 1$.  Thus $|\cdot|$ is discrete on $k$ if
and only if $\rho_1 < 1$.  It is easy to see that $\rho_1 = 0$ if and
only if $|\cdot|$ is the trivial absolute value function on $k$.  Let
us suppose from now on in this section that $|\cdot|$ is nontrivial
and discrete on $k$, so that $0 < \rho_1 < 1$.  Under these
conditions, one can check that the supremum in (\ref{rho_1 = sup {|x|
    : x in k, |x| < 1}}) is attained, which is to say that there is an
$x_1 \in k$ such that
\begin{equation}
\label{|x_1| = rho_1}
        |x_1| = \rho_1.
\end{equation}
This uses the fact that (\ref{{|x| : x in k, x ne 0}}) is a discrete
subgroup of ${\bf R}_+$.  It follows that
\begin{equation}
\label{|x_1^j| = |x_1|^j = rho_1^j}
        |x_1^j| = |x_1|^j = \rho_1^j
\end{equation}
for each $j \in {\bf Z}$, so that $\rho_1^j$ is an element of
(\ref{{|x| : x in k, x ne 0}}) for each $j \in {\bf Z}$.

        Conversely, let $w \in k$ with $w \ne 0$ be given, and let us
show that
\begin{equation}
\label{|w| = rho^j}
        |w| = \rho^j
\end{equation}
for some $j \in {\bf Z}$.  Of course, there is a $j \in {\bf Z}$ such
that
\begin{equation}
\label{rho_1^{j+ 1} < |w| le rho_1^j}
        \rho_1^{j + 1} < |w| \le \rho_1^j,
\end{equation}
because $|w| > 0$.  Suppose for the sake of a contradiction that
$|w| < \rho_1^j$, and put $u = w / x_1^j$, so that
\begin{equation}
\label{|u| = |w| / |x_1|^j = |w| / rho_1^j < 1}
        |u| = |w| / |x_1|^j = |w| / \rho_1^j < 1.
\end{equation}
The definition (\ref{rho_1 = sup {|x| : x in k, |x| < 1}}) of $\rho_1$
implies that $|u| \le \rho_1$, and hence that
\begin{equation}
\label{|w| = |u x_1^j| le rho_1^{j + 1}}
        |w| = |u \, x_1^j| \le \rho_1^{j + 1},
\end{equation}
contradicting (\ref{rho_1^{j+ 1} < |w| le rho_1^j}).  Thus we get
(\ref{|w| = rho^j}), which means that (\ref{{|x| : x in k, x ne 0}})
consists of exactly the integer powers of $\rho_1$ under these
conditions.

\section{The ultrametric case}
\label{ultrametric case}

        Let $k$ be a field, and let $|\cdot|$ be an ultrametric absolute
value function on $k$.  It is easy to see that the closed unit ball
\begin{equation}
\label{overline{B}(0, 1) = {x in k : |x| le 1}}
        \overline{B}(0, 1) = \{x \in k : |x| \le 1\}
\end{equation}
in $k$ with respect to $|\cdot|$ is a subring of $k$ under these
conditions.  Note that $\overline{B}(0, 1)$ contains the
multiplicative identity element $1$ in $k$, by (\ref{|1| = 1}).  If $x
\in k$ satisfies $|x| = 1$, then $x$ and $x^{-1}$ are in
$\overline{B}(0, 1)$, so that $x$ is invertible in $\overline{B}(0,
1)$ as a ring.  Conversely, if $x$ is an element of $\overline{B}(0,
1)$ that is invertible in $\overline{B}(0, 1)$, then $|x| \le 1$, $x
\ne 0$, and $|x^{-1}| \le 1$, which implies that $|x| = 1$.

        Similarly, the open unit ball
\begin{equation}
\label{B(0, 1) = {x in k : |x| < 1}}
        B(0, 1) = \{x \in k : |x| < 1\}
\end{equation}
in $k$ with respect to $|\cdot|$ is an ideal in $\overline{B}(0, 1)$.
More precisely, $B(0, 1)$ is a maximal ideal in $\overline{B}(0, 1)$,
by the remarks in the previous paragraph.  This implies that the
quotient ring
\begin{equation}
\label{overline{B}(0, 1) / B(0, 1)}
        \overline{B}(0, 1) / B(0, 1)
\end{equation}
is a field, which also follows more directly from the characterization
of invertible elements in $\overline{B}(0, 1)$.  Of course, the
multiplicative identity element $1$ in $k$ is in $\overline{B}(0, 1)$
and not in $B(0, 1)$, by (\ref{|1| = 1}), so that its image in the
quotient (\ref{overline{B}(0, 1) / B(0, 1)}) is nonzero.  The quotient
(\ref{overline{B}(0, 1) / B(0, 1)}) is known as the \emph{residue
  field}\index{residue field} associated to $|\cdot|$ on $k$.

        If $|\cdot|$ is the trivial absolute value function on any
field $k$, then $\overline{B}(0, 1)$ is equal to $k$, $B(0, 1)$ is the
trivial ideal $\{0\}$, and hence (\ref{overline{B}(0, 1) / B(0, 1)})
is isomorphic to $k$.  If $k = {\bf Q}$ equipped with the $p$-adic
absolute value function for some prime number $p$, then
$\overline{B}(0, 1)$ is the same as the ring ${\bf Z}_p$ of $p$-adic
integers, and $B(0, 1)$ reduces to $p \, {\bf Z}_p$.  It follows that
the residue field (\ref{overline{B}(0, 1) / B(0, 1)}) is isomorphic to
${\bf Z} / p {\bf Z}$ in this case, as in Section \ref{p-adic
  integers}.  If $|\cdot|$ is an ultrametric absolute value function
on a field $k$ and $a$ is a positive real number, then $|\cdot|^a$ is
also an ultrametric absolute value function on $k$, as in Section
\ref{definitions, examples}.  Clearly $\overline{B}(0, 1)$ and $B(0,
1)$ are the same for $|\cdot|^a$ as for $|\cdot|$, so that the residue
field (\ref{overline{B}(0, 1) / B(0, 1)}) is the same for $|\cdot|^a$
as for $|\cdot|$ too.

        Suppose that $k$ is a field of characteristic $p$ for some
prime number $p$, so that $p \cdot 1 = 0$ in $k$.  If $|\cdot|$ is an
ultrametric absolute value function on $k$, then $n \cdot 1$ is an
element of $\overline{B}(0, 1)$ for every positive integer $n$, whose
image in (\ref{overline{B}(0, 1) / B(0, 1)}) is the same as $n \cdot
1$ in (\ref{overline{B}(0, 1) / B(0, 1)}).  It follows that $p \cdot 1
= 0$ in (\ref{overline{B}(0, 1) / B(0, 1)}) as well, so that the
residue field also has characteristic $p$ under these conditions.

        Now let $k_0$ be a field, let $T$ be an indeterminate, and let
$k$ be the corresponding field $k_0((T))$ of formal series in $T$ with
coefficients in $k_0$ discussed in Section \ref{formal series}.  Also
let $|f(T)|$ be the absolute value function on $k_0((T))$ associated
to some $r \in (0, 1)$ as in (\ref{|f(T)| = r^n}).  In this situation,
$\overline{B}(0, 1)$ is the ring $k_0[[T]]$ of formal power series in
$T$ with coefficients in $k_0$, and $B(0, 1)$ is equal to $T \,
k_0[[T]]$.  Note that there is a natural homomorphism from $k_0[[T]]$
onto $k_0$, which sends a formal power series to its constant term.
The kernel of this homomorphism is equal to $T \, k_0[[T]]$, which
leads to an isomorphism from the residue field (\ref{overline{B}(0, 1)
  / B(0, 1)}) onto $k_0$.

        Let $|\cdot|$ be an ultrametric absolute value function on a
field $k$ again, and let $k_1$ be a subfield of $k$.  Thus the
restriction of $|\cdot|$ to $k_1$ is an ultrametric absolute value
function on $k_1$, and the corresponding open and closed unit balls in
$k_1$ are the same as the intersections of their counterparts in $k$
with $k_1$.  This leads to a natural injective homomorphism from the
residue field associated to $k_1$ into the residue field associated to
$k$.  If $k_1$ is dense in $k$ with respect to the ultrametric
corresponding to $|\cdot|$, then one can check that this homomorphism
is surjective.  In particular, the residue field associated to the
completion of $k$ with respect to $|\cdot|$ is isomorphic to the
residue field associated to $k$.

        Note that $\overline{B}(0, 1)$ can be expressed as the union of
a family of pairwise-disjoint open balls in $k$ of radius $1$, which
are the cosets of $B(0, 1)$ in $\overline{B}(0, 1)$.  If the residue
field (\ref{overline{B}(0, 1) / B(0, 1)}) has only finitely many
elements, then $\overline{B}(0, 1)$ can be expressed as the union of
finitely many pairwise-disjoint open balls of radius $1$.  Conversely,
if $\overline{B}(0, 1)$ can be covered by finitely many open balls of
radius $1$, then one can check that the residue field
(\ref{overline{B}(0, 1) / B(0, 1)}) has only finitely many elements.
Of course, any open ball of radius $1$ in $k$ that intersects
$\overline{B}(0, 1)$ is contained in $\overline{B}(0, 1)$, because of
the ultrametric version of the triangle inequality.  Thus any open
ball in $k$ of radius $1$ that intersects $\overline{B}(0, 1)$ is
actually a coset of $B(0, 1)$ in $\overline{B}(0, 1)$.

        Suppose that $|\cdot|$ is discrete on $k$, as in the previous
section.  Let $\rho_1$ be as in (\ref{rho_1 = sup {|x| : x in k, |x| <
    1}}), so that $0 \le \rho_1 < 1$, and $B(0, 1)$ is the same as
\begin{equation}
\label{overline{B}(0, rho_1) = {x in k : |x| le rho_1}}
        \overline{B}(0, \rho_1) = \{x \in k : |x| \le \rho_1\}.
\end{equation}
Suppose also that $|\cdot|$ is nontrivial on $k$, so that $\rho_1 >
0$.  Remember that the nonzero values of $|\cdot|$ on $k$ are the same
as the integer powers of $\rho_1$ under these conditions, as in the
previous section.  If the residue field (\ref{overline{B}(0, 1) / B(0,
  1)}) has only finitely many elements, then the closed unit ball in
$k$ can be covered by finitely many closed balls of radius $\rho_1$.
This implies that every closed ball in $k$ of radius $\rho_1^j$ for
some $j \in {\bf Z}$ can be covered by finitely many closed balls of
radius $\rho_1^{j + 1}$, using translations and dilations on $k$.
Repeating the process, it follows that every closed ball in $k$ can be
covered by finite many closed balls of arbitrarily small radius.  Thus
closed balls in $k$ are totally bounded under these conditions, with
respect to the ultrametric on $k$ corresponding to $|\cdot|$.  If $k$
is also complete, then it follows that closed balls in $k$ are
compact.  In this case, closed and bounded subsets of $k$ are compact,
because closed subsets of compact sets are compact.

        Conversely, let $|\cdot|$ be an ultrametric absolute value function
on a field $k$, and suppose that $\overline{B}(0, 1)$ is totally
bounded with respect to the corresponding ultrametric on $k$.  In
particular, this implies that $\overline{B}(0, 1)$ can be covered by
finitely many open balls of radius $1$, so that the residue field
(\ref{overline{B}(0, 1) / B(0, 1)}) has only finitely many elements.
We also get that $B(0, 1)$ can be covered by finitely many closed
balls of radius $1/2$, say, so that there are finitely many elements
$x_1, \ldots, x_n$ of $k$ such that
\begin{equation}
\label{B(0, 1) subseteq bigcup_{j = 1}^n overline{B}(x_j, 1/2)}
        B(0, 1) \subseteq \bigcup_{j = 1}^n \overline{B}(x_j, 1/2).
\end{equation}
We may as well suppose that
\begin{equation}
\label{overline{B}(x_j, 1/2) cap B(0, 1) ne emptyset}
        \overline{B}(x_j, 1/2) \cap B(0, 1) \ne \emptyset
\end{equation}
for each $j = 1, \ldots, n$, since otherwise $\overline{B}(x_j, 1/2)$
is not needed in the covering of $B(0, 1)$.  It follows that $x_j \in
B(0, 1)$ for each $j$, by the ultrametric version of the triangle
inequality, which is to say that $|x_j| < 1$ for each $j$.  If $x \in
k$ and $|x| < 1$, so that $x \in B(0, 1)$, then $x \in
\overline{B}(x_j, 1/2)$ for some $j$, by (\ref{B(0, 1) subseteq
  bigcup_{j = 1}^n overline{B}(x_j, 1/2)}).  Thus
\begin{equation}
\label{|x| le max(|x_1|, ldots, |x_n|, 1/2) < 1}
        |x| \le \max(|x_1|, \ldots, |x_n|, 1/2) < 1,
\end{equation}
by the ultrametric version of the triangle inequality again, which
shows that $|\cdot|$ is discrete on $k$ under these conditions.

\section{Haar measure, 2}
\label{haar measure, 2}

        If $p$ is a prime number, then the $p$-adic numbers ${\bf Q}_p$
form a locally compact commutative topological group with respect to
addition.  If $H$ is a choice of Haar measure on ${\bf Q}_p$, then
$H({\bf Z}_p)$ should be positive and finite, because ${\bf Z}_p$ is
nonempty, open, and compact.  It is convenient to normalize Haar
measure on ${\bf Q}_p$ so that
\begin{equation}
\label{H({bf Z}_p) = 1}
        H({\bf Z}_p) = 1,
\end{equation}
in which case $H$ is unique.  Let us check that
\begin{equation}
\label{H(p^j {bf Z}_p) = p^{-j}}
        H(p^j \, {\bf Z}_p) = p^{-j}
\end{equation}
for each integer $j$, using (\ref{H({bf Z}_p) = 1}) and the
translation-invariance of $H$.  If $j \ge 0$, then this follows from
the fact that ${\bf Z}_p$ can be expressed as the union of $p^j$
pairwise-disjoint translates of $p^j \, {\bf Z}_p$, because (\ref{{bf
    Z}_p / p^j {bf Z}_p}) has exactly $p^j$ elements.  If $j < 0$,
then $p^{-j} \, {\bf Z}_p$ can be expressed as the union of $p^{-j}$
pairwise-disjoint translates of ${\bf Z}_p$, since one can multiply by
$p^j$ to reduce to the previous situation.  Of course, (\ref{H(p^j {bf
    Z}_p) = p^{-j}}) implies that the Haar measure of every translate
of $p^j \, {\bf Z}_p$ should be equal to $p^{-j}$ too.  One way to get
the existence of Haar measure on ${\bf Q}_p$ is to first define a Haar
integral on $C_{com}({\bf Q}_p)$, as a Riemann integral.

        Let $k$ be a field, and let $|\cdot|$ be a quasimetric
absolute value function on $k$.  In particular, $k$ is a commutative
topological group with respect to addition, and using the topology on
$k$ determined by the quasimetric corresponding to $|\cdot|$.  If
$|\cdot|$ is the trivial absolute value function on $k$, then the
corresponding topology on $k$ is discrete, and counting measure on $k$
satisfies the requirements of Haar measure, as before.  If $k$ is
locally compact with respect to the topology determined by the
quasimetric associated to $|\cdot|$, then we have seen that $k$ is
complete, as in Section \ref{local compactness}.  

        If $k$ is complete and archimedian, then $k$ is isomorphic to
the real or complex numbers, and $|\cdot|$ corresponds to a positive
power of the standard absolute value function on ${\bf R}$ or ${\bf
  C}$, as in Sections \ref{complex numbers} and \ref{complex numbers,
  2}.  Thus $k$ is isomorphic to ${\bf R}$ or ${\bf C}$ as a
topological group with respect to addition, and where ${\bf R}$ and
${\bf C}$ are equipped with their standard topologies.  As in Section
\ref{haar measure}, Lebesgue measure on ${\bf R}$ satisfies the
requirements of Haar measure.  Similarly, ${\bf C}$ is isomorphic to
${\bf R}^2$ as a topological group with respect to addition, and
two-dimensional Lebesgue measure satisfies the requirements of Haar
measure on ${\bf C}$.

        Let us suppose from now on in this section that $|\cdot|$ is a
nontrivial ultrametric absolute value function on $k$.  If $k$ is
locally compact, then the closed unit ball in $k$ is compact and hence
totally bounded, as in Section \ref{local compactness}.  This implies
that the residue field (\ref{overline{B}(0, 1) / B(0, 1)}) is finite,
and that $|\cdot|$ is discrete on $k$, by the remarks at the end of
the previous section.  Conversely, if the residue field is finite, and
if $|\cdot|$ is discrete on $k$, then we have seen that closed balls
in $k$ are totally bounded, as in the previous section.  If $k$ is
also complete with respect to the ultrametric corresponding to
$|\cdot|$, then it follows that closed balls in $k$ are compact, so
that $k$ is locally compact.

        Continuing with these hypotheses, let $N$ be the number of
elements of the residue field (\ref{overline{B}(0, 1) / B(0, 1)}), and
let $\rho_1$ be as in (\ref{rho_1 = sup {|x| : x in k, |x| < 1}}).
Thus $\rho_1 < 1$, because $|\cdot|$ is discrete on $k$, and $\rho_1 >
0$, because $|\cdot|$ is nontrivial on $k$.  Remember that the nonzero
values of $|\cdot|$ on $k$ are the same as the integer powers of
$\rho_1$, as in Section \ref{discrete absolute value functions}.  We
have also seen that $\overline{B}(0, 1)$ can be expressed as the union
of $N$ pairwise-disjoint open balls of radius $1$, as in the previous
section.  Equivalently, $\overline{B}(0, 1)$ can be expressed as the
union of $N$ pairwise-disjoint closed balls of radius $\rho_1$, by the
definition of $\rho_1$.  This implies that every closed ball in $k$ of
radius $\rho_1^j$ for some $j \in {\bf Z}$ can be expressed as the
union of $N$ pairwise-disjoint closed balls of radius $\rho_1^{j +
  1}$.  It follows that the Haar measure of a closed ball in $k$ of
radius $\rho_1^j$ is equal to $N$ times the Haar measure of a closed
ball of radius $\rho_1^{j + 1}$, for any choice of Haar measure on
$k$, using invariance under translations.  If Haar measure on $k$ is
normalized so that
\begin{equation}
\label{H(overline{B}(0, 1)) = 1}
        H(\overline{B}(0, 1)) = 1,
\end{equation}
then we get that
\begin{equation}
\label{H(overline{B}(x, rho_1^j)) = N^{-j}}
        H(\overline{B}(x, \rho_1^j)) = N^{-j}
\end{equation}
for every $x \in k$ and $j \in {\bf Z}$.  As before, one can first
define a Haar integral on $C_{com}(k)$ as a Riemann integral, and then
get Haar measure on $k$ using the Riesz representation theorem.

        Let $k_0$ be a field, let $T$ be an indeterminate, and
consider the corresponding field $k_0((T))$ of formal series, as in
Section \ref{formal series}.  Also let $|f(T)|$ be the absolute value
function on $k_0((T))$ associated to some $r \in (0, 1)$ as in
(\ref{|f(T)| = r^n}).  Thus $|f(T)|$ is nontrivial and discrete on
$k$, and we have seen that $k_0((T))$ is complete with respect to
$|f(T)|$.  The corresponding residue field is isomorphic to $k_0$, as
in the previous section.  If $k_0$ has only finitely many elements,
then we have already mentioned in Section \ref{formal series} that
$k_0[[T]]$ is a compact subset of $k_0((T))$, and hence that closed
balls in $k_0((T))$ are compact.

        In particular, $k_0((T))$ is locally compact when $k_0$ has
only finitely many elements, in which case Haar measure on $k_0((T))$
can be analyzed as before.  Alternatively, $k_0[[T]]$ is a compact
commutative topological group with respect to addition when $k_0$
has only finitely many elements, using the topology induced by
the one on $k_0((T))$.  As a topological group with respect to addition,
$k_0[[T]]$ is isomorphic to a product of a sequence of copies of $k_0$,
where $k_0$ is considered as a topological group with respect to addition
and the discrete topology.  It is convenient to take Haar measure on $k_0$
to be counting measure divided by the total number of elements of $k_0$,
so that the Haar measure of $k_0$ is equal to $1$.  With this normalization,
Haar measure on $k_0[[T]]$ corresponds to a product measure on a product
of a sequence of copies of $k_0$, and one can get Haar measure on $k_0((T))$
from Haar measure on $k_0[[T]]$.

\section{Norms and ultranorms}
\label{norms, ultranorms}

        Let $k$ be a field, let $|\cdot|$ be an absolute value function
on $k$, and let $V$ be a vector space over $k$.  A nonnegative
real-valued function $N$ on $V$ is said to be a
\emph{norm}\index{norms} if it satisfies the following three
conditions: first, $N(v) = 0$ if and only if $v = 0$; second,
\begin{equation}
\label{N(t v) = |t| N(v)}
        N(t \, v) = |t| \, N(v)
\end{equation}
for every $t \in k$ and $v \in V$; and third,
\begin{equation}
\label{N(v + w) le N(v) + N(w)}
        N(v + w) \le N(v) + N(w)
\end{equation}
for every $v, w \in V$.  Under these conditions,
\begin{equation}
\label{d(v, w) = N(v - w)}
        d(v, w) = N(v - w)
\end{equation}
defines a metric on $V$.  More precisely, the fact that (\ref{d(v, w)
  = N(v - w)}) is symmetric in $v$ and $w$ follows from (\ref{N(t v) =
  |t| N(v)}) with $t = -1$, since $|-1| = 1$, as in (\ref{|-1| = 1}).

        A norm $N$ on $V$ is said to be an
\emph{ultranorm}\index{ultranorms} if
\begin{equation}
\label{N(v + w) le max(N(v), N(w))}
        N(v + w) \le \max(N(v), N(w))
\end{equation}
for every $v, w \in V$, which implies (\ref{N(v + w) le N(v) + N(w)}).
In this case, (\ref{d(v, w) = N(v - w)}) is an ultrametric on $V$.  If
$V \ne \{0\}$ and $N$ is an ultranorm on $V$, then it is easy to see
that $|\cdot|$ is an ultrametric absolute value function on $k$,
because of (\ref{N(t v) = |t| N(v)}).  If $|\cdot|$ is the trivial
absolute value function on $k$, then one can get an ultranorm $N$ on
$V$ by putting $N(v) = 1$ for every $v \in V$ with $v \ne 0$, and
$N(0) = 0$.  Let us call this the \emph{trivial
  ultranorm}\index{trivial ultranorm} on $V$, for which the
corresponding metric on $V$ is the discrete metric.

        If $N$ is any norm on $V$, then one can check that
\begin{equation}
\label{max(N(v) - N(w), N(w) - N(v)) le N(v - w)}
        \max(N(v) - N(w), N(w) - N(v)) \le N(v - w)
\end{equation}
for every $v, w \in V$, using (\ref{N(v + w) le N(v) + N(w)}).  This
is a special case of (\ref{|f(x) - f(y)| = max(f(x) - f(y), f(y) -
  f(x)) le C d_1(x, y)}) in Section \ref{lipschitz mappings}.  It
follows that $N$ is continuous as a real-valued function on $V$, with
respect to the topology on $V$ determined by the metric (\ref{d(v, w)
  = N(v - w)}) corresponding to $N$, and with respect to the standard
topology on ${\bf R}$.  If $N$ is an ultrametric on $V$, then
\begin{equation}
\label{N(v) = N(w)}
        N(v) = N(w)
\end{equation}
for every $v, w \in V$ such that $N(v - w) < N(v)$.  This is a special
case of (\ref{d(x, y) = d(x, z)}) in Section \ref{metrics,
  ultrametrics}.

        Let $n$ be a positive integer, and let $k^n$ be the set of $n$-tuples
of elements of $k$.  This is a vector space over $k$ with respect to
coordinatewise addition and scalar multiplication.  Put
\begin{equation}
\label{N_0(v) = max(|v_1|, ldots, |v_n|)}
        N_0(v) = \max(|v_1|, \ldots, |v_n|)
\end{equation}
for each $v = (v_1, \ldots, v_n) \in k^n$, which defines a norm on
$k^n$.  If $|\cdot|$ is an ultrametric absolute value function on $k$,
then $N_0$ is an ultranorm on $k^n$.  The topology on $k^n$ determined
by the corresponding metric as in (\ref{d(v, w) = N(v - w)}) is the
same as the product topology associated to the topology on $k$
determined by the metric (\ref{d(x, y) = |x - y|, 2}) corresponding to
$|\cdot|$.

        Let $e(1), \ldots, e(n)$ be the standard basis vectors in $k^n$,
so that the $l$th coordinate of $e(j)$ is equal to $1$ when $j = l$,
and to $0$ otherwise.  Thus each $v \in k^n$ can be expressed as
\begin{equation}
\label{v = sum_{j = 1}^n v_j e(j)}
        v = \sum_{j = 1}^n v_j \, e(j).
\end{equation}
If $N$ is any norm on $k^n$, then it follows that
\begin{equation}
\label{N(v) le ... le (sum_{j = 1}^n N(e(j))) N_0(v)}
        N(v) \le \sum_{j = 1}^n N(v_j \, e(j)) = \sum_{j = 1}^n |v_j| \, N(e(j))
              \le \Big(\sum_{j = 1}^n N(e(j))\Big) \, N_0(v)
\end{equation}
for every $v \in k^n$.  If $N$ is an ultranorm on $k^n$, then we get that
\begin{eqnarray}
\label{N(v) le ... le (max_{1 le j le n} N(e(j))) N_0(v)}
         N(v) \le \max_{1 \le j \le n} N(v_j \, e(j))
              & = & \max_{1 \le j \le n} (|v_j| \, N(e(j))) \\
              & \le & \Big(\max_{1 \le j \le n} N(e(j))\Big) \, N_0(v) \nonumber
\end{eqnarray}
for every $v \in k^n$.

        Let $V$ be any vector space over $k$, and let $N_1$, $N_2$
be norms on $V$.  Suppose that there is a positive real number $C$
such that
\begin{equation}
\label{N_1(v) le C N_2(v)}
        N_1(v) \le C \, N_2(v)
\end{equation}
for every $v \in V$.  This implies that the corresponding metrics as
in (\ref{d(v, w) = N(v - w)}) satisfy the analogous condition.  It
follows that every open set in $V$ with respect to the topology
determined by the metric associated to $N_1$ is also an open set with
respect to the topology determined by the metric associated to $N_2$.

        Conversely, suppose that every open set in $V$ with respect to
the topology determined by the metric associated to $N_1$ is also an
open set with respect to the topology determined by the metric
associated to $N_2$.  Let $r_1 > 0$ be given, and note that the open ball
\begin{equation}
\label{{v in V : N_1(v) < r_1}}
        \{v \in V : N_1(v) < r_1\}
\end{equation}
with respect to $N_1$ centered at $0$ with radius $r_1$ is an open set
with respect to the topology determined by the metric associated to
$N_1$.  The hypothesis that this also be an open set in $V$ with
respect to the topology determined by the metric associated to $N_2$
implies that there is an $r_2 > 0$ such that
\begin{equation}
\label{{v in V : N_2(v) < r_2}}
        \{v \in V : N_2(v) < r_2\}
\end{equation}
is contained in (\ref{{v in V : N_1(v) < r_1}}).  If $|\cdot|$ is not
the trivial absolute value function on $k$, then one can use this with
$r_1 = 1$ and the homogeneity property (\ref{N(t v) = |t| N(v)}) of
norms to get that there is a $C > 0$ such that (\ref{N_1(v) le C
  N_2(v)}) holds for every $v \in V$.

        Let $V$ be a vector space over $k$ again, and let $N$ be a
norm on $V$.  If $V$ is not complete as a metric space with respect to
the corresponding metric (\ref{d(v, w) = N(v - w)}), then one can
define its completion in the usual way.  The vector space operations
and norm $N$ can be extended to the completion, in such a way that the
completion becomes a vector space over $k$, the extension of the norm
to the completion is also a norm, and $V$ is a dense linear subspace
of the completion.  If $V$ is complete and $k$ is not complete with
respect to the metric associated to $|\cdot|$, then one can extend
scalar multiplication on $V$ to the completion of $k$, so that $V$
becomes a vector space over the completion of $k$.  One can also check
that $N$ will still be a norm on $V$ as a vector over the completion
of $k$, which is to say that (\ref{N(t v) = |t| N(v)}) holds for $t$
in the completion of $k$.

\section{Finite-dimensional vector spaces}
\label{finite-dimensional vector spaces}

        Let $k$ be a field, and let $|\cdot|$ be an absolute value
function on $k$.  Also let $n$ be a positive integer, and let
$N$ be a norm on $k^n$.  If $N_0$ is the norm (\ref{N_0(v) = max(|v_1|, 
ldots, |v_n|)}) on $k^n$, then we have seen that
\begin{equation}
\label{N(v) le C_1 N_0(v)}
        N(v) \le C_1 \, N_0(v)
\end{equation}
for some $C_1 > 0$ and every $v \in k^n$, as in (\ref{N(v) le ... le
  (sum_{j = 1}^n N(e(j))) N_0(v)}) and (\ref{N(v) le ... le (max_{1 le
    j le n} N(e(j))) N_0(v)}).  Under certain conditions, we would like
to show that there is a $C_2 > 0$ such that
\begin{equation}
\label{N_0(v) le C_2 N(v)}
        N_0(v) \le C_2 \, N(v)
\end{equation}
for every $v \in k^n$.  This would imply that the topologies on $k^n$
determined by the metrics associated to $N$ and $N_0$ as in (\ref{d(v,
  w) = N(v - w)}) are the same, as in the previous section.

        Suppose for the moment that $k$ is locally compact with respect
to the metric associated to $|\cdot|$, and that $|\cdot|$ is not the
trivial absolute value function on $k$.  This implies that every
closed ball in $k$ is compact, as in Section \ref{local compactness}.
Note that a closed ball in $k^n$ with respect to $N_0$ with radius $r >
0$ is the same a Cartesian product of $n$ closed balls in $k$ with
respect to $|\cdot|$ with radius $r$, by the definition of $N_0$.
Thus closed balls in $k^n$ with respect to $N_0$ are also compact with
respect to the product topology on $k^n$, because products of compact
sets are compact.

        Remember that $|\cdot|$ is continuous as a real-valued function
on $k$, with respect to the topology on $k$ determined by the metric
corresponding to $|\cdot|$, and the standard topology on ${\bf R}$.
This implies that $N_0$ is continuous as a real-valued function on
$k^n$, with respect to the product topology on $k^n$, and the standard
topology on ${\bf R}$.  This could also be derived from the analogue
of (\ref{max(N(v) - N(w), N(w) - N(v)) le N(v - w)}) for $N_0$.  It
follows that
\begin{equation}
\label{{v in k^n : N_0(v) = 1}}
        \{v \in k^n : N_0(v) = 1\}
\end{equation}
is a closed set in $k^n$ with respect to the product topology.  This
implies that (\ref{{v in k^n : N_0(v) = 1}}) is a compact subset of
$k^n$ with respect to the product topology, because closed subsets of
compact sets are compact.

        It is easy to see that $N$ is also continuous as a real-valued
function on $k^n$, with respect to the product topology on $k^n$ and
the standard topology on ${\bf R}$, using (\ref{max(N(v) - N(w), N(w)
  - N(v)) le N(v - w)}) and (\ref{N(v) le C_1 N_0(v)}).  This implies
that $N$ attains its minimum on (\ref{{v in k^n : N_0(v) = 1}}),
because (\ref{{v in k^n : N_0(v) = 1}}) is nonempty and compact.  Let
$c$ be the minimum value of $N$ on (\ref{{v in k^n : N_0(v) = 1}}), so
that $c > 0$, because $N(v) > 0$ when $v \ne 0$.  We would like to
check that
\begin{equation}
\label{c N_0(v) le N(v)}
        c \, N_0(v) \le N(v)
\end{equation}
for every $v \in k^n$ under these conditions.  More precisely, (\ref{c
  N_0(v) le N(v)}) is trivial when $v = 0$, and (\ref{c N_0(v) le
  N(v)}) holds by definition of $c$ when $N_0(v) = 1$.  If $v \in k^n$
and $v \ne 0$, then there is a $t \in k$ such that $|t| = N_0(v) > 0$,
by the definition (\ref{N_0(v) = max(|v_1|, ldots, |v_n|)}) of
$N_0(v)$.  Thus $N_0(t^{-1} \, v) = 1$, which implies that $N(t^{-1}
\, v) \ge c$, and hence (\ref{c N_0(v) le N(v)}), by (\ref{N(t v) =
  |t| N(v)}).  Of course, (\ref{N_0(v) le C_2 N(v)}) is the same as
(\ref{c N_0(v) le N(v)}), with $C_2 = 1/c$.

        Now suppose that $k$ is complete, and let us show that there
is a $C_2 > 0$ such that (\ref{N_0(v) le C_2 N(v)}) holds for every
$v \in k^n$.  To do this, we use induction on $n$.  The $n = 1$ case
is very easy, using the homogeneity property (\ref{N(t v) = |t| N(v)})
of norms.  Thus we let an integer $n \ge 2$ be given such that the
statement holds for $n - 1$, and we would like to prove the analogous
statement for $n$.

        Let $N$ be a norm on $k^n$, as before.  Put
\begin{equation}
\label{L = {v in k^n : v_n = 0}}
        L = \{v \in k^n : v_n = 0\},
\end{equation}
which is an $(n - 1)$-dimensional linear subspace of $k^n$ that can be
identified with $k^{n - 1}$ in an obvious way.  Note that the
restriction of $N_0$ to $L$ corresponds exactly to the analogue of
$N_0$ on $k^{n - 1}$.  Thus the induction hypothesis can be applied to
the restriction of $N$ to $L$, to get that there is a $C_2' > 0$ such
that
\begin{equation}
\label{N_0(v) le C_2' N(v)}
        N_0(v) \le C_2' \, N(v)
\end{equation}
for every $v \in L$.

        It is easy to see that $k^n$ is complete with respect to the
metric associated to $N_0$ when $k$ is complete.  Similarly, $L$ is
complete with respect to the metric associated to the restriction of
$N_0$ to $L$.  One can check that $L$ is also complete with respect to
the metric associated to the restriction of $N$ to $L$, using
(\ref{N(v) le C_1 N_0(v)}) and (\ref{N_0(v) le C_2' N(v)}).  It
follows from this that $L$ is a closed set in $k^n$ with respect to
the metric associated to $N$.  More precisely, if a sequence
$\{v(l)\}_{l = 1}^\infty$ of elements of $L$ converges to an element
$v$ of $k^n$ with respect to this metric, then $\{v(l)\}_{l =
  1}^\infty$ is a Cauchy sequence in $L$.  The completeness of $L$
implies that $\{v(l)\}_{l = 1}^\infty$ already converges to an element
$v'$ of $L$, so that $v = v'$, because the limit of a convergent
sequence in a metric space is unique.  In particular, $v \in L$, as
desired.

        Let $e(n)$ be the element of $k^n$ whose $j$th coordinate is
equal to $0$ when $j < n$, and whose $n$th coordinate is equal to $1$,
as in the previous section.  Thus $e(n) \not\in L$, and hence there is
a positive real number $c$ such that
\begin{equation}
\label{N(e(n) - w) ge c}
        N(e(n) - w) \ge c
\end{equation}
for every $w \in L$, because $L$ is a closed set in $k^n$ with respect
to the metric associated to $N$.  Equivalently, $N(v) \ge c$ for every
$v \in k^n$ such that $v_n = 1$.  This implies that
\begin{equation}
\label{N(v) ge c |v_n|}
        N(v) \ge c \, |v_n|
\end{equation}
for every $v \in k^n$ with $v_n \ne 0$, by the homogeneity property
(\ref{N(t v) = |t| N(v)}) of norms.  Of course, (\ref{N(v) ge c
  |v_n|}) is trivial when $v_n = 0$.

        Let $v \in k^n$ be given, and observe that
\begin{equation}
\label{N(v - v_n e(n)) le N(v) + |v_n| N(e(n)) le (1 + c^{-1} N(e(n))) N(v)}
        N(v - v_n \, e(n)) \le N(v) + |v_n| \, N(e(n))
                            \le (1 + c^{-1} \, N(e(n))) \, N(v),
\end{equation}
by (\ref{N(v) ge c |v_n|}).  By construction, $v - v_n \, e(n) \in L$,
so that
\begin{equation}
\label{N_0(v - v_n e(n)) le ... le C_2' (1 + c^{-1} N(e(n)) N(v)}
        N_0(v - v_n \, e(n)) \le C_2' \, N(v - v_n \, e(n))
                              \le C_2' \, (1 + c^{-1} \, N(e(n)) \, N(v),
\end{equation}
by (\ref{N_0(v) le C_2' N(v)}).  Note that
\begin{equation}
\label{N_0(v - v_n e(n)) = max(|v_1|, ldots, |v_{n - 1}|)}
        N_0(v - v_n \, e(n)) = \max(|v_1|, \ldots, |v_{n - 1}|),
\end{equation}
which implies that $N_0(v) = \max(N_0(v - v_n \, e(n)), |v_n|)$.  It
follows that
\begin{equation}
\label{N_0(v) le max(C_2' (1 + c^{-1} N(e(n))), c^{-1}) N(v)}
        N_0(v) \le \max(C_2' \, (1 + c^{-1} \, N(e(n))), c^{-1}) \, N(v)
\end{equation}
for every $v \in k^n$, as desired, by (\ref{N(v) ge c |v_n|}) and
(\ref{N_0(v - v_n e(n)) le ... le C_2' (1 + c^{-1} N(e(n)) N(v)}).

\section{Some related facts}
\label{some related facts}

        Let $k$ be a field, and let $|\cdot|$ be an absolute value function
on $k$.  Also let $V$ be a vector space over $k$, and let $N$ be a
norm on $V$.  If $k$ is complete with respect to the metric associated
to $|\cdot|$, and if $V$ has finite dimension as a vector space over
$k$, then $V$ is also complete with respect to the metric associated
to $N$.  Of course, this is trivial when $V = \{0\}$, and otherwise
one can reduce to the case where $V = k^n$ for some positive integer
$n$.  It is easy to see that $k^n$ is complete with respect to the metric
associated to the norm $N_0$ in (\ref{N_0(v) = max(|v_1|, ldots, |v_n|)}),
as mentioned in the previous section, and completeness with respect to
the metric associated to any other norm $N$ can be derived from
(\ref{N(v) le C_1 N_0(v)}) and (\ref{N_0(v) le C_2 N(v)}).

        Suppose that $|\cdot|$ is nontrivial on $k$, and that $k$ is
locally compact.  If $V$ is finite-dimensional, then $V$ is locally
compact with respect to the topology determined by the metric
associated to $N$, and in fact closed balls in $V$ with respect to $N$
are compact with respect to the corresponding topology.  To see this,
it suffices to consider the case where $V = k^n$, as before.  We have
already seen that closed balls in $k^n$ with respect to $N_0$ are
compact, as in the previous section.  This implies that closed balls
with respect to $N$ are compact, by (\ref{N(v) le C_1 N_0(v)}) and
(\ref{N_0(v) le C_2 N(v)}), and because closed subsets of compact sets
are compact.

        If $|\cdot|$ is the trivial absolute value function on $k$,
then $N_0$ is the trivial ultranorm on $k^n$ for each positive integer
$n$.  If $N$ is any other norm on $k^n$, then we still have (\ref{N(v)
  le C_1 N_0(v)}) and (\ref{N_0(v) le C_2 N(v)}), because $k$ is
complete.  Of course, the proof in the previous section could be
simplified in this case.  It follows that the topology on $k^n$
determined by the metric associated to $N$ is still the discrete
topology.  There is an analogous statement for arbitrary
finite-dimensional vector spaces over $k$, as before.

        Let $|\cdot|$ be an absolute value function on a field $k$ again,
and let $N$ be a norm on a vector space $V$ over $k$, which is not
necessarily finite-dimensional.  If $k$ is complete and $W$ is a
finite-dimensional linear subspace of $V$, then $W$ is complete with
respect to the restriction to $W$ of the metric associated to $N$, as
before.  This implies that $W$ is a closed set in $V$ with respect to
the topology determined by the metric associated to $N$.  This follows
from a standard argument that was also used in the previous section.

        Let $W$ be a linear subspace of $V$ again.  Suppose that there
is a positive real number $c < 1$ such that for each $v \in V$ there
is a $w \in W$ such that
\begin{equation}
\label{N(v - w) le c N(v)}
        N(v - w) \le c \, N(v).
\end{equation}
Applying this property to $v - w$ and repeating the process, one can
check that $W$ is dense in $V$ with respect to the metric associated
to $N$.  If $W$ is a closed set in $V$ with respect to the topology
determined by this metric, then it follows that $V = W$.  In
particular, this holds when $k$ is complete and $W$ is
finite-dimensional, as in the preceding paragraph.

        Let us now restrict our attention to the case where $|\cdot|$
is nontrivial on $k$.  Thus there is a $t_0 \in k$ such that $t_ 0 \ne 0$
and $|t_0| \ne 1$.  As usual, we may as well ask that $|t_0| > 1$,
since otherwise we could replace $t_0$ with $1/t_0$.

        Suppose that $V$ is locally compact with respect to the topology
determined by the metric associated to $N$.  This implies that every
closed ball in $V$ with respect to $N$ is compact, because $|\cdot|$
is nontrivial on $k$, as in Section \ref{local compactness}.  In
particular, the closed unit ball in $V$ is compact, and hence totally
bounded.  Thus for each $\epsilon > 0$ the closed unit ball in $V$
with respect to $N$ can be covered by finitely many closed balls of radius
$\epsilon$.  Let us apply this with $\epsilon = 1/(2 \, |t_0|)$, to get
finitely many vectors $w_1, \ldots, w_n$ in $V$ such that for each $u \in V$
with $N(u) \le 1$ there is a $j \in \{1, \ldots, n\}$ that satisfies
\begin{equation}
\label{N(u - w_j) le 1/(2 |t_0|)}
        N(u - w_j) \le 1/(2 \, |t_0|).
\end{equation}

        Let $W$ be the linear span of $w_1, \ldots, w_n$ in $V$, and let
$v \in V$ be given.  We would like to show that there is a $w \in W$
that satisfies (\ref{N(v - w) le c N(v)}) with $c = 1/2$.  Of course,
this is trivial when $v = 0$, and so we may suppose that $v \ne 0$.
Let $l$ be the unique integer such that $|t_0|^{l - 1} < N(v) \le
|t_0|^l$.  If $u = t_0^{-l} \, v$, then $N(u) \le 1$, and hence there
is a $j \in \{1, \ldots, n\}$ that satisfies (\ref{N(u - w_j) le 1/(2
  |t_0|)}).  It follows that
\begin{equation}
\label{N(v - t_0^l w_j) = |t_0|^l N(u - w_j) le |t_0|^{l - 1}/2 < N(v)/2}
        N(v - t_0^l \, w_j) = |t_0|^l \, N(u - w_j) \le |t_0|^{l - 1}/2 < N(v)/2.
\end{equation}
This is exactly what we wanted, since $t_0^l \, w_j \in W$.

        This implies that $W$ is dense in $V$, as before.  If $k$ is complete,
then $V = W$, because $W$ is finite-dimensional, by construction.  This
shows that $V$ is finite-dimensional when $V$ is locally compact, $|\cdot|$
is nontrivial on $k$, and $k$ is complete.

\section{Haar measure, 3}
\label{haar measure, 3}

        Let $X$ be a locally compact Hausdorff topological space,
and let $W$ be an open subset of $X$.  It is well known that $X$ is
regular as a topological space under these conditions, so that for
each $p \in W$ there is an open subset $V(p)$ of $X$ such that $p \in
V(p)$ and the closure $\overline{V(p)}$ of $V(p)$ in $X$ is contained
in $W$.  Because $X$ is locally compact, one can also choose $V(p)$ so
that $\overline{V}(p)$ is compact for each $p \in W$.  If there is
also a base for the topology of $X$ with only finitely or countably
many elements, then Lindel\"of's theorem implies that $W$ can be
expressed as the union of finitely or countably many of the sets
$V(p)$.  In this case, it follows that $W$ is
\emph{$\sigma$-compact},\index{sigma-compact sets@$\sigma$-compact sets}
which means that $W$ can be expressed as the union of finitely or countably
many compact sets, because $W$ is the union of finitely or countably
many of the sets $\overline{V(p)}$.    

        Now let $n$ be a positive integer, and let $A_1, A_2, \ldots, A_n$
be $n$ commutative topological groups.  The Cartesian product
\begin{equation}
\label{A = prod_{j = 1}^n A_j}
        A = \prod_{j = 1}^n A_j
\end{equation}
is also a commutative group, where the group operations on $A$ are
defined coordinatewise.  If $A$ is equipped with the product topology
corresponding to the given topologies on the $A_j$'s, then it is easy
to see that $A$ is a topological group as well.  Suppose from now on
that $A_j$ is locally compact for each $j$, which implies that $A$ is
locally compact too.  This uses the fact that Cartesian products of
compact subsets of the $A_j$'s are compact subsets of $A$, with
respect to the product topology.

        Under these conditions, there is a Haar measure $H_j$ on
$A_j$ for each $j$.  To get a Haar measure $H$ on $A$, one can
basically use a suitable product measure construction, but there are
some details involved with this.  Suppose first that there is a base
$\mathcal{B}_j$ for the topology of $A_j$ with only finitely or
countably many elements for each $j$.  Let $\mathcal{B}$ be the
collection of sets $U \subseteq A$ of the form
\begin{equation}
\label{U = prod_{j = 1}^n U_j}
        U = \prod_{j = 1}^n U_j,
\end{equation}
where $U_j \in \mathcal{B}_j$ for each $j$.  It is easy to see that
$\mathcal{B}$ is a base for the product topology on $A$, and that
$\mathcal{B}$ has only finitely or countably many elements.

        In this situation, $A_j$ is $\sigma$-compact for each $j$,
by the remarks at the beginning of the section.  This implies that
$A_j$ is $\sigma$-finite with respect to $H_j$ for each $j$, because
compact sets have finite Haar measure.  This permits one to apply the
standard construction of product measures to the $H_j$'s, as Borel
measures on the $A_j$'s.  This leads to a product measure $H$, defined
on a suitable $\sigma$-algebra of subsets of $A$.  In particular, if
$E \subseteq A$ is a product of Borel subsets of the $A_j$'s, then $E$
is measurable with respect to the product product measure
construction.

        Of course, open subsets of the $A_j$'s are Borel sets, by
definition, so that products of open subsets of the $A_j$'s are
measurable with respect to the product measure construction.  Thus the
elements of the base $\mathcal{B}$ for the topology of $A$ mentioned
earlier are all measurable with respect to the product measure
construction.  Remember that every open set in $A$ can be expressed as
a union of elements of $\mathcal{B}$, and that $\mathcal{B}$ has only
finitely or countably many elements.  It follows that every open set
in $A$ is measurable with respect to the product measure construction,
since it can be expressed as a union of finitely or countably many
measurable sets.  This implies that Borel subsets of $A$ are measurable
with respect to the product measure construction under these conditions.

        Thus the product measure $H$ may be considered as a Borel
measure on $A$.  One can also check that $H$ is invariant under
translations on $A$, because $H_j$ is invariant under translations on
$A_j$ for each $j$.  If $V$ is a nonempty open subset of $A$, then $V$
contains a product of nonempty open subsets of the $A_j$'s, by
definition of the product topology.  This implies that $H(V) > 0$,
because of the corresponding property of the $H_j$'s, which is one of
the requirements of a Haar measure.  If $K$ is a compact subset of
$A$, then the projection of $K$ in $A_j$ is compact for each $j$,
because the projection mappings are continuous.  Of course, $K$ is
contained in the Cartesian product of its projections in the $A_j$'s.
This implies that $H(K)$ is finite, because of the corresponding
property of the $H_j$'s.  Note that open subsets of $A$ are
$\sigma$-compact in this situation, by the remarks at the beginning of
the section.  Because of this, it is well known that $H$ satisfies the
regularity conditions required of a Haar measure.

        Alternatively, there is another product measure construction
for nonnegative Borel measures on locally compact Hausdorff topological
spaces, where the Borel measures satisfy suitable regularity conditions,
but the topological spaces are not required to have countable bases.
Because of the Riesz representation theorem, one can look at this
in terms of the corresponding nonnegative linear functionals on
spaces of continuous functions with compact support.  In order to define
such a linear functional on the product space, one can basically apply
the given linear functionals on the factors.  In our situation, this
means that one can get a Haar integral on $A$ using Haar integrals
on the $A_j$'s.

        Let $k$ be a field, and let $|\cdot|$ be an absolute value function
on $k$.  Also let $n$ be a positive integer, so that $k^n$ is an
$n$-dimensional vector space over $k$, with respect to coordinatewise
addition and scalar multiplication.  As usual, we take $k^n$ to be
equipped with the product topology corresponding to the topology on
$k$ determined by the metric associated to $|\cdot|$.  Note that $k^n$
is a commutative topological group with respect to addition and this
topology, because $k$ is a commutative topological group with respect
to addition.

        Suppose that $k$ is locally compact, so that $k^n$ is locally compact 
too, as before.  Of course, if $|\cdot|$ is the trivial absolute value
function on $k$, then the corresponding topology on $k$ is the
discrete topology, and the associated product topology on $k^n$ is the
discrete topology too.  In this case, one can simply use counting
measure as Haar measure, on $k$ and $k^n$.  Let us suppose from now on
that $|\cdot|$ is not the trivial absolute value function on $k$.
This implies that closed and bounded subsets of $k$ are compact, as in
Section \ref{local compactness}.

        It follows that $k$ is $\sigma$-compact under these conditions,
because $k$ can be expressed as the union of a sequence of closed
balls, each of which is compact.  It is well known that compact metric
spaces are separable, and hence that $k$ is separable.  This implies
that there is a countable base for the topology of $k$, because $k$ is
a separable metric space.  This permits one to look at Haar measure on
$k^n$ in terms of the usual product measure construction applied to
$n$ copies of Haar measure on $k$, as before.  Otherwise, one can
get a Haar integral on $k^n$ using a Haar integral on $k$.

\section{Summable functions}
\label{summable functions}

        Let $X$ be a nonempty set, and let $f(x)$ be a nonnegative
real-valued function on $X$.  Thus
\begin{equation}
\label{sum_{x in A} f(x)}
        \sum_{x \in A} f(x)
\end{equation}
is defined for every nonempty finite subset $A$ of $X$, and the sum
\begin{equation}
\label{sum_{x in X} f(x)}
        \sum_{x \in X} f(x)
\end{equation}
can be defined as the supremum of (\ref{sum_{x in A} f(x)}) over all
nonempty finite subsets $A$ of $X$.  More precisely, the supremum is
finite when there is a finite upper bound for the subsums (\ref{sum_{x
    in A} f(x)}), and otherwise (\ref{sum_{x in X} f(x)}) is equal to
$+\infty$.  If $g$ is another nonnegative real-valued function on $X$,
then one can check that
\begin{equation}
\label{sum_{x in X} (f(x) + g(x)) = sum_{x in X} f(x) + sum_{x in X} g(x)}
        \sum_{x \in X} (f(x) + g(x)) = \sum_{x \in X} f(x) + \sum_{x \in X} g(x),
\end{equation}
where the right side of (\ref{sum_{x in X} (f(x) + g(x)) = sum_{x in
    X} f(x) + sum_{x in X} g(x)}) is interpreted as being $+\infty$
when either of the individual sums is equal to $+\infty$.  Similarly,
if $a$ is a nonnegative real number, then
\begin{equation}
\label{sum_{x in X} a f(x) = a sum_{x in X} f(x)}
        \sum_{x \in X} a \, f(x) = a \, \sum_{x \in X} f(x),
\end{equation}
where the right side of (\ref{sum_{x in X} a f(x) = a sum_{x in X}
  f(x)}) is interpreted as being $+\infty$ when (\ref{sum_{x in X}
  f(x)}) is $+\infty$ and $a > 0$, and equal to $0$ when $a = 0$, even
when (\ref{sum_{x in X} f(x)}) is infinite.

        A nonnegative real-valued function $f(x)$ on $X$ is said to be
\emph{summable}\index{summable functions} on $X$ if (\ref{sum_{x in X} f(x)})
is finite.  Let $\epsilon > 0$ be given, and put
\begin{equation}
\label{E(f, epsilon) = {x in X : f(x) ge epsilon}}
        E(f, \epsilon) = \{x \in X : f(x) \ge \epsilon\}.
\end{equation}
If $A$ is a nonempty finite subset of $E(f, \epsilon)$, then
\begin{equation}
\label{epsilon (A) le sum_{x in A} f(x) le sum_{x in X} f(x)}
        \epsilon \, (\# A) \le \sum_{x \in A} f(x) \le \sum_{x \in X} f(x),
\end{equation}
where $\# A$ denotes the number of elements of $A$.  If $f$ is
summable on $X$, then it follows that $E(f, \epsilon)$ has only
finitely many elements, and that
\begin{equation}
\label{E(f, epsilon) le epsilon^{-1} sum_{x in X} f(x)}
        \# E(f, \epsilon) \le \epsilon^{-1} \, \sum_{x \in X} f(x)
\end{equation}
for each $\epsilon > 0$.  In particular, if $f$ is summable on $X$,
then the set of $x \in X$ such that $f(x) \ne 0$ has only finitely
or countably many elements, since it can be expressed as the
union of $E(f, 1/n)$ over all positive integers $n$.

        Now let $k$ be a field, and let $V$ be a vector space over $k$.
Also let $|\cdot|$ be an absolute value function on $k$, and let
$N$ be a norm on $V$ with respect to $|\cdot|$, as in Section \ref{norms,
ultranorms}.  Under these conditions, a $V$-valued function $f(x)$ on
$X$ is said to be \emph{summable}\index{summable functions} if $N(f(x))$
is summable on $X$ as a nonnegative real-valued function on $X$,
as in the preceding paragraph.  Of course, one can consider $V = k$
as a one-dimensional vector space over $k$, and $|\cdot|$ as a norm
on this vector space, so that this definition applies to $k$-valued
functions on $X$ in particular.

        Let $\ell^1(X, V)$\index{l^1(X, V)@$\ell^1(X, V)$} be the space
of summable $V$-valued functions on $X$, in the sense just described,
and put
\begin{equation}
\label{||f||_1 = ||f||_{ell^1(X, V)} = sum_{x in X} N(f(x))}
        \|f\|_1 = \|f\|_{\ell^1(X, V)} = \sum_{x \in X} N(f(x))
\end{equation}
for every $f \in \ell^1(X, V)$.  If $f, g \in \ell^1(X, V)$, then it
is easy to see that $f + g$ is summable on $X$ too, and that
\begin{eqnarray}
\label{||f + g||_1 = sum_{x in X} N(f(x) + g(x)) le ... = ||f||_1 + ||g||_1}
        \|f + g\|_1 & = & \sum_{x \in X} N(f(x) + g(x)) \\
 & \le & \sum_{x \in X} (N(f(x)) + N(g(x))) = \|f\|_1 + \|g\|_1. \nonumber
\end{eqnarray}
Similarly, if $f \in \ell^1(X, V)$ and $t \in k$, then $t \, f(x)$ is
summable on $X$ as well, and
\begin{equation}
\label{||t f||_1 = sum_{x in X} N(t f(x)) = ... = |t| ||f||_1}
 \|t \, f\|_1 = \sum_{x \in X} N(t \, f(x)) = \sum_{x \in X} |t| \, N(f(x))
              = |t| \, \|f\|_1.
\end{equation}
This shows that $\ell^1(X, V)$ is a vector space over $k$ with respect
to pointwise addition and scalar multiplication, and that $\|f\|_1$ is
a norm on $\ell^1(X, V)$.

        If $f$ is any $V$-valued function on $X$, then the 
\emph{support}\index{supports of functions} of $f$ is defined by
\begin{equation}
\label{supp f = {x in X : f(x) ne 0}}
        \supp f = \{x \in X : f(x) \ne 0\}.
\end{equation}
Let $c_{00}(X, V)$\index{c_{00}(X, V)@$c_{00}(X, V)$} be the
collection of $V$-valued functions $f$ on $X$ such that $\supp f$ has
only finitely many elements.  This is a vector space over $k$ with
respect to pointwise addition and scalar multiplication, and a linear
subspace of $\ell^1(X, V)$.  If $f$ is summable on $X$, then we have
seen that the support of $N(f(x))$ has only finitely or countably many
elements, which means that $\supp f$ has the same property.  If
$|\cdot|$ is the trivial absolute value function on $k$, and $N$ is
the trivial ultranorm on $V$, then $f$ is summable on $V$ if and only
if $\supp f$ has only finitely many elements.

        Suppose that $f$ is a summable $V$-valued function on $X$, and
let $\epsilon > 0$ be given.  Thus $N(f(x))$ is summable as a
nonnegative real-valued function on $X$, and hence there is a finite
set $A(\epsilon) \subseteq X$ such that
\begin{equation}
\label{sum_{x in X} N(f(x)) < sum_{x in A(epsilon)} N(f(x)) + epsilon}
        \sum_{x \in X} N(f(x)) < \sum_{x \in A(\epsilon)} N(f(x)) + \epsilon,
\end{equation}
by the definition of the sum on the left side of (\ref{sum_{x in X}
  N(f(x)) < sum_{x in A(epsilon)} N(f(x)) + epsilon}).  If we put
$f_\epsilon(x) = f(x)$ when $x \in A(\epsilon)$, and $f_\epsilon(x) =
0$ otherwise, then $f_\epsilon \in c_{00}(X, V)$, and
\begin{equation}
\label{||f - f_epsilon||_1 = sum_{x in X setminus A(epsilon)} N(f(x)) < epsilon}
 \|f - f_\epsilon\|_1 = \sum_{x \in X \setminus A(\epsilon)} N(f(x)) < \epsilon.
\end{equation}
This uses (\ref{sum_{x in X} N(f(x)) < sum_{x in A(epsilon)} N(f(x)) +
  epsilon}) in the second step, and it follows that $c_{00}(X, V)$ is
dense in $\ell^1(X, V)$ with respect to the metric associated to the
$\ell^1$ norm.

        If $f \in c_{00}(X, V)$, then
\begin{equation}
\label{sum_{x in X} f(x), 2}
        \sum_{x \in X} f(x)
\end{equation}
can be defined as an element of $V$ in the obvious way, and satisfies
\begin{equation}
\label{N(sum_{x in X} f(x)) le sum_{x in X} N(f(x)) = ||f||_1}
        N\Big(\sum_{x \in X} f(x)\Big) \le \sum_{x \in X} N(f(x)) = \|f\|_1.
\end{equation}
Of course,
\begin{equation}
\label{f mapsto sum_{x in X} f(x)}
        f \mapsto \sum_{x \in X} f(x)
\end{equation}
defines a linear mapping from $c_{00}(X, V)$ into $V$, so that
\begin{equation}
\label{N(sum_{x in X} f(x) - sum_{x in X} g(x)) = ... le ||f - g||_1}
        N\Big(\sum_{x \in X} f(x) - \sum_{x \in X} g(x)\Big)
            = N\Big(\sum_{x \in X} (f(x) - g(x))\Big) \le \|f - g\|_1
\end{equation}
for every $f, g \in c_{00}(X, V)$, as in (\ref{N(sum_{x in X} f(x)) le
  sum_{x in X} N(f(x)) = ||f||_1}).  This implies that (\ref{f mapsto
  sum_{x in X} f(x)}) is uniformly continuous with respect to the
metric associated to the $\ell^1$ norm on $c_{00}(X, V)$, and the
metric associated to $N$ on $V$.

        If $V$ is complete with respect to the metric associated to $N$,
then there is a unique extension of (\ref{f mapsto sum_{x in X} f(x)})
to a uniformly continuous mapping from $\ell^1(X, V)$ into $V$, with
respect to the metric associated to the $\ell^1$ norm on $\ell^1(X,
V)$, and the metric associated to $N$ on $V$.  This follows from the
discussion in Section \ref{continuous extensions}, and uses the fact
that $c_{00}(X, V)$ is dense in $\ell^1(X, V)$.  In this situation, it
is easy to see that this extension of (\ref{f mapsto sum_{x in X}
  f(x)}) to $\ell^1(X, V)$ is also linear, and satisfies
(\ref{N(sum_{x in X} f(x)) le sum_{x in X} N(f(x)) = ||f||_1}).  At
any rate, one can use this to define (\ref{sum_{x in X} f(x), 2}) as
an element of $V$ when $f \in \ell^1(X, V)$ and $V$ is complete.

        Alternatively, if $f \in \ell^1(X, V)$, then we have seen that
$\supp f$ has only finitely or countably many elements.  If $\supp f$
has only finitely many elements, then we already know how to define
(\ref{sum_{x in X} f(x), 2}), and so we may as well suppose that
$\supp f$ is countably infinite.  Thus the elements of $\supp f$ may be
enumerated by a sequence, and (\ref{sum_{x in X} f(x), 2}) may be identified
with an infinite series with terms in $V$.  The summability of $f$ implies
that this series is absolutely convergent, and hence that the partial sums
form a Cauchy sequence in $V$, as in Section \ref{infinite series}.
If $V$ is complete, then it follows that the sequence of partial sums
converges in $V$.  Any other enumeration of the elements of $\supp f$
would lead to a rearrangement of the same series, and one can check that
the corresponding sums would be the same, using absolute convergence.
Of course, this approach to the definition of (\ref{sum_{x in X} f(x), 2})
for $f \in \ell^1(X, V)$ is equivalent to the one described in the
preceding paragraph.

        Let us continue to suppose that $V$ is complete, and verify that
$\ell^1(X, V)$ is complete with respect to the $\ell^1$ metric.  Let
$\{f_j\}_{j = 1}^\infty$ be a Cauchy sequence in $\ell^1(X, V)$, so that
for each $\epsilon > 0$ there is an $L(\epsilon) \ge 1$ such that
\begin{equation}
\label{||f_j - f_l||_1 < epsilon}
        \|f_j - f_l\|_1 < \epsilon
\end{equation}
for every $j, l \ge L(\epsilon)$.  In particular,
\begin{equation}
\label{N(f_j(x) - f_l(x)) le ||f_j - f_l||_1 < epsilon}
        N(f_j(x) - f_l(x)) \le \|f_j - f_l\|_1 < \epsilon
\end{equation}
for every $x \in X$ and $j, l \ge L(\epsilon)$, which implies that
$\{f_j(x)\}_{j = 1}^\infty$ is a Cauchy sequence in $V$ for every $x
\in X$.  If $V$ is complete, then it follows that $\{f_j(x)\}_{j =
  1}^\infty$ converges to an element of $V$ for every $x \in X$, and
we let $f(x)$ denote the limit of this sequence.  We would like to
show that $f \in \ell^1(X, V)$, and that $\{f_j\}_{j = 1}^\infty$ converges
to $f$ with respect to the $\ell^1$ norm.

        If $A$ is any nonempty finite subset of $X$, then
\begin{equation}
\label{sum_{x in A} N(f_j(x) - f(x)) = ...}
        \sum_{x \in A} N(f_j(x) - f(x))
               = \lim_{l \to \infty} \sum_{x \in A} N(f_j(x) - f_l(x))
\end{equation}
for every $j$, because $N$ is continuous on $V$ with respect to the
metric associated to itself.  Combining this with (\ref{||f_j -
  f_l||_1 < epsilon}), we get that
\begin{equation}
\label{sum_{x in A} N(f_j(x) - f(x)) le epsilon}
        \sum_{x \in A} N(f_j(x) - f(x)) \le \epsilon
\end{equation}
for every $j \ge L(\epsilon)$.  This implies that
\begin{equation}
\label{sum_{x in X} N(f_j(x) - f(x)) le epsilon}
        \sum_{x \in X} N(f_j(x) - f(x)) \le \epsilon
\end{equation}
for every $j \ge L(\epsilon)$, by taking the supremum over all finite
subsets $A$ of $X$ in (\ref{sum_{x in A} N(f_j(x) - f(x)) le
  epsilon}).  In particular, one can use this to get that $f$ is
summable on $X$, because $f_j$ is summable on $X$ for every $j$, by
hypothesis.  It follows easily that $\{f_j\}_{j = 1}^\infty$ converges
to $f$ with respect to the $\ell^1$ norm, as desired, using
(\ref{sum_{x in X} N(f_j(x) - f(x)) le epsilon}) again.

\section{Vanishing at infinity}
\label{vanishing at infinity}

        As in the previous section, we let $X$ be a nonempty set, $k$
be a field, and $V$ be a vector space over $k$.  We also let $|\cdot|$
be an absolute value function on $k$, and $N$ be a norm on $V$ with
respect to $|\cdot|$.  Under these conditions, a $V$-valued function
$f$ on $X$ is said to be \emph{bounded}\index{bounded functions} if
$N(f(x))$ is bounded on $X$, as a nonnegative real-valued function on
$X$.  It is easy to see that the space $\ell^\infty(X,
V)$\index{l^infty(X, V)@$\ell^\infty(X, V)$} of bounded $V$-valued
functions on $X$ is a vector space over $k$ with respect to pointwise
addition and scalar multiplication.  If $f \in \ell^\infty(X, V)$,
then we put
\begin{equation}
\label{||f||_infty = ||f||_{ell^infty(X, V)} = sup_{x in X} N(f(x))}
        \|f\|_\infty = \|f\|_{\ell^\infty(X, V)} = \sup_{x \in X} N(f(x)),
\end{equation}
which is easily seen to be a norm on $\ell^\infty(X, V)$.

        A $V$-valued function $f$ on $X$ is said to \emph{vanish at 
infinity}\index{vanishing at infinity} on $X$ if for each $\epsilon > 0$,
\begin{equation}
\label{N(f(x)) < epsilon}
        N(f(x)) < \epsilon
\end{equation}
for all but finitely many $x \in X$.  Let $c_0(X, V)$\index{c_0(X,
  V)@$c_0(X, V)$} be the collection of $V$-valued functions on $X$
that vanish at infinity.  It is easy to see that $c_0(X, V)$ is a
linear subspace of $\ell^\infty(X, V)$, and that the space $c_{00}(X,
V)$ of $V$-valued functions on $X$ with finite support is a linear
subspace of $c_0(X, V)$.  One can also check that $c_0(X, V)$ is a
closed set in $\ell^\infty(X, V)$ with respect to the metric
associated to the $\ell^\infty$ norm.  If $|\cdot|$ is the trivial
absolute value function on $k$, and $N$ is the trivial ultranorm on
$V$, then a $V$-valued function $f$ on $X$ vanishes at infinity
if and only if $\supp f$ has only finitely many elements.

        Let $f$ be a $V$-valued function on $X$ that vanishes at
infinity, and let $\epsilon > 0$ be given.  Put $f_\epsilon(x) = f(x)$
when $N(f(x)) \ge \epsilon$ and $f_\epsilon(x) = 0$ otherwise, so that
$\supp f_\epsilon$ has only finitely many elements, by hypothesis.
Thus $f_\epsilon \in c_{00}(X, V)$, and it is easy to see that
\begin{equation}
\label{||f - f_epsilon||_infty le epsilon}
        \|f - f_\epsilon\|_\infty \le \epsilon,
\end{equation}
by construction.  This shows that $c_{00}(X, V)$ is dense in $c_0(X,
V)$, with respect to the metric associated to the $\ell^\infty$ metric.

        If $V$ is complete with respect to the metric associated to $N$,
then $\ell^\infty(X, V)$ is complete with respect to the metric
corresponding to the $\ell^\infty$ norm.  In fact, it is well known
that the space of bounded mappings from $X$ into any complete metric
space is itself a complete metric space with respect to the supremum
metric.  It follows that $c_0(X, V)$ is complete with respect to the
$\ell^\infty$ norm when $V$ is complete, because a closed subset of a
complete metric space is also complete, with respect to the restriction
of the metric to the subset.

        Let us suppose from now on in this section that $|\cdot|$
is an ultrametric absolute value function on $k$, and that $N$ is an
ultranorm on $V$ with respect to $|\cdot|$.  It is easy to see that
the $\ell^\infty$ norm is an ultranorm on $\ell^\infty(X, V)$ under these
conditions.  If $f \in c_{00}(X, V)$, then
\begin{equation}
\label{sum_{x in X} f(x), 3}
        \sum_{x \in X} f(x)
\end{equation}
can be defined in the usual way, and satisfies
\begin{equation}
\label{N(sum_{x in X} f(x)) le max_{x in X} N(f(x)) = ||f||_infty}
        N\Big(\sum_{x \in X} f(x)\Big) \le \max_{x \in X} N(f(x)) = \|f\|_\infty,
\end{equation}
by the ultrametric version of the triangle inequality.  As before,
\begin{equation}
\label{f mapsto sum_{x in X} f(x), 2}
        f \mapsto \sum_{x \in X} f(x)
\end{equation}
is a linear mapping from $c_{00}(X, V)$ into $V$, and in this case we
have that
\begin{equation}
\label{N(sum_{x in X} f(x) - sum_{x in X} g(x)) = ... le ||f - g||_infty}
        N\Big(\sum_{x \in X} f(x) - \sum_{x \in X} g(x)\Big)
           = N\Big(\sum_{x \in X} (f(x) - g(x))\Big) \le \|f - g\|_\infty
\end{equation}
for every $f, g \in c_{00}(X, V)$, as in (\ref{N(sum_{x in X} f(x)) le
  max_{x in X} N(f(x)) = ||f||_infty}).  It follows that (\ref{f
  mapsto sum_{x in X} f(x), 2}) is uniformly continuous with respect
to the metric associated to the $\ell^\infty$ norm on $c_{00}(X, V)$,
and the metric associated to $N$ on $V$.

        If $V$ is complete with respect to the metric associated to $N$,
then there is a unique extension of (\ref{f mapsto sum_{x in X} f(x),
  2}) to a uniformly continuous mapping from $c_0(X, V)$ into $V$,
with respect to the metric on $c_0(X, V)$ associated to the
$\ell^\infty$ norm, and the metric on $V$ associated to $N$.  This
follows from the discussion in Section \ref{continuous extensions},
and uses the fact that $c_{00}(X, V)$ is dense in $c_0(X, V)$.  As in
the previous section, it is easy to see that this extension of (\ref{f
  mapsto sum_{x in X} f(x), 2}) to $c_0(X, V)$ is linear, and that it
also satisfies (\ref{N(sum_{x in X} f(x)) le max_{x in X} N(f(x)) =
  ||f||_infty}).  As before, this can be used to define (\ref{sum_{x
    in X} f(x), 3}) when $f \in c_0(X, V)$ and $V$ is complete.

        If $f \in c_0(X, V)$, then for each positive integer $n$,
$N(f(x)) \ge 1/n$ for at most finitely many $x \in X$.  This implies
that $f(x) \ne 0$ for at most finitely or countably many $x \in X$, by
taking the union over $n \ge 1$.  If $f(x) \ne 0$ for only finitely
many $x \in X$, then (\ref{sum_{x in X} f(x), 3}) can be defined in
the usual way.  Otherwise, the support of $f$ is countably infinite,
so that its elements may be enumerated by an infinite sequence.  Thus
(\ref{sum_{x in X} f(x), 3}) can be identified with an infinite
series, as in the previous section.  In this case, the terms of the
series converge to $0$, which implies that the corresponding sequence
of partial sums is a Cauchy sequence in $V$, because of the
ultrametric version of the triangle inequality.  This is similar to
the ultrametric case in Section \ref{infinite series}, and it follows
that the sequence of partial sums converges in $V$ when $V$ is
complete.  One can also check that rearrangements of the series have
the same sum under these conditions, so that the value of the sum
does not depend on the particular enumeration of the elements of
the support of $f$.  As before, this approach to the definition of
(\ref{sum_{x in X} f(x), 3}) for $f \in c_0(X, V)$ is equivalent to
the one described in the preceding paragraph.

\section{Double sums}
\label{double sums}

        Let $X_1$, $X_2$ be nonempty sets, and let $X = X_1 \times X_2$
be their Cartesian product.  Also let $f(x_1, x_2)$ be a nonnegative
real-valued function on $X$, so that the sum
\begin{equation}
\label{sum_{(x_1, x_2) in X} f(x_1, x_2)}
        \sum_{(x_1, x_2) \in X} f(x_1, x_2)
\end{equation}
can be defined as in Section \ref{summable functions}.  Similarly,
\begin{equation}
\label{sum_{x_2 in X_2} f(x_1, x_2)}
        \sum_{x_2 \in X_2} f(x_1, x_2)
\end{equation}
is defined as a nonnegative extended real number for each $x_1 \in X_1$, and
\begin{equation}
\label{sum_{x_1 in X_1} f(x_1, x_2)}
        \sum_{x_1 \in X_1} f(x_1, x_2)
\end{equation}
is defined as a nonnegative extended real number for each $x_2 \in X_2$.
This permits the iterated sums
\begin{equation}
\label{sum_{x_1 in X_1} (sum_{x_2 in X_2} f(x_1, x_2))}
        \sum_{x_1 \in X_1} \Big(\sum_{x_2 \in X_2} f(x_1, x_2)\Big)
\end{equation}
and
\begin{equation}
\label{sum_{x_2 in X_2} (sum_{x_1 in X_1} f(x_1, x_2))}
        \sum_{x_2 \in X_2} \Big(\sum_{x_1 \in X_1} f(x_1, x_2)\Big)
\end{equation}
to be defined in essentially the same way as before.  More precisely,
if (\ref{sum_{x_2 in X_2} f(x_1, x_2)}) is equal to $+\infty$ for any
$x_1 \in X_1$, then (\ref{sum_{x_1 in X_1} (sum_{x_2 in X_2} f(x_1,
  x_2))}) is interpreted as being $+\infty$, and otherwise
(\ref{sum_{x_1 in X_1} (sum_{x_2 in X_2} f(x_1, x_2))}) may be defined
as in Section \ref{summable functions}.  The sum (\ref{sum_{x_2 in
    X_2} (sum_{x_1 in X_1} f(x_1, x_2))}) is defined analogously, and
one can show that (\ref{sum_{(x_1, x_2) in X} f(x_1, x_2)}),
(\ref{sum_{x_1 in X_1} (sum_{x_2 in X_2} f(x_1, x_2))}), and
(\ref{sum_{x_2 in X_2} (sum_{x_1 in X_1} f(x_1, x_2))}) are equal to
each other under these conditions.  One can start by checking that
(\ref{sum_{(x_1, x_2) in X} f(x_1, x_2)}) and (\ref{sum_{x_1 in X_1}
  (sum_{x_2 in X_2} f(x_1, x_2))}) are each less than or equal to the
other, by considering approximations to these sums by finite subsums.
This implies that (\ref{sum_{(x_1, x_2) in X} f(x_1, x_2)}) and
(\ref{sum_{x_1 in X_1} (sum_{x_2 in X_2} f(x_1, x_2))}) are equal to
each other, and the equality of (\ref{sum_{(x_1, x_2) in X} f(x_1,
  x_2)}) and (\ref{sum_{x_2 in X_2} (sum_{x_1 in X_1} f(x_1, x_2))})
is similar.

        Now let $k$ be a field, and let $V$ be a vector space over $k$.
If $f(x_1, x_2)$ is a $V$-valued function on $X$ with finite support,
then all of the sums mentioned in the preceding paragraph can be
defined as elements of $V$.  More precisely, (\ref{sum_{x_2 in X_2}
  f(x_1, x_2)}) is equal to $0$ for all but finitely many $x_1 \in
X_1$, and (\ref{sum_{x_1 in X_1} f(x_1, x_2)}) is equal to $0$ for all
but finitely many $x_2 \in X_2$, so that (\ref{sum_{x_1 in X_1}
  (sum_{x_2 in X_2} f(x_1, x_2))}) and (\ref{sum_{x_2 in X_2}
  (sum_{x_1 in X_1} f(x_1, x_2))}) are defined.  Of course,
(\ref{sum_{(x_1, x_2) in X} f(x_1, x_2)}), (\ref{sum_{x_1 in X_1}
  (sum_{x_2 in X_2} f(x_1, x_2))}), and (\ref{sum_{x_2 in X_2}
  (sum_{x_1 in X_1} f(x_1, x_2))}) are equal under these conditions.
This is basically the same as the case where $X_1$ and $X_2$ are both
finite sets.

        Suppose that $|\cdot|$ is an absolute value function on $k$,
$N$ is a norm on $V$ with respect to $|\cdot|$, and that $V$ is complete
with respect to the metric corresponding to $N$.  Let $f(x_1, x_2)$ be
a summable $V$-valued function on $X$, so that $N(f(x_1, x_2))$ is
summable as a nonnegative real-valued function on $X$.  Thus
(\ref{sum_{(x_1, x_2) in X} f(x_1, x_2)}) is defined as an element of
$V$, as in Section \ref{summable functions}.  Note that $f(x_1, x_2)$
is also summable as a $V$-valued function of $x_2 \in X_2$ for every
$x_1 \in X_1$, and as a $V$-valued function of $x_1 \in X_1$ for every
$x_2 \in X_2$.  This implies that (\ref{sum_{x_2 in X_2} f(x_1, x_2)})
is defined as an element of $V$ for every $x_1 \in X_1$, and that
(\ref{sum_{x_1 in X_1} f(x_1, x_2)}) is defined as an element of $V$
for every $x_2 \in X_2$.  We also have that
\begin{equation}
\label{N(sum_{x_2 in X_2} f(x_1, x_2)) le sum_{x_2 in X_2} N(f(x_1, x_2))}
 N\Big(\sum_{x_2 \in X_2} f(x_1, x_2)\Big) \le \sum_{x_2 \in X_2} N(f(x_1, x_2))
\end{equation}
for every $x_1 \in X_1$, and that
\begin{equation}
\label{N(sum_{x_1 in X_1} f(x_1, x_2)) le sum_{x_1 in X_1} N(f(x_1, x_2))}
 N\Big(\sum_{x_1 \in X_1} f(x_1, x_2)\Big) \le \sum_{x_1 \in X_1} N(f(x_1, x_2))
\end{equation}
for every $x_2 \in X_2$, as in (\ref{N(sum_{x in X} f(x)) le sum_{x in
    X} N(f(x)) = ||f||_1}).  It follows that
\begin{equation}
\label{sum_{x_1 in X_1} N(sum_{x_2 in X_2} f(x_1, x_2)) le ...}
        \sum_{x_1 \in X_1} N\Big(\sum_{x_2 \in X_2} f(x_1, x_2)\Big)
         \le \sum_{x_1 \in X_1} \Big(\sum_{x_2 \in X_2} N(f(x_1, x_2))\Big),
\end{equation}
and that
\begin{equation}
\label{sum_{x_2 in X_2} N(sum_{x_1 in X_1} f(x_1, x_2)) le ...}
        \sum_{x_2 \in X_2} N\Big(\sum_{x_1 \in X_1} f(x_1, x_2)\Big)
          \le \sum_{x_2 \in X_2} \Big(\sum_{x_1 \in X_1} N(f(x_1, x_2))\Big).
\end{equation}
Thus (\ref{sum_{x_2 in X_2} f(x_1, x_2)}) is summable as a $V$-valued
function of $x_1 \in X_1$, and (\ref{sum_{x_1 in X_1} f(x_1, x_2)}) is
summable as a $V$-valued function of $x_2 \in X_2$, by the remarks
about nonnegative real-valued functions on $X$ at the beginning of the
section.  This implies that (\ref{sum_{x_1 in X_1} (sum_{x_2 in X_2}
  f(x_1, x_2))}) and (\ref{sum_{x_2 in X_2} (sum_{x_1 in X_1} f(x_1,
  x_2))}) are defined as elements of $V$, as in Section \ref{summable
  functions}.  One can check that (\ref{sum_{(x_1, x_2) in X} f(x_1,
  x_2)}), (\ref{sum_{x_1 in X_1} (sum_{x_2 in X_2} f(x_1, x_2))}), and
(\ref{sum_{x_2 in X_2} (sum_{x_1 in X_1} f(x_1, x_2))}) are equal to
each other under these conditions.  This follows from the remarks in
the preceding paragraph when $f$ has finite support in $X$, and
otherwise one can approximate $f$ by functions with finite support in
$X$ with respect to the $\ell^1$ norm, as in Section \ref{summable
  functions}.

        Let us now restrict our attention to the case where $|\cdot|$
is an ultrametric absolute value function on $k$, and $N$ is an
ultranorm on $V$.  We continue to ask that $V$ be complete with
respect to the ultrametric corresponding to $N$.  If $f(x_1, x_2)$ is
a $V$-valued function on $X$ that vanishes at infinity, then
(\ref{sum_{(x_1, x_2) in X} f(x_1, x_2)}) can be defined as an element
of $V$, as in the previous section.  It is easy to see that $f(x_1,
x_2)$ also vanishes at infinity as a $V$-valued function of $x_2 \in
X_2$ for every $x_1 \in X_1$, and as a $V$-valued function of $x_1 \in
X_1$ for every $x_2 \in X_2$.  This implies that (\ref{sum_{x_2 in
    X_2} f(x_1, x_2)}) is defined as an element of $V$ for every $x_1
\in X_1$, and that (\ref{sum_{x_1 in X_1} f(x_1, x_2)}) is defined as
an element of $V$ for every $x_2 \in X_2$, as in the previous section.
We also have that
\begin{equation}
\label{N(sum_{x_2 in X_2} f(x_1, x_2)) le max_{x_2 in X_2} N(f(x_1, x_2))}
 N\Big(\sum_{x_2 \in X_2} f(x_1, x_2)\Big) \le \max_{x_2 \in X_2} N(f(x_1, x_2))
\end{equation}
for every $x_1 \in X_1$, and that
\begin{equation}
\label{N(sum_{x_1 in X_1} f(x_1, x_2)) le max_{x_1 in X_1} N(f(x_1, x_2))}
 N\Big(\sum_{x_1 \in X_1} f(x_1, x_2)\Big) \le \max_{x_1 \in X_1} N(f(x_1, x_2))
\end{equation}
for every $x_2 \in X_2$, by (\ref{N(sum_{x in X} f(x)) le max_{x in X}
  N(f(x)) = ||f||_infty}).  Using this, one can check that
(\ref{sum_{x_2 in X_2} f(x_1, x_2)}) vanishes at infinity as a
$V$-valued function of $x_1 \in X_1$, and that (\ref{sum_{x_1 in X_1}
  f(x_1, x_2)}) vanishes at infinity as a function of $x_2 \in X_2$,
because $f(x_1, x_2)$ vanishes at infinity on $X$.  This implies that
(\ref{sum_{x_1 in X_1} (sum_{x_2 in X_2} f(x_1, x_2))}) and
(\ref{sum_{x_2 in X_2} (sum_{x_1 in X_1} f(x_1, x_2))}) are defined as
elements of $V$, as in the previous section.  Moreover,
\begin{eqnarray}
\label{N(sum_{x_1 in X_1} (sum_{x_2 in X_2} f(x_1, x_2))) le ...}
        N\Big(\sum_{x_1 \in X_1} \Big(\sum_{x_2 \in X_2} f(x_1, x_2)\Big)\Big)
        & \le & \max_{x_1 \in X_1} N\Big(\sum_{x_2 \in X_2} f(x_1, x_2)\Big) \\
        & \le & \max_{(x_1, x_2) \in X} N(f(x_1, x_2)) \nonumber
\end{eqnarray}
and
\begin{eqnarray}
\label{N(sum_{x_2 in X_2} (sum_{x_1 in X_1} f(x_1, x_2))) le ...}
        N\Big(\sum_{x_2 \in X_2} \Big(\sum_{x_1 \in X_1} f(x_1, x_2)\Big)\Big)
        & \le & \max_{x_2 \in X_2} N\Big(\sum_{x_1 \in X_1} f(x_1, x_2)\Big) \\
        & \le & \max_{(x_1, x_2) \in X} N(f(x_1, x_2)), \nonumber
\end{eqnarray}
using (\ref{N(sum_{x in X} f(x)) le max_{x in X} N(f(x)) =
  ||f||_infty}) in the first steps in (\ref{N(sum_{x_1 in X_1}
  (sum_{x_2 in X_2} f(x_1, x_2))) le ...}) and (\ref{N(sum_{x_2 in
    X_2} (sum_{x_1 in X_1} f(x_1, x_2))) le ...}), and
(\ref{N(sum_{x_2 in X_2} f(x_1, x_2)) le max_{x_2 in X_2} N(f(x_1,
  x_2))}) and (\ref{N(sum_{x_1 in X_1} f(x_1, x_2)) le max_{x_1 in
    X_1} N(f(x_1, x_2))}) in the second steps, respectively.  As
before, one can show that (\ref{sum_{(x_1, x_2) in X} f(x_1, x_2)}),
(\ref{sum_{x_1 in X_1} (sum_{x_2 in X_2} f(x_1, x_2))}), and
(\ref{sum_{x_2 in X_2} (sum_{x_1 in X_1} f(x_1, x_2))}) are the same
under these conditions, by approximating $f$ by functions with finite
support in $X$ with respect to the $\ell^\infty$ norm.

\section{Generalized convergence}
\label{generalized convergence}

        Let $X$ be a nonempty set, let $k$ be a field, and let $V$
be a vector space over $k$.  Also let $|\cdot|$ be an absolute value
function on $k$, and let $N$ be a norm on $V$ with respect to
$|\cdot|$.  If $f$ is a $V$-valued function on $X$ and $A$ is a finite
subset of $X$, then the sum
\begin{equation}
\label{sum_{x in A} f(x), 2}
        \sum_{x \in A} f(x)
\end{equation}
can be defined as an element of $V$ in the usual way, which is
interpreted to be $0$ when $A = \emptyset$.  The family of these
finite sums may be considered as a net of elements of $V$, indexed by
the collection of finite subsets of $X$.  More precisely, the
collection of finite subsets of $X$ is partially ordered by inclusion.
If $A_1$, $A_2$ are finite subsets of $X$, then their union $A_1 \cup
A_2$ is a finite subset of $X$ that contains $A_1$ and $A_2$, which
implies that the collection of finite subsets of $X$ is a directed
system.  The convergence\index{generalized convergence} of the sum
\begin{equation}
\label{sum_{x in X} f(x), 4}
        \sum_{x \in X} f(x)
\end{equation}
in $V$ can be defined in terms of the convergence of the corresponding
net of finite subsums (\ref{sum_{x in A} f(x), 2}) in $V$.  This means
that there is a $v \in V$ with the property that for every $\epsilon >
0$ there is a finite set $A(\epsilon) \subseteq X$ such that
\begin{equation}
\label{N(sum_{x in A} f(x) - v) < epsilon}
        N\Big(\sum_{x \in A} f(x) - v\Big) < \epsilon
\end{equation}
for every finite set $A \subseteq X$ with $A(\epsilon) \subseteq A$.
It is easy to see that the limit $v$ of this net is unique when it
exists, in which case the value of the sum (\ref{sum_{x in X} f(x),
  4}) is defined to be $v$.

        Similarly, let us say that the sum (\ref{sum_{x in X} f(x), 4})
satisfies the \emph{generalized Cauchy criterion}\index{generalized Cauchy
criterion} if for each $\epsilon > 0$ there is a finite set $A_0(\epsilon)
\subseteq X$ such that
\begin{equation}
\label{N(sum_{x in B} f(x)) < epsilon}
        N\Big(\sum_{x \in B} f(x)\Big) < \epsilon
\end{equation}
for every finite set $B \subseteq X$ with $A_0(\epsilon) \cap B =
\emptyset$.  If the sum (\ref{sum_{x in X} f(x), 4}) converges in the
sense described in the preceding paragraph, then one can check that it
satisfies the generalized Cauchy criterion, in essentially the same
way as for ordinary infinite series.  More precisely, let
$A(\epsilon)$ be a finite subset of $X$ for which (\ref{N(sum_{x in A}
  f(x) - v) < epsilon}) holds, and let $B \subseteq X$ be a finite set
that is disjoint from $A(\epsilon)$.  Thus (\ref{N(sum_{x in A} f(x) -
  v) < epsilon}) can be applied to $A = A(\epsilon)$ and $A =
A(\epsilon) \cup B$, to get that
\begin{eqnarray}
\label{N(sum_{x in B} f(x)) = ... < epsilon + epsilon = 2 epsilon}
  \quad N\Big(\sum_{x \in B} f(x)\Big) & = &
 N\Big(\sum_{x \in A(\epsilon) \cup B} f(x) - \sum_{x \in A(\epsilon)} f(x)\Big) \\
         & \le & N\Big(\sum_{x \in A(\epsilon) \cup B} f(x) - v\Big)
                 + N\Big(\sum_{x \in A(\epsilon)} f(x) - v\Big) \nonumber \\
           & < & \epsilon + \epsilon = 2 \, \epsilon. \nonumber
\end{eqnarray}
This shows that (\ref{N(sum_{x in B} f(x)) < epsilon}) holds with
$A_0(\epsilon) = A(\epsilon/2)$, and one can take $A_0(\epsilon) =
A(\epsilon)$ when $N$ is an ultranorm on $V$.

        Suppose that (\ref{sum_{x in X} f(x), 4}) satisfies the generalized
Cauchy criterion, and let $A_0(\epsilon)$ be as in the preceding paragraph.
If $x \in X \setminus A_0(\epsilon)$, then we can take $B = \{x\}$ in
(\ref{N(sum_{x in B} f(x)) < epsilon}), to get that
\begin{equation}
\label{N(f(x)) < epsilon, 2}
        N(f(x)) < \epsilon.
\end{equation}
This implies that $f$ vanishes at infinity on $X$.  Conversely, if $N$
is an ultranorm on $V$, and if $f$ vanishes at infinity on $X$, then
it is easy to see that (\ref{sum_{x in X} f(x), 4}) satisfies the
generalized Cauchy criterion.

        Let $N$ be any norm on $V$ again, and suppose that $f$ is a
summable $V$-valued function on $X$, as in Section \ref{summable functions}.
Thus for each $\epsilon > 0$ there is a finite set $A_0(\epsilon) \subseteq X$
such that
\begin{equation}
\label{sum_{x in X setminus A_0(epsilon)} N(f(x)) < epsilon}
        \sum_{x \in X \setminus A_0(\epsilon)} N(f(x)) < \epsilon,
\end{equation}
as in (\ref{sum_{x in X} N(f(x)) < sum_{x in A(epsilon)} N(f(x)) + epsilon}).
If $B \subseteq X$ is a finite set that is disjoint from $A_0(\epsilon)$,
then we get that
\begin{equation}
\label{N(sum_{x in B} f(x)) le ... < epsilon}
        N\Big(\sum_{x \in B} f(x)\Big) \le \sum_{x \in B} N(f(x))
           \le \sum_{x \in X \setminus A_0(\epsilon)} N(f(x)) < \epsilon.
\end{equation}
This shows that (\ref{sum_{x in X} f(x), 4}) satisfies the generalized
Cauchy criterion under these conditions.

        If (\ref{sum_{x in X} f(x), 4}) satisfies the generalized Cauchy
criterion, then the norms of the finite subsums (\ref{sum_{x in A}
  f(x), 2}) are uniformly bounded, over all finite sets $A \subseteq
X$.  More precisely, let $A_0(1)$ be a finite subset of $X$ such that
(\ref{N(sum_{x in B} f(x)) < epsilon}) holds with $\epsilon = 1$ for
all finite sets $B \subseteq X \setminus A_0(1)$.  If $A$ is any
finite subset of $X$, then
\begin{eqnarray}
\label{N(sum_{x in A} f(x)) le ... le 1 + sum_{x in A_0(1)} N(f(x))}
        N\Big(\sum_{x \in A} f(x)\Big)
              & \le & N\Big(\sum_{x \in A \setminus A_0(1)} f(x)\Big)
                      + N\Big(\sum_{x \in A \cap A_0(1)} f(x)\Big) \\
              & \le & 1 + \sum_{x \in A_0(1)} N(f(x)), \nonumber
\end{eqnarray}
using (\ref{N(sum_{x in B} f(x)) < epsilon}) applied to $B = A
\setminus A_0(1)$ in the second step.

        Suppose that $k = {\bf R}$ with the standard absolute value
function, and that $V = {\bf R}$ with the standard absolute value
function as the norm.  If $f$ is a real-valued function on $X$ such
that (\ref{sum_{x in X} f(x), 4}) satisfies the generalized Cauchy
criterion, then the finite subsums (\ref{sum_{x in A} f(x), 2}) are
uniformly bounded, as in the previous paragraph.  In particular, this
can be applied to finite sets $A \subseteq X$ such that $f(x) \ge 0$
for every $x \in X$, or such that $f(x) \le 0$ for every $x \in X$.
Using this, one can check that $f$ is actually summable on $X$ in this
case.  The analogous statement also holds when $k = V = {\bf C}$,
with the standard absolute value function, by applying the previous
argument to the real and imaginary parts of a complex-valued function
on $X$.

        Let $k$, $V$ be arbitary again, and suppose that
(\ref{sum_{x in X} f(x), 4}) satisfies the generalized Cauchy criterion.
Thus $f$ vanishes at infinity on $X$, as before, which implies that
the support of $f$ has only finitely or countably many elements.  If
the support of $f$ is finite, then (\ref{sum_{x in X} f(x), 4})
obviously converges.  Otherwise, suppose that the support of $f$ is
countably infinite, and let $\{x_j\}_{j = 1}^\infty$ be a sequence of
distinct elements of $X$ that includes every element in the support of
$f$.  It is easy to see that the partial sums
\begin{equation}
\label{sum_{j = 1}^n f(x_j)}
        \sum_{j = 1}^n f(x_j)
\end{equation}
of the infinite series
\begin{equation}
\label{sum_{j = 1}^infty f(x_j)}
        \sum_{j = 1}^\infty f(x_j)
\end{equation}
form a Cauchy sequence in $V$ under these conditions, because
(\ref{sum_{x in X} f(x), 4}) satisfies the generalized Cauchy
criterion.  If $V$ is complete with respect to the metric associated
to $N$, then it follows that this sequence of partial sums converges
in $V$, which is to say that (\ref{sum_{j = 1}^infty f(x_j)})
converges in $V$ in the usual sense.  In this case, it is easy to see
that (\ref{sum_{x in X} f(x), 4}) also converges in $V$ in the sense
described at the beginning of the section, with the same value of the
sum.  This uses the generalized Cauchy criterion again, to get that
finite sums of the form (\ref{sum_{x in A} f(x), 2}) are close to
partial sums of the form (\ref{sum_{j = 1}^n f(x_j)}) under suitable
conditions.

\chapter{Power series}
\label{power series}

\section{Complex coefficients}
\label{complex coefficients}

        Let $a_0, a_1, a_2, a_3, \ldots$ be a sequence of complex numbers,
and consider the power series
\begin{equation}
\label{sum_{j = 0}^infty a_j z^j}
        \sum_{j = 0}^\infty a_j \, z^j,
\end{equation}
where $z \in {\bf C}$.  As usual, $z^j$ is interpreted as being equal
to $1$ for every $z \in {\bf C}$ when $j = 0$.  If (\ref{sum_{j =
    0}^infty a_j z^j}) converges for some $z \in {\bf C}$, then $\{a_j
\, z^j\}_{j = 0}^\infty$ converges to $0$ as a sequence in ${\bf C}$,
which implies in particular that $\{a_j \, z^j\}_{j = 0}^\infty$ is
bounded.  Using this, one can check that $\sum_{j = 0}^\infty a_j \,
w^j$ converges absolutely for every $w \in {\bf C}$ such that $|w| <
|z|$, by comparison with a convergent geometric series.

        The \emph{radius of convergence}\index{radius of convergence} $\rho$
of (\ref{sum_{j = 0}^infty a_j z^j}) may be defined by
\begin{equation}
\label{rho = sup {|z| : z in {bf C}, sum_{j = 0}^infty a_j z^j converges}}
 \rho = \sup \bigg\{|z| : z \in {\bf C}, \ \sum_{j = 0}^\infty a_j \, z^j
                                            \hbox{ converges}\bigg\}.
\end{equation}
Of course, (\ref{sum_{j = 0}^infty a_j z^j}) converges trivially when
$z = 0$, so that (\ref{rho = sup {|z| : z in {bf C}, sum_{j = 0}^infty
    a_j z^j converges}}) is the supremum of a nonempty set.  If
(\ref{sum_{j = 0}^infty a_j z^j}) converges for $z \in {\bf C}$ with
arbitrarily large modulus, then (\ref{rho = sup {|z| : z in {bf C},
    sum_{j = 0}^infty a_j z^j converges}}) is interpreted as being
$+\infty$, as usual.  If $w \in {\bf C}$ satisfies $|w| < \rho$, then
it is easy to see that $\sum_{j = 0}^\infty a_j \, w^j$ converges
absolutely, using the remarks in the previous paragraph.  If $w \in
{\bf C}$ satisfies $|w| > \rho$, then $\sum_{j = 0}^\infty a_j \, w^j$
does not converge, by the definition (\ref{rho = sup {|z| : z in {bf
      C}, sum_{j = 0}^infty a_j z^j converges}}) of $\rho$.  Note that
$\rho$ is uniquely determined by these two properties.  It is well
known that
\begin{equation}
\label{rho = (limsup_{j to infty} |a_j|^{1/j})^{-1}}
        \rho = \Big(\limsup_{j \to \infty} |a_j|^{1/j}\Big)^{-1},
\end{equation}
by the root test, with the standard conventions that $1/0 = +\infty$
and $1/+\infty = 0$.  

        Suppose that
\begin{equation}
\label{sum_{j = 0}^infty |a_j| r^j}
        \sum_{j = 0}^\infty |a_j| \, r^j
\end{equation}
converges for some nonnegative real number $r$.  This implies that
(\ref{sum_{j = 0}^infty a_j z^j}) converges absolutely for every $z
\in {\bf C}$ with $|z| \le r$, by the comparison test.  Moreover,
\begin{equation}
\label{|sum_{j = 0}^infty a_j z^j - sum_{j = 0}^n a_j z^j| = ...}
 \biggl|\sum_{j = 0}^\infty a_j \, z^j - \sum_{j = 0}^n a_j \, z^j\biggr|
 = \biggl|\sum_{j = n + 1}^\infty a_j \, z^j\biggr|
 \le \sum_{j = n + 1}^\infty |a_j| \, |z|^j \le \sum_{j = n + 1}^\infty |a_j| \, r^j
\end{equation}
for each $z \in {\bf C}$ with $|z| \le r$ and nonnegative integer $n$.
It follows that the partial sums
\begin{equation}
\label{sum_{j = 0}^n a_j z^j}
        \sum_{j = 0}^n a_j \, z^j
\end{equation}
converge to (\ref{sum_{j = 0}^infty a_j z^j}) uniformly on the closed
disk
\begin{equation}
\label{{z in {bf C} : |z| le r}}
        \{z \in {\bf C} : |z| \le r\}
\end{equation}
under these conditions.  This shows that (\ref{sum_{j = 0}^infty a_j
  z^j}) defines a continuous function on (\ref{{z in {bf C} : |z| le
    r}}), since the partial sums (\ref{sum_{j = 0}^n a_j z^j}) are
continuous.

        Similarly, if the radius of convergence $\rho$ of
(\ref{sum_{j = 0}^infty a_j z^j}) is positive, then
(\ref{sum_{j = 0}^infty a_j z^j}) defines a continuous function
on the open disk
\begin{equation}
\label{{z in {bf C} : |z| < rho}}
        \{z \in {\bf C} : |z| < \rho\}.
\end{equation}
More precisely, let $z_0 \in {\bf C}$ with $|z_0| < \rho$ be given, and
let $r$ be a positive real number such that
\begin{equation}
\label{|z_0| < r < rho}
        |z_0| < r < \rho.
\end{equation}
Thus (\ref{sum_{j = 0}^infty |a_j| r^j}) converges, because $r <
\rho$, so that (\ref{sum_{j = 0}^infty a_j z^j}) defines a continuous
function on (\ref{{z in {bf C} : |z| le r}}), as in the previous
paragraph.  In particular, (\ref{sum_{j = 0}^infty a_j z^j}) is
continuous at $z_0$ as a function on (\ref{{z in {bf C} : |z| le r}}),
which implies that (\ref{sum_{j = 0}^infty a_j z^j}) is continuous at
$z_0$ as a function on (\ref{{z in {bf C} : |z| < rho}}), because
$|z_0| < r$.  It follows that (\ref{sum_{j = 0}^infty a_j z^j}) is
continuous as a function on (\ref{{z in {bf C} : |z| < rho}}) at every
point in (\ref{{z in {bf C} : |z| < rho}}), as desired.

        Consider the power series
\begin{equation}
\label{sum_{j = 1}^infty j a_j z^{j - 1}}
        \sum_{j = 1}^\infty j \, a_j \, z^{j - 1},
\end{equation}
obtained by differentiating (\ref{sum_{j = 0}^infty a_j z^j}) term by
term.  It is well known that the radius of convergence of (\ref{sum_{j
    = 1}^infty j a_j z^{j - 1}}) is equal to the radius of convergence
of (\ref{sum_{j = 0}^infty a_j z^j}), which can be derived from
(\ref{rho = (limsup_{j to infty} |a_j|^{1/j})^{-1}}), for instance.
Alternatively, if (\ref{sum_{j = 1}^infty j a_j z^{j - 1}}) converges
absolutely for some $z \in {\bf C}$, then (\ref{sum_{j = 0}^infty a_j
  z^j}) converges absolutely as well, by the comparison test.
Conversely, if $z \in {\bf C}$ and $|z| < \rho$, then (\ref{sum_{j =
    0}^infty a_j z^j}) converges by comparison with a convergent
geometric series, in which case (\ref{sum_{j = 1}^infty j a_j z^{j -
    1}}) converges absolutely too.  It is also well known that
(\ref{sum_{j = 0}^infty a_j z^j}) is a holomorphic function on
(\ref{{z in {bf C} : |z| < rho}}), whose complex derivative is given
by (\ref{sum_{j = 1}^infty j a_j z^{j - 1}}).

        The complex exponential function\index{exponential function}
may be defined for $z \in {\bf C}$ by
\begin{equation}
\label{E(z) = sum_{j = 0}^infty frac{z^j}{j!}}
        E(z) = \sum_{j = 0}^\infty \frac{z^j}{j!},
\end{equation}
where $j!$ is $j$ factorial, the product of the positive integers from
$1$ to $j$, which is interpreted as being equal to $1$ when $j = 0$.
It is easy to see that this series converges absolutely for every $z
\in {\bf C}$, using the ratio test, for instance.  If $w, z \in {\bf C}$,
then
\begin{equation}
\label{E(w + z) = sum_{l = 0}^infty frac{(w + z)^l}{l!} = ...}
        E(w + z) = \sum_{l = 0}^\infty \frac{(w + z)^l}{l!}
 = \sum_{l = 0}^\infty \sum_{j = 0}^l \frac{w^j \, z^{l - j}}{j! \, (l - j)!},
\end{equation}
using the binomial theorem in the second step.  This implies that
\begin{equation}
\label{E(w + z) = E(w) E(z)}
        E(w + z) = E(w) \, E(z)
\end{equation}
for every $w, z \in {\bf C}$, because the right side of (\ref{E(w + z)
  = sum_{l = 0}^infty frac{(w + z)^l}{l!} = ...}) is the same as the
Cauchy product of the series representing $E(w)$ and $E(z)$.  More
precisely, this uses the absolute convergence of the series
representing $E(w)$ and $E(z)$, to ensure that the right side of
(\ref{E(w + z) = sum_{l = 0}^infty frac{(w + z)^l}{l!} = ...})  is
equal to the right side of (\ref{E(w + z) = E(w) E(z)}).

        If we take $w = -z$ in (\ref{E(w + z) = E(w) E(z)}), then we get
that
\begin{equation}
\label{E(z) E(-z) = E(0) = 1}
        E(z) \, E(-z) = E(0) = 1
\end{equation}
for every $z \in {\bf C}$.  Equivalently, this means that $E(z) \ne 0$
for every $z \in {\bf C}$, and that $E(-z) = 1/E(z)$.  Observe that
\begin{equation}
\label{overline{E(z)} = E(overline{z})}
        \overline{E(z)} = E(\overline{z})
\end{equation}
for every $z \in {\bf C}$, where $\overline{z}$ is the complex
conjugate of $z$, since one can take the complex-conjugate of
(\ref{E(z) = sum_{j = 0}^infty frac{z^j}{j!}}) term by term.
This implies that
\begin{equation}
\label{E(2 re z) = E(z + overline{z}) = E(z) E(overline{z}) = ... = |E(z)|^2}
        E(2 \re z) = E(z + \overline{z}) = E(z) \, E(\overline{z})
                              = E(z) \, \overline{E(z)} = |E(z)|^2
\end{equation}
for every $z \in {\bf C}$, where $\re z$ denotes the real part of $z$.
In particular,
\begin{equation}
\label{|E(i y)| = 1}
        |E(i \, y)| = 1
\end{equation}
for every $y \in {\bf R}$.

        Note that $E(x) \in {\bf R}$ for every $x \in {\bf R}$,
by (\ref{E(z) = sum_{j = 0}^infty frac{z^j}{j!}}).  More precisely,
$E(x) \ge 1$ when $x \ge 0$, which implies that $0 < E(x) \le 1$ when
$x \le 0$, because $E(x) = 1/E(-x)$.  Of course, $E(x)$ is the same as
the usual real exponential function on ${\bf R}$ when $x \in {\bf R}$.
It is easy to see that $E(x)$ is strictly increasing on $[0,
  +\infty)$, and that $E(x) \to +\infty$ as $x \to +\infty$, directly
  from (\ref{E(z) = sum_{j = 0}^infty frac{z^j}{j!}}).  This implies
that $E(x)$ is also strictly increasing on $(-\infty, 0]$, and
hence on ${\bf R}$, and that $E(x) \to 0$ as $x \to -\infty$,
using $E(x) = 1/E(-x)$ again.

        If the coefficients $a_j$ of (\ref{sum_{j = 0}^infty a_j z^j})
are real numbers, then one can think of (\ref{sum_{j = 0}^infty a_j
  z^j}) as a power series on ${\bf R}$, with many of the same
properties as before.  In particular, the radius of convergence of
(\ref{sum_{j = 0}^infty a_j z^j}) as a power series on ${\bf R}$ is
the same as the radius of convergence of (\ref{sum_{j = 0}^infty a_j
  z^j}) as a power series on ${\bf C}$.  If the radius of convergence
$\rho$ is positive, then (\ref{sum_{j = 0}^infty a_j z^j}) defines a
continuous real-valued function on the open interval $(-\rho, \rho)$
in ${\bf R}$.  This function is also differentiable on $(-\rho,
\rho)$, with the derivative given by the power series (\ref{sum_{j =
    1}^infty j a_j z^{j - 1}}).

\section{Ultrametric absolute value functions}
\label{ultrametric absolute value functions}

        Let $k$ be a field, and let $a_0, a_1, a_2, a_3, \ldots$ be a
sequence of elements of $k$.  Consider the corresponding formal power
series
\begin{equation}
\label{f(X) = sum_{j = 0}^infty a_j X^j}
        f(X) = \sum_{j = 0}^\infty a_j \, X^j,
\end{equation}
where $X$ is an indeterminate.  As in \cite{cas, fg}, we shall use
upper-case letters like $X$ for indeterminates, and lower-case letters
like $x$ for elements of $k$ or other fields.  Let $|\cdot|$ be an
ultrametric absolute value function on $k$, and suppose that $k$ is
complete with respect to the ultrametric that corresponds to
$|\cdot|$, as in (\ref{d(x, y) = |x - y|, 2}).  If $x \in k$, then
\begin{equation}
\label{sum_{j = 0}^infty a_j x^j}
        \sum_{j = 0}^\infty a_j \, x^j,
\end{equation}
converges in $k$ exactly when $\{a_j \, x^j\}_{j = 0}^\infty$
converges to $0$ in $k$, as in Section \ref{infinite series}.
Equivalently, this means that
\begin{equation}
\label{|a_j x^j| = |a_j| |x|^j to 0 as j to infty}
        |a_j \, x^j| = |a_j| \, |x|^j \to 0 \quad\hbox{as } j \to \infty
\end{equation}
as a sequence of nonnegative real numbers.  In this case, the value of
(\ref{sum_{j = 0}^infty a_j x^j}) may be denoted $f(x)$.

        The radius of convergence\index{radius of convergence} $\rho$
of this power series may be defined as the supremum of the set of
nonnegative real numbers $r$ such that
\begin{equation}
\label{lim_{j to infty} |a_j| r^j = 0}
        \lim_{j \to \infty} |a_j| \, r^j = 0.
\end{equation}
Of course, $r = 0$ automatically has this property, so that this is
the supremum of a nonempty set.  If (\ref{lim_{j to infty} |a_j| r^j =
  0}) holds for some arbitrarily large real numbers $r$, then the
supremum is interpreted as being equal to $+\infty$, as usual.  If
$r$, $t$ are nonnegative real numbers such that $t \le r$ and $r$
satisfies (\ref{lim_{j to infty} |a_j| r^j = 0}), then
\begin{equation}
\label{lim_{j to infty} |a_j| t^j = 0}
        \lim_{j \to \infty} |a_j| \, t^j = 0
\end{equation}
as well.  This implies that (\ref{lim_{j to infty} |a_j| t^j = 0})
holds when $t < \rho$, by the definition of $\rho$.  If $t > \rho$,
then (\ref{lim_{j to infty} |a_j| t^j = 0}) does not hold, again by
the definition of $\rho$.  It is easy to see that $\rho$ is uniquely
determined by these two properties.  One can also check that
\begin{equation}
\label{rho = (limsup_{j to infty} |a_j|^{1/j})^{-1}, 2}
        \rho = \Big(\limsup_{j \to \infty} |a_j|^{1/j}\Big)^{-1},
\end{equation}
using standard properties of the limsup, and with the usual
conventions for $1/0$ and $1/\infty$.  It follows from these
properties of $\rho$ that (\ref{sum_{j = 0}^infty a_j x^j}) converges
in $k$ when $x \in k$ satisfies $|x| < \rho$, and not when $|x| >
\rho$.  However, this may not determine $\rho$ uniquely, depending on
the possible values of $|\cdot|$ on $k$.

        If (\ref{sum_{j = 0}^infty a_j x^j}) converges for some $x \in k$,
then
\begin{equation}
\label{|sum_{j = 0}^infty a_j x^j| le max_{j ge 0} |a_j x^j| = ...}
 \biggl|\sum_{j = 0}^\infty a_j \, x^j\biggr| \le \max_{j \ge 0} |a_j \, x^j|
                                              = \max_{j \ge 0} (|a_j| \, |x|^j),
\end{equation}
as in (\ref{|sum_{j = 1}^infty a_j| le max_{j ge 1} |a_j|}) in Section
\ref{infinite series}.  Suppose now that $r$ is a positive real number
that satisfies (\ref{lim_{j to infty} |a_j| r^j = 0}), which implies
that (\ref{|a_j x^j| = |a_j| |x|^j to 0 as j to infty}) holds for
every $x \in k$ with $|x| \le r$.  Thus (\ref{sum_{j = 0}^infty a_j
  x^j}) converges when $|x| \le r$, and
\begin{eqnarray}
\label{|sum_{j = 0}^infty a_j x^j - sum_{j = 0}^n a_j x^j| = ...}
 \quad \biggl|\sum_{j = 0}^\infty a_j \, x^j - \sum_{j = 0}^n a_j \, x^j\biggr|
       & = & \biggl|\sum_{j = n + 1}^\infty a_j \, x^j\biggr| \\
 & \le & \max_{j \ge n + 1} (|a_j| \, |x|^j) \le \max_{j \ge n + 1} (|a_j| \, r^j)
                                                            \nonumber
\end{eqnarray}
for every nonnegative integer $n$, by (\ref{|sum_{j = 0}^infty a_j
  x^j| le max_{j ge 0} |a_j x^j| = ...}).  It follows that the partial sums
\begin{equation}
\label{sum_{j = 0}^n a_j x^j}
        \sum_{j = 0}^n a_j \, x^j
\end{equation}
converge to (\ref{sum_{j = 0}^infty a_j x^j}) uniformly on the set
\begin{equation}
\label{{x in k : |x| le r}}
        \{x \in k : |x| \le r\},
\end{equation}
because the right side of (\ref{|sum_{j = 0}^infty a_j x^j - sum_{j =
    0}^n a_j x^j| = ...})  tends to $0$ as $n \to \infty$, by
(\ref{lim_{j to infty} |a_j| r^j = 0}).  This implies that
(\ref{sum_{j = 0}^infty a_j x^j}) defines a continuous function on
(\ref{{x in k : |x| le r}}) under these conditions, because the
partial sums (\ref{sum_{j = 0}^n a_j x^j}) are continuous.

        Using this, one can check that (\ref{sum_{j = 0}^infty a_j x^j})
defines a continuous function on
\begin{equation}
\label{{x in k : |x| < rho}}
        \{x \in k : |x| < \rho\}
\end{equation}
when $\rho > 0$.  Of course, if $\rho < +\infty$ and
\begin{equation}
\label{lim_{j to infty} |a_j| rho^j = 0}
        \lim_{j \to \infty} |a_j| \, \rho^j = 0,
\end{equation}
then one might as well apply the previous discussion to $r = \rho$,
to get that (\ref{sum_{j = 0}^infty a_j x^j}) is continuous on
\begin{equation}
\label{{x in k : |x| le rho}}
        \{x \in k : |x| \le \rho\}.
\end{equation}
Otherwise, one can apply the previous discussion to each $r > 0$ such
that $r < \rho$.  More precisely, if $x_0 \in k$ satisfies $|x_0| <
\rho$, then one can choose $r > 0$ such that $|x_0| \le r < \rho$.
One can also choose $r$ so that $r > |x_0|$ when $x_0 \ne 0$, but this
is not really necessary here, because (\ref{{x in k : |x| le r}}) is
an open set in $k$ when $r > 0$ and $|\cdot|$ is an ultrametric
absolute value function.

\section{Differentiation and Lipschitz conditions}
\label{differentiation, lipschitz conditions}

        Let $k$ be a field, and let $|\cdot|$ be an ultrametric absolute
value function on $k$.  Also let $r$ be a positive real number, and
suppose that $x, y \in k$ satisfy $|x|, |y| \le r$, so that $|x - y|
\le r$ too.  If $j$ is a positive integer, then
\begin{equation}
\label{x^j - y^j = ((x - y) + y)^j - y^j = ...}
        x^j - y^j = ((x - y) + y)^j - y^j 
               = \sum_{l = 1}^j {j \choose l} \cdot (x - y)^l \, y^{j - l},
\end{equation}
by the binomial theorem.  This implies that
\begin{equation}
\label{|x^j - y^j| le |x - y| r^{j - 1}}
        |x^j - y^j| \le |x - y| \, r^{j - 1},
\end{equation}
because of the ultrametric version of the triangle inequality, which
implies in particular that $|n \cdot 1| \le 1$ for every positive
integer $n$.  Similarly, if $j \ge 2$, then
\begin{eqnarray}
\label{x^j - y^j - j y^{j - 1} (x - y) = ...}
        x^j - y^j - j \, y^{j - 1} \, (x - y) 
          & = & ((x - y) + y)^j - y^j - j \, y^{j - 1} \, (x - y) \\
 & = & \sum_{l = 2}^j {j \choose l} \cdot (x - y)^l \, y^{j - l}, \nonumber
\end{eqnarray}
and hence
\begin{equation}
\label{|x^j - y^j - j y^{j - 1} (x - y)| le |x - y|^2 r^{j - 2}}
        |x^j - y^j - j \, y^{j - 1} \, (x - y)| \le |x - y|^2 \, r^{j - 2}.
\end{equation}

        Suppose that $k$ is complete with respect to the ultrametric
that corresponds to $|\cdot|$.  Let $a_0, a_1, a_2, a_3, \ldots$ be a
sequence of elements of $k$ that satisfies
\begin{equation}
\label{lim_{j to infty} |a_j| r^j = 0, 2}
        \lim_{j \to \infty} |a_j| \, r^j = 0,
\end{equation}
and put
\begin{equation}
\label{f(x) = sum_{j = 0}^infty a_j x^j}
        f(x) = \sum_{j = 0}^\infty a_j \, x^j
\end{equation}
for each $x \in k$ with $|x| \le r$.  As before, the convergence of
this series follows from the completeness of $k$, although this is not
needed when $a_j = 0$ for all but finitely many $j$, in which case
$f(x)$ is a polynomial function.  Note that
\begin{equation}
\label{|f(x)| le max_{j ge 0} (|a_j| r^j)}
        |f(x)| \le \max_{j \ge 0} (|a_j| \, r^j)
\end{equation}
for every $x \in k$ with $|x| \le r$, as in (\ref{|sum_{j = 0}^infty
  a_j x^j| le max_{j ge 0} |a_j x^j| = ...}).

        Now let $x, y \in k$ be given, with $|x|. |y| \le r$, so that the
series expansions for both $f(x)$ and $f(y)$ converge.  Thus
\begin{equation}
\label{f(x) - f(y) = sum_{j = 1}^infty a_j (x^j - y^j)}
        f(x) - f(y) = \sum_{j = 1}^\infty a_j \, (x^j - y^j),
\end{equation}
and hence
\begin{equation}
\label{|f(x) - f(y)| le max_{j ge 1} (|a_j| |x^j - y^j|)}
        |f(x) - f(y)| \le \max_{j \ge 1} (|a_j| \, |x^j - y^j|),
\end{equation}
as in (\ref{|sum_{j = 1}^infty a_j| le max_{j ge 1} |a_j|}) in Section
\ref{infinite series}.  Combining this with (\ref{|x^j - y^j| le |x -
  y| r^{j - 1}}), we get that
\begin{equation}
\label{|f(x) - f(y)| le max_{j ge 1} (|a_j| r^{j - 1}) |x - y|}
        |f(x) - f(y)| \le \max_{j \ge 1} (|a_j| \, r^{j - 1}) \, |x - y|.
\end{equation}

        If $f(X)$ is the formal power series associated to this
sequence of coefficients, as in (\ref{f(X) = sum_{j = 0}^infty a_j X^j}),
then the formal derivative of $f(X)$ is the formal power series
\begin{equation}
\label{f'(X) = sum_{j = 1}^infty j cdot a_j X^{j - 1}}
        f'(X) = \sum_{j = 1}^\infty j \cdot a_j \, X^{j - 1}.
\end{equation}
Of course,
\begin{equation}
\label{|j cdot a_j x^{j - 1}| le |a_j| |x|^{j - 1}}
        |j \cdot a_j \, x^{j - 1}| \le |a_j| \, |x|^{j - 1}
\end{equation}
for every $x \in k$ and positive integer $j$, because $|j \cdot 1| \le
1$.  If $|x| \le r$, then
\begin{equation}
\label{lim_{j to infty} j cdot a_j x^{j -1} = 0}
        \lim_{j \to \infty} j \cdot a_j \, x^{j -1} = 0
\end{equation}
in $k$, by (\ref{lim_{j to infty} |a_j| r^j = 0, 2}), so that 
\begin{equation}
\label{sum_{j = 1}^infty j cdot a_j x^{j - 1}}
        \sum_{j = 1}^\infty j \cdot a_j \, x^{j - 1}.
\end{equation}
converges in $k$.  Let $f'(x)$ denote the value of the sum
(\ref{sum_{j = 1}^infty j cdot a_j x^{j - 1}}), so that
\begin{equation}
\label{|f'(x)| le max_{j ge 1} (|j cdot a_j| r^{j - 1}) le ...}
        |f'(x)| \le \max_{j \ge 1} (|j \cdot a_j| \, r^{j - 1})
                  \le \max_{j \ge 1} (|a_j| \, r^{j - 1})
\end{equation}
for every $x \in k$ with $|x| \le r$, which is the analogue of
(\ref{|f(x)| le max_{j ge 0} (|a_j| r^j)}) for $f'(x)$.

        If $x, y \in k$ satisfy $|x|, |y| \le r$, then
\begin{eqnarray}
\label{f(x) - f(y) - f'(y) (x - y) = ...}
\lefteqn{f(x) - f(y) - f'(y) \, (x - y)}                         \\
 & = & \sum_{j = 1}^\infty a_j (x^j - y^j)
             - \sum_{j = 1}^\infty j \cdot a_j \, y^{j - 1} \, (x - y) \nonumber \\
 & = & \sum_{j = 2}^\infty a_j \, (x^j - y^j - j \cdot y^{j - 1} \, (x - y)).
                                                                    \nonumber
\end{eqnarray}
This implies that
\begin{equation}
\label{|f(x) - f(y) - f'(y) (x - y)| le ...}
        \quad |f(x) - f(y) - f'(y) \, (x - y)|
         \le \max_{j \ge 2} (|a_j| \, |x^j - y^j - j \cdot y^{j - 1} \, (x - y)|),
\end{equation}
as in (\ref{|sum_{j = 1}^infty a_j| le max_{j ge 1} |a_j|}) in Section
\ref{infinite series}.  It follows that
\begin{equation}
\label{|f(x) - f(y) - f'(y) (x - y)| le ..., 2}
        |f(x) - f(y) - f'(y) \, (x - y)|
               \le \max_{j \ge 2} (|a_j| \, r^{j - 2}) \, |x - y|^2,
\end{equation}
by (\ref{|x^j - y^j - j y^{j - 1} (x - y)| le |x - y|^2 r^{j - 2}}).
In particular, this shows that $f'(y)$ is the derivative of $f(y)$ in
the usual sense when $|\cdot|$ is not the trivial absolute value
function on $k$.

        Observe that
\begin{eqnarray}
\label{|f(x) - f(y)| le ...}
 |f(x) - f(y)| & \le & \max\Big(|f'(x)| \, |x - y|,
                      \, \max_{j \ge 2} (|a_j| \, r^{j - 2}) \, |x - y|^2\Big) \\
& = & \max\Big(|f'(y)|, \, \max_{j \ge 2} (|a_j| \, r^{j - 2}) \, |x - y|\Big)
                                                    \, |x - y| \nonumber
\end{eqnarray}
for every $x, y \in k$ with $|x|, |y| \le r$, by (\ref{|f(x) - f(y) -
  f'(y) (x - y)| le ...}) and the ultrametric version of the triangle
inequality.  Similarly,
\begin{equation}
\label{|f'(y)| |x - y| le ...}
        |f'(y)| \, |x - y| \le \max\Big(|f(x) - f(y)|,
                          \max_{j \ge 2} (|a_j| \, r^{j - 2}) \, |x - y|^2\Big)
\end{equation}
for every $x, y \in k$ with $|x|, |y| \le r$.  If we also have that
\begin{equation}
\label{max_{j ge 2} (|a_j| r^{j - 2}) |x - y| < |f'(y)|}
        \max_{j \ge 2} (|a_j| \, r^{j - 2}) \, |x - y| < |f'(y)|,
\end{equation}
then it follows that
\begin{equation}
\label{|f(x) - f(y)| = |f'(y)| |x - y|}
        |f(x) - f(y)| = |f'(y)| \, |x - y|.
\end{equation}
More precisely, this is trivial when $x = y$, and otherwise one can
multiply both sides of (\ref{max_{j ge 2} (|a_j| r^{j - 2}) |x - y| <
  |f'(y)|}) by $|x - y|$, and still have a strict inequality.

        If $x, y \in k$ satisfy $|x|, |y| \le r$ again, then
\begin{equation}
\label{|f'(x) - f'(y)| le max_{j ge 2} (|j cdot a_j| r^{j - 2}) |x - y| le ...}
 |f'(x) - f'(y)| \le \max_{j \ge 2} (|j \cdot a_j| \, r^{j - 2}) \, |x - y|
                  \le \max_{j \ge 2} (|a_j| \, r^{j - 2}) \, |x - y|,
\end{equation}
which is the analogue of (\ref{|f(x) - f(y)| le max_{j ge 1} (|a_j|
  r^{j - 1}) |x - y|}) for $f'(x)$ in place of $f(x)$.  This implies that
\begin{equation}
\label{|f'(x)| = |f'(y)|}
        |f'(x)| = |f'(y)|
\end{equation}
when (\ref{max_{j ge 2} (|a_j| r^{j - 2}) |x - y| < |f'(y)|}) holds,
by the ultrametric version of the triangle inequality.

        As mentioned earlier, the completeness of $k$ is only needed
in this section to ensure the convergence of the various infinite series.
If $a_j = 0$ for all but finitely many $j$, then completeness of $k$
is not needed, and the various estimates for $f(x)$ and $f'(x)$ still hold.

\section{Hensel's lemma}
\label{hensel's lemma}

        Let $k$ be a field, let $|\cdot|$ be an ultrametric absolute
value function on $k$, and suppose that $k$ is complete with respect
to the corresponding ultrametric.  Also let $r$ be a positive real
number, and let $a_0, a_1, a_2, a_3, \ldots$ be a sequence of elements
of $k$ that satisfies (\ref{lim_{j to infty} |a_j| r^j = 0, 2}).  Thus
$f(x)$ may be defined for $x \in k$ with $|x| \le r$ as in (\ref{f(x)
  = sum_{j = 0}^infty a_j x^j}).  Let $x_0 \in k$ be given, with
$|x_0| \le r$ and $f'(x_0) \ne 0$.  If $z \in k$ is sufficiently close
to $f(x_0)$, then we would like to find an $x \in k$ close to $x_0$
that satisfies
\begin{equation}
\label{f(x) = z}
        f(x) = z.
\end{equation}
In particular, if $|x - x_0| \le r$, then $|x| \le r$, so that $f(x)$
is defined.  In order to find $x$, we consider an appropriate sequence
of approximations.

        Suppose that the $l$th approximation $x_l \in k$ has been
chosen for some nonnegative integer $l$, in such a way that $|x_l| \le r$
and $f'(x_l) \ne 0$, where $f'$ is as defined in the previous section.
Let us choose $x_{l + 1} \in k$ so that
\begin{equation}
\label{f(x_l) + f'(x_l) (x_{l + 1} - x_l) = z}
        f(x_l) + f'(x_l) \, (x_{l + 1} - x_l) = z,
\end{equation}
which is to say that
\begin{equation}
\label{x_{l + 1} = x_l + f'(x_l)^{-1} (z - f(x_l))}
        x_{l + 1} = x_l + f'(x_l)^{-1} \, (z - f(x_l)).
\end{equation}
Thus
\begin{equation}
\label{|x_{l + 1} - x_l| = |f'(x_l)|^{-1} |z - f(x_l)|}
        |x_{l + 1} - x_l| = |f'(x_l)|^{-1} \, |z - f(x_l)|.
\end{equation}

        Suppose also that
\begin{equation}
\label{|x_{l + 1} - x_l| le r}
        |x_{l + 1} - x_l| \le r,
\end{equation}
which implies that $|x_{l + 1}| \le r$, so that $f(x_{l + 1})$ is
defined.  In this case,
\begin{equation}
\label{f(x_{l + 1}) - z = f(x_{l + 1}) - f(x_l) - f'(x_l) (x_{l + 1} - x_l)}
        f(x_{l + 1}) - z = f(x_{l + 1}) - f(x_l) - f'(x_l) \, (x_{l + 1} - x_l),
\end{equation}
by (\ref{f(x_l) + f'(x_l) (x_{l + 1} - x_l) = z}), and hence
\begin{eqnarray}
\label{|f(x_{l + 1}) - z| = ...}
        |f(x_{l + 1}) - z| & = &
            |f(x_{l + 1}) - f(x_l) - f'(x_l) \, (x_{l + 1} - x_l)| \\
 & \le & \max_{j \ge 2} (|a_j| \, r^{j - 2}) \, |x_{l + 1} - x_l|^2, \nonumber
\end{eqnarray}
using (\ref{|f(x) - f(y) - f'(y) (x - y)| le ..., 2}) in the second
step.  Plugging (\ref{|x_{l + 1} - x_l| = |f'(x_l)|^{-1} |z - f(x_l)|})
into the right side of (\ref{|f(x_{l + 1}) - z| = ...}), we get that
\begin{equation}
\label{|f(x_{l + 1}) - z| le ...}
        |f(x_{l + 1}) - z| \le \max_{j \ge 2} (|a_j| \, r^{j - 2})
                               \, |f'(x_l)|^{-2} \, |f(x_l) - z|^2.
\end{equation}

        Similarly, if (\ref{|x_{l + 1} - x_l| le r}) holds, and thus
$|x_{l + 1}| \le r$, then
\begin{equation}
\label{|f'(x_{l + 1}) - f'(x_l)| le ...}
        |f'(x_{l + 1}) - f'(x_l)| \le \max_{j \ge 2} (|a_j| \, r^{j - 2})
                                                       \, |x_{l + 1} - x_l|,
\end{equation}
by (\ref{|f'(x) - f'(y)| le max_{j ge 2} (|j cdot a_j| r^{j - 2}) |x -
  y| le ...}).  As before, we can combine this with (\ref{|x_{l + 1} -
  x_l| = |f'(x_l)|^{-1} |z - f(x_l)|}), to get that
\begin{equation}
\label{|f'(x_{l + 1}) - f'(x_l)| le ..., 2}
        |f'(x_{l + 1}) - f'(x_l)| \le \max_{j \ge 2} (|a_j| \, r^{j - 2})
                                      \, |f'(x_l)|^{-1} \, |f(x_l) - z|.
\end{equation}
If
\begin{equation}
\label{max_{j ge 2} (|a_j| r^{j - 2}) |f'(x_l)|^{-1} |f(x_l) - z| < |f'(x_l)|}
 \max_{j \ge 2} (|a_j| \, r^{j - 2}) \, |f'(x_l)|^{-1} \, |f(x_l) - z| < |f'(x_l)|,
\end{equation}
then (\ref{|f'(x_{l + 1}) - f'(x_l)| le ..., 2}) implies that
\begin{equation}
\label{|f'(x_{l + 1}) - f'(x_l)| < |f'(x_l)|}
        |f'(x_{l + 1}) - f'(x_l)| < |f'(x_l)|,
\end{equation}
and hence
\begin{equation}
\label{|f'(x_{l + 1})| = |f'(x_l)|}
        |f'(x_{l + 1})| = |f'(x_l)|,
\end{equation}
by the ultrametric version of the triangle inequality.  Of course, we
would like (\ref{|f'(x_{l + 1})| = |f'(x_l)|}) to hold for each $l \ge
0$, so that $|f'(x_l)| = |f'(x_0)|$ for every $l$.

        Put
\begin{equation}
\label{b_l = max_{j ge 2} (|a_j| r^{j - 2}) |f'(x_l)|^{-2} |f(x_l) - z|}
 b_l = \max_{j \ge 2} (|a_j| \, r^{j - 2}) \, |f'(x_l)|^{-2} \, |f(x_l) - z|,
\end{equation}
so that (\ref{max_{j ge 2} (|a_j| r^{j - 2}) |f'(x_l)|^{-1} |f(x_l) -
  z| < |f'(x_l)|}) is the same as saying that $b_l < 1$.  Using this
notation, (\ref{|f(x_{l + 1}) - z| le ...}) can be reexpressed as
\begin{equation}
\label{|f(x_{l + 1}) - z| le b_l |f(x_l) - z|}
        |f(x_{l + 1}) - z| \le b_l \, |f(x_l) - z|.
\end{equation}
If (\ref{|x_{l + 1} - x_l| le r}) holds and $b_l < 1$, then
\begin{eqnarray}
\label{b_{l + 1} = ...}
 b_{l + 1} & = & \max_{j \ge 2} (|a_j| \, r^{j - 2}) \, |f'(x_{l + 1})|^{-2}
                                                   \, |f(x_{l + 1}) - z| \\
 & = & \max_{j \ge 2} (|a_j| \, r^{j - 2}) \, |f'(x_l)|^{-2} \, |f(x_{l + 1}) - z|,
                                                           \nonumber
\end{eqnarray}
by (\ref{|f'(x_{l + 1})| = |f'(x_l)|}).  Combining this with
(\ref{|f(x_{l + 1}) - z| le b_l |f(x_l) - z|}), we get that
\begin{equation}
\label{b_{l + 1} le ... = b_l^2}
        b_{l + 1} \le b_l \, \max_{j \ge 2} (|a_j| \, r^{j - 2})
                          \, |f'(x_l)|^{-2} \, |f(x_l) - z| = b_l^2,
\end{equation}
using the definition (\ref{b_l = max_{j ge 2} (|a_j| r^{j - 2})
  |f'(x_l)|^{-2} |f(x_l) - z|}) of $b_l$ in the second step.  In particular,
this implies that $b_{l + 1} < 1$ under these conditions.

        If (\ref{|x_{l + 1} - x_l| le r}) holds and $b_l < 1$, then
$f'(x_{l + 1}) \ne 0$, by (\ref{|f'(x_{l + 1})| = |f'(x_l)|}), so that
we can repeat the process.  If $x_{l + 2}$ is obtained from $x_{l + 1}$
as in (\ref{x_{l + 1} = x_l + f'(x_l)^{-1} (z - f(x_l))}), then we get that
\begin{equation}
\label{|x_{l + 2} - x_{l + 1}| = |f'(x_{l + 1})|^{-1} |z - f(x_{l + 1})|}
        |x_{l + 2} - x_{l + 1}| = |f'(x_{l + 1})|^{-1} \, |z - f(x_{l + 1})|,
\end{equation}
as in (\ref{|x_{l + 1} - x_l| = |f'(x_l)|^{-1} |z - f(x_l)|}).  This
implies that
\begin{equation}
\label{|x_{l + 2} - x_{l + 1}| le b_l |f'(x_l)|^{-1} |z - f(x_l)|}
        |x_{l + 2} - x_{l + 1}| \le b_l \, |f'(x_l)|^{-1} \, |z - f(x_l)|,
\end{equation}
by (\ref{|f'(x_{l + 1})| = |f'(x_l)|}) and (\ref{|f(x_{l + 1}) - z| le
  b_l |f(x_l) - z|}).  Equivalently,
\begin{equation}
\label{|x_{l + 2} - x_{l + 1}| le b_l |x_{l + 1} - x_l|}
        |x_{l + 2} - x_{l + 1}| \le b_l \, |x_{l + 1} - x_l|
\end{equation}
under these conditions, by (\ref{|x_{l + 1} - x_l| = |f'(x_l)|^{-1} |z
  - f(x_l)|}).

        In order for all of this to work, we need $z$ to be sufficiently
close to $f(x_0)$, as mentioned at the beginning of the section.
More precisely, let us suppose that
\begin{equation}
\label{|z - f(x_0)| le |f'(x_0)| r}
        |z - f(x_0)| \le |f'(x_0)| \, r
\end{equation}
and
\begin{equation}
\label{max_{j ge 2} (|a_j| r^{j - 2}) |f(x_0) - z| < |f'(x_0)|^2}
        \max_{j \ge 2} (|a_j| \, r^{j - 2}) \, |f(x_0) - z| < |f'(x_0)|^2.
\end{equation}
Note that (\ref{|z - f(x_0)| le |f'(x_0)| r}) is the same as
(\ref{|x_{l + 1} - x_l| le r}) with $l = 0$, by (\ref{|x_{l + 1} -
  x_l| = |f'(x_l)|^{-1} |z - f(x_l)|}).  Similarly, (\ref{max_{j ge 2}
  (|a_j| r^{j - 2}) |f(x_0) - z| < |f'(x_0)|^2}) is the same as
(\ref{max_{j ge 2} (|a_j| r^{j - 2}) |f'(x_l)|^{-1} |f(x_l) - z| <
  |f'(x_l)|}) with $l = 0$, which is the same as saying that $b_0 <
1$.  Under these conditions, we can repeat the process described above
to get a sequence $\{x_l\}_{l = 0}^\infty$ of elements of $k$ defined
recursively by (\ref{x_{l + 1} = x_l + f'(x_l)^{-1} (z - f(x_l))}),
and which satisfies (\ref{|x_{l + 1} - x_l| le r}) and $b_l < 1$ for
every $l$.

        It is easy to see that
\begin{equation}
\label{lim_{l to infty} b_l = 0}
        \lim_{l \to \infty} b_l = 0,
\end{equation}
by (\ref{b_{l + 1} le ... = b_l^2}) and the hypothesis that $b_0 < 1$.
This implies that
\begin{equation}
\label{lim_{l to infty} |x_{l + 1} - x_l| = 0}
        \lim_{l \to \infty} |x_{l + 1} - x_l| = 0,
\end{equation}
because of (\ref{|x_{l + 2} - x_{l + 1}| le b_l |x_{l + 1} - x_l|}).
Thus $\{x_l\}_{l = 0}^\infty$ is a Cauchy sequence of elements of $k$,
by the ultrametric version of the triangle inequality.  It follows
that $\{x_l\}_{l = 0}^\infty$ converges to an element $x$ of $k$,
since $k$ is supposed to be complete.  If $a_j = 0$ for all but
finitely many $j$, then the completeness of $k$ is not needed to
define $f(x)$ as in (\ref{f(x) = sum_{j = 0}^infty a_j x^j}), but it
is still needed here.

        Of course, $|x_{l + 1} - x_l|$ decreases monotonically, by
(\ref{|x_{l + 2} - x_{l + 1}| le b_l |x_{l + 1} - x_l|}) and the
fact that $b_l < 1$ for each $l$.  This implies that
\begin{equation}
\label{|x_l - x_0| le |x_1 - x_0|}
        |x_l - x_0| \le |x_1 - x_0|
\end{equation}
for each $l \ge 0$, by the ultrametric version of the triangle
inequality.  Equivalently,
\begin{equation}
\label{|x_l - x_0| le |f'(x_0)|^{-1} |z - f(x_0)| le r}
        |x_l - x_0| \le |f'(x_0)|^{-1} \, |z - f(x_0)| \le r
\end{equation}
for every $l \ge 0$, using (\ref{|x_{l + 1} - x_l| = |f'(x_l)|^{-1} |z
  - f(x_l)|}) with $l = 0$ in the first step, and (\ref{|z - f(x_0)|
  le |f'(x_0)| r}) in the second step.  It follows that
\begin{equation}
\label{|x - x_0| le |f'(x_0)|^{-1} |z - f(x_0)| le r}
        |x - x_0| \le |f'(x_0)|^{-1} \, |z - f(x_0)| \le r,
\end{equation}
where $x$ is the limit of $\{x_l\}_{l = 1}^\infty$, as in the preceding
paragraph.  In particular,
\begin{equation}
\label{|x| le max(|x - x_0|, |x_0|) le r}
        |x| \le \max(|x - x_0|, |x_0|) \le r,
\end{equation}
since $|x_0| \le r$ by hypothesis, so that $f(x)$ is defined.  It is
easy to see that
\begin{equation}
\label{lim_{l to infty} |f(x_l) - z| = 0}
        \lim_{l \to \infty} |f(x_l) - z| = 0,
\end{equation}
by (\ref{|f(x_{l + 1}) - z| le b_l |f(x_l) - z|}) and (\ref{lim_{l to
    infty} b_l = 0}).  Thus $f(x) = z$, as desired, because $f$ is
continuous on the closed ball in $k$ centered at $0$ with radius $r$,
as in Section \ref{ultrametric absolute value functions}.

        Remember that (\ref{|f'(x_{l + 1})| = |f'(x_l)|}) holds for each
$l \ge 0$ in this situation, which implies that
\begin{equation}
\label{|f'(x_l)| = |f'(x_0)|}
        |f'(x_l)| = |f'(x_0)|
\end{equation}
for each $l \ge 0$.  Taking the limit as $l \to \infty$, we get that
\begin{equation}
\label{|f'(x)| = |f'(x_0)|}
        |f'(x)| = |f'(x_0)|,
\end{equation}
because $f'$ is also continuous on the closed ball in $k$ centered at
$0$ with radius $r$.  Alternatively, (\ref{|f'(x)| = |f'(x_0)|}) could
be derived from (\ref{max_{j ge 2} (|a_j| r^{j - 2}) |f(x_0) - z| <
  |f'(x_0)|^2}) and the first inequality in (\ref{|x - x_0| le
  |f'(x_0)|^{-1} |z - f(x_0)| le r}), as in (\ref{|f'(x)| = |f'(y)|}).
More precisely, let $w$ be any element of $k$ that satisfies
\begin{equation}
\label{|w - x_0| le |f'(x_0)|^{-1} |z - f(x_0)|}
        |w - x_0| \le |f'(x_0)|^{-1} \, |z - f(x_0)|.
\end{equation}
Thus $|w - x_0| \le r$, by (\ref{|z - f(x_0)| le |f'(x_0)| r}), which
implies that
\begin{equation}
\label{|w| le max(|w - x_0|, |x_0|) le r}
        |w| \le \max(|w - x_0|, |x_0|) \le r.
\end{equation}
as before.  Under these conditions, we have that
\begin{equation}
\label{|f'(w)| = |f'(x_0)|}
        |f'(w)| = |f'(x_0)|,
\end{equation}
as in (\ref{|f'(x)| = |f'(y)|}), where $x$ and $y$ in (\ref{|f'(x)| =
  |f'(y)|}) correspond to $w$ and $x_0$ here, respectively.  This uses
(\ref{max_{j ge 2} (|a_j| r^{j - 2}) |f(x_0) - z| < |f'(x_0)|^2}) and
(\ref{|w - x_0| le |f'(x_0)|^{-1} |z - f(x_0)|}) to get the hypothesis
(\ref{max_{j ge 2} (|a_j| r^{j - 2}) |x - y| < |f'(y)|}) for
(\ref{|f'(x)| = |f'(y)|}).

        If $w$ is any element of $k$ that satisfies
(\ref{|w - x_0| le |f'(x_0)|^{-1} |z - f(x_0)|}) and hence
(\ref{|w| le max(|w - x_0|, |x_0|) le r}), then
\begin{equation}
\label{|x - w| le max(|x - x_0|, |x_0 - w|) le |f'(x_0)|^{-1} |z - f(x_0)|}
 |x - w| \le \max(|x - x_0|, |x_0 - w|) \le |f'(x_0)|^{-1} \, |z - f(x_0)|,
\end{equation}
by (\ref{|x - x_0| le |f'(x_0)|^{-1} |z - f(x_0)| le r}).  In this
case, we get that
\begin{equation}
\label{|f(x) - f(w)| = |f'(x)| |x - w| = |f'(x_0)| |x - w|}
        |f(x) - f(w)| = |f'(x)| \, |x - w| = |f'(x_0)| \, |x - w|,
\end{equation}
as in (\ref{|f(x) - f(y)| = |f'(y)| |x - y|}), where $x$ and $y$ in
(\ref{|f(x) - f(y)| = |f'(y)| |x - y|}) correspond to $w$ and $x$
here, respectively.  The hypothesis (\ref{max_{j ge 2} (|a_j| r^{j -
    2}) |x - y| < |f'(y)|}) for (\ref{|f(x) - f(y)| = |f'(y)| |x -
  y|}) follows from (\ref{max_{j ge 2} (|a_j| r^{j - 2}) |f(x_0) - z|
  < |f'(x_0)|^2}), (\ref{|f'(x)| = |f'(x_0)|}), and (\ref{|x - w| le
  max(|x - x_0|, |x_0 - w|) le |f'(x_0)|^{-1} |z - f(x_0)|}) here, and
(\ref{|x - w| le max(|x - x_0|, |x_0 - w|) le |f'(x_0)|^{-1} |z -
  f(x_0)|}) is also used in the second step in (\ref{|f(x) - f(w)| =
  |f'(x)| |x - w| = |f'(x_0)| |x - w|}).  This implies that $w = x$
when $f(w) = f(x)$, since $f'(x_0) \ne 0$ by hypothesis.  This shows
that $x$ is the unique element of $k$ that satisfies (\ref{|x - x_0|
  le |f'(x_0)|^{-1} |z - f(x_0)| le r}) and $f(x) = z$ under these
conditions.

\section{Some additional remarks}
\label{some additional remarks}

        Let us continue with the notation in the previous section.
Note that
\begin{equation}
\label{|f'(x_0)| le max_{j ge 1} (|a_j| r^{j - 1}) = ...}
        |f'(x_0)| \le \max_{j \ge 1} (|a_j| \, r^{j - 1})
                   = r \, \max_{j \ge 1} (|a_j| r^{j - 2}),
\end{equation}
by (\ref{|f'(x)| le max_{j ge 1} (|j cdot a_j| r^{j - 1}) le ...})  in
Section \ref{differentiation, lipschitz conditions} applied to $x_0$.
Suppose that
\begin{equation}
\label{max_{j ge 1} (|a_j| r^{j - 2}) |f(x_0) - z| < |f'(x_0)|^2}
        \max_{j \ge 1} (|a_j| \, r^{j - 2}) \, |f(x_0) - z| < |f'(x_0)|^2.
\end{equation}
This automatically implies (\ref{max_{j ge 2} (|a_j| r^{j - 2})
  |f(x_0) - z| < |f'(x_0)|^2}), since the maximum is taken over all $j
\ge 1$ instead of $j \ge 2$.  Combining (\ref{|f'(x_0)| le max_{j ge
    1} (|a_j| r^{j - 1}) = ...}) and (\ref{max_{j ge 1} (|a_j| r^{j -
    2}) |f(x_0) - z| < |f'(x_0)|^2}), we also get that
\begin{equation}
\label{|f'(x_0)| |f(x_0) - z| le ... < r |f'(x_0)|^2}
        |f'(x_0)| \, |f(x_0) - z|
             \le r \, \max_{j \ge 1} (|a_j| \, r^{j - 2}) \, |f(x_0) - z|
                < r \, |f'(x_0)|^2,
\end{equation}
which implies that
\begin{equation}
\label{|z - f(x_0)| < |f'(x_0)| r}
        |z - f(x_0)| < |f'(x_0)| \, r.
\end{equation}
This is a bit stronger than (\ref{|z - f(x_0)| le |f'(x_0)| r}), so
that (\ref{max_{j ge 1} (|a_j| r^{j - 2}) |f(x_0) - z| < |f'(x_0)|^2})
implies both (\ref{|z - f(x_0)| le |f'(x_0)| r}) and (\ref{max_{j ge
    2} (|a_j| r^{j - 2}) |f(x_0) - z| < |f'(x_0)|^2}).  Thus it
suffices to ask that (\ref{max_{j ge 1} (|a_j| r^{j - 2}) |f(x_0) - z|
  < |f'(x_0)|^2}) hold, in order to get a unique point $x \in k$ that
satisfies (\ref{|x - x_0| le |f'(x_0)|^{-1} |z - f(x_0)| le r}) and
$f(x) = z$, as in the previous section.

        Suppose for the moment that
\begin{equation}
\label{r max_{j ge 2} (|a_j| r^{j - 2}) < |f'(x_0)|}
        r \, \max_{j \ge 2} (|a_j| \, r^{j - 2}) < |f'(x_0)|.
\end{equation}
If (\ref{|z - f(x_0)| le |f'(x_0)| r}) also holds, then we get that
\begin{equation}
\label{r max_{j ge 2} (|a_j| r^{j - 2}) |f(x_0) - z| < |f'(x_0)|^2 r}
 r \, \max_{j \ge 2} (|a_j| \, r^{j - 2}) \, |f(x_0) - z| < |f'(x_0)|^2 \, r,
\end{equation}
by multiplying the left and right sides of (\ref{|z - f(x_0)| le
  |f'(x_0)| r}) and (\ref{r max_{j ge 2} (|a_j| r^{j - 2}) <
  |f'(x_0)|}) together.  Of course, (\ref{r max_{j ge 2} (|a_j| r^{j -
    2}) |f(x_0) - z| < |f'(x_0)|^2 r}) is equivalent to (\ref{max_{j
    ge 2} (|a_j| r^{j - 2}) |f(x_0) - z| < |f'(x_0)|^2}), by dividing
by $r$, so that (\ref{|z - f(x_0)| le |f'(x_0)| r}) automatically
implies (\ref{max_{j ge 2} (|a_j| r^{j - 2}) |f(x_0) - z| <
  |f'(x_0)|^2}) when (\ref{r max_{j ge 2} (|a_j| r^{j - 2}) <
  |f'(x_0)|}) holds.  Note that
\begin{eqnarray}
\label{|f'(x_0) - a_1| = ... = r max_{j ge 2} (|a_j| r^{j - 2})}
 |f'(x_0) - a_1| & = & \biggl|\sum_{j = 2}^\infty j \cdot a_j \, x_0^j\biggr|
                 \le \max_{j \ge 2} (|j \cdot a_j| \, r^{j - 1}) \\
                & \le & \max_{j \ge 2} (|a_j| \, r^{j - 1}) 
                   = r \, \max_{j \ge 2} (|a_j| \, r^{j - 2}), \nonumber
\end{eqnarray}
using the definition of $f'(x_0)$ in the first step, and (\ref{|sum_{j
    = 1}^infty a_j| le max_{j ge 1} |a_j|}) in Section \ref{infinite
  series} in the second step.  This implies that
\begin{equation}
\label{|f'(x_0)| = |a_1|}
        |f'(x_0)| = |a_1|
\end{equation}
when (\ref{r max_{j ge 2} (|a_j| r^{j - 2}) < |f'(x_0)|}) holds,
by the ultrametric version of the triangle inequality.

        Similarly, if
\begin{equation}
\label{r max_{j ge 2} (|a_j| r^{j - 2}) < |a_1| = |f'(0)|}
        r \, \max_{j \ge 2} (|a_j| \, r^{j - 2}) < |a_1| = |f'(0)|,
\end{equation}
then (\ref{|f'(x_0)| = |a_1|}) holds again, by (\ref{|f'(x_0) - a_1| =
  ... = r max_{j ge 2} (|a_j| r^{j - 2})}) and the ultrametric version
of the triangle inequality.  In this case, we also have that
\begin{eqnarray}
\label{|f(x_0) - f(0)| = ... le |a_1| r = |f'(0)| r}
  \quad |f(x_0) - f(0)| = \biggl|\sum_{j = 1}^\infty a_j \, x_0^j\biggr|
                     & \le & \max_{j \ge 1} (|a_j| \, |x_0|^j)  \\
 & \le & \max_{j \ge 1} (|a_j| \, r^j) \le |a_1| \, r = |f'(0)| \, r, \nonumber
\end{eqnarray}
since $|x_0| \le r$.  Consider the inequality
\begin{equation}
\label{|z - f(0)| le |a_1| r = |f'(0)| r}
        |z - f(0)| \le |a_1| \, r = |f'(0)| \, r,
\end{equation}
which is the analogue of (\ref{|z - f(x_0)| le |f'(x_0)| r}) with
$x_0$ replaced by $0$.  If (\ref{r max_{j ge 2} (|a_j| r^{j - 2}) <
  |a_1| = |f'(0)|}) holds, then (\ref{|z - f(x_0)| le |f'(x_0)| r}) is
equivalent to (\ref{|z - f(0)| le |a_1| r = |f'(0)| r}), by
(\ref{|f'(x_0)| = |a_1|}), (\ref{|f(x_0) - f(0)| = ... le |a_1| r =
  |f'(0)| r}), and the ultrametric version of the triangle inequality.
If (\ref{r max_{j ge 2} (|a_j| r^{j - 2}) < |a_1| = |f'(0)|}) and
(\ref{|z - f(0)| le |a_1| r = |f'(0)| r}) both hold, then the analogue
of (\ref{max_{j ge 2} (|a_j| r^{j - 2}) |f(x_0) - z| < |f'(x_0)|^2})
with $x_0$ replaced by $0$ holds too, as in the preceding paragraph.

        Now suppose that        
\begin{equation}
\label{|f'(x_0)| le r max_{j ge 2} (|a_j| r^{j - 2})}
        |f'(x_0)| \le r \, \max_{j \ge 2} (|a_j| \, r^{j - 2}),
\end{equation}
which is the opposite of (\ref{r max_{j ge 2} (|a_j| r^{j - 2}) <
  |f'(x_0)|}).  If (\ref{max_{j ge 2} (|a_j| r^{j - 2}) |f(x_0) - z| <
  |f'(x_0)|^2}) also holds, then we get that
\begin{equation}
\label{|f'(x_0)| |f(x_0) - z| le ... < r |f'(x_0)|^2, 2}
        |f'(x_0)| \, |f(x_0) - z|
             \le r \, \max_{j \ge 2} (|a_j| \, r^{j - 2}) \, |f(x_0) - z|
               < r \, |f'(x_0)|^2,
\end{equation}
as in (\ref{|f'(x_0)| |f(x_0) - z| le ... < r |f'(x_0)|^2}).  This
implies (\ref{|z - f(x_0)| < |f'(x_0)| r}), as before, so that
(\ref{max_{j ge 2} (|a_j| r^{j - 2}) |f(x_0) - z| < |f'(x_0)|^2})
automatically implies (\ref{|z - f(x_0)| le |f'(x_0)| r}) when
(\ref{|f'(x_0)| le r max_{j ge 2} (|a_j| r^{j - 2})}) holds.  Observe
that (\ref{|f'(x_0)| le r max_{j ge 2} (|a_j| r^{j - 2})}) holds when
\begin{equation}
\label{|a_1| le r max_{j ge 2} (|a_j| r^{j - 2})}
        |a_1| \le r \, \max_{j \ge 2} (|a_j| \, r^{j - 2}),
\end{equation}
by (\ref{|f'(x_0)| le max_{j ge 1} (|a_j| r^{j - 1}) = ...}).  Of
course, (\ref{|a_1| le r max_{j ge 2} (|a_j| r^{j - 2})}) is the
opposite of (\ref{r max_{j ge 2} (|a_j| r^{j - 2}) < |a_1| =
  |f'(0)|}).

\section{Another look at regularity}
\label{another look at regularity}

        Let $k$ be a field, and let $|\cdot|$ be an ultrametric absolute
value function on $k$.  If $x, y \in k$ and $j$ is a positive integer,
then
\begin{eqnarray}
\label{(x - y) sum_{l = 0}^{j - 1} x^l y^{j - l - 1} = ... = x^j - y^j}
        (x - y) \, \sum_{l = 0}^{j - 1} x^l \, y^{j - l - 1} 
         & = & x \, \sum_{l = 0}^{j - 1} x^l \, y^{j - l - 1} 
                - y \, \sum_{l = 0}^{j - 1} x^l \, y^{j - l - 1} \\
 & = & \sum_{l = 0}^{j - 1} x^{l + 1} \, y^{j - l - 1} 
                        - \sum_{l = 0}^{j - 1} x^l \, y^{j - l} \nonumber \\
 & = & \sum_{l = 1}^j x^l \, y^{j - l} - \sum_{l = 0}^{j - 1} x^l \, y^{j - l}
           = x^j - y^j.                                         \nonumber
\end{eqnarray}
This implies that
\begin{eqnarray}
\label{|x^j - y^j| = ... le |x - y| (max (|x|, |y|))^{j - 1}}
 \quad |x^j - y^j|
 & = & |x - y| \, \biggl|\sum_{l = 0}^{j - 1} x^l \, y^{j - l - 1}\biggr| \\
 & \le & |x - y| \, \Big(\max_{0 \le l \le j - 1} |x^l \, y^{j - l - 1}|\Big)
                                                                 \nonumber \\
 & \le & |x - y| \, \Big(\max (|x|, |y|)\Big)^{j - 1}, \nonumber
\end{eqnarray}
by the ultrametric version of the triangle inequality.  This gives
another way to look at (\ref{|x^j - y^j| le |x - y| r^{j - 1}}) in
Section \ref{differentiation, lipschitz conditions}.

        Suppose that $k$ is complete with respect to the ultrametric
associated to $|\cdot|$, and let $a_0, a_1, a_2, a_3, \ldots$ be a sequence
of elements of $k$ that satisfies
\begin{equation}
\label{lim_{j to infty} |a_j| r^j = 0, 3}
        \lim_{j \to \infty} |a_j| \, r^j = 0
\end{equation}
for some positive real number $r$.  Thus the corresponding power
series
\begin{equation}
\label{f(x) = sum_{j = 0}^infty a_j x^j, 2}
        f(x) = \sum_{j = 0}^\infty a_j \, x^j
\end{equation}
is defined for each $x \in k$ with $|x| \le r$, as usual.  Note that
the Lipschitz condition (\ref{|f(x) - f(y)| le max_{j ge 1} (|a_j|
  r^{j - 1}) |x - y|}) in Section \ref{differentiation, lipschitz
  conditions} could be derived from (\ref{|f(x) - f(y)| le max_{j ge
    1} (|a_j| |x^j - y^j|)}) using (\ref{|x^j - y^j| = ... le |x - y|
  (max (|x|, |y|))^{j - 1}}) instead of (\ref{|x^j - y^j| le |x - y|
  r^{j - 1}}).

        Let $x, x_0 \in k$ be given, with $|x|, |x_0| \le r$, and observe
that
\begin{equation}
\label{f(x) - f(x_0) = sum_{j = 1}^infty a_j (x^j - x_0^j) = ...}
 \quad  f(x) - f(x_0) = \sum_{j = 1}^\infty a_j \, (x^j - x_0^j)
                      = \sum_{j = 1}^\infty a_j \, (x - x_0) \,
                          \Big(\sum_{l = 0}^{j - 1} x^l \, x_0^{j - l - 1}\Big),
\end{equation}
using (\ref{(x - y) sum_{l = 0}^{j - 1} x^l y^{j - l - 1} = ... = x^j - y^j})
with $y = x_0$ in the second step.  Put
\begin{equation}
\label{b_l = sum_{j = l + 1}^infty a_j x_0^{j - l - 1}}
        b_l = \sum_{j = l + 1}^\infty a_j \, x_0^{j - l - 1}
\end{equation}
for each nonnegative integer $l$, where the convergence of the series
follows from (\ref{lim_{j to infty} |a_j| r^j = 0, 3}).  We also have that
\begin{equation}
\label{|b_l| le ... le max_{j ge l + 1} (|a_j| r^{j - l - 1})}
        |b_l| \le \max_{j \ge l + 1} |a_j \, x_0^{j - l - 1}|
               \le \max_{j \ge l + 1} (|a_j| \, r^{j - l - 1})
\end{equation}
for each $l \ge 0$, using (\ref{|sum_{j = 1}^infty a_j| le max_{j ge
    1} |a_j|}) in Section \ref{infinite series} in the first step.
Thus
\begin{equation}
\label{|b_l| r^l le max_{j ge l + 1} (|a_j| r^{j - 1})}
        |b_l| \, r^l \le \max_{j \ge l + 1} (|a_j| \, r^{j - 1})
\end{equation}
for each $l \ge 0$, which tends to $0$ as $l \to \infty$, by
(\ref{lim_{j to infty} |a_j| r^j = 0, 3}).  This implies that the
power series
\begin{equation}
\label{g_0(x) = sum_{l = 0}^infty b_l x^l}
        g_0(x) = \sum_{l = 0}^\infty b_l \, x^l
\end{equation}
also converges for every $x \in k$ with $|x| \le r$.  Of course,
\begin{equation}
\label{|a_j| |x|^l |x_0|^{j - l - 1} le |a_j| r^{j - 1} to 0 as j to infty}
        |a_j| \, |x|^l \, |x_0|^{j - l - 1} \le |a_j| \, r^{j - 1} \to 0
                                           \quad\hbox{as } j \to \infty,
\end{equation}
by (\ref{lim_{j to infty} |a_j| r^j = 0, 3}), which permits us to
interchange the order of summation in (\ref{f(x) - f(x_0) = sum_{j =
    1}^infty a_j (x^j - x_0^j) = ...}).  It follows that
\begin{equation}
\label{f(x) - f(x_0) = (x - x_0) g_0(x)}
        f(x) - f(x_0) = (x - x_0) \, g_0(x)
\end{equation}
for every $x \in k$ with $|x| \le r$.

        By construction,
\begin{eqnarray}
\label{g_0(x_0) = ... = sum_{j = 1}^infty j a_j x_0^{j - 1} = f'(x_0)}
        g_0(x_0) = \sum_{l = 0}^\infty \Big(\sum_{j = l + 1}^\infty
                                     a_j \, x_0^{j - l - 1}\Big) \, x_0^l
 & = & \sum_{j = 1}^\infty \Big(\sum_{l = 0}^{j - 1} a_j \, x_0^{j - 1}\Big) \\
 & = & \sum_{j = 1}^\infty j \, a_j \, x_0^{j - 1} = f'(x_0).        \nonumber
\end{eqnarray}
Here the order of summation has been interchaged in the second step,
using (\ref{|a_j| |x|^l |x_0|^{j - l - 1} le |a_j| r^{j - 1} to 0 as j
  to infty}) with $x = x_0$.  We also have that
\begin{equation}
\label{|g_0(x) - g_0(y)| le max_{l ge 1} (|b_l| r^{l - 1}) |x - y|}
        |g_0(x) - g_0(y)| \le \max_{l \ge 1} (|b_l| \, r^{l - 1}) \, |x - y|
\end{equation}
for every $x, y \in k$ with $|x|, |y| \le r$, by (\ref{|f(x) - f(y)|
  le max_{j ge 1} (|a_j| r^{j - 1}) |x - y|}) in Section
\ref{differentiation, lipschitz conditions} applied to $g_0$.
Combining this with (\ref{|b_l| le ... le max_{j ge l + 1} (|a_j| r^{j
    - l - 1})}), we get that
\begin{equation}
\label{|g_0(x) - g_0(y)| le max_{j ge 2} (|a_j| r^{j - 2}) |x - y|}
        |g_0(x) - g_0(y)| \le \max_{j \ge 2} (|a_j| \, r^{j - 2}) \, |x - y|
\end{equation}
for every $x, y \in k$ with $|x|, |y| \le r$.  Note that
\begin{equation}
\label{f(x) - f(x_0) - f'(x_0) (x - x_0) = (x - x_0) (g_0(x) - g_0(x_0))}
        f(x) - f(x_0) - f'(x_0) \, (x - x_0) = (x - x_0) \, (g_0(x) - g_0(x_0))
\end{equation}
for every $x \in k$ with $|x| \le r$, by (\ref{f(x) - f(x_0) = (x -
  x_0) g_0(x)}) and (\ref{g_0(x_0) = ... = sum_{j = 1}^infty j a_j
  x_0^{j - 1} = f'(x_0)}).  Thus
\begin{eqnarray}
 |f(x) - f(x_0) - f'(x_0) \, (x - x_0)| & = & |x - x_0| \, |g_0(x) - g_0(x_0)| \\
         & \le & \max_{j \ge 2} (|a_j| \, r^{j - 2}) \, |x - x_0|^2 \nonumber
\end{eqnarray}
for every $x \in k$ with $|x| \le r$, using (\ref{|g_0(x) - g_0(y)| le
  max_{j ge 2} (|a_j| r^{j - 2}) |x - y|}) with $y = x_0$ in the
second step.  This gives another way to look at (\ref{|f(x) - f(y) -
  f'(y) (x - y)| le ..., 2}) in Section \ref{differentiation,
  lipschitz conditions}.

\section{Changing centers}
\label{changing centers}

        Let $k$ be a field, and let
\begin{equation}
\label{f(X) = sum_{j = 0}^infty a_j X^j, 2}
        f(X) = \sum_{j = 0}^\infty a_j \, X^j
\end{equation}
be a formal power series with coefficients in $k$.  If $X$ and $Y$ are
commuting indeterminates, then we have that
\begin{equation}
\label{f(X + Y) = sum_{l = 0}^infty a_l (X + Y)^l = ...}
        f(X + Y) = \sum_{l = 0}^\infty a_l \, (X + Y)^l
  = \sum_{l = 0}^\infty a_l \, \Big(\sum_{j = 0}^l {l \choose j} \cdot X^j
                                                             \, Y^{l - j}\Big)
\end{equation}
as a formal power series in $X$ and $Y$.  Equivalently,
\begin{equation}
\label{f(X + Y) = ...}
        f(X + Y) = \sum_{j = 0}^\infty \Big(\sum_{l = j}^\infty {l \choose j}
                                      \cdot a_l \, Y^{l - j}\Big) \, X^j,
\end{equation}
by interchanging the order of summation in (\ref{f(X + Y) = sum_{l =
    0}^infty a_l (X + Y)^l = ...}).

        Suppose that $k = {\bf C}$, with the standard absolute value
function, and that
\begin{equation}
\label{sum_{j = 0}^infty |a_j| r^j, 2}
        \sum_{j = 0}^\infty |a_j| \, r^j
\end{equation}
converges for some $r > 0$.  This implies that
\begin{equation}
\label{sum_{j = 0}^infty a_j z^j, 2}
        \sum_{j = 0}^\infty a_j \, z^j
\end{equation}
converges absolutely for every $z \in {\bf C}$ with $|z| \le r$, and
we denote the value of the sum by $f(z)$.  Let $z_0 \in {\bf C}$ be
given, with $|z_0| < r$, and put
\begin{equation}
\label{r_0 = r - |z_0| > 0}
        r_0 = r - |z_0| > 0.
\end{equation}
Also let $w \in {\bf C}$ be given, with $|w| \le r_0$, so that
\begin{equation}
\label{|w + z_0| le |w| + |z_0| le r_0 + |z_0| = r}
        |w + z_0| \le |w| + |z_0| \le r_0 + |z_0| = r.
\end{equation}
Thus $f(w + z_0)$ is defined, and can be expressed as
\begin{equation}
\label{f(w + z_0) = sum_{l = 0}^infty a_l (w + z_0)^l = ...}
        f(w + z_0) = \sum_{l = 0}^\infty a_l \, (w + z_0)^l
 = \sum_{l = 0}^\infty a_l \, \Big(\sum_{j = 0}^l {l \choose j} \, w^j \,
                                                           z_0^{l - j}\Big).
\end{equation}
Note that
\begin{equation}
\label{sum_l |a_l| (sum_{j = 0}^l {l choose j} |w|^j |z_0|^{l -j}) = ...}
 \quad \sum_{l = 0}^\infty |a_l| \, \Big(\sum_{j = 0}^l {l \choose j} \, |w|^j
                                                        \, |z_0|^{l -j}\Big)
          = \sum_{l = 0}^\infty |a_l| \, (|z_0| + |w|)^l
           \le \sum_{l = 0}^\infty |a_l| \, r^l,
\end{equation}
which is finite, by hypothesis.  This permits us to interchange the
order of summation in (\ref{f(w + z_0) = sum_{l = 0}^infty a_l (w +
  z_0)^l = ...}), to get that
\begin{equation}
\label{f(w + z_0) = ...}
 f(w + z_0) = \sum_{j = 0}^\infty \Big(\sum_{l = j}^\infty
                          {l \choose j} \, a_l \, z_0^{l - j}\Big) \, w^j.
\end{equation}
More precisely,
\begin{equation}
\label{sum_j (sum_{l = j}^infty {l choose j} |a_l| |z_0|^{l - j}) r_0^j = ...}
 \sum_{j = 0}^\infty \Big(\sum_{l = j}^\infty {l \choose j} \, |a_l| \,
                                                |z_0|^{l - j}\Big) \, |w|^j
   = \sum_{l = 0}^\infty \Big(\sum_{j = 0}^l {l \choose j} \, |a_l| \,
                                                |z_0|^{l - j} \, |w|^j\Big)
\end{equation}
is the same as (\ref{sum_l |a_l| (sum_{j = 0}^l {l choose j} |w|^j
  |z_0|^{l -j}) = ...}), and hence is finite, by hypothesis.  This
implies that the sum in $l$ on the right side of (\ref{f(w + z_0) =
  ...})  converges absolutely for each $j \ge 0$, and that the
resulting sum in $j$ converges absolutely too.

        Let $k$ be an arbitrary field again, and let $|\cdot|$ be an
ultrametric absolute value function on $k$.  Suppose that $k$ is
complete with respect to the metric corresponding to $|\cdot|$, and
that the coefficients $a_j$ of (\ref{f(X) = sum_{j = 0}^infty a_j X^j,
  2}) satisfy
\begin{equation}
\label{lim_{j to infty} |a_j| r^j = 0, 4}
        \lim_{j \to \infty} |a_j| \, r^j = 0
\end{equation}
for some $r > 0$.  This implies that
\begin{equation}
\label{f(x) = sum_{j = 0}^infty a_j x^j, 3}
        f(x) = \sum_{j = 0}^\infty a_j \, x^j
\end{equation}
is defined for $x \in k$ with $|x| \le r$, as before.  If $w, x_0 \in k$
satisfy $|w|, |x_0| \le r$, then
\begin{equation}
\label{|w + x_0| le max(|w|, |x_0|) le r}
        |w + x_0| \le \max(|w|, |x_0|) \le r,
\end{equation}
by the ultrametric version of the triangle inequality.  Thus $f(w +
x_0)$ is defined, and can be expressed as
\begin{equation}
\label{f(w + x_0) = sum_{l = 0}^infty a_l (w + x_0)^l = ...}
        f(w + x_0) = \sum_{l = 0}^\infty a_l \, (w + x_0)^l
 = \sum_{l = 0}^\infty a_l \, \Big(\sum_{j = 0}^l {l \choose j} \cdot w^j \,
                                                         x_0^{l - j}\Big).
\end{equation}
Remember that $|N \cdot 1| \le 1$ for every positive integer $N$,
because $|\cdot|$ is an ultrametric absolute value function on $k$.
This implies that
\begin{equation}
\label{|a_l| |{l choose j} cdot w^j x_0^{l - j}| le ... le |a_l| r^l}
        |a_l| \, \biggl|{l \choose j} \cdot w^j \, x_0^{l - j}\biggr|
                      \le |a_l| \, |w|^j \, |x_0|^{l - j} \le |a_l| \, r^l
\end{equation}
for all $l \ge j \ge 0$, which tends to $0$ as $l \to \infty$, by
(\ref{lim_{j to infty} |a_j| r^j = 0, 4}).  This permits us to
interchange the order of summation in (\ref{f(w + x_0) = sum_{l =
    0}^infty a_l (w + x_0)^l = ...}), as follows.

        Put
\begin{equation}
\label{widetilde{a}_j = sum_{l = j}^infty {l choose j} cdot a_l x_0^{l - j}}
 \widetilde{a}_j = \sum_{l = j}^\infty {l \choose j} \cdot a_l \, x_0^{l - j}
\end{equation}
for each nonnegative integer $j$, where the convergence of the series
in $k$ follows from (\ref{lim_{j to infty} |a_j| r^j = 0, 4}) and the
hypothesis that $|x_0| \le r$.  We also have that
\begin{equation}
\label{|widetilde{a}_j| le ... le max_{l ge j} |a_l| r^{l - j}}
 |\widetilde{a}_j|
       \le \max_{l \ge j} \biggl|{l \choose j} \cdot a_l \, x_0^{l - j}\biggr|
        \le \max_{l \ge j} (|a_l| \, |x_0|^{l - j})
         \le \max_{l \ge j} (|a_l| \, r^{l - j})
\end{equation}
for each $j \ge 0$, using (\ref{|sum_{j = 1}^infty a_j| le max_{j ge
    1} |a_j|}) in Section \ref{infinite series} in the first step.
Thus
\begin{equation}
\label{|widetilde{a}_j| r^j le max_{l ge j} (|a_l| r^l)}
        |\widetilde{a}_j| \, r^j \le \max_{l \ge j} (|a_l| \, r^l)
\end{equation}
for each $j \ge 0$, which tends to $0$ as $j \to \infty$, by
(\ref{lim_{j to infty} |a_j| r^j = 0, 4}).  If $w \in k$ satisfies
$|w| \le r$, as before, then we get that
\begin{equation}
\label{f(w + x_0) = sum_{j = 0}^infty widetilde{a}_j w^j}
        f(w + x_0) = \sum_{j = 0}^\infty \widetilde{a}_j \, w^j,
\end{equation}
by interchanging the order of summation in (\ref{f(w + x_0) = sum_{l =
    0}^infty a_l (w + x_0)^l = ...}).  Note that the series on the
right side of (\ref{f(w + x_0) = sum_{j = 0}^infty widetilde{a}_j
  w^j}) converges when $|w| \le r$, because (\ref{|widetilde{a}_j| r^j
  le max_{l ge j} (|a_l| r^l)}) tends to $0$ as $j \to \infty$.

\section{Contraction mappings}
\label{contraction mappings}

        Let $k$ be a field, and let $|\cdot|$ be an ultrametric absolute
value function on $k$.  Also let $x_0 \in k$ and $t > 0$ be given,
and consider the closed ball
\begin{equation}
\label{overline{B}(x_0, t) = {x in k : |x - x_0| le t}}
        \overline{B}(x_0, t) = \{x \in k : |x - x_0| \le t\}
\end{equation}
centered at $x_0$ with radius $t$ in $k$.  Suppose that $g$ is a
$k$-valued function on $\overline{B}(x_0, t)$ which is Lipschitz with
constant $c \ge 0$ with respect to the ultrametric on $k$ associated
to $|\cdot|$, so that
\begin{equation}
\label{|g(x) - g(y)| le c |x - y|}
        |g(x) - g(y)| \le c \, |x - y|
\end{equation}
for every $x, y \in \overline{B}(x_0, t)$.  Let $\alpha \in k$ be
given, and put
\begin{equation}
\label{f(x) = alpha x + g(x)}
        f(x) = \alpha \, x + g(x)
\end{equation}
for each $x \in \overline{B}(x_0, t)$.  Thus
\begin{eqnarray}
\label{|f(x) - f(y)| le ... le max(|alpha|, c) |x - y|}
        |f(x) - f(y)| & \le & \max(|\alpha| \, |x - y|, |g(x) - g(y)|) \\
                      & \le & \max(|\alpha|, c) \, |x - y| \nonumber
\end{eqnarray}
for every $x, y \in \overline{B}(x_0, t)$, by the ultrametric version of
the triangle inequality.

        Similarly,
\begin{equation}
\label{|alpha| |x - y| le max(|f(x) - f(y)|, |g(x) - g(y)|)}
        |\alpha| \, |x - y| \le \max(|f(x) - f(y)|, |g(x) - g(y)|)
\end{equation}
for each $x, y \in \overline{B}(x_0, t)$.  If
\begin{equation}
\label{|g(x) - g(y)| < |alpha| |x - y|}
        |g(x) - g(y)| < |\alpha| \, |x - y|
\end{equation}
for each $x, y \in \overline{B}(x_0, t)$ with $x \ne y$, then we get
that
\begin{equation}
\label{|alpha| |x - y| le |f(x) - f(y)|}
        |\alpha| \, |x - y| \le |f(x) - f(y)|
\end{equation}
for every $x, y \in \overline{B}(x_0, t)$, which is trivial when $x =
y$.  It follows that
\begin{equation}
\label{|f(x) - f(y)| = |alpha| |x - y|}
        |f(x) - f(y)| = |\alpha| \, |x - y|
\end{equation}
for every $x, y \in \overline{B}(x_0, t)$ under these conditions,
since (\ref{|f(x) - f(y)| le ... le max(|alpha|, c) |x - y|}) holds
with $c = |\alpha|$.

        Let $z \in k$ be given, with
\begin{equation}
\label{|f(x_0) - z| le |alpha| t}
        |f(x_0) - z| \le |\alpha| \, t.
\end{equation}
Suppose that $\alpha \ne 0$, and put
\begin{equation}
\label{h(x) = alpha^{-1} (z - g(x))}
        h(x) = \alpha^{-1} \, (z - g(x))
\end{equation}
for every $x \in \overline{B}(x_0, t)$.  Thus
\begin{equation}
\label{alpha (x - h(x)) = f(x) - z}
        \alpha \, (x - h(x)) = f(x) - z
\end{equation}
for every $x \in \overline{B}(x_0, t)$, so that (\ref{|f(x_0) - z| le
  |alpha| t}) is the same as saying that
\begin{equation}
\label{|h(x_0) - x_0| le t}
        |h(x_0) - x_0| \le t.
\end{equation}
Note that $h$ is Lipschitz with constant $c/|\alpha|$, because $g$ is
Lipschitz with constant $c$.  In particular, if $c \le |\alpha|$, then
$h$ is Lipschitz with constant $1$, and (\ref{|h(x_0) - x_0| le t})
implies that
\begin{equation}
\label{h(overline{B}(x_0, t)) subseteq overline{B}(x_0, t)}
        h(\overline{B}(x_0, t)) \subseteq \overline{B}(x_0, t),
\end{equation}
by the ultrametric version of the triangle inequality.

        If $k$ is complete with respect to the ultrametric associated to
$|\cdot|$, then every closed set $E \subseteq k$ is also complete with
respect to the restriction of this ultrametric to $E$, including $E =
\overline{B}(x_0, r)$.  Suppose that
\begin{equation}
\label{c < |alpha|}
        c < |\alpha|,
\end{equation}
so that $h$ is Lipschitz with constant $c/|\alpha| < 1$.  Under these
conditions, the contraction mapping theorem implies that there is a
unique $x \in \overline{B}(x_0, t)$ such that
\begin{equation}
\label{h(x) = x}
        h(x) = x.
\end{equation}
Equivalently, (\ref{h(x) = x}) means that
\begin{equation}
\label{f(x) = z, 2}
        f(x) = z,
\end{equation}
by (\ref{alpha (x - h(x)) = f(x) - z}).

\section{Contraction mappings, 2}
\label{contraction mappings, 2}

        Let $k$ be a field with an ultrametric absolute value function
$|\cdot|$ again, and suppose that $k$ is complete with respect to
the associated ultrametric.  Also let $a_0, a_1, a_2, a_3, \ldots$
be a sequence of elements of $k$ such that
\begin{equation}
\label{lim_{j to infty} |a_j| r^j = 0, 5}
        \lim_{j \to \infty} |a_j| \, r^j = 0
\end{equation}
for some $r > 0$, so that the corresponding power series
\begin{equation}
\label{f(x) = sum_{j = 0}^infty a_j x^j, 4}
        f(x) = \sum_{j = 0}^\infty a_j \, x^j
\end{equation}
converges for every $x \in k$ with $|x| \le r$.  If $x_0 \in k$ and $t
> 0$ satisfy $|x_0| \le r$ and
\begin{equation}
\label{t le r}
        t \le r,
\end{equation}
then
\begin{equation}
\label{overline{B}(x_0, t) subseteq ... subseteq overline{B}(0, r)}
         \overline{B}(x_0, t) \subseteq \overline{B}(x_0, r)
                               \subseteq \overline{B}(0, r),
\end{equation}
by the ultrametric version of the triangle inequality.  In particular,
this implies that $f$ is defined on $\overline{B}(x_0, t)$.

        Put
\begin{equation}
\label{g(x) = f(x) - f'(x_0) x}
        g(x) = f(x) - f'(x_0) \, x
\end{equation}
for each $x \in \overline{B}(0, r)$.  Thus
\begin{eqnarray}
\label{g(x) - g(y) = f(x) - f(y) - f'(x_0) (x - y) = ...}
        g(x) - g(y) & = & f(x) - f(y) - f'(x_0) \, (x - y) \nonumber \\
 & = & f(x) - f(y) - f'(y) \, (x - y) + (f'(y) - f'(x_0)) \, (x - y) 
\end{eqnarray}
for every $x, y \in \overline{B}(0, r)$, and hence
\begin{eqnarray}
\label{|g(x) - g(y)| le ...}
\lefteqn{|g(x) - g(y)|} \\
 & \le & \max(|f(x) - f(y) - f'(y) \, (x - y)|, |f'(y) - f'(x_0)| \, |x - y|).
                                                   \nonumber
\end{eqnarray}
The right side of (\ref{|g(x) - g(y)| le ...}) can be estimated using
(\ref{|f(x) - f(y) - f'(y) (x - y)| le ..., 2}) and (\ref{|f'(x) -
  f'(y)| le max_{j ge 2} (|j cdot a_j| r^{j - 2}) |x - y| le ...}) in
Section \ref{differentiation, lipschitz conditions}, to get that
\begin{equation}
\label{|g(x) - g(y)| le ..., 2}
        |g(x) - g(y)| \le \Big(\max_{j \ge 2} (|a_j| \, r^{j - 2})\Big) \,
                                \max(|x - y|, |y - x_0|) \, |x - y|
\end{equation}
for every $x, y \in \overline{B}(0, r)$.  If $x, y \in
\overline{B}(x_0, t)$, then $|x - y| \le t$, by the ultrametric
version of the triangle inequality, and so
\begin{equation}
\label{|g(x) - g(y)| le (max_{j ge 2} (|a_j| r^{j - 2})) t |x - y|}
 |g(x) - g(y)| \le \Big(\max_{j \ge 2} (|a_j| \, r^{j - 2})\Big) \, t \, |x - y|.
\end{equation}

        Suppose that $f'(x_0) \ne 0$, and let us take
\begin{equation}
\label{alpha = f'(x_0)}
        \alpha = f'(x_0),
\end{equation}
so that (\ref{f(x) = alpha x + g(x)}) is the same as (\ref{g(x) = f(x)
  - f'(x_0) x}).  The restriction of $g$ to $\overline{B}(x_0, t)$ is
Lipschitz with constant
\begin{equation}
\label{c = t max_{j ge 2} (|a_j| r^{j - 2})}
        c = t \, \max_{j \ge 2} (|a_j| \, r^{j - 2}),
\end{equation}
by (\ref{|g(x) - g(y)| le (max_{j ge 2} (|a_j| r^{j - 2})) t |x - y|}).
Thus (\ref{c < |alpha|}) holds when
\begin{equation}
\label{t max_{j ge 2} (|a_j| r^{j - 2}) < |f'(x_0)|}
        t \, \max_{j \ge 2} (|a_j| \, r^{j - 2}) < |f'(x_0)|.
\end{equation}
The hypotheses (\ref{|z - f(x_0)| le |f'(x_0)| r}) and (\ref{max_{j ge
    2} (|a_j| r^{j - 2}) |f(x_0) - z| < |f'(x_0)|^2}) in Section
\ref{hensel's lemma} basically correspond to (\ref{|f(x_0) - z| le
  |alpha| t}), (\ref{t le r}), and (\ref{t max_{j ge 2} (|a_j| r^{j -
    2}) < |f'(x_0)|}) here.  More precisely, if we take
\begin{equation}
\label{t = |f'(x_0)|^{-1} |f(x_0) - z|}
        t = |f'(x_0)|^{-1} \, |f(x_0) - z|,
\end{equation}
then (\ref{|f(x_0) - z| le |alpha| t}) becomes an equality, by
(\ref{alpha = f'(x_0)}).  With this choice of $t$, (\ref{|z - f(x_0)|
  le |f'(x_0)| r}) is equivalent to (\ref{t le r}), and (\ref{max_{j
    ge 2} (|a_j| r^{j - 2}) |f(x_0) - z| < |f'(x_0)|^2}) is equivalent
to (\ref{t max_{j ge 2} (|a_j| r^{j - 2}) < |f'(x_0)|}).  This shows
that the conclusions of Section \ref{hensel's lemma} can also be
derived from the discussion in the previous section.

        Note that
\begin{equation}
\label{|f'(x)| = |f'(x_0)|, 2}
        |f'(x)| = |f'(x_0)|
\end{equation}
for every $x \in \overline{B}(x_0, t)$ when $f'(x_0) \ne 0$ and $t >
0$ satisfies (\ref{t le r}) and (\ref{t max_{j ge 2} (|a_j| r^{j - 2})
  < |f'(x_0)|}).  This is basically the same as (\ref{|f'(x)| =
  |f'(y)|}) in Section \ref{differentiation, lipschitz conditions}.

\section{Strassmann's theorem}
\label{strassmann's theorem}

        Let $k$ be a field, and let $|\cdot|$ be an ultrametric absolute
value function on $k$.  Suppose that $k$ is complete with respect to the
ultrametric associated to $|\cdot|$, and let $a_0, a_1, a_2, a_3, \ldots$
be a sequence of elements of $k$ that satisfies
\begin{equation}
\label{lim_{j to infty} |a_j| r^j = 0, 6}
        \lim_{j \to \infty} |a_j| \, r^j = 0
\end{equation}
for some positive real number $r$.  Thus the power series
\begin{equation}
\label{f(x) = sum_{j = 0}^infty a_j x^j, 5}
        f(x) = \sum_{j = 0}^\infty a_j \, x^j
\end{equation}
converges for every $x \in k$ with $|x| \le r$, as usual.  Suppose
also that $a_j \ne 0$ for some $j$, and let $N$ be a nonnegative
integer such that
\begin{equation}
\label{|a_j| r^j < |a_N| r^N for every j > N}
        |a_j| \, r^j < |a_N| \, r^N \quad\hbox{for every } j > N.
\end{equation}
More precisely, one can take $N$ to be the largest nonnegative integer
such that
\begin{equation}
\label{|a_N| r^N = max_{j ge 0} (|a_j| r^j)}
        |a_N| \, r^N = \max_{j \ge 0} (|a_j| \, r^j),
\end{equation}
and this will be the smallest $N$ that satisfies (\ref{|a_j| r^j <
  |a_N| r^N for every j > N}).

        Under these conditions, \emph{Strassmann's 
theorem}\index{Strassmann's theorem} implies that $f$ can have at most
$N$ zeros in the closed ball $\overline{B}(0, r)$ in $k$ centered at $0$
and with radius $r$.  We follow the very nice proof by induction on
$N$ in \cite{cas}, which is also discussed in \cite{fg}.  To deal with
the base case $N = 0$, observe that
\begin{equation}
\label{|f(x) - a_0| = ... le max_{j ge 1} |a_j| r^j}
        |f(x) - a_0| = \biggl|\sum_{j = 1}^\infty a_j \, x^j\biggr|
           \le \max_{j \ge 1} (|a_j| \, |x|^j) \le \max_{j \ge 1} (|a_j| \, r^j)
\end{equation}
for every $x \in k$ with $|x| \le r$, using (\ref{|sum_{j = 1}^infty
  a_j| le max_{j ge 1} |a_j|}) in Section \ref{infinite series} in the
second step.  If (\ref{|a_j| r^j < |a_N| r^N for every j > N}) holds
with $N = 0$, then we get that
\begin{equation}
\label{|f(x) - a_0| < |a_0|}
        |f(x) - a_0| < |a_0|
\end{equation}
for every $x \in k$ with $|x| \le r$, which implies that
\begin{equation}
\label{|f(x)| = |a_0|}
        |f(x)| = |a_0|,
\end{equation}
by the ultrametric version of the triangle inequality.  In particular,
this means that $f(x) \ne 0$ for every $x \in k$ with $|x| \le r$,
because $a_0 \ne 0$ in this situation.

        Now suppose that $N \ge 1$, and that Strassmann's theorem holds
for $N - 1$.  If $f(x) \ne 0$ for every $x \in k$ with $|x| \le r$,
then there is nothing to do, and so we suppose also that there is an
$x_0 \in k$ such that $|x_0| \le r$ and $f(x_0) = 0$.  Let
\begin{equation}
\label{g_0(x) = sum_{l = 0}^infty b_l x^l, 2}
        g_0(x) = \sum_{l = 0}^\infty b_l \, x^l
\end{equation}
be as in (\ref{g_0(x) = sum_{l = 0}^infty b_l x^l}) in Section
\ref{another look at regularity}, where
\begin{equation}
\label{b_l = sum_{j = l + 1}^infty a_j x_0^{j - l - 1}, 2}
        b_l = \sum_{j = l + 1}^\infty a_j \, x_0^{j - l - 1}
\end{equation}
for each nonnegative integer $l$, as in (\ref{b_l = sum_{j = l +
    1}^infty a_j x_0^{j - l - 1}}).  Thus
\begin{equation}
\label{|b_l| r^l le max_{j ge l + 1} (|a_j| r^{j - 1}), 2}
        |b_l| \, r^l \le \max_{j \ge l + 1} (|a_j| \, r^{j - 1})
\end{equation}
for each $l \ge 0$, as in (\ref{|b_l| r^l le max_{j ge l + 1} (|a_j|
  r^{j - 1})}).  By construction,
\begin{eqnarray}
\label{|b_l - a_{l + 1}| = ... le max_{j ge l + 2} |a_j| r^{j - l - 1}}
 |b_l - a_{l + 1}| = \biggl|\sum_{j = l + 2}^\infty a_j \, x_0^{j - l - 1}\biggr|
                   & \le & \max_{j \ge l + 2} (|a_j| \, |x_0|^{j - l - 1}) \\
                    & \le & \max_{j \ge l + 2} (|a_j| \, r^{j - l - 1}) \nonumber
\end{eqnarray}
for each $l \ge 0$, using (\ref{|sum_{j = 1}^infty a_j| le max_{j ge
    1} |a_j|}) in Section \ref{infinite series} in the second step.
If $N$ satisfies (\ref{|a_j| r^j < |a_N| r^N for every j > N}), then
it is easy to see that
\begin{equation}
\label{|b_{N - 1} - a_N| < |a_N|}
        |b_{N - 1} - a_N| < |a_N|,
\end{equation}
and hence
\begin{equation}
\label{|b_{N - 1}| = |a_N|}
        |b_{N - 1}| = |a_N|,
\end{equation}
by the ultrametric version of the triangle inequality.  It follows that
\begin{equation}
\label{|b_l| r^l < |a_N| r^{N - 1} = |b_{N - 1}| r^{N - 1}}
        |b_l| \, r^l < |a_N| \, r^{N - 1} = |b_{N - 1}| \, r^{N - 1}
\end{equation}
for every $l > N - 1$, by (\ref{|a_j| r^j < |a_N| r^N for every j >
  N}) and (\ref{|b_l| r^l le max_{j ge l + 1} (|a_j| r^{j - 1}), 2}).
This means that $g_0$ satisfies the analogous condition with $N - 1$
instead of $N$, so that $g_0$ has at most $N - 1$ zeros in
$\overline{B}(0, r)$, by the induction hypothesis.  We also have that
\begin{equation}
\label{f(x) = (x - x_0) g_0(x)}
        f(x) = (x - x_0) \, g_0(x)
\end{equation}
for every $x \in k$ with $|x| \le r$, by (\ref{f(x) - f(x_0) = (x -
  x_0) g_0(x)}) in Section \ref{another look at regularity}, and the
condition that $f(x_0) = 0$.  This implies that $f$ can have at most
$N$ zeros in $\overline{B}(0, r)$, as desired.

        As an application of Strassmann's theorem, let $N$ be a
positive integer that satisfies (\ref{|a_N| r^N = max_{j ge 0} (|a_j| r^j)}).
As before, one can take $N$ to be the largest positive integer such that
\begin{equation}
\label{|a_N| r^N = max_{j ge 1} (|a_j| r^j)}
        |a_N| \, r^N = \max_{j \ge 1} (|a_j| \, r^j).
\end{equation}
If $z$ is any element of $k$, then $f(x) - z$ can be expressed by a
power series in $x$ on $\overline{B}(0, r)$, where the constant term
is equal to $a_0 - z$, and the other coefficients are the same as the
coefficients $a_j$ of $f$, with $j \ge 1$.  This means that $N$ also
satisfies the analogue of (\ref{|a_j| r^j < |a_N| r^N for every j >
  N}) for $f(x) - z$ instead of $f(x)$, because $N \ge 1$ by
hypothesis.  Thus Strassmann's theorem implies that $f(x) - z$ can
have at most $N$ zeros in $\overline{B}(0, r)$, which is to say that
there are at most $N$ points $x \in \overline{B}(0, r)$ at which $f(x)
= z$.

\section{The exponential function}
\label{exponential function}

        Let $k$ be a field of characteristic $0$, and let
\begin{equation}
\label{E(X) = sum_{j = 0}^infty frac{X^j}{j!}}
        E(X) = \sum_{j = 0}^\infty \frac{X^j}{j!}
\end{equation}
be the exponential function\index{exponential function} on $k$, which
is initially considered as a formal power series in the indeterminate
$X$.  Of course, this uses the natural embedding of ${\bf Q}$ into
$k$, that results from $k$ having characteristic $0$.  If $X$ and $Y$
are commuting indeterminates, then it is easy to see that
\begin{equation}
\label{E(X + Y) = E(X) E(Y)}
        E(X + Y) = E(X) \, E(Y)
\end{equation}
as formal power series in $X$ and $Y$, as in (\ref{E(w + z) = E(w)
  E(z)}) in Section \ref{complex coefficients}.  Note that the
formal derivative of $E(X)$ is equal to $E(X)$, as usual.

        Let $|\cdot|$ be an absolute value function on $k$, and let
us consider the convergence properties of
\begin{equation}
\label{E(x) = sum_{j = 0}^infty frac{x^j}{j!}}
        E(x) = \sum_{j = 0}^\infty \frac{x^j}{j!}
\end{equation}
for $x \in k$.  To do this, we also
suppose from now on in this section that $k$ is complete with respect
to the metric that corresponds to $|\cdot|$.  Remember that $|\cdot|$
induces an absolute value function on ${\bf Q}$, using the natural
embedding of ${\bf Q}$ into $k$.  The behavior of this absolute value
function on ${\bf Q}$ is characterized by Ostrowksi's theorem,
as discussed at the end of Section \ref{some more refinements}.
This leads to various cases for $k$, as follows.

        Suppose first that the induced absolute value function on
${\bf Q}$ is archimedian, and hence that $|\cdot|$ is archimedian on $k$.
In this case, Ostrowski's theorem implies that the induced absolute
value function on ${\bf Q}$ is a positive power of the standard
absolute value function.  We may as well suppose that the induced
absolute value function on ${\bf Q}$ is equal to the standard absolute
value function, by replacing $|\cdot|$ on $k$ by a suitable power of
itself, as in Section \ref{another refinement}.  Under these
conditions, $k$ is isomorphic to ${\bf R}$ or ${\bf C}$, and $|\cdot|$
corresponds to the standard absolute value function, as in Sections
\ref{complex numbers} and \ref{complex numbers, 2}.  Similarly,
(\ref{E(x) = sum_{j = 0}^infty frac{x^j}{j!}}) corresponds to the
usual exponential function on ${\bf R}$ or ${\bf C}$, as in Section
\ref{complex coefficients}.

        Otherwise, if the induced absolute value function on ${\bf Q}$
is non-archimedian, then $|\cdot|$ is non-archimedian on $k$ too.
This implies that $|\cdot|$ is an ultrametric absolute value function
on $k$, as in Section \ref{some more refinements}.  In this situation,
the right side of (\ref{E(x) = sum_{j = 0}^infty frac{x^j}{j!}})
converges in $k$ exactly when
\begin{equation}
\label{x^j/j! to 0 as j to infty}
        x^j/j! \to 0 \quad\hbox{as } j \to \infty,
\end{equation}
as in Section \ref{infinite series}.  If $x, y \in k$ both have this
property, then the Cauchy product of the series used to define $E(x)$
and $E(y)$ also converges in $k$, and has sum equal to $E(x) \, E(y)$,
as in Section \ref{infinite series} again.  This means that the series
used to define $E(x + y)$ converges in $k$, and that
\begin{equation}
\label{E(x + y) = E(x) E(y)}
        E(x + y) = E(x) \, E(y),
\end{equation}
since the series defining $E(x + y)$ is the same as the Cauchy product
of the series defining $E(x)$ and $E(y)$, as before.

        In particular, the induced absolute value function on ${\bf Q}$
is non-archimedian when it is trivial, in which case (\ref{x^j/j! to 0
  as j to infty}) holds exactly when $|x| < 1$.  Suppose now that the
induced absolute value function on ${\bf Q}$ is non-archimedian and
not trivial.  By Ostrowski's theorem, there is a prime number $p$ such
that the induced absolute value function on ${\bf Q}$ is equal to a
positive power of the $p$-adic absolute value.  As before, we may as
well suppose that the induced absolute value function on ${\bf Q}$ is
equal to the $p$-adic absolute value, by replacing $|\cdot|$ with a
suitable positive power of itself on $k$, if necessary.  We shall
restrict our attention to this case for the rest of the section.

        Note that
\begin{equation}
\label{|x^j / j!| = |x|^j / |j!|_p}
        |x^j / j!| = |x|^j / |j!|_p
\end{equation}
for every $x \in k$ and nonnegative integer $j$, where $|j!|_p$ is the
$p$-adic absolute value of $j!$.  This follows from the hypothesis
that the absolute value function induced on ${\bf Q}$ by $|\cdot|$ on
$k$ is equal to $|\cdot|_p$.  In order to find out when $x \in k$
satisfies (\ref{x^j/j! to 0 as j to infty}), we would like to estimate
$|j!|_p$, which is determined by the total number of factors of $p$ in
$j!$.  Let $[r]$ denote the integer part of a nonnegative real number
$r$, which is the largest nonnegative integer less than or equal to
$r$.  It is well known that the total number of factors of $p$ in $j!$
is equal to
\begin{equation}
\label{sum_{l = 1}^infty [j/p^l]}
        \sum_{l = 1}^\infty [j/p^l]
\end{equation}
for each nonnegative integer $j$.  More precisely, $[j/p^l]$ is the
number of positive integers less than or equal to $j$ that are
divisible by $p^l$.  Thus the first term in the series (\ref{sum_{l =
    1}^infty [j/p^l]}) counts a factor of $p$ in $j!$ for each
positive integer less than or equal to $j$ that is divisible by $p$.
The second term in (\ref{sum_{l = 1}^infty [j/p^l]}) counts another
factor of $p$ for each positive integer less than or equal to $j$ that
is divisible by $p^2$, and so on.

        Of course, the infinite series (\ref{sum_{l = 1}^infty [j/p^l]})
is really a finite sum, since $[j/p^l] = 0$ when $j < p^l$.
We also have that
\begin{equation}
\label{sum_{l = 1}^infty [j/p^l] < ... = frac{j}{p - 1}}
        \sum_{l = 1}^\infty [j/p^l] < \sum_{l = 1}^\infty j/p^l
                      = j \, \frac{(1/p)}{1 - (1/p)} = \frac{j}{p - 1}
\end{equation}
for each positive integer $j$.  It follows that
\begin{equation}
\label{1/|j!|_p le p^{j/(p - 1)}}
        1/|j!|_p \le p^{j/(p - 1)}
\end{equation}
for every nonnegative integer $j$, so that
\begin{equation}
\label{|x^j / j!| = |x|^j/|j!|_p le (|x| p^{1/(p - 1)})^j}
        |x^j / j!| = |x|^j/|j!|_p \le (|x| \, p^{1/(p - 1)})^j
\end{equation}
for every $x \in k$.  If $x \in k$ satisfies
\begin{equation}
\label{|x| < p^{-1/(p - 1)}}
        |x| < p^{-1/(p - 1)},
\end{equation}
then $|x| \, p^{1/(p - 1)} < 1$, and hence $(|x| \, p^{1/(p - 1)})^ j
\to 0$ as $j \to \infty$.  Combining this with (\ref{|x^j / j!| =
  |x|^j/|j!|_p le (|x| p^{1/(p - 1)})^j}), we get that (\ref{x^j/j!
  to 0 as j to infty}) holds when $x$ satisfies (\ref{|x| < p^{-1/(p -
    1)}}).

        If $j = p^n$ for some positive integer $n$, then we have that
\begin{equation}
\label{sum_{l = 1}^infty [j/p^l] = ... = frac{j - 1}{p - 1}}
 \sum_{l = 1}^\infty [j/p^l] = \sum_{l = 1}^n p^{n - l} = \sum_{l = 0}^{n - 1} p^l
                                = \frac{p^n - 1}{p - 1} = \frac{j - 1}{p - 1}.
\end{equation}
This means that
\begin{equation}
\label{1/|j!|_p = p^{(j - 1)/(p - 1)}}
        1/|j!|_p = p^{(j - 1)/(p - 1)}
\end{equation}
when $j$ is a positive power of $p$, so that
\begin{equation}
\label{|x^j / j!| = |x|^j / |j!|_p = (|x| p^{1/(p - 1)})^j p^{-1/(p - 1)}}
        |x^j / j!| = |x|^j / |j!|_p = (|x| \, p^{1/(p - 1)})^j \, p^{-1/(p - 1)}
\end{equation}
for every $x \in k$ in this case.  If $x \in k$ satisfies
\begin{equation}
\label{|x| ge p^{-1/(p - 1)}}
        |x| \ge p^{-1/(p - 1)},
\end{equation}
then $|x| \, p^{1/(p - 1)} \ge 1$, and the right side of (\ref{|x^j /
  j!| = |x|^j / |j!|_p = (|x| p^{1/(p - 1)})^j p^{-1/(p - 1)}}) is
greater than or equal to $p^{-1/(p - 1)}$.  This implies that $x^j/j!$
does not converge to $0$ in $k$ as $j \to \infty$ when $x$ satisfies
(\ref{|x| ge p^{-1/(p - 1)}}), since (\ref{|x^j / j!| = |x|^j / |j!|_p
  = (|x| p^{1/(p - 1)})^j p^{-1/(p - 1)}}) holds for infinitely many
$j$.  It follows that (\ref{x^j/j! to 0 as j to infty}) holds if and
only if $x \in k$ satisfies (\ref{|x| < p^{-1/(p - 1)}}) under these
conditions.

\section{Some additional properties}
\label{some additional properties}

          Let $k$ be a field with characteristic $0$ again, and let
$|\cdot|$ be an ultrametric absolute value function on $k$.  Suppose
from now on in this section that $k$ is complete with respect to the
ultrametric associated to $|\cdot|$.  Let us also suppose for the
first part of the section that the induced absolute value function on
${\bf Q}$ is the $p$-adic absolute value for some prime number $p$.
Thus the domain of the exponential function (\ref{E(x) = sum_{j =
    0}^infty frac{x^j}{j!}}) is given by
\begin{equation}
\label{D = B(0, p^{-1/(p - 1)}) = {x in k : |x| < p^{-1/(p - 1)}}}
        D = B(0, p^{-1/(p - 1)}) = \{x \in k : |x| < p^{-1/(p - 1)}\},
\end{equation}
as in the previous section.

        Note that
\begin{equation}
\label{|E(x) - 1| le max_{j ge 1} |x^j/j!|}
        |E(x) - 1| \le \max_{j \ge 1} |x^j/j!|
\end{equation}
for every $x \in D$, as in (\ref{|sum_{j = 1}^infty a_j| le max_{j ge
    1} |a_j|}) in Section \ref{infinite series}.  This implies that
\begin{equation}
\label{|E(x) - 1| < 1}
        |E(x) - 1| < 1
\end{equation}
for every $x \in D$, by (\ref{|x^j / j!| = |x|^j/|j!|_p le (|x|
  p^{1/(p - 1)})^j}) and the definition (\ref{D = B(0, p^{-1/(p - 1)})
  = {x in k : |x| < p^{-1/(p - 1)}}}) of $D$.  It follows that
\begin{equation}
\label{|E(x)| = 1}
        |E(x)| = 1
\end{equation}
for every $x \in D$, by the ultrametric version of the triangle
inequality.  

        Let $j$ be a positive integer, and remember that
\begin{equation}
\label{(p - 1) sum_{l = 1}^infty [j/p^l] < j}
        (p - 1) \, \sum_{l = 1}^\infty [j/p^l] < j,
\end{equation}
by (\ref{sum_{l = 1}^infty [j/p^l] < ... = frac{j}{p - 1}}).  The left
side of (\ref{(p - 1) sum_{l = 1}^infty [j/p^l] < j}) is an integer,
and so we get that
\begin{equation}
\label{(p - 1) sum_{l = 1}^infty [j/p^l] le j - 1}
        (p - 1) \, \sum_{l = 1}^\infty [j/p^l] \le j - 1.
\end{equation}
Equivalently,
\begin{equation}
\label{sum_{l = 1}^infty [j/p^l] le frac{j - 1}{p - 1}}
        \sum_{l = 1}^\infty [j/p^l] \le \frac{j - 1}{p - 1},
\end{equation}
which implies that
\begin{equation}
\label{1/|j!|_p le p^{(j - 1)/(p - 1)}}
        1/|j!|_p \le p^{(j - 1)/(p - 1)}
\end{equation}
for every positive integer $j$.  A more precise analysis of
(\ref{sum_{l = 1}^infty [j/p^l]}) and hence $|j!|_p$ is given in
Problem 165 on p113 of \cite{fg}.  Note that the sum in the statement
of Lemma 4.5.4 on p112 of \cite{fg} should start at $i = 1$, as well
as the sum in Problem 164 on p113 of \cite{fg}.

        Using (\ref{1/|j!|_p le p^{(j - 1)/(p - 1)}}), we get that
\begin{equation}
\label{|x^{j - 1} / j!| = |x|^{j - 1} / |j!|_p le (|x| p^{1/(p - 1)})^{j - 1}}
        |x^{j - 1} / j!| = |x|^{j - 1} / |j!|_p \le (|x| \, p^{1/(p - 1)})^{j - 1}
\end{equation}
for every $x \in k$ and positive integer $j$.  Equivalently,
\begin{equation}
\label{|x^j / j!| le (|x| p^{1/(p - 1)})^{j - 1} |x|}
        |x^j / j!| \le (|x| \, p^{1/(p - 1)})^{j - 1} \, |x|
\end{equation}
for every $x \in k$ and positive integer $j$, which is a refinement
of (\ref{|x^j / j!| = |x|^j/|j!|_p le (|x| p^{1/(p - 1)})^j}).  If
$x$ satisfies (\ref{|x| < p^{-1/(p - 1)}}), then we get that
\begin{equation}
\label{|x^j / j!| le |x|}
        |x^j / j!| \le |x|
\end{equation}
for every positive integer $j$.  Thus
\begin{equation}
\label{|E(x) - 1| le |x|}
        |E(x) - 1| \le |x|
\end{equation}
for every $x \in D$, by (\ref{|E(x) - 1| le max_{j ge 1} |x^j/j!|}).
In particular,
\begin{equation}
\label{|E(x) - 1| < p^{-1/(p - 1)}}
        |E(x) - 1| < p^{-1/(p - 1)}
\end{equation}
for every $x \in D$, by (\ref{|E(x) - 1| le |x|}) and the definition
(\ref{D = B(0, p^{-1/(p - 1)}) = {x in k : |x| < p^{-1/(p - 1)}}}) of
$D$.  

        More precisely, (\ref{|x^j / j!| le (|x| p^{1/(p - 1)})^{j - 1} |x|})
implies that
\begin{equation}
\label{|x^j / j!| < |x|}
        |x^j / j!| < |x|
\end{equation}
for every $x \in D$ with $x \ne 0$, and every integer $j \ge 2$.  We
also have that
\begin{equation}
\label{|E(x) - 1 - x| le max_{j ge 2} |x^j / j!|}
        |E(x) - 1 - x| \le \max_{j \ge 2} |x^j / j!|
\end{equation}
for every $x \in D$, as in (\ref{|sum_{j = 1}^infty a_j| le max_{j ge
    1} |a_j|}) in Section \ref{infinite series}.  Combining (\ref{|x^j
  / j!| < |x|}) and (\ref{|E(x) - 1 - x| le max_{j ge 2} |x^j / j!|}),
we get that
\begin{equation}
\label{|E(x) - 1 - x| < |x|}
        |E(x) - 1 - x| < |x|
\end{equation}
for every $x \in D$ with $x \ne 0$.  This implies that
\begin{equation}
\label{|E(x) - 1| = |x|}
        |E(x) - 1| = |x|
\end{equation}
for every $x \in D$ with $x \ne 0$, by the ultrametric version of the
triangle inequality.  Of course, (\ref{|E(x) - 1| = |x|}) also holds
when $x = 0$, so that (\ref{|E(x) - 1| = |x|}) holds for every $x \in
D$.

        Let $z \in k$ be given, with
\begin{equation}
\label{|z - 1| < p^{-1/(p - 1)}}
        |z - 1| < p^{-1/(p - 1)}.
\end{equation}
We would like to find $x \in D$ such that $E(x) = z$.  If $z = 1$,
then we can take $x = 0$, and so we suppose from now on that $z \ne
1$.  Put $r = |z - 1| > 0$ and $x_0 = 0$, and let us apply the
discussion in Section \ref{hensel's lemma} to $f = E$.  Thus $f(x_0) =
E(0) = 1$ and $f'(x_0) = E'(0) = 1$ in this case.  It is easy to see
that (\ref{|z - f(x_0)| le |f'(x_0)| r}) holds under these conditions,
by definition of $r$.  One can also check that (\ref{max_{j ge 2}
  (|a_j| r^{j - 2}) |f(x_0) - z| < |f'(x_0)|^2}) holds in this
situation, because of (\ref{1/|j!|_p le p^{(j - 1)/(p - 1)}}).  The
discussion in Section \ref{hensel's lemma} leads to a point $x \in k$
such that $|x| \le r$ and $E(x) = f(x) = z$, as desired.

        Let us now consider the case where the induced absolute value
function on ${\bf Q}$ is trivial.  Thus the domain of the exponential
function (\ref{E(x) = sum_{j = 0}^infty frac{x^j}{j!}}) is
\begin{equation}
\label{D = B(0, 1) = {x in k : |x| < 1}}
        D = B(0, 1) = \{x \in k : |x| < 1\},
\end{equation}
as in the previous section.  It is easy to see that (\ref{|E(x) - 1| <
  1}) still holds for every $x \in D$ in this situation, using
(\ref{|E(x) - 1| le max_{j ge 1} |x^j/j!|}) and the fact that
\begin{equation}
\label{|x^j / j!| = |x|^j}
        |x^j / j!| = |x|^j
\end{equation}
for each $x \in k$ and positive integer $j$.  This implies that
(\ref{|E(x)| = 1}) also holds for every $x \in D$, as before.
Similarly, (\ref{|E(x) - 1| le |x|}) holds for every $x \in D$, by
(\ref{|E(x) - 1| le max_{j ge 1} |x^j/j!|}), (\ref{|x^j / j!| =
  |x|^j}), and the definition (\ref{D = B(0, 1) = {x in k : |x| < 1}})
of $D$.  We also have (\ref{|E(x) - 1| = |x|}) for every $x \in D$,
for essentially the same reasons as before.  If $z \in k$ satisfies
\begin{equation}
\label{|z - 1| < 1}
        |z - 1| < 1,
\end{equation}
then one can find an $x \in D$ such that $E(x) = z$, using the
discussion in Section \ref{hensel's lemma} again.

\chapter{Geometry and measure}
\label{geometry, measure}

\section{Diameters of sets}
\label{diameters of sets}

        Let $(M, d(x, y))$ be a metric space.  As usual, a subset $A$
of $M$ is said to be \emph{bounded}\index{bounded sets} if $A$ is
contained in a ball of finite radius in $M$.  If $A \subseteq M$ is
bounded and nonempty, then the \emph{diameter}\index{diameters of
  sets} of $A$ is defined by
\begin{equation}
\label{diam A = sup{d(x, y) : x, y in A}}
        \diam A = \sup\{d(x, y) : x, y \in A\}.
\end{equation}
It will be convenient to put $\diam A = +\infty$ when $A$ is
unbounded, and $\diam A = 0$ when $A = \emptyset$.  Note that
\begin{equation}
\label{diam overline{A} = diam A}
        \diam \overline{A} = \diam A
\end{equation}
for every $A \subseteq M$, where $\overline{A}$ denotes the closure of
$A$ in $M$.

        Let $A \subseteq M$ and a positive real number $r$ be given, and put
\begin{equation}
\label{A_r = bigcup_{x in A} B(x, r) = {y in M : d(x, y) < r for some x in A}}
        A_r = \bigcup_{x \in A} B(x, r)
            = \{y \in M : d(x, y) < r \hbox{ for some } x \in A\}.
\end{equation}
Thus $A \subseteq A_r$, and $A_r$ is an open set in $M$, because $A_r$
is a union of open sets.  It is easy to see that
\begin{equation}
\label{diam A_r le diam A + 2 r}
        \diam A_r \le \diam A + 2 \, r,
\end{equation}
by the triangle inequality.  If $d(\cdot, \cdot)$ is an ultrametric on
$M$, then we have that
\begin{equation}
\label{diam A_r le max(diam A, r)}
        \diam A_r \le \max(\diam A, r).
\end{equation}

        If $M$ is any metric space, $A \subseteq M$ is bounded,
and $p \in A$, then
\begin{equation}
\label{A subseteq overline{B}(p, diam A)}
        A \subseteq \overline{B}(p, \diam A).
\end{equation}
Of course,
\begin{equation}
\label{diam overline{B}(x, r) le 2 r}
        \diam \overline{B}(x, r) \le 2 \, r
\end{equation}
for every $x \in M$ and $r \ge 0$, by the triangle inequality.  If
$d(\cdot, \cdot)$ is an ultrametric on $M$, then
\begin{equation}
\label{diam overline{B}(x, r) le r}
        \diam \overline{B}(x, r) \le r
\end{equation}
for every $x \in M$ and $r \ge 0$.  If $M$ is the real line with
the standard metric, and $A \subseteq {\bf R}$ is bounded, then
there is a closed interval $I \subseteq {\bf R}$ such that
\begin{equation}
\label{diam A = diam I}
        \diam A = \diam I.
\end{equation}
In this case, $I$ is the same as the closed ball centered at its
midpoint, with radius equal to one-half the diameter of $I$,
which is the same as the length of $I$.

        Let $n$ be a positive integer, and suppose that $M = {\bf R}^n$,
equipped with the metric
\begin{equation}
\label{d(x, y) = max_{1 le j le n} |x_j - y_j|}
        d(x, y) = \max_{1 \le j \le n} |x_j - y_j|.
\end{equation}
Of course, the topology on ${\bf R}^n$ determined by this metric is
the same as the standard topology on ${\bf R}^n$, which is the product
topology on ${\bf R}^n$ associated to the standard topology on ${\bf
  R}$.  Note that open and closed balls in ${\bf R}^n$ with respect to
(\ref{d(x, y) = max_{1 le j le n} |x_j - y_j|}) are open and closed
cubes with sides parallel to the axes, respectively.  If $A \subseteq
{\bf R}^n$ is a bounded set with diameter less than or equal to $r$
for some $r \ge 0$, then the $n$ coordinate projections of $A$ into
${\bf R}$ each have diameter less than or equal to $r$ as well.  This
implies that $A$ is contained in a closed cube in ${\bf R}^n$ with
sides parallel to the axes and with side length $r$, which is the same
as a closed ball in ${\bf R}^n$ with radius $r/2$ with respect to
the metric (\ref{d(x, y) = max_{1 le j le n} |x_j - y_j|}).

\section{Hausdorff content}
\label{hausdorff content}

        Let $(M, d(x, y))$ be a metric space, and let $\alpha$ be a 
nonnegative real number.  If $A \subseteq M$ is bounded and nonempty,
then $\diam A$ is defined as in the previous section, and so
\begin{equation}
\label{(diam A)^alpha}
        (\diam A)^\alpha
\end{equation}
is defined for each $\alpha > 0$.  If $\alpha = 0$, then (\ref{(diam
  A)^alpha}) is interpreted as being equal to $1$ when $A$ is bounded
and nonempty, even when $A$ has only one element, so that $\diam A =
0$.  It will be convenient to interpret (\ref{(diam A)^alpha}) as
being equal to $+\infty$ for every $\alpha \ge 0$ when $A$ is
unbounded, including $\alpha = 0$.  Similarly, we interpret
(\ref{(diam A)^alpha}) as being equal to $0$ for every $\alpha \ge 0$
when $A = \emptyset$, including $\alpha = 0$.

        The \emph{$\alpha$-dimensional Hausdorff 
content}\index{Hausdorff content} $H^\alpha_{con}(E)$ of a set
$E \subseteq M$ is defined to be the infimum of
\begin{equation}
\label{sum_j (diam A_j)^alpha}
        \sum_j (\diam A_j)^\alpha
\end{equation}
over all collections $\{A_j\}_j$ of finitely or countably many subsets
of $M$ such that $E \subseteq \bigcup_j A_j$.  If $A_j$ is unbounded
for any $j$, then $(\diam A_j)^\alpha = +\infty$, and (\ref{sum_j
  (diam A_j)^alpha}) is infinite.  If there are infinitely many
$A_j$'s, then (\ref{sum_j (diam A_j)^alpha}) is interpreted as the
supremum over all finite subsums, as in Section \ref{summable
  functions}, which may be infinite even when $A_j$ is bounded for
each $j$.

        By construction,
\begin{equation}
\label{H^alpha_{con}(E) le (diam E)^alpha}
        H^\alpha_{con}(E) \le (\diam E)^\alpha
\end{equation}
for every $E \subseteq M$ and $\alpha \ge 0$, since we can cover $E$
by itself.  In particular,
\begin{equation}
\label{H^alpha_{con}(emptyset) = 0}
        H^\alpha_{con}(\emptyset) = 0.
\end{equation}
Alternatively, one might consider the empty set to be covered by the
empty family of subsets of $M$, for which the corresponding empty sum
(\ref{sum_j (diam A_j)^alpha}) is equal to $0$.  If $E \subseteq
\widetilde{E} \subseteq M$, then it is easy to see that
\begin{equation}
\label{H^alpha_{con}(E) le H^alpha_{con}(widetilde{E})}
        H^\alpha_{con}(E) \le H^\alpha_{con}(\widetilde{E})
\end{equation}
for every $\alpha \ge 0$, because every covering of $\widetilde{E}$ in
$M$ is a covering of $E$ too.  If $E_1, E_2, E_3, \ldots$ is any
sequence of subsets of $M$, then one can check that
\begin{equation}
\label{H^alpha_{con}(bigcup_{k = 1}^infty E_l) le ...}
        H^\alpha_{con}\Big(\bigcup_{l = 1}^\infty E_l\Big)
                      \le \sum_{l = 1}^\infty H^\alpha_{con}(E_l),
\end{equation}
for every $\alpha \ge 0$, by combining coverings of the $E_l$'s to get
coverings of their union.  This is a bit simpler in the case of finite
unions, and of course (\ref{H^alpha_{con}(bigcup_{k = 1}^infty E_l) le
  ...})  is trivial when $H^\alpha_{con}(E_l) = +\infty$ for any $l$.
Thus $H^\alpha_{con}$ defines an outer measure on $M$ for every
$\alpha \ge 0$.

        Because of (\ref{diam overline{A} = diam A}), one might as well
restrict one's attention to coverings of $E \subseteq M$ by closed
subsets of $M$ in the definition of $H^\alpha_{con}(E)$.
Alternatively, one can restrict one's attention to coverings of $E$ by
open subsets of $M$ and get the same result for $H^\alpha_{con}(E)$,
because of (\ref{diam A_r le diam A + 2 r}).  This implies that one
can restrict one's attention to coverings of $E$ by finitely many open
subsets of $M$ when $E$ is compact.  If $M$ is the real line with the
standard metric, then one can restrict one's attention to coverings of
$E$ by intervals, using (\ref{diam A = diam I}).  Similarly, if
$d(\cdot, \cdot)$ is an ultrametric on $M$, then one can restrict
one's attention to coverings of $E \subseteq M$ by closed balls in the
definition of $H^\alpha_{con}(E)$, by the remarks in the previous
section.

\section{Restricting the diameters}
\label{restricting the diameters}

        Let $(M, d(x, y))$ be a metric space again, and let $\alpha \ge 0$
and $0 < \delta \le \infty$ be given.  If $E \subseteq M$, then 
$H^\alpha_\delta(E)$ is defined to be the infimum of
\begin{equation}
\label{sum_j (diam A_j)^alpha, 2}
        \sum_j (\diam A_j)^\alpha
\end{equation}
over all collections $\{A_j\}_j$ of finitely or countably many subsets
of $M$ such that $E \subseteq \bigcup_j A_j$ and
\begin{equation}
\label{diam A_j < delta}
        \diam A_j < \delta
\end{equation}
for each $j$, if such a covering exists, and otherwise we put
$H^\alpha_\delta(E) = +\infty$.  If $M$ is separable, then $M$ can be
covered by finitely or countably many balls of radius $r$ for each $r
> 0$, which implies that these coverings always exist.  If $0 < \delta
\le \eta \le \infty$, then it is easy to see that
\begin{equation}
\label{H^alpha_eta(E) le H^alpha_delta(E)}
        H^\alpha_\eta(E) \le H^\alpha_\delta(E),
\end{equation}
because the coverings used to define $H^\alpha_\delta(E)$ can also be
used in the definition of $H^\alpha_\eta(E)$ in this case.

        Note that $H^\alpha_\infty(E)$ is the infimum of
(\ref{sum_j (diam A_j)^alpha, 2}) over all collections $\{A_j\}_j$ of
finitely or countably many bounded subsets of $M$ such that $E
\subseteq \bigcup_j A_j$, and that these coverings always exist.  The
only difference between this and the definition of $H^\alpha_{con}(E)$
is the restriction to coverings of $E$ by bounded subsets of $M$.
This does not affect the infimum, because (\ref{sum_j (diam
  A_j)^alpha, 2}) is infinite as soon as $A_j$ is unbounded for any
$j$.  It follows that
\begin{equation}
\label{H^alpha_infty(E) = H^alpha_{con}(E)}
        H^\alpha_\infty(E) = H^\alpha_{con}(E)
\end{equation}
for every $E \subseteq M$.

        It is easy to see that $H^\alpha_\delta$ is an outer measure on $M$
for each $\delta > 0$, for essentially the same reasons as for
$H^\alpha_{con}$.  As before, one can restrict one's attention to
coverings of $E$ by open subsets of $M$ or closed subsets of $M$ in
the definition of $H^\alpha_\delta(E)$, and to finite coverings of $E$
when $E$ is compact.  If $M = {\bf R}$ with the standard metric,
then one can restrict one's attention to coverings of $E \subseteq {\bf R}$
by intervals.  If $d(x, y)$ is an ultrametric on a set $M$, then one
can restrict one's attention to coverings of $E \subseteq M$ by closed
balls.

        Alternatively, let $\widetilde{H}^\alpha_\delta(E)$ be the infimum
of (\ref{sum_j (diam A_j)^alpha, 2}) over all collections $\{A_j\}_j$
of finitely or countably many subsets of $M$ such that $E \subseteq
\bigcup_j A_j$ and
\begin{equation}
\label{diam A_j le delta}
        \diam A_j \le \delta
\end{equation}
for each $j$, if such a covering exists, and otherwise put
$\widetilde{H}^\alpha_\delta(E) = +\infty$.  The only difference
between this and $H^\alpha_\delta(E)$ is that we replace (\ref{diam
  A_j < delta}) with (\ref{diam A_j le delta}).  Of course, (\ref{diam
  A_j le delta}) is vacuous when $\delta = \infty$, so that
$\widetilde{H}^\alpha_\infty(E) = H^\alpha_{con}(E)$ by definition.
If $M$ is separable, then these coverings exist for every $\delta >
0$, as before.

        If $0 < \delta \le \eta \le \infty$, then
\begin{equation}
\label{widetilde{H}^alpha_eta(E) le widetilde{H}^alpha_delta(E)}
        \widetilde{H}^\alpha_\eta(E) \le \widetilde{H}^\alpha_\delta(E),
\end{equation}
for the same reasons as in (\ref{H^alpha_eta(E) le H^alpha_delta(E)}).
Similarly, it is easy to see that
\begin{equation}
\label{widetilde{H}^alpha_delta(E) le H^alpha_delta(E)}
        \widetilde{H}^\alpha_\delta(E) \le H^\alpha_\delta(E),
\end{equation}
because any of the coverings of $E$ that can be used in the definition
of $H^\alpha_\delta(E)$ can also be used in the definition of
$\widetilde{H}^\alpha_\delta(E)$.  If $0 < \delta < \eta \le \infty$,
then
\begin{equation}
\label{H^alpha_eta(E) le widetilde{H}^alpha_delta(E)}
        H^\alpha_\eta(E) \le \widetilde{H}^\alpha_\delta(E),
\end{equation}
because the coverings of $E$ used in the definition of
$\widetilde{H}^\alpha_\delta(E)$ can also be used in the definition of
$H^\alpha_\eta(E)$.

        As before, one can check that  $\widetilde{H}^\alpha_\delta$
is an outer measure on $M$ for each $\delta > 0$.  One can also
restrict one's attention to coverings of $E$ by closed subsets of $M$
in the definition of $\widetilde{H}^\alpha_\delta(E)$, because of
(\ref{diam overline{A} = diam A}).  If $M = {\bf R}$ with the standard
metric, then one can restrict one's attention to coverings of $E$ by
closed intervals in the definition of
$\widetilde{H}^\alpha_\delta(E)$, by (\ref{diam A = diam I}).  If
$d(x, y)$ is an ultrametric on $M$, then one can restrict one's
attention to coverings of $E$ by closed balls in the definition of
$\widetilde{H}^\alpha_\delta(E)$, using the remarks in Section
\ref{diameters of sets}.

        The condition (\ref{diam A_j le delta}) seems to be used more
commonly, but there are some advantages to (\ref{diam A_j < delta}).
In particular, the argument for restricting to coverings of $E$ by
open subsets of $M$ does not work for the condition (\ref{diam A_j le
  delta}).  However, the two conditions are essentially equivalent in
the limit as $\delta \to 0$, by (\ref{widetilde{H}^alpha_delta(E) le
  H^alpha_delta(E)}) and (\ref{H^alpha_eta(E) le
  widetilde{H}^alpha_delta(E)}).

\section{Hausdorff measure}
\label{hausdorff measure}

        Let $(M, d(x, y))$ be a metric space again, and let $\alpha \ge 0$
be given.  The \emph{$\alpha$-dimensional Hausdorff
  measure}\index{Hausdorff measure} of $E \subseteq M$ is defined by
\begin{equation}
\label{H^alpha(E) = sup_{delta > 0} H^alpha_delta(E)}
        H^\alpha(E) = \sup_{\delta > 0} H^\alpha_\delta(E).
\end{equation}
This can also be considered as a limit of $H^\alpha_\delta(E)$ as
$\delta \to 0$, because of (\ref{H^alpha_eta(E) le H^alpha_delta(E)}).
It is easy to see that $H^\alpha$ is an outer measure on $M$, because
of the corresponding property of $H^\alpha_\delta$ for each $\delta >
0$.  If $\alpha = 0$, then one can check that $H^\alpha$ reduces to
counting measure on $M$.

        Suppose for the moment that $M = {\bf R}$ with the standard metric.
In this case, we have that
\begin{equation}
\label{H^1_delta(E) = H^1_{con}(E)}
        H^1_\delta(E) = H^1_{con}(E)
\end{equation}
for each $\delta > 0$ and $E \subseteq {\bf R}$.  This is because
every interval in ${\bf R}$ can be expressed as a union of finitely
many subintervals with arbitrarily small diameter, where the sum of
the lengths of the smaller intervals is equal to the length of the
initial interval.  It follows that
\begin{equation}
\label{H^1(E) = H^1_{con}(E)}
        H^1(E) = H^1_{con}(E)
\end{equation}
for every $E \subseteq {\bf R}$, by taking the supremum over $\delta >
0$ of (\ref{H^1_delta(E) = H^1_{con}(E)}).  Note that $H^1_{con}$ is
basically the same as Lebesgue outer measure on ${\bf R}$.  In
particular,
\begin{equation}
\label{H^1_{con}([a, b]) = b - a}
        H^1_{con}([a, b]) = b - a
\end{equation}
for every $a, b \in {\bf R}$ with $a \le b$, where $[a, b]$ is the
usual closed interval in ${\bf R}$ from $a$ to $b$, consisting of
$x \in {\bf R}$ with $a \le x \le b$.  More precisely,
\begin{equation}
\label{H^1_{con}([a, b]) le diam [a, b] = b - a}
        H^1_{con}([a, b]) \le \diam [a, b] = b - a
\end{equation}
by (\ref{H^alpha_{con}(E) le (diam E)^alpha}) in Section
\ref{hausdorff content}.  In order to get the opposite inequality,
it suffices to show that
\begin{equation}
\label{b - a le sum_j diam A_j}
        b - a \le \sum_j \diam A_j
\end{equation}
for any collection $\{A_j\}_j$ of finitely or countably many subsets
of ${\bf R}$ such that $[a, b] \subseteq \bigcup_j A_j$.  As in
Section \ref{hausdorff content}, one can reduce to the case where the
$A_j$'s are open intervals in ${\bf R}$, and it is enough to consider
coverings of $[a, b]$ by finitely many open intervals, because $[a,
  b]$ is compact.

        Let $(M, d(x, y))$ be an arbitrary metric space again, and
suppose that $E_1$, $E_2$ are subsets of $M$ such that
\begin{equation}
\label{d(x, y) ge eta}
        d(x, y) \ge \eta
\end{equation}
for some $\eta > 0$ and every $x \in E_1$, $y \in E_2$.  This implies
that any subset of $M$ with diameter less than $\delta \le \eta$
cannot intersect both $E_1$ and $E_2$.  Thus coverings of $E_1 \cup
E_2$ by such subsets of $M$ can be split into coverings of $E_1$ and
$E_2$ separately.  It follows that
\begin{equation}
\label{H^alpha_delta(E_1 cup E_2) ge H^alpha_delta(E_1) + H^alpha_delta(E_2)}
 H^\alpha_\delta(E_1 \cup E_2) \ge H^\alpha_\delta(E_1) + H^\alpha_\delta(E_2)
\end{equation}
for $0 < \delta \le \eta$, by splitting the corresponding sums
(\ref{sum_j (diam A_j)^alpha, 2}).  Taking the limit as $\delta \to
0$, we get that
\begin{equation}
\label{H^alpha(E_1 cup E_2) ge H^alpha(E_1) + H^alpha(E_2)}
        H^\alpha(E_1 \cup E_2) \ge H^\alpha(E_1) + H^\alpha(E_2)
\end{equation}
under these conditions.  Because $H^\alpha$ is an outer measure on
$M$, there is a standard way to define a collection of measurable
subsets of $M$ with respect to $H^\alpha$.  More precisely, this
collection of measurable sets is a $\sigma$-algebra of subsets of $M$,
and $H^\alpha$ is countably additive on this $\sigma$-algebra.  It is
also well known that the Borel subsets of $M$ are measurable with
respect to $H^\alpha$, because of (\ref{H^alpha(E_1 cup E_2) ge
  H^alpha(E_1) + H^alpha(E_2)}).

        Of course,
\begin{equation}
\label{H^alpha_{con}(E) le H^alpha(E)}
        H^\alpha_{con}(E) \le H^\alpha(E)
\end{equation}
for every $E \subseteq M$, by (\ref{H^alpha_infty(E) =
  H^alpha_{con}(E)}).  If $H^\alpha_{con}(E) = 0$, then the coverings
of $E$ for which the corresponding sums (\ref{sum_j (diam A_j)^alpha})
are small automatically involve subsets of $M$ with small diameter.
This implies that $H^\alpha_\delta(E) = 0$ for every $\delta > 0$, and
hence that $H^\alpha(E) = 0$.  It follows that $H^\alpha(E) = 0$ if
and only if $H^\alpha_{con}(E) = 0$, using (\ref{H^alpha_{con}(E) le
  H^alpha(E)}) for the ``only if'' part.

\section{Hausdorff dimension}
\label{Hausdorff dimension}

         Let $(M, d(x, y))$ be a metric space, and let $E \subseteq M$,
$0 \le \alpha < \beta < \infty$, and $0 < \delta < \infty$ be given.
If $\{A_j\}_j$ is a collection of finitely or countably many subsets
of $M$ such that $E \subseteq \bigcup_j A_j$ and $\diam A_j < \delta$
for each $j$, then
\begin{equation}
\label{H^beta_delta(E) le ... le delta^{beta - alpha} sum_j (diam A_j)^alpha}
        H^\beta_\delta(E) \le \sum_j (\diam A_j)^\beta
                         \le \delta^{\beta - \alpha} \sum_j (\diam A_j)^\alpha,
\end{equation}
using the definition of $H^\beta_\delta(E)$ in the first step.  Taking
the infimum over all such collections $\{A_j\}_j$, we get that
\begin{equation}
\label{H^beta_delta(E) le ... le delta^{beta - alpha} H^alpha(E)}
        H^\beta_\delta(E) \le \delta^{\beta - \alpha} \, H^\alpha_\delta(E)
                         \le \delta^{\beta - \alpha} \, H^\alpha(E),
\end{equation}
using the definition (\ref{H^alpha(E) = sup_{delta > 0}
  H^alpha_delta(E)}) of $H^\alpha(E)$ in the second step.  If
$H^\alpha(E) < \infty$, then it follows that
\begin{equation}
\label{H^beta(E) = 0}
        H^\beta(E) = 0
\end{equation}
for every $\beta > \alpha$, by taking the limit as $\delta \to 0$ in
(\ref{H^beta_delta(E) le ... le delta^{beta - alpha} H^alpha(E)}).

        The \emph{Hausdorff dimension}\index{Hausdorff dimension}
of $E$ may be defined as
\begin{equation}
\label{{dim}_H E = inf {alpha ge 0 : H^alpha(E) < infty}}
        {\dim}_H E = \inf \{\alpha \ge 0 : H^\alpha(E) < \infty\},
\end{equation}
which is interpreted as being $+\infty$ when $H^\alpha(E) = \infty$ for
every $\alpha \ge 0$.  If $\beta \ge 0$ satisfies $\beta > \dim_H E$,
then there is an $\alpha \ge 0$ such that $\alpha < \beta$ and
$H^\alpha(E) < \infty$.  This implies that $H^\beta(E) = 0$, as in
(\ref{H^beta(E) = 0}).  Thus $\dim_H E$ can also be given by
\begin{equation}
\label{{dim}_H E = inf {beta ge 0 : H^beta(E) = 0}}
        {\dim}_H E = \inf \{\beta \ge 0 : H^\beta(E) = 0\},
\end{equation}
where the infimum is interpreted as being $+\infty$ when $H^\beta(E) >
0$ for every $\beta \ge 0$.  Combining this with the remarks at the
end of the preceding section, we get that $\dim_H E$ may be defined
equivalently by
\begin{equation}
\label{{dim}_H(E) = inf {beta ge 0 : H^beta_{con}(E) = 0}}
        {\dim}_H(E) = \inf \{\beta \ge 0 : H^\beta_{con}(E) = 0\},
\end{equation}
which is also interpreted as being $+ \infty$ when $H^\beta_{con}(E) >
0$ for every $\beta \ge 0$.

        If $E \subseteq \widetilde{E} \subseteq M$, then it is easy to
see that
\begin{equation}
\label{{dim}_H E le {dim}_H widetilde{E}}
        {\dim}_H E \le {\dim}_H \widetilde{E},
\end{equation}
because of the analogous property of the Hausdorff measures.  Let
$\alpha \ge 0$ be given, and suppose that $E_1, E_2, E_3, \ldots$ is a
sequence of subsets of $M$ such that
\begin{equation}
\label{{dim}_H E_l le alpha}
        {\dim}_H E_l \le \alpha
\end{equation}
for each $l \ge 1$.  This implies that
\begin{equation}
\label{H^beta(E_l) = 0}
        H^\beta(E_l) = 0
\end{equation}
for every $\beta > \alpha$ and $l \ge 1$, as in (\ref{H^beta(E) = 0}).
It follows that
\begin{equation}
\label{H^beta(bigcup_{l = 1}^infty E_l) = 0}
        H^\beta\Big(\bigcup_{l = 1}^\infty E_l\Big) = 0
\end{equation}
for every $\beta > \alpha$, and hence that
\begin{equation}
\label{{dim}_H (bigcup_{l = 1}^infty E_l) le alpha}
        {\dim}_H \Big(\bigcup_{l = 1}^\infty E_l\Big) \le \alpha.
\end{equation}
In particular, a subset of $M$ with only finitely or countably many
elements has Hausdorff dimension $0$.

        Let $a$ be a positive real number, and suppose that $d(x, y)^a$
is also a metric on $M$.  As in Section \ref{metrics, ultrametrics},
this condition holds automatically when $0 < a \le 1$, and for every
$a > 0$ when $d(x, y)$ is an ultrametric on $M$.  If $A \subseteq M$,
then it is easy to see that the diameter of $A$ with respect to $d(x,
y)^a$ is the same as the diameter of $A$ with respect to $d(x, y)$ to
the power $a$.  Using this, one can check that the
$\alpha$-dimensional Hausdorff content of $E \subseteq M$ with respect
to $d(x, y)^a$ is the same as the $(\alpha \, a)$-dimensional content
of $E$ with respect to $d(x, y)$.  Similarly, the $\alpha$-dimensional
Hausdorff measure of $E \subseteq M$ with respect to $d(x, y)^a$ is
equal to the $(\alpha \, a)$-dimensional Hausdorff measure of $E$ with
respect to $d(x, y)$.  There is an analogous statement for
$H^\alpha_\delta$ for each $\delta > 0$, where the $\delta$ used for
$d(x, y)$ corresponds to $\delta^a$ for $d(x, y)^a$.  It follows that
the Hausdorff dimension of $E \subseteq M$ with respect to $d(x, y)^a$
is equal to the Hausdorff dimension of $E$ with respect to $d(x, y)$
divided by $a$.

\section{Some regularity properties}
\label{some regularity properties}

        Let $(M, d(x, y))$ be a metric space, and let $E$ be a subset
of $M$.  Suppose that the $\alpha$-dimensional Hausdorff measure of
$E$ is finite for some $\alpha \ge 0$.  In particular, this means that
$H^\alpha_\delta(E) < \infty$ for each $\delta > 0$, where
$H^\alpha_\delta$ is as in Section \ref{restricting the diameters}.
Let $l$ be a positive integer, and consider $H^\alpha_\delta(E)$ with
$\delta = 1/l$.  By definition of $H^\alpha_\delta(E)$, there is a
collection $\{A_{j, l}\}_j$ of finitely or countable many subsets of
$M$ such that $E \subseteq \bigcup_j A_{j, l}$,
\begin{equation}
\label{diam A_{j, n} < 1/l}
        \diam A_{j, l} < 1/l
\end{equation}
for each $j$, and
\begin{equation}
\label{sum_j (diam A_{j,l})^alpha < H^alpha_{1/l}(E) + 1/l le H^alpha(E) + 1/l}
 \sum_j (\diam A_{j,l})^\alpha < H^\alpha_{1/l}(E) + 1/l \le H^\alpha(E) + 1/l.
\end{equation}

        As in Section \ref{restricting the diameters}, we may as well
ask that $A_{j, l}$ be an open set in $M$ for each $j$ as well.  This
implies that
\begin{equation}
\label{U_l = bigcup_j A_{j, l}}
        U_l = \bigcup_j A_{j, n}
\end{equation}
is an open set in $M$ for each $l$, so that
\begin{equation}
\label{widetilde{E} = bigcap_{l = 1}^infty U_l}
        \widetilde{E} = \bigcap_{l = 1}^\infty U_l
\end{equation}
is a Borel set in $M$, and more precisely a $G_\delta$ set in $M$.  Of
course, $E \subseteq U_l$ for each $l$, by construction, and hence $E
\subseteq \widetilde{E}$.  Similarly, if $A_{j, l}$ is a Borel set for
each $j$ and $l$, then $U_l$ is a Borel set for each $l$, and hence
$\widetilde{E}$ is a Borel set.  In particular, if we had used the
outer measures $\widetilde{H}^\alpha_\delta$ in Section
\ref{restricting the diameters} instead of $H^\alpha_\delta$, then we
could take $A_{j,l}$ to be a closed set for each $j$ and $l$, if not
an open set.

        Let us check that
\begin{equation}
\label{H^alpha(widetilde{E}) = H^alpha(E)}
        H^\alpha(\widetilde{E}) = H^\alpha(E).
\end{equation}
To do this, it suffices to show that
\begin{equation}
\label{H^alpha(widetilde{E}) le H^alpha(E)}
        H^\alpha(\widetilde{E}) \le H^\alpha(E),
\end{equation}
because the opposite inequality follows from the fact that $E
\subseteq \widetilde{E}$.  Observe that
\begin{equation}
\label{H^alpha_{1/l}(widetilde{E}) le sum_j (diam A_{j, l})^alpha}
        H^\alpha_{1/l}(\widetilde{E}) \le \sum_j (\diam A_{j, l})^\alpha
\end{equation}
for each $l$, because $\{A_{j, l}\}_j$ can also be used as a covering
of $\widetilde{E}$ in the definition of
$H^\alpha_{1/l}(\widetilde{E})$ for each $l$.  Combining this with
(\ref{sum_j (diam A_{j,l})^alpha < H^alpha_{1/l}(E) + 1/l le
  H^alpha(E) + 1/l}), we get that
\begin{equation}
\label{H^alpha_{1/l}(widetilde{E}) < H^alpha(E) + 1/l}
        H^\alpha_{1/l}(\widetilde{E}) < H^\alpha(E) + 1/l
\end{equation}
for each $l$.  If we take the limit as $l \to \infty$ in
(\ref{H^alpha_{1/l}(widetilde{E}) < H^alpha(E) + 1/l}), then we get
(\ref{H^alpha(widetilde{E}) le H^alpha(E)}), as desired.

        It would be nice if for every $\epsilon > 0$ there is an
open set $U \subseteq M$ such that $E \subseteq U$ and
\begin{equation}
\label{H^alpha(U) < H^alpha(E) + epsilon}
        H^\alpha(U) < H^\alpha(E) + \epsilon,
\end{equation}
but this does not always work.  If $M$ is the real line with the
standard metric and $\alpha < 1$, for instance, then $H^\alpha(U) =
+\infty$ for every nonempty open set $U \subseteq {\bf R}$.  It is
well known that this type of approximation does work when $E$ is
contained in a countable union of open sets with finite measure.  More
precisely, if $E$ is a Borel set, then this follows from Theorem 2.2.2
on p60 of \cite{fed}, or Theorem 1.10 on p11 of \cite{pm}.  If $E$ is
not a Borel set, then one can reduce to this case by replacing $E$
with the intersection of $\widetilde{E}$ from (\ref{widetilde{E} =
  bigcap_{l = 1}^infty U_l}) with the countable union of open sets
just mentioned.

        Alternatively, put
\begin{equation}
\label{V_n = bigcap_{l = 1}^n U_l}
        V_n = \bigcap_{l = 1}^n U_l
\end{equation}
for each positive integer $n$, so that $V_n$ is an open set, $V_{n + 1}
\subseteq V_n$, and
\begin{equation}
\label{bigcap_{n = 1}^infty V_n = bigcap_{l = 1}^infty U_l = widetilde{E}}
        \bigcap_{n = 1}^\infty V_n = \bigcap_{l = 1}^\infty U_l = \widetilde{E}.
\end{equation}
If $H^\alpha(V_n) < \infty$ for any $n$, then
\begin{equation}
\label{lim_{n to infty} H^alpha(V_n) = H^alpha(widetilde{E}) = H^alpha(E)}
        \lim_{n \to \infty} H^\alpha(V_n) = H^\alpha(\widetilde{E}) = H^\alpha(E),
\end{equation}
by standard results from measure theory, since $H^\alpha$ is countably
additive on Borel sets.  This gives another way to look at
(\ref{H^alpha(U) < H^alpha(E) + epsilon}), using the definition of
$H^\alpha$.  Similarly, if $E$ is contained in an open set $W
\subseteq M$ with $H^\alpha(W) < \infty$, then one can apply the same
argument to $V_n \cap W$.  A slightly more complicated version of this
argument can also be used when $E$ is contained in a countable union
of open sets with finite $H^\alpha$ measure.

        In some situations, we may have that
\begin{equation}
\label{H^alpha(A) le (diam A)^alpha}
        H^\alpha(A) \le (\diam A)^\alpha
\end{equation}
for every $A \subseteq M$, and for a fixed $\alpha \ge 0$.  This holds
when $A = {\bf R}$ with the standard metric and $\alpha = 1$, for
instance.  In this case, if $U_l$ is as in (\ref{U_l = bigcup_j A_{j,
    l}}), then
\begin{equation}
\label{H^alpha(U_l) le sum_j H^alpha(A_{j, l}) le sum_j (diam A_{j, l})^alpha}
 H^\alpha(U_l) \le \sum_j H^\alpha(A_{j, l}) \le \sum_j (\diam A_{j, l})^\alpha
\end{equation}
for each $l$.  Combining this with (\ref{sum_j (diam A_{j,l})^alpha <
  H^alpha_{1/l}(E) + 1/l le H^alpha(E) + 1/l}), we get that
\begin{equation}
\label{H^alpha(U_l) < H^alpha(E) + 1/l}
        H^\alpha(U_l) < H^\alpha(E) + 1/l
\end{equation}
for each $l$.

        Suppose instead that
\begin{equation}
\label{H^alpha(A) le C (diam A)^alpha}
        H^\alpha(A) \le C \, (\diam A)^\alpha
\end{equation}
for some $C \ge 1$ and every $A \subseteq M$.  This implies that
\begin{equation}
\label{H^alpha(U_l) le C sum_j (diam A_{j, l})^alpha}
        H^\alpha(U_l) \le C \, \sum_j (\diam A_{j, l})^\alpha
\end{equation}
for each $l$, as in (\ref{H^alpha(U_l) le sum_j H^alpha(A_{j, l}) le
  sum_j (diam A_{j, l})^alpha}).  It follows that $H^\alpha(U_l)$ is
finite for each $l$, by (\ref{sum_j (diam A_{j,l})^alpha <
  H^alpha_{1/l}(E) + 1/l le H^alpha(E) + 1/l}) and the hypothesis that
$H^\alpha(E) < \infty$.  In particular, $H^\alpha(V_n)$ is finite for
each $n$, where $V_n$ is as in (\ref{V_n = bigcap_{l = 1}^n U_l}), so
that (\ref{lim_{n to infty} H^alpha(V_n) = H^alpha(widetilde{E}) =
  H^alpha(E)}) holds, as before.

\section{Lipschitz mappings, 2}
\label{lipschitz mappings, 2}

        Let $(M_1, d_1(x, y))$ and $(M_2, d_2(u, v))$ be metric spaces,
and suppose that $f$ is a Lipschitz mapping from $M_1$ into $M_2$ of
order $a > 0$ with constant $C \ge 0$, as in Section \ref{lipschitz
mappings}.  If $A$ is a bounded subset of $M_1$, then $f(A)$ is a
bounded subset of $M_2$, and
\begin{equation}
\label{diam f(A) le C (diam A)^a}
        \diam f(A) \le C \, (\diam A)^a.
\end{equation}
More precisely, $\diam A$ refers to the diameter of $A$ as a subset of
$M_1$, with respect to $d_1(x, y)$, and $\diam f(A)$ refers to the
diameter of $f(A)$ in $M_2$, with respect to $d_2(u, v)$.  This also
works when $A$ is unbounded, in which case we interpret $C \, (\diam
A)^a$ to be $+\infty$ when $C > 0$, and to be $0$ when $C = 0$.

        Let $E \subseteq M_1$ and a nonnegative real number $\alpha$
be given.  If $\{A_j\}_j$ is any collection of finitely or countably
many subsets of $M_1$ such that $E \subseteq \bigcup_j A_j$, then
$\{f(A_j)\}_j$ is a collection of finitely or countably many subsets
of $M_2$ such that $f(E) \subseteq \bigcup_j f(A_j)$.  Using this and
(\ref{diam f(A) le C (diam A)^a}), it is easy to see that
\begin{equation}
\label{H^alpha_{con}(f(E)) le C^alpha H^{alpha a}_{con}(E)}
        H^\alpha_{con}(f(E)) \le C^\alpha \, H^{\alpha \, a}_{con}(E).
\end{equation}
Similarly, if $\delta, \eta > 0$ satisfy $C \, \delta^a \le \eta$,
then we get that
\begin{equation}
\label{H^alpha_eta(f(E)) le C^alpha H^{alpha a}_delta(E)}
        H^\alpha_\eta(f(E)) \le C^\alpha \, H^{\alpha \, a}_\delta(E).
\end{equation}
It follows that
\begin{equation}
\label{H^alpha_eta(f(E)) le C^alpha H^{alpha a}(E)}
        H^\alpha_\eta(f(E)) \le C^\alpha \, H^{\alpha \, a}(E)
\end{equation}
for every $\eta > 0$, and hence that
\begin{equation}
\label{H^alpha(f(E)) le C^alpha H^{alpha a}(E)}
        H^\alpha(f(E)) \le C^\alpha \, H^{\alpha \, a}(E).
\end{equation}
As before, the various measures of $E$ are defined using $d_1(x, y)$
on $M_1$, and the various measures of $f(E)$ are defined using $d_2(u,
v)$ on $M_2$.  If $\alpha = 0$, then $C^\alpha$ is interpreted as
being equal to $1$ for every $C \ge 0$.  If $\alpha > 0$ and $C = 0$,
then the right sides of these inequalities may be interpreted as being
equal to $0$, even when the corresponding measure of $E$ is infinite.
In particular,
\begin{equation}
\label{{dim}_H f(E) le a^{-1} {dim}_H E}
        {\dim}_H f(E) \le a^{-1} \, {\dim}_H E
\end{equation}
for every $E \subseteq M_1$ under these conditions.

        A mapping $f : M_1 \to M_2$ is said to be
\emph{bilipschitz}\index{bilipschitz mappings} with constant $C \ge 1$ if
\begin{equation}
\label{C^{-1} d_1(x, y) le d_2(f(x), f(y)) le C d_1(x, y)}
        C^{-1} \, d_1(x, y) \le d_2(f(x), f(y)) \le C \, d_1(x, y)
\end{equation}
for every $x, y \in M_1$.  In this case, we have that
\begin{equation}
\label{C^{-alpha} H^alpha(E) le H^alpha(f(E)) le C^alpha H^alpha(E)}
        C^{-\alpha} \, H^\alpha(E) \le H^\alpha(f(E)) \le C^\alpha \, H^\alpha(E)
\end{equation}
for every $\alpha \ge 0$ and $E \subseteq M_1$.  More precisely, the
second inequality in (\ref{C^{-alpha} H^alpha(E) le H^alpha(f(E)) le
  C^alpha H^alpha(E)}) follows from (\ref{H^alpha(f(E)) le C^alpha
  H^{alpha a}(E)}), with $a = 1$.  The first inequality in
(\ref{C^{-alpha} H^alpha(E) le H^alpha(f(E)) le C^alpha H^alpha(E)})
is essentially the same as the second inequality, applied to the
inverse of $f$.  This also uses the fact that we can restrict our
attention to coverings of $f(E)$ by subsets of $f(E)$ in the
definition of $H^\alpha(f(E))$.  As before, it follows that
\begin{equation}
\label{{dim}_H f(E) = {dim}_H E}
        {\dim}_H f(E) = {\dim}_H E
\end{equation}
for every $E \subseteq M_1$ under these conditions.  Note that $f$ is
bilipschitz with constant $C = 1$ if and only if $f$ is an isometric
embedding.

        Let us now take $M_2 = {\bf R}$, equipped with the standard metric.
Remember that
\begin{equation}
\label{f_p(x) = d_1(x, p), 2}
        f_p(x) = d_1(x, p)
\end{equation}
is a Lipschitz mapping of order $1$ from $M_1$ into ${\bf R}$ with
constant $C = 1$ for every $p \in M_1$, as in Section \ref{lipschitz
  mappings}.  If $p, q \in E \subseteq M_1$, then
\begin{equation}
\label{f_p(p) = d_1(p, p) = 0 and f_p(q) = d_1(p, q)}
        f_p(p) = d_1(p, p) = 0 \quad\hbox{and}\quad f_p(q) = d_1(p, q)
\end{equation}
are elements of $f_p(E)$, so that
\begin{equation}
\label{d_1(p, q) le diam f_p(E)}
        d_1(p, q) \le \diam f_p(E).
\end{equation}
If $E$ is connected, then $f_p(E)$ is a connected subset of ${\bf R}$,
because $f$ is continuous.  Of course, $H^1_{con}$ is the same as
Lebesgue outer measure on ${\bf R}$, as in Section \ref{hausdorff
  measure}, and the Lebesgue measure of a connected subset of ${\bf
  R}$ is equal to its diameter.  Thus
\begin{equation}
\label{diam f_p(E) = H^1_{con}(f_p(E))}
        \diam f_p(E) = H^1_{con}(f_p(E))
\end{equation}
when $E$ is connected.  We also have that
\begin{equation}
\label{H^1_{con}(f_p(E)) le H^1_{con}(E)}
        H^1_{con}(f_p(E)) \le H^1_{con}(E),
\end{equation}
by (\ref{H^alpha_{con}(f(E)) le C^alpha H^{alpha a}_{con}(E)}).  It
follows that
\begin{equation}
\label{d_1(p, q) le H^1_{con}(E)}
        d_1(p, q) \le H^1_{con}(E)
\end{equation}
for every $p, q \in E$ when $E$ is connected, and hence that
\begin{equation}
\label{diam E le H^1_{con}(E) le H^1(E)}
        \diam E \le H^1_{con}(E) \le H^1(E).
\end{equation}

\section{Spherical measure}
\label{spherical measure}

        Let $(M, d(x, y))$ be a metric space, and let $\alpha$ be a
nonnegative real number.  Suppose that in the definition of
$\alpha$-dimensional Hausdorff content of $E \subseteq M$, we restrict
our attention to collections $\{A_j\}_j$ of finitely or countably many
\emph{closed balls} in $M$, instead of arbitary subsets of $M$.  Let
us call the infimum of (\ref{sum_j (diam A_j)^alpha}) in Section
\ref{hausdorff content} over all such coverings of $E$ the
\emph{$\alpha$-dimensional spherical content}\index{spherical content}
of $E$.  Thus $H^\alpha_{con}(E)$ is automatically less than or equal
to the $\alpha$-dimensional spherical content of $E$, since all of the
coverings of $E$ used in the definition of the $\alpha$-dimensional
spherical content of $E$ can also be used in the definition of
$H^\alpha_{con}(E)$.  In the other direction, the $\alpha$-dimensional
spherical content of $E$ is less than or equal to $2^\alpha$ times
$H^\alpha_{con}(E)$.  This follows from the fact that every nonempty
bounded set $A \subseteq M$ is contained in a closed ball with radius
$\diam A$ and diameter less than or equal to $2 \, \diam A$, as in
Section \ref{diameters of sets}.  If $A$ is unbounded, then one might
interpret $M$ as a ball of infinite radius, but this case does not
really matter for estimating the $\alpha$-dimensional spherical content
of $E$ in terms of $H^\alpha_{con}(E)$.

        Similarly, we can restrict our attention to coverings of $E$ by
collections $\{A_j\}_j$ of finitely or countably many closed balls
$A_j$ in $M$ with diameter less than $\delta$ in Section
\ref{restricting the diameters}, to get spherical versions of the
outer measures $H^\alpha_\delta$.  As before, $H^\alpha_\delta(E)$ is
automatically less than or equal to its spherical version, because
every covering of $E$ used for the spherical version can also be used
for $H^\alpha_\delta(E)$.  One can also check that the spherical
version of $H^\alpha_\delta(E)$ is less than or equal to $2^\alpha$
times $H^\alpha_{2 \, \delta}(E)$, using the same fact about nonempty
bounded subsets $A$ of $M$ being contained in closed balls as in the
preceding paragraph.  Note that $\delta$ is replaced by $2 \, \delta$
in $H^\alpha_{2 \, \delta}(E)$, because of the extra factor of $2$ in
the estimate for the diameter of the closed ball that contains $A$.

        The supremum of the spherical version of $H^\alpha_\delta(E)$
over $\delta > 0$ is known as the \emph{$\alpha$-dimensional spherical
  measure}\index{spherical measure} of $E$.  Thus $H^\alpha(E)$ is
less than or equal to the $\alpha$-dimensional spherical measure of
$E$, because of the analogous inequality for $H^\alpha_\delta(E)$.
The $\alpha$-dimensional spherical measure of $E$ is less than or
equal to $2^\alpha$ times $H^\alpha(E)$, because of the corresponding
estimate for the spherical version of $H^\alpha_\delta(E)$.  Note that
the analogue of Hausdorff dimension for spherical measures is the same
as the ordinary Hausdorff dimension.

        If $M$ is the real line with the standard metric, then
spherical measures are the same as Hausdorff measures.  This is because
Hausdorff measures on ${\bf R}$ can already be defined in terms of
coverings by closed intervals, which are the same as closed balls in
this case.  Similarly, if $d(x, y)$ is an ultrametric on any set $M$,
then we have seen that the corresponding Hausdorff measures can be
defined in terms of coverings by closed balls, so that spherical measures
are the same as Hausdorff measures in this situation as well.

        Let $(M, d(x, y))$ be an arbitrary metric space again, and let
$Y$ be a subset of $M$.  Thus $Y$ can also be considered as a metric space,
using the restriction of $d(x, y)$ to $Y$.  If $E \subseteq Y$, then
the $\alpha$-dimensional Hausdorff measure of $E$ as a subset of $M$
is equal to the $\alpha$-dimensional Hausdorff measure of $E$ as a
subset of $Y$.  Of course, coverings of $E$ as a subset of $Y$ are already
coverings of $E$ as a subset of $M$.  In the other direction, if
$\{A_j\}_j$ is any collection of subsets of $M$ such that $E \subseteq
\bigcup_j A_j$, then $\{A_j \cap Y\}_j$ is a collection of subsets of
$Y$ such that $E \subseteq \bigcup_j A_j \cap Y$, and
\begin{equation}
\label{diam (A_j cap Y) le diam A_j}
        \diam (A_j \cap Y) \le \diam A_j
\end{equation}
for each $j$.  This fact was implicitly used in the first inequality
in (\ref{C^{-alpha} H^alpha(E) le H^alpha(f(E)) le C^alpha
  H^alpha(E)}).  More precisely, this is essentially the same as
(\ref{C^{-alpha} H^alpha(E) le H^alpha(f(E)) le C^alpha H^alpha(E)}),
with $C = 1$.

        However, this type of argument does not work for spherical measures.
One problem is that the intersection of a ball in $M$ with $Y$ may not
be a ball in $Y$.  Another problem is that although every ball in $Y$
can be expressed as the intersection of $Y$ with a ball in $M$ with
the same radius, the ball in $M$ may have larger diameter.  One way to
try to avoid the first problem is to restrict one's attention further
to coverings of a set $E$ by balls centered on $E$.  This leads to another
problem, which is that such a covering of $E$ might not be admissible
as a covering of a subset of $E$.

        Hausdorff and spherical measures do have a number of properties
in common, and indeed both are examples of a well-known construction
of Carath\'eodory.  In particular, $\alpha$-dimensional spherical
measure defines an outer measure on $M$ for each $\alpha \ge 0$.  This
leads to a corresponding $\sigma$-algebra of measurable subsets of $M$
for each $\alpha \ge 0$, on which $\alpha$-dimensional spherical
measure is countably additive.  If $E_1, E_2 \subseteq M$ satisfy
(\ref{d(x, y) ge eta}) in Section \ref{hausdorff measure} for some
$\eta > 0$, then it is easy to see that the analogues of
(\ref{H^alpha_delta(E_1 cup E_2) ge H^alpha_delta(E_1) +
  H^alpha_delta(E_2)}) and (\ref{H^alpha(E_1 cup E_2) ge H^alpha(E_1)
  + H^alpha(E_2)}) for spherical measures also hold.  This implies
that Borel subsets of $M$ are measurable with respect to
$\alpha$-dimensional spherical measure for each $\alpha \ge 0$.  If $E
\subseteq M$ has finite $\alpha$-dimensional spherical measure for
some $\alpha \ge 0$, then there is a Borel set $\widetilde{E}
\subseteq M$ that contains $E$ and has the same $\alpha$-dimensional
spherical measure as $E$.  This can be derived from the same type of
argument as in Section \ref{some regularity properties}, although in
this case the sets $U_l$ in (\ref{U_l = bigcup_j A_{j, l}}) are
$F_\sigma$ sets.

        One could also consider coverings by open balls instead of
closed balls, and get similar conclusions.  In some situations,
the analogues of spherical measures using coverings by open balls
are the same as for coverings by closed balls.  This happens when
the diameter of any open or closed ball of radius $r$ in $M$ is equal
to $2 \, r$, for instance.  One can check that this happens when
$d(x, y)$ is an ultrametric on $M$ too, using the fact that the
diameter of a ball of radius $r$ is less than or equal to $r$.
Of course, in some examples of ultrametric spaces, every closed ball
of positive radius can be expressed as an open ball, and vice-versa.

        Suppose that $f$ is a Lipschitz mapping of order $a > 0$
from a metric space $M_1$ into another metric space $M_2$, with
constant $C \ge 0$.  If $B$ is an open or closed ball in $M_1$
centered at a point $x$ and with radius $r$, then $f(B)$ is contained
in the open or closed ball in $M_2$ centered at $f(x)$ with radius $C
\, r^a$, respectively.  In order to estimate the effect on spherical
measures, one should look at the diameter of a ball in $M_2$ that
contains $f(B)$, such as the one just mentioned.  In some situations,
it may be possible to represent $B$ as a ball in $M_1$ centered at $x$
with more than one radius $r$.  In this case, it is better to take $r$
to be as small as possible when working with closed balls, and at
least approximately minimal when working with open balls.

        If $f$ is a bilipschitz mapping from $M_1$ onto $M_2$, then
$f$ is a Lipschitz mapping of order $1$ from $M_1$ into $M_2$, and
the inverse mapping $f^{-1}$ defines a Lipschitz mapping of order $1$
from $M_2$ into $M_1$.  The behavior of these Lipschitz mappings can
be analyzed as in the previous paragraph.  However, if $f$ is not
surjective, then the inverse mapping $f^{-1}$ is only defined as a
Lipschitz mapping from $f(E)$ into $M_1$, and we are back to some of
the problems mentioned earlier in the section.

\section{Euclidean spaces}
\label{euclidean spaces}

        Let $(M, d(x, y))$ be a metric space, and let $\mathcal{A}$ be
a $\sigma$-algebra of subsets of $M$ that contains the Borel subsets
of $M$.  Also let $\mu$ be an outer measure defined on $\mathcal{A}$,
and suppose that
\begin{equation}
\label{mu(A) le C (diam A)^alpha}
        \mu(A) \le C \, (\diam A)^\alpha
\end{equation}
for some $\alpha, C \ge 0$ and every $A \in \mathcal{A}$.  If $E \in
\mathcal{A}$, and if $\{A_j\}_j$ is a collection of finitely or
countably many elements of $\mathcal{A}$ such that $E \subseteq
\bigcup_j A_j$, then
\begin{equation}
\label{mu(E) le sum_j mu(A_j) le sum_j C (diam A_j)^alpha}
        \mu(E) \le \sum_j \mu(A_j) \le \sum_j C \, (\diam A_j)^\alpha.
\end{equation}
This implies that
\begin{equation}
\label{mu(E) le C H^alpha_{con}(E)}
        \mu(E) \le C \, H^\alpha_{con}(E),
\end{equation}
by taking the infimum over all such collections $\{A_j\}_j$ in
(\ref{mu(E) le sum_j mu(A_j) le sum_j C (diam A_j)^alpha}).  More
precisely, this also uses the fact that one can restrict one's
attention to coverings of $E$ by Borel subsets of $M$ in the
definition of $H^\alpha_{con}(E)$, as Section \ref{hausdorff content}.

        Similarly, suppose that $\mu$ satisfies
(\ref{mu(A) le C (diam A)^alpha}) for some $\alpha, C \ge 0$ and all
closed balls $A$ in $M$.  If $E \in \mathcal{A}$, then the same type
of argument as in the preceding paragraph implies that $\mu(E)$ is
less than or equal to $C$ times the spherical version of $H^\alpha_{con}(E)$.
Of course, one may be able to use a smaller constant $C$ in this case,
since one is considering a smaller collection of sets $A$ in
(\ref{mu(A) le C (diam A)^alpha}).

        Let $n$ be a positive integer, and suppose that $M = {\bf R}^n$,
with the standard metric.  If $\mu$ is Lebesgue measure on ${\bf
  R}^n$, then it is easy to see that (\ref{mu(A) le C (diam A)^alpha})
holds with $\alpha = n$ for some $C \ge 0$ and all Borel sets $A
\subseteq {\bf R}^n$.  If $A$ is a ball in ${\bf R}^n$, then
(\ref{mu(A) le C (diam A)^alpha}) holds with $\alpha = n$ and the same
constant $C$ as for the unit ball.  It is well known that (\ref{mu(A)
  le C (diam A)^alpha}) actually holds with $\alpha = n$ for all Borel
sets $A \subseteq {\bf R}^n$, and with the same constant $C$ as for
balls.  This is known as the \emph{isodiametric inequality}.  In fact,
$n$-dimensional Hausdorff measure is equal to a constant multiple of
Lebesgue measure on ${\bf R}^n$, where the constant corresponds to the
one in the isodiametric inequality.  It can also be shown that
$n$-dimensional Hausdorff measure is equal to $n$-dimensional
spherical measure on ${\bf R}^n$.  These measures are often defined
with additional constant factors included, so that they agree with
Lebesgue measure on ${\bf R}^n$.

        It is much easier to check directly that
\begin{equation}
\label{H^n(A) le C' (diam A)^n}
        H^n(A) \le C' \, (\diam A)^n
\end{equation}
for some $C' \ge 0$ and every $A \subseteq {\bf R}^n$, where $C'$
depends only on $n$.  It suffices to consider the case where $A$ is a
cube, which can be covered by smaller cubes in the usual way.  This
implies that $H^n(U)$ is bounded by a constant multiple of the
Lebesgue measure of $U$ for every open set $U \subseteq {\bf R}^n$, by
expressing $U$ as a union of cubes with disjoint interiors.  The same
inequality holds for all Lebesgue measurable subsets of ${\bf R}^n$,
by approximation by open sets.

        Remember that ${\bf R}^n$ is a locally compact commutative
topological group with respect to addition and the standard topology.
Of course, Lebesgue measure satisfies the requirements of Haar measure
on ${\bf R}^n$, and one can check that $H^n$ satisfies the
requirements of Haar measure on ${\bf R}^n$ too.  This implies that
$H^n$ is a constant multiple of Lebesgue measure on the Borel subsets
of ${\bf R}^n$, by the uniqueness of Haar measure.  More precisely,
Hausdorff measure of any dimension is invariant under translations on
${\bf R}^n$, because the standard Euclidean metric on ${\bf R}^n$ is
invariant under translations.  The $n$-dimensional Hausdorff measure
of a bounded subset of ${\bf R}^n$ is finite, as in the preceding
paragraph.  If $U$ is a nonempty open subset of ${\bf R}^n$, then the
Lebesgue measure of $U$ is positive, and hence $H^n(U) > 0$, since we
have already seen that the Lebesgue measure of $U$ is bounded by a
constant times $H^n(U)$.  One can also verify that $H^n(U) > 0$ more
directly.

        Let $Q$ be a cube in ${\bf R}^n$ with sides parallel to the axes,
with sidelength equal to $1$.  Thus $H^n(Q)$ and the Lebesgue measure
are both positive and finite, and $H^n(Q)$ is a constant multiple of
the Lebesgue measure of $Q$, by invariance under translations.  One
can extend this to cubes with sidelength equal to $2^l$ for some $l
\in {\bf Z}$, using invariance under translations again.  This implies
that $H^n(U)$ is equal to the same constant multiple of the Lebesgue
measure of $U$ when $U \subseteq {\bf R}^n$ is an open set, by
expressing $U$ as a union of such cubes with disjoint interiors.  It
follows that $H^n(E)$ is equal to the same constant multiple of the
Lebesgue measure of $E$ for every Borel set $E \subseteq {\bf R}^n$,
by the outer regularity properties of both measures.

        If $N$ is any norm on ${\bf R}^n$, then there is a 
translation-invariant metric on ${\bf R}^n$ associated to $N$, as in
(\ref{d(v, w) = N(v - w)}) in Section \ref{norms, ultranorms}.  This
norm is equivalent to the standard Euclidean norm on ${\bf R}^n$, in
the sense that each is bounded by a constant multiple of the other, as
in Section \ref{finite-dimensional vector spaces}.  Of course, this
implies that the corresponding metrics satisfy the same property, and
in particular that they determine the same topology on ${\bf R}^n$.
Similarly, it is easy to see that $n$-dimensional Hausdorff measure on
${\bf R}^n$ with respect to the metric associated to $N$ is comparable
to $n$-dimensional Hausdorff measure with respect to the standard
metric on ${\bf R}^n$.  As before, one can check that $n$-dimensional
Hausdorff measure on ${\bf R}^n$ with respect to the metric associated
to $N$ satisfies the requirements of Haar measure, and hence is equal
to a constant multiple of Lebesgue measure on the Borel subsets of
${\bf R}^n$.  This is all much simpler when $N(v)$ is the maximum of
the absolute values of the coordinates of $v \in {\bf R}^n$, as in
(\ref{N_0(v) = max(|v_1|, ldots, |v_n|)}).  In this case, the
corresponding metric on ${\bf R}^n$ is the same as (\ref{d(x, y) =
  max_{1 le j le n} |x_j - y_j|}) in Section \ref{diameters of sets}.

\section{A simple covering argument}
\label{a simple covering argument}

        Let $(M, d(x, y))$ be an ultrametric space, and let $B$ be
an open or closed ball in $M$ with radius $r > 0$.  If $x \in B$, then
$B$ may be considered as the open or closed ball in $M$ centered at
$x$ with radius $r$, as appropriate.  This follows from (\ref{B(x, r)
  = B(y, r)}) in Section \ref{metrics, ultrametrics} in the case of
open balls, and (\ref{overline{B}(x, r) = overline{B}(y, r)}) in the
case of closed balls.

        Let $B$, $B'$ be open or closed balls in $M$ with radii
$r, r' > 0$, respectively, such that
\begin{equation}
\label{B cap B' ne emptyset}
        B \cap B' \ne \emptyset.
\end{equation}
If $x \in B \cap B'$, then $B$ and $B'$ can both be considered
as balls centered at $x$ in $M$, as in the preceding paragraph.  Using
this, one can check that either
\begin{equation}
\label{B subseteq B'}
        B \subseteq B'
\end{equation}
or
\begin{equation}
\label{B' subseteq B}
        B' \subseteq B.
\end{equation}
More precisely, (\ref{B subseteq B'}) holds when $r < r'$, and
(\ref{B' subseteq B}) holds when $r' < r$.  If $r = r'$, then (\ref{B
  subseteq B'}) holds when $B = B(x, r)$, and (\ref{B' subseteq B})
holds when $B' = \overline{B}(x, r')$.  If $r = r'$ and $B$, $B'$ are
both open or both closed, then $B = B'$.  Of course, $B$ and $B'$ may
be the same as subsets of $M$, even if they are initially defined as
balls of different radii, or one is initially defined as an open ball
and the other is initially defined as a closed ball.  This also works
when $B$ or $B'$ is a closed ball of radius $0$ in $M$, which is to
say a subset of $M$ with exactly one element.

        Now let $\mathcal{B}$ be a collection of open or closed balls
in $M$.  An element $B$ of $\mathcal{B}$ is said to be \emph{maximal}
in $\mathcal{B}$ if it is maximal with respect to inclusion, which is
to say that for each $B' \in \mathcal{B}$ with $B \subseteq B'$ we
have that $B = B'$.  It is easy to see that the maximal elements of
$\mathcal{B}$ are pairwise disjoint, by the remarks in the preceding
paragraph.

        Let $\mathcal{B}_0$ be the collection of maximal elements of
$\mathcal{B}$, so that
\begin{equation}
\label{bigcup {B : B in mathcal{B}_0} subseteq bigcup {B : B in mathcal{B}}}
        \bigcup \{B : B \in \mathcal{B}_0\}
          \subseteq \bigcup \{B : B \in \mathcal{B}\},
\end{equation}
since $\mathcal{B}_0 \subseteq \mathcal{B}$.  Let us say that
$\mathcal{B}$ is a nice collection of balls in $M$ if every element of
$\mathcal{B}$ is contained in a maximal element of $\mathcal{B}$.
In this case, we have that
\begin{equation}
\label{bigcup {B : B in mathcal{B}} subseteq bigcup {B : B in mathcal{B}_0}}
        \bigcup \{B : B \in \mathcal{B}\}
          \subseteq \bigcup \{B : B \in \mathcal{B}_0\},
\end{equation}
and hence
\begin{equation}
\label{bigcup {B in mathcal{B}_0} = bigcup {B in mathcal{B}}}
        \bigcup \{B \in \mathcal{B}_0\} = \bigcup \{B \in \mathcal{B}\}.
\end{equation}
Actually, in order to get (\ref{bigcup {B : B in mathcal{B}} subseteq
  bigcup {B : B in mathcal{B}_0}}), it suffices to know that for each
$B \in \mathcal{B}$ and $x \in B$ there is a $B_0 \in \mathcal{B}_0$
such that $x \in B_0$.  However, this would imply that $B \subseteq
B_0$ or $B_0 \subseteq B$, since $B \cap B_0 \ne \emptyset$.  If $B_0
\subseteq B$, then $B = B_0$, because $B_0$ is supposed to be maximal
in $\mathcal{B}$.  Thus $B \subseteq B_0$ under these conditions.  In
particular, for each $B \in \mathcal{B}$, this could be applied to any
$x \in B$, to get that $B$ is contained in some $B_0 \in
\mathcal{B}_0$.  It follows that (\ref{bigcup {B : B in mathcal{B}}
  subseteq bigcup {B : B in mathcal{B}_0}}) holds if and only if
$\mathcal{B}$ is nice.

        Let $B_1 \in \mathcal{B}$ be given, and put
\begin{equation}
\label{mathcal{C}(B_1) = {B in mathcal{B} : B_1 subseteq B}}
        \mathcal{C}(B_1) = \{B \in \mathcal{B} : B_1 \subseteq B\}.
\end{equation}
Note that $B_1$ is automatically an element of $\mathcal{C}(B_1)$.  If
$x_1$ is any element of $B_1$, then the elements of $\mathcal{C}(B_1)$
may be considered as balls centered at $x_1$, by the remarks at the
beginning of the section.  In particular, the elements of
$\mathcal{C}(B_1)$ are linearly ordered by inclusion.  If $B \in
\mathcal{C}(B_1)$ and $B$ is a maximal element of $\mathcal{B}$, then
$B$ is a maximal element of $\mathcal{C}(B_1)$.  Conversely, if $B \in
\mathcal{C}(B_1)$ is a maximal element of $\mathcal{C}(B_1)$, then it
is easy to see that $B$ is a maximal element of $\mathcal{B}$ too.  Of
course, if $\mathcal{C}(B_1)$ has only finitely many elements, then it
has a maximal element.  If $\mathcal{B}$ has only finitely many
elements, then $\mathcal{C}(B_1)$ has only finitely many elements for
each $B_1 \in \mathcal{B}$, and $\mathcal{B}$ is nice.

        Put
\begin{equation}
\label{U(B_1) = bigcup {B : B in mathcal{C}(B_1)}}
        U(B_1) = \bigcup \{B : B \in \mathcal{C}(B_1)\},
\end{equation}
so that $B_1 \subseteq U(B_1)$, by construction.  Let $x_1$ be an
element of $B_1$ again, so that the elements of $\mathcal{C}(B_1)$ may
be considered as balls centered at $x_1$, as before.  If
$\mathcal{C}(B_1)$ has a maximal element, then it is the same as
$U(B_1)$ as a subset of $M$.  If $\mathcal{C}(B_1)$ contains balls of
arbitrarily large radius, then $U(B_1) = M$.  Of course, if $M$ is
bounded, then $M$ may be considered as a ball centered at $x_1$ with
sufficiently large radius.  Otherwise, suppose that $\mathcal{C}(B_1)$
does not have a maximal element, and that the elements of
$\mathcal{C}(B_1)$ can be expressed as open or closed balls with
bounded radius.  In this case, $U(B_1)$ is an open ball centered at
$x_1$ with positive finite radius, which is the supremum of the radii
of the elements of $\mathcal{C}(B_1)$.

        Observe that
\begin{equation}
\label{bigcup {U(B_1) : B_1 in mathcal{B}} = bigcup {B : B in mathcal{B}}}
 \bigcup \{U(B_1) : B_1 \in \mathcal{B}\} = \bigcup \{B : B \in \mathcal{B}\}.
\end{equation}
More precisely, the left side of (\ref{bigcup {U(B_1) : B_1 in
    mathcal{B}} = bigcup {B : B in mathcal{B}}}) is contained in the
right side of (\ref{bigcup {U(B_1) : B_1 in mathcal{B}} = bigcup {B :
    B in mathcal{B}}}), because $U(B_1)$ is a union of elements
$\mathcal{B}$ for every $B_1 \in \mathcal{B}$, by construction.
Similarly, the right side of (\ref{bigcup {U(B_1) : B_1 in mathcal{B}}
  = bigcup {B : B in mathcal{B}}}) is contained in the left side of
(\ref{bigcup {U(B_1) : B_1 in mathcal{B}} = bigcup {B : B in
    mathcal{B}}}), because $B_1 \subseteq U(B_1)$ for every $B_1 \in
\mathcal{B}$.

        Suppose that $B_1, B_2 \in \mathcal{B}$ satisfy $B_1 \subseteq
B_2$, so that $B_2 \in \mathcal{C}(B_1)$, and in fact
$\mathcal{C}(B_2) \subseteq \mathcal{C}(B_1)$.  Let us check that
\begin{equation}
\label{U(B_1) = U(B_2)}
        U(B_1) = U(B_2).
\end{equation}
Of course, $U(B_2) \subseteq U(B_1)$, because $\mathcal{C}(B_2)
\subseteq \mathcal{C}(B_1)$.  To get the opposite inclusion, let $B
\in \mathcal{C}(B_1)$ be given, and let us verify that $B \subseteq
U(B_2)$.  Note that either either $B \subseteq B_2$ or $B_2 \subseteq
B$, because $B_2 \in \mathcal{C}(B_1)$ and the elements of
$\mathcal{C}(B_1)$ are linearly ordered by inclusion.  If $B \subseteq
B_2$, then $B \subseteq U(B_2)$, since $B_2 \subseteq U(B_2)$, by
construction.  Otherwise, if $B_2 \subseteq B$, then $B \in
\mathcal{C}(B_2)$, and hence $B \subseteq U(B_2)$, as desired.

        If $B_1, B_1' \in \mathcal{B}$ are both contained in $B_2 \in
\mathcal{B}$, then it follows that
\begin{equation}
\label{U(B_1) = U(B_1')}
        U(B_1) = U(B_1').
\end{equation}
Conversely, suppose that $B_1, B_1' \in \mathcal{B}$ have the property that
\begin{equation}
\label{U(B_1) cap U(B_1') ne emptyset}
        U(B_1) \cap U(B_1') \ne \emptyset.
\end{equation}
This implies that there are $B \in \mathcal{C}(B_1)$, $B' \in
\mathcal{C}(B_1')$ such that $B \cap B' \ne \emptyset$, by definition
of $U(B_1)$, $U(B_1')$.  In this case, $B \subseteq B'$ or $B'
\subseteq B$, as mentioned earlier in the section.  Put $B_2 = B'$
when $B \subseteq B'$, and $B_2 = B$ when $B' \subseteq B$.  Thus $B_2
\in \mathcal{B}$, because $B, B' \in \mathcal{B}$, by definition of
$\mathcal{C}(B_1)$, $\mathcal{C}(B_1')$.  We also have that $B_1
\subseteq B$ and $B_1' \subseteq B'$, by definition of
$\mathcal{C}(B_1)$, $\mathcal{C}(B_1')$, and hence that $B_1, B_1'
\subseteq B_2$.  This brings us back to the situation described at the
beginning of the paragraph, so that (\ref{U(B_1) = U(B_1')}) holds
under these conditions.  It follows that for each $B_1, B_1' \in
\mathcal{B}$, either
\begin{equation}
\label{U(B_1) cap U(B_1') = emptyset}
        U(B_1) \cap U(B_1') = \emptyset,
\end{equation}
or (\ref{U(B_1) = U(B_1')}) holds.

        As a basic class of examples, let $V$ be a nonempty open subset
of $M$, and let $\mathcal{B}$ be the collection of open or closed balls
in $M$ with positive radius that are contained in $V$.  Thus
\begin{equation}
\label{bigcup {B : B in mathcal{B}} = V}
        \bigcup \{B : B \in \mathcal{B}\} = V
\end{equation}
in this case.  If $V = M$, then $\mathcal{B}$ consists of all open or
closed balls in $M$ with positive radius, and $U(B_1) = M$ for every
$B_1 \in \mathcal{B}$.  Otherwise, if $V \ne M$, then for each $B_1
\in \mathcal{B}$, the elements of $\mathcal{C}(B_1)$ have bounded
radius.  This implies that $U(B_1)$ is an open or closed ball of
finite radius for each $B_1 \in \mathcal{B}$.  By construction,
\begin{equation}
\label{U(B_1) subseteq V}
        U(B_1) \subseteq V
\end{equation}
for every $B_1 \in \mathcal{B}$, because $U(B_1)$ is a union of
elements of $\mathcal{B}$, each of which is contained in $V$.  It
follows that $U(B_1) \in \mathcal{B}$ for every $B_1 \in \mathcal{B}$,
and these are the maximal elements of $\mathcal{B}$.

        Now let $\alpha$ be a nonnegative real number, and let
$\mathcal{B}$ be a nonempty collection of open or closed balls in $M$
such that
\begin{equation}
\label{sum_{B in mathcal{B}} (diam B)^alpha < infty}
        \sum_{B \in \mathcal{B}} (\diam B)^\alpha < \infty.
\end{equation}
If $B_1 \in \mathcal{B}$ and $\diam B_1 > 0$, then it is easy to see
that $\mathcal{C}(B_1)$ has only finitely many elements, because of
(\ref{sum_{B in mathcal{B}} (diam B)^alpha < infty}).  In particular,
this implies that $\mathcal{C}(B_1)$ has a maximal element, as before.
Otherwise, if $\diam B_1 = 0$, then $B_1$ consists of a single point.
This happens for closed balls of radius $0$, but it can also happen
for open or closed balls of positive radius around an isolated point
in $M$.  If $B_1$ is not maximal, then there is a $B_2 \in
\mathcal{B}$ such that $B_1 \subseteq B_2$ and $B_1 \ne B_2$.  This
implies that $\diam B_2 > 0$, and hence that $B_2$ is contained in a
maximal element of $\mathcal{B}$.  This shows that $\mathcal{B}$ is
nice under these conditions.

\section{Haar and Hausdorff measures}
\label{haar, hausdorff measures}

        Let $k$ be a field, and let $|\cdot|$ be an ultrametric absolute
value function on $k$.  In this section, we suppose that $|\cdot|$ is
nontrivial on $k$, and that $k$ is locally compact with respect to the
corresponding ultrametric.  This implies that closed balls in $k$ are
compact, as in Section \ref{local compactness}.  It follows that
$|\cdot|$ is a discrete absolute value function on $k$, in the sense
described in Section \ref{discrete absolute value functions}, and as
discussed at the end of Section \ref{ultrametric case}.  In this
situation, we also have that the corresponding residue field
(\ref{overline{B}(0, 1) / B(0, 1)}) is finite, as mentioned near the
end of Section \ref{ultrametric case}.

        Let $\rho_1$ be as in (\ref{rho_1 = sup {|x| : x in k, |x| < 1}})
in Section \ref{discrete absolute value functions}, so that $0 <
\rho_1 < 1$, because $|\cdot|$ is discrete and nontrivial.  Thus the
positive values of $|x|$ on $k$ are integer powers of $\rho_1$, as in
Section \ref{discrete absolute value functions}.  In particular,
\begin{equation}
\label{diam overline{B}(x, rho_1^j) = rho_1^j}
        \diam \overline{B}(x, \rho_1^j) = \rho_1^j
\end{equation}
for each $x \in k$ and $j \in {\bf Z}$, where $\overline{B}(x, r)$ is
the closed ball in $k$ centered at $x \in k$ with radius $r > 0$ with
respect to the ultrametric associated to $|\cdot|$.  Remember that
$\overline{B}(0, 1)$ is a subring of $k$, $B(0, 1)$ is an ideal in
$\overline{B}(0, 1)$, and that the quotient ring
\begin{equation}
\label{overline{B}(0, 1) / B(0, 1), 2}
        \overline{B}(0, 1) / B(0, 1)
\end{equation}
is the residue field corresponding to $|\cdot|$ on $k$, as in Section
\ref{ultrametric case}.  Let $\alpha$ be the positive real number that
satisfies
\begin{equation}
\label{rho_1^{-alpha} = num (overline{B}(0, 1) / B(0, 1))}
        \rho_1^{-\alpha} = \#(\overline{B}(0, 1) / B(0, 1)),
\end{equation}
where $\#(\overline{B}(0, 1) / B(0, 1))$ is the number of elements of
the residue field (\ref{overline{B}(0, 1) / B(0, 1), 2}).

        Let $H$ be Haar measure on $k$, normalized so that
\begin{equation}
\label{H(overline{B}(0, 1)) = 1, 2}
        H(\overline{B}(0, 1)) = 1.
\end{equation}
Observe that
\begin{equation}
\label{H(overline{B}(x, rho_1^j)) = rho_1^{alpha j}}
        H(\overline{B}(x, \rho_1^j)) = \rho_1^{\alpha \, j}
\end{equation}
for every $x \in k$ and $j \in {\bf Z}$, by (\ref{H(overline{B}(x,
  rho_1^j)) = N^{-j}}) in Section \ref{haar measure, 2}.  Equivalently,
\begin{equation}
\label{H(overline{B}(x, rho_1^j)) = (diam overline{B}(x, rho_1^j))^alpha}
        H(\overline{B}(x, \rho_1^j)) = (\diam \overline{B}(x, \rho_1^j))^\alpha
\end{equation}
for every $x \in k$ and $j \in {\bf Z}$, by (\ref{diam overline{B}(x,
  rho_1^j) = rho_1^j}).  It follows that
\begin{equation}
\label{H(A) le (diam A)^alpha}
        H(A) \le (\diam A)^\alpha
\end{equation}
for every bounded Borel set $A \subseteq k$.  More precisely, if
$\diam A = 0$, then $A$ contains at most one element, and (\ref{H(A)
  le (diam A)^alpha}) says that $H(A) = 0$.  Otherwise, if the
diameter of $A$ is positive, then it is equal to $\rho_1^j$ for some
$j \in {\bf Z}$.  In this case, $A$ is contained in a closed ball of
radius $\rho_1^j$, so that (\ref{H(A) le (diam A)^alpha}) follows from
(\ref{H(overline{B}(x, rho_1^j)) = rho_1^{alpha j}}).  Using
(\ref{H(A) le (diam A)^alpha}), we get that
\begin{equation}
\label{H(E) le H^alpha_{con}(E)}
        H(E) \le H^\alpha_{con}(E)
\end{equation}
for every Borel set $E \subseteq k$, where $H^\alpha_{con}(E)$ is the
$\alpha$-dimensional Hausdorff content of $E$ with respect to the
ultrametric associated to $|\cdot|$ on $k$.  This is the same as
(\ref{mu(E) le C H^alpha_{con}(E)}) in Section \ref{euclidean spaces},
with $\mu = H$ and $C = 1$.

        Remember that $\overline{B}(0, 1)$ can be expressed as the
union of finitely many pairwise disjoint open balls of radius $1$,
where the number of these open balls of radius $1$ is equal to
(\ref{rho_1^{-alpha} = num (overline{B}(0, 1) / B(0, 1))}).
Equivalently, $\overline{B}(0, 1)$ can be expressed as the union of
the same number of pairwise-disjoint closed balls of radius $\rho_1$,
because an open ball in $k$ of radius $1$ is the same as a closed ball
in $k$ centered at the same point with radius $\rho_1$ in this
situation.  It follows that for each $x \in k$ and $j \in {\bf Z}$,
$\overline{B}(x, \rho_1^j)$ can be expressed as the union of the same
number of pairwise-disjoint closed balls of radius $\rho_1^{j + 1}$.
If $l$ is any positive integer, then one can repeat the process to
get that $\overline{B}(x, \rho_1^j)$ can be expressed as the union of
\begin{equation}
\label{rho_1^{- alpha l}}
        \rho_1^{- \alpha \, l}
\end{equation}
pairwise-disjoint closed balls of radius $\rho_1^{j + l}$.

        Using this and (\ref{diam overline{B}(x, rho_1^j) = rho_1^j}),
one can check directly that
\begin{equation}
\label{H^alpha_delta(overline{B}(x, rho_1^j)) le rho_1^{alpha j}}
        H^\alpha_\delta(\overline{B}(x, \rho_1^j)) \le \rho_1^{\alpha \, j}
\end{equation}
for every $x \in k$, $j \in {\bf Z}$, and $\delta > 0$.  This implies that
\begin{equation}
\label{H^alpha(overline{B}(x, rho_1^j)) le rho_1^{alpha j}}
        H^\alpha(\overline{B}(x, \rho_1^j)) \le \rho_1^{\alpha \, j}
\end{equation}
for every $x \in k$ and $j \in {\bf Z}$, by taking the supremum of the
left side of (\ref{H^alpha_delta(overline{B}(x, rho_1^j)) le
  rho_1^{alpha j}}) over $\delta > 0$.  It follows that
\begin{equation}
\label{H^alpha(A) le (diam A)^alpha, 2}
        H^\alpha(A) \le (\diam A)^\alpha
\end{equation}
for every $A \subseteq k$, as before.  We also get that $H^\alpha(E) \le
H^\alpha_{con}(E)$ for every $E \subseteq k$, and hence that
\begin{equation}
\label{H^alpha(E) = H^alpha_{con}(E)}
        H^\alpha(E) = H^\alpha_{con}(E),
\end{equation}
since $H^\alpha_{con}(E) \le H^\alpha(E)$ automatically.

        Observe that
\begin{equation}
\label{H^alpha_{con}(overline{B}(x, rho_1^j)) ge rho_1^{alpha j}}
        H^\alpha_{con}(\overline{B}(x, \rho_1^j)) \ge \rho_1^{\alpha \, j}
\end{equation}
for every $x \in k$ and $j \in {\bf Z}$, by (\ref{H(overline{B}(x,
  rho_1^j)) = rho_1^{alpha j}}) and (\ref{H(E) le H^alpha_{con}(E)}).
This can also be verified more directly, as follows.  In order to
estimate $H^\alpha_{con}(E)$, it suffices to consider coverings of
$\overline{B}(x, \rho_1^j)$ by finitely or countably many closed balls
of positive radius in $k$, as in Section \ref{hausdorff content}.  Of
course, this uses the fact that the metric associated to $|\cdot|$ on
$k$ is an ultrametric.  More precisely, it suffices to consider
coverings of $\overline{B}(x, \rho_1^j)$ by finitely many closed balls
of positive radius, because closed balls of positive radius are open
sets in $k$, and $\overline{B}(x, \rho_1^k)$ is compact.  In this
situation, one can reduce further to coverings of $\overline{B}(x,
\rho_1^j)$ by finitely many closed balls of equal radius, using the
earlier arguments about expressing closed balls as unions of balls of
smaller radius.  If $l$ is a nonnegative integer, then (\ref{rho_1^{-
    alpha l}}) is the minimum number of closed balls of radius
$\rho_1^{j + l}$ needed to cover $\overline{B}(x, \rho_1^j)$.  This
implies (\ref{H^alpha_{con}(overline{B}(x, rho_1^j)) ge rho_1^{alpha
    j}}), as desired.

        Combining (\ref{H^alpha(overline{B}(x, rho_1^j)) le rho_1^{alpha j}})
and (\ref{H^alpha_{con}(overline{B}(x, rho_1^j)) ge rho_1^{alpha j}}),
we get that
\begin{equation}
\label{H^alpha(overline{B}(x, rho_1^j)) = ... = rho_1^j}
 H^\alpha(\overline{B}(x, \rho_1^j)) = H^\alpha_{con}(\overline{B}(x, \rho_1^j))
                                    = \rho_1^j 
\end{equation}
for every $x \in k$ and $j \in {\bf Z}$.  It follows that
\begin{equation}
\label{H^alpha(E) = H(E)}
        H^\alpha(E) = H(E)
\end{equation}
when $E$ is a ball in $k$, by (\ref{H(overline{B}(x, rho_1^j)) =
  rho_1^{alpha j}}) and (\ref{H^alpha(overline{B}(x, rho_1^j)) = ... =
  rho_1^j}).  If $E$ is a nonempty proper open subset of $k$, then $E$
can be expressed as a union of pairwise-disjoint balls of positive
radius in $k$, as in the previous section.  Note that $k$ is
separable, because closed balls in $k$ are separable.  This implies
that any collection of pairwise-disjoint balls in $k$ of positive
radius can have only finitely or countably many elements.  In
particular, if $E$ is a nonempty proper open subset of $k$, then $E$
can be expressed as the union of finitely or countably many
pairwise-disjoint balls in $k$.  This permits (\ref{H^alpha(E) =
  H(E)}) to be derived from the corresponding statement for balls in
$k$ in this case.  One can check that (\ref{H^alpha(E) = H(E)}) holds
for Borel sets $E \subseteq k$, using this and outer regularity
properties of $H$ and $H^\alpha$.

        Of course, the metric on $k$ associated to $|\cdot|$ is
invariant under translations, by construction.  This implies that the
corresponding Hausdorff measure of any dimension is invariant under
translations as well.  One can check that $H^\alpha$ satisfies the
other requirements of Haar measure on $k$, and indeed this may be
considered as a way to construct Haar measure on $k$.  Similarly,
the argument in the preceding paragraph may be considered as a way
to deal with uniqueness of Haar measure in this situation.

        Let $p$ be a prime number, and suppose that $k = {\bf Q}_p$,
equipped with the $p$-adic metric.  This satisfies the conditions
mentioned at the beginning of the section, with $\rho_1 = 1/p$.  In
this case, (\ref{overline{B}(0, 1) / B(0, 1), 2}) has exactly $p$
elements, and (\ref{rho_1^{-alpha} = num (overline{B}(0, 1) / B(0, 1))})
corresponds to taking $\alpha = 1$.  Thus $1$-dimensional Hausdorff
measure is the same as Haar measure on ${\bf Q}_p$.

        Let $k$ be any field with an ultrametric absolute value function
$|\cdot|$ that satisfies the conditions mentioned at the beginning of
the section again.  Also let $n$ be a positive integer, so that $k^n$
may be considered as an $n$-dimensional vector space over $k$.
Similarly, $k^n$ may be considered as a locally compact commutative
group with respect to addition, using the product topology associated
to the topology on $k$ determined by the ultrametric corresponding to
$|\cdot|$.  Remember that Haar measure on $k^n$ basically corresponds
to a product of $n$ copies of Haar measure on $k$, as in Section
\ref{haar measure, 3}.

        Let $N_0$ be the ultranorm on $k^n$ defined in 
(\ref{N_0(v) = max(|v_1|, ldots, |v_n|)}) in Section \ref{norms, ultranorms},
so that closed balls in $k^n$ with respect to $N_0$ are the same as
products of $n$ closed balls in $k$ with respect to $|\cdot|$, with
the same radius.  If $\alpha$ is as in (\ref{rho_1^{-alpha} = num
  (overline{B}(0, 1) / B(0, 1))}), then $(\alpha \, n)$-dimensional
Hausdorff measure on $k^n$ with respect to the ultrametric associated
to $N_0$ satisfies the requirements of Haar measure, and has other
properties like those in the $n = 1$ case.  Any other norm on $k^n$ is
equivalent to $N_0$, as in Section \ref{finite-dimensional vector
  spaces}, which implies that the corresponding Hausdorff measures are
bounded by constant multiples of each other too.  Using this, one can
check that $(\alpha \, n)$-dimensional Hausdorff measure on $k^n$ with
respect to the metric associated to $N$ also satisfies the
requirements of Haar measure.

\section{Similarities}
\label{similarities}

        Let $k$ be a field, let $|\cdot|$ be an ultrametric
absolute value function on $k$, and suppose that $k$ is complete with
respect to the ultrametric corresponding to $|\cdot|$.  Also let $a_0,
a_1, a_2, a_3, \ldots$ be a sequence of elements of $k$, and consider
the corresponding power series
\begin{equation}
\label{f(x) = sum_{j = 0}^infty a_j x^j, 6}
        f(x) = \sum_{j = 0}^\infty a_j \, x^j.
\end{equation}
If
\begin{equation}
\label{lim_{j to infty} |a_j| r^j = 0, 7}
        \lim_{j \to \infty} |a_j| \, r^j = 0
\end{equation}
for some $r > 0$, then $f(x)$ may be considered as a $k$-valued
function on
\begin{equation}
\label{D = overline{B}(0, r)}
        D = \overline{B}(0, r).
\end{equation}
Alternatively, if $\rho > 0$, and (\ref{lim_{j to infty} |a_j| r^j =
  0, 7}) holds for every $r \in (0, \rho)$, then $f(x)$ may be
considered as a $k$-valued function on
\begin{equation}
\label{D = B(0, rho)}
        D = B(0, \rho).
\end{equation}
If (\ref{lim_{j to infty} |a_j| r^j = 0, 7}) holds for every $r > 0$,
then we can take $\rho = +\infty$ and $D = k$.  Note that the series
expansion (\ref{sum_{j = 1}^infty j cdot a_j x^{j - 1}}) in Section
\ref{differentiation, lipschitz conditions} for $f'(x)$ converges for
every $x \in D$.  This uses the fact that $|j \cdot a_j| \le |a_j|$
for each $j$, by the ultrametric version of the triangle inequality.

        Let $B$ be an open or closed ball centered at a point $x \in k$
with radius $t > 0$ such that $B \subseteq D$.  Let us say that $B$ is
``admissible'' if it satisfies the following four conditions.  First,
$f'(x) \ne 0$.  Second,
\begin{equation}
\label{|f'(y)| = |f'(x)|}
        |f'(y)| = |f'(x)|
\end{equation}
for every $y \in B$.  Third,
\begin{equation}
\label{|f(y) - f(w)| = |f'(x)| |y - w|}
        |f(y) - f(w)| = |f'(x)| \, |y - w|
\end{equation}
for every $y, w \in B$.  To state the fourth condition, let $B'$ be
the ball in $k$ centered at $f(x)$, with radius equal to $t \,
|f'(x)|$, and which is an open ball when $B$ is an open ball, and a
closed ball when $B$ is a closed ball.  Thus (\ref{|f(y) - f(w)| =
  |f'(x)| |y - w|}) implies that $f(B) \subseteq B'$, and the fourth
condition is that
\begin{equation}
\label{f(B) = B'}
        f(B) = B'.
\end{equation}
Of course, $B$ may be considered as being centered at any of its
elements, as mentioned at the beginning of Section \ref{a simple
  covering argument}.  However, it is easy to see that the
admissibility of $B$ does not depend on the choice $x$ of center of
$B$.  This also uses the fact that $B'$ may be considered as being
centered at any of its elements, which include the elements of $f(B)$.
If $x \in D$ and $f'(x) \ne 0$, then the discussion in Sections
\ref{contraction mappings} and \ref{contraction mappings, 2} gives a
criterion for the admissibility of balls centered at $x$.  In
particular, this criterion implies that balls centered at $x$ with
sufficiently small radius are admissible.

        If $\mathcal{B}$ is the collection of these admissible balls, then
\begin{equation}
\label{bigcup {B : B in mathcal{B}} = {x in D : f'(x) ne 0}}
        \bigcup \{B : B \in \mathcal{B}\} = \{x \in D : f'(x) \ne 0\}.
\end{equation}
More precisely, the left side of (\ref{bigcup {B : B in mathcal{B}} =
  {x in D : f'(x) ne 0}}) is contained in the right side of
(\ref{bigcup {B : B in mathcal{B}} = {x in D : f'(x) ne 0}}), because
of the first and second conditions in the definition of admissibility.
Similarly, the right side of (\ref{bigcup {B : B in mathcal{B}} = {x
    in D : f'(x) ne 0}}) is contained in the left side of (\ref{bigcup
  {B : B in mathcal{B}} = {x in D : f'(x) ne 0}}), because every $x \in
D$ with $f'(x) \ne 0$ is contained in an admissible ball, as mentioned
at the end of the preceding paragraph.  If $D$ is as in (\ref{D =
  overline{B}(0, r)}) and $j \cdot a_j \ne 0$ for some positive
integer $j$, then Strassmann's theorem implies that $f'(x) = 0$ for
only finitely many $x \in D$.  Otherwise, if $D$ is as in (\ref{D =
  B(0, rho)}) or $D = k$, and if $j \cdot a_j \ne 0$ for some $j \ge
1$, then one can use Strassmann's theorem to get that $f'(x) = 0$ for
at most finitely or countably many $x \in D$.

        Let $B_1$ be an admissible ball, let $\mathcal{C}(B_1)$ be as
in (\ref{mathcal{C}(B_1) = {B in mathcal{B} : B_1 subseteq B}}), and
let $U(B_1)$ be as in (\ref{U(B_1) = bigcup {B : B in
    mathcal{C}(B_1)}}).  If the elements of $\mathcal{C}(B_1)$ have
bounded radii, then we have seen that $U(B_1)$ is an open or closed
ball with finite radius.  In this situation, one can also check that
$U(B_1)$ is automatically admissible, because it is the union of a
collection of admissible balls centered at the same point.  This
implies that $U(B_1)$ is an element of $\mathcal{C}(B_1)$, which is
automatically maximal, and hence that $U(B_1)$ is a maximal element of
$\mathcal{B}$.  Note that $\mathcal{C}(B_1)$ only has finitely many
elements when the elements of $\mathcal{C}(B_1)$ have bounded radii
and $|\cdot|$ is a discrete absolute value function on $k$, as in
Section \ref{discrete absolute value functions}.

        Let us suppose from now on in this section that $|\cdot|$
is nontrivial on $k$, and that $k$ is locally compact.  As in the
previous section, we let $H$ be Haar measure on $k$, normalized as in
(\ref{H(overline{B}(0, 1)) = 1, 2}).  We also let $\rho_1$ be as in
(\ref{rho_1 = sup {|x| : x in k, |x| < 1}}) in Section \ref{discrete
  absolute value functions}, and we let $\alpha > 0$ be as in
(\ref{rho_1^{-alpha} = num (overline{B}(0, 1) / B(0, 1))}).  Put
\begin{equation}
\label{x E = {x y : y in E}}
        x \, E = \{x \, y : y \in E\}
\end{equation}
for each $x \in k$ and $E \subseteq k$.  Of course, (\ref{x E = {x y :
    y in E}}) reduces to $\{0\}$ when $x = 0$ and $E \ne \emptyset$,
and (\ref{x E = {x y : y in E}}) is the empty set for every $x \in k$
when $E = \emptyset$.  If $x \ne 0$, then multiplication by $x$
defines a homeomorphism on $k$, and the mapping from $E \subseteq k$
to (\ref{x E = {x y : y in E}}) preserves open sets, closed sets,
compact sets, and Borel sets.  Let us check that
\begin{equation}
\label{H(x E) = |x|^alpha H(E)}
        H(x \, E) = |x|^\alpha \, H(E)
\end{equation}
for every $x \in k$ and Borel set $E \subseteq k$.  This is trivial
when $x = 0$, and so we may as well restrict our attention to $x \ne
0$.  In this case, it is easy to see that $H(x \, E)$ also satisfies
the requirements of Haar measure on $k$.  The uniqueness of Haar
measure implies that $H(x \, E)$ is equal to a positive real number
times $H(E)$, where this positive real number depends on $x$ but not
$E$.  In order to determine this positive real number, it suffices
to take $E = \overline{B}(0, 1)$, for which we have that
\begin{equation}
\label{H(x overline{B}(0, 1)) = H(overline{B}(0, |x|)) = |x|^alpha}
        H(x \, \overline{B}(0, 1)) = H(\overline{B}(0, |x|)) = |x|^\alpha,
\end{equation}
by (\ref{H(overline{B}(x, rho_1^j)) = rho_1^{alpha j}}).  This implies
(\ref{H(x E) = |x|^alpha H(E)}), because of the normalization
(\ref{H(overline{B}(0, 1)) = 1, 2}) for $H$.

        Let $x_0 \in k$ and $j_0 \in {\bf Z}$ be given, and let $g$ be
a mapping from $\overline{B}(x_0, \rho_1^{j_0})$ into $k$.  Let $l_0$
be another integer, and suppose that
\begin{equation}
\label{|g(x) - g(y)| = rho_1^{l_0} |x - y|}
        |g(x) - g(y)| = \rho_1^{l_0} \, |x - y|
\end{equation}
for every $x, y \in \overline{B}(x_0, \rho_1^{j_0})$.  This implies
that $g(\overline{B}(x_0, \rho_1^{j_0})) \subseteq
\overline{B}(g(x_0), \rho_1^{j_0 + l_0})$, and in fact one can show that
\begin{equation}
\label{g(overline{B}(x_0, rho_1^{j_0})) = ...}
 g(\overline{B}(x_0, \rho_1^{j_0})) = \overline{B}(g(x_0), \rho_1^{j_0 + l_0})
\end{equation}
in this situation.  Let us simply take this as an additional condition
here, for the sake of convenience.  If $x \in \overline{B}(x_0,
\rho_1^{j_0})$, $j \in {\bf Z}$, and $j \ge j_0$, then $\rho_1^j \le
\rho_1^{j_0}$, and hence
\begin{equation}
\label{overline{B}(x, rho_1^j) subseteq overline{B}(x_0, rho_1^{j_0})}
        \overline{B}(x, \rho_1^j) \subseteq \overline{B}(x_0, \rho_1^{j_0}),
\end{equation}
by the ultrametric version of the triangle inequality.  In this case,
it is easy to see that
\begin{equation}
\label{g(overline{B}(x, rho_1^j)) = overline{B}(g(x), rho_1^{j + l_0})}
        g(\overline{B}(x, \rho_1^j)) = \overline{B}(g(x), \rho_1^{j + l_0}),
\end{equation}
using (\ref{|g(x) - g(y)| = rho_1^{l_0} |x - y|}) and
(\ref{g(overline{B}(x_0, rho_1^{j_0})) = ...}).  This can also be
shown directly, as for (\ref{g(overline{B}(x_0, rho_1^{j_0})) = ...}).

        If $E \subseteq \overline{B}(x_0, \rho_1^{j_0})$ is a Borel set,
then one can check that
\begin{equation}
\label{H(g(E)) = rho_1^{alpha l_0} H(E)}
        H(g(E)) = \rho_1^{\alpha \, l_0} \, H(E).
\end{equation}
If $H$ is identified with $\alpha$-dimensional Hausdorff measure on
$k$, as in the previous section, then (\ref{H(g(E)) = rho_1^{alpha
    l_0} H(E)}) can be obtained directly from (\ref{|g(x) - g(y)| =
  rho_1^{l_0} |x - y|}).  Otherwise, one can first check that
(\ref{H(g(E)) = rho_1^{alpha l_0} H(E)}) holds when $E$ is a ball
contained in $\overline{B}(x_0, \rho_1^{j_0})$, using
(\ref{H(overline{B}(x, rho_1^j)) = rho_1^{alpha j}}) and
(\ref{g(overline{B}(x, rho_1^j)) = overline{B}(g(x), rho_1^{j +
    l_0})}).  If $E$ is an open set contained in $\overline{B}(x_0,
\rho_1^{j_0})$, then $E$ can be expressed as the union of the maximal
balls contained in $E$, as in Section \ref{a simple covering
  argument}.  The maximal balls contained in $E$ are pairwise
disjoint, as before, and there can be only finitely or countably many
of these maximal balls contained in $E$, because $\overline{B}(x_0,
\rho_1^{j_0})$ is compact and hence separable.  This permits
(\ref{H(g(E)) = rho_1^{alpha l_0} H(E)}) to be derived from the
analogous statement for balls contained in $\overline{B}(x_0,
\rho_1^{j_0})$, by countable additivity of Haar measure.  Once one has
(\ref{H(g(E)) = rho_1^{alpha l_0} H(E)}) for open subsets of
$\overline{B}(x_0, \rho_1^{j_0})$, the analogous statement for Borel
sets can be derived from the outer regularity of Haar measure.

\backmatter

\newpage

\addcontentsline{toc}{chapter}{Index}

\printindex

\end{document}